\let\noarrow = t

\input eplain


\magnification=\magstep1

\topskip1truecm
\def\pagewidth#1{
  \hsize=#1
}

\def\pageheight#1{
  \vsize=#1
}

\pageheight{23.5truecm} \pagewidth{16truecm}

\abovedisplayskip=3mm \belowdisplayskip=3mm
\abovedisplayshortskip=0mm \belowdisplayshortskip=2mm
\def\spacing{{\smallskip}}
\frenchspacing
\parindent1pc

\normalbaselineskip=13pt \baselineskip=13pt
\def\spacing{{\smallskip}}


\newdimen\abstractmargin
\abstractmargin=3pc


\newdimen\footnotemargin
\footnotemargin=1pc


\font\eightrm=cmr8 \relax 
\font\sixrm=cmr6 \relax 
\font\eighti=cmmi8 \relax     \skewchar\eighti='177 
\font\sixi=cmmi6 \relax       \skewchar\sixi='177   
\font\eightsy=cmsy8 \relax    \skewchar\eightsy='60 
\font\sixsy=cmsy6 \relax      \skewchar\sixsy='60   
\font\eightbf=cmbx8 \relax 
\font\sixbf=cmbx6 \relax   
\font\eightit=cmti8 \relax 
\font\eightsl=cmsl8 \relax 
\font\eighttt=cmtt8 \relax 

\catcode`\@=11
\newskip\ttglue

\def\eightpoint{\def\rm{\fam0\eightrm}%
 \textfont0=\eightrm \scriptfont0=\sixrm
 \scriptscriptfont0=\fiverm
 \textfont1=\eighti \scriptfont1=\sixi
 \scriptscriptfont0=\fivei
 \textfont2=\eightsy \scriptfont2=\sixsy
 \scriptscriptfont2=\fivesy
 \textfont3=\tenex \scriptfont3=\tenex
 \scriptscriptfont3=\tenex
 \textfont\itfam\eightit \def\it{\fam\itfam\eightit}%
 \textfont\slfam\eightsl \def\sl{\fam\slfam\eightsl}%
 \textfont\ttfam\eighttt \def\tt{\fam\ttfam\eighttt}%
 \textfont\bffam\eightbf \scriptfont\bffam\sixbf
   \scriptscriptfont\bffam\fivebf \def\bf{\fam\bffam\eightit}%
 \tt \ttglue=.5em plus.25em minus.15em
 \normalbaselineskip=9pt
 \setbox\strutbox\hbox{\vrule height7pt depth3pt width0pt}%
 \let\sc=\sixrm \let\big=\eifgtbig \normalbaselines\rm}


 \font\titlefont=cmbx12 scaled\magstep1
 \font\sectionfont=cmbx12
 \font\ssectionfont=cmsl10
 \font\claimfont=cmsl10

 \font\normalfont=cmr10

\catcode`\@=11 \font\teneusm=eusm10 
\font\seveneusm=eusm7  \font\fiveeusm=eusm5
\newfam\eusmfam \textfont\eusmfam=\teneusm
\scriptfont\eusmfam=\seveneusm \scriptscriptfont\eusmfam=\fiveeusm
\def\hexnumber@#1{\ifcase#1
0\or1\or2\or3\or4\or5\or6\or7\or8\or9\or         A\or B\or C\or
D\or E\or F\fi } \edef\eusm@{\hexnumber@\eusmfam}
\def\euscr{\ifmmode\let\next\euscr@\else
\def\next{\errmessage{Use \string\euscr\space only in math mode}}\fi\next}
\def\euscr@#1{{\euscr@@{#1}}} \def\euscr@@#1{\fam\eusmfam#1} \catcode`\@=12

\catcode`\@=11 \font\teneuex=euex10 
 \font\seveneuex=euex7  \newfam\euexfam
\textfont\euexfam=\teneuex  \scriptfont\euexfam=\seveneuex
 \def\hexnumber@#1{\ifcase#1
0\or1\or2\or3\or4\or5\or6\or7\or8\or9\or         A\or B\or C\or
D\or E\or F\fi } \edef\euex@{\hexnumber@\euexfam}
\def\euscrex{\ifmmode\let\next\euscrex@\else
\def\next{\errmessage{Use \string\euscrex\space only in math mode}}\fi\next}
\def\euscrex@#1{{\euscrex@@{#1}}} \def\euscrex@@#1{\fam\euexfam#1}
\catcode`\@=12

\catcode`\@=11 \font\teneufb=eufb10 
\font\seveneufb=eufb7  \font\fiveeufb=eufb5
\newfam\eufbfam \textfont\eufbfam=\teneufb
\scriptfont\eufbfam=\seveneufb \scriptscriptfont\eufbfam=\fiveeufb
\def\hexnumber@#1{\ifcase#1
0\or1\or2\or3\or4\or5\or6\or7\or8\or9\or         A\or B\or C\or
D\or E\or F\fi } \edef\eufb@{\hexnumber@\eufbfam}
\def\euscrfb{\ifmmode\let\next\euscrfb@\else
\def\next{\errmessage{Use \string\euscrfb\space only in math mode}}\fi\next}
\def\euscrfb@#1{{\euscrfb@@{#1}}} \def\euscrfb@@#1{\fam\eufbfam#1}
\catcode`\@=12

\catcode`\@=11 \font\teneufm=eufm10 
\font\seveneufm=eufm7  \font\fiveeufm=eufm5
\newfam\eufmfam \textfont\eufmfam=\teneufm
\scriptfont\eufmfam=\seveneufm \scriptscriptfont\eufmfam=\fiveeufm
\def\hexnumber@#1{\ifcase#1
0\or1\or2\or3\or4\or5\or6\or7\or8\or9\or         A\or B\or C\or
D\or E\or F\fi } \edef\eufm@{\hexnumber@\eufmfam}
\def\euscrfm{\ifmmode\let\next\euscrfm@\else
\def\next{\errmessage{Use \string\euscrfm\space only in math mode}}\fi\next}
\def\euscrfm@#1{{\euscrfm@@{#1}}} \def\euscrfm@@#1{\fam\eufmfam#1}
\catcode`\@=12

\catcode`\@=11 \font\teneusb=eusb10 
\font\seveneusb=eusb7  \font\fiveeusb=eusb5
\newfam\eusbfam \textfont\eusbfam=\teneusb
\scriptfont\eusbfam=\seveneusb \scriptscriptfont\eusbfam=\fiveeusb
\def\hexnumber@#1{\ifcase#1
0\or1\or2\or3\or4\or5\or6\or7\or8\or9\or         A\or B\or C\or
D\or E\or F\fi } \edef\eusb@{\hexnumber@\eusbfam}
\def\euscrsb{\ifmmode\let\next\euscrsb@\else
\def\next{\errmessage{Use \string\euscrsb\space only in math mode}}\fi\next}
\def\euscrsb@#1{{\euscrsb@@{#1}}} \def\euscrsb@@#1{\fam\eusbfam#1}
\catcode`\@=12

\catcode`\@=11 \font\tenmsa=msam10 
\font\sevenmsa=msam7  \font\fivemsa=msam5
\font\tenmsb=msbm10  \font\sevenmsb=msbm7
 \font\fivemsb=msbm5 \newfam\msafam
\newfam\msbfam \textfont\msafam=\tenmsa
\scriptfont\msafam=\sevenmsa
  \scriptscriptfont\msafam=\fivemsa
\textfont\msbfam=\tenmsb  \scriptfont\msbfam=\sevenmsb
  \scriptscriptfont\msbfam=\fivemsb
\def\hexnumber@#1{\ifcase#1 0\or1\or2\or3\or4\or5\or6\or7\or8\or9\or
        A\or B\or C\or D\or E\or F\fi }
\edef\msa@{\hexnumber@\msafam} \edef\msb@{\hexnumber@\msbfam}
\mathchardef\square="0\msa@03 \mathchardef\subsetneq="3\msb@28
\mathchardef\supsetneq="3\msb@29 \mathchardef\ltimes="2\msb@6E
\mathchardef\rtimes="2\msb@6F \mathchardef\dabar="0\msa@39
\mathchardef\daright="0\msa@4B \mathchardef\daleft="0\msa@4C

\def\Bbb{\ifmmode\let\next\Bbb@\else
        \def\next{\errmessage{Use \string\Bbb\space only in math mode}}\fi\next}
\def\Bbb@#1{{\Bbb@@{#1}}}
\def\Bbb@@#1{\fam\msbfam#1}
\catcode`\@=12



\newcount\senu
\def\senum{\number\senu}
\newcount\ssnu
\def\ssnum{\number\ssnu}
\newcount\sssn
\def\sssnu{\number\sssn}
\newcount\fonu
\def\fonum{\number\fonu}

\def\num{{\senum.\ssnum}}
\def\snum{{\senum.\ssnum.\sssnu}}
\def\numfo{{\senum.\ssnum.\fonum}}
\def\snumfo{{\senum.\ssnum.\sssnu.\fonum}}


\outer\def\section#1\par{\vskip0pt
  plus.3\vsize\penalty-20\vskip0pt
  plus-.3\vsize\bigskip\vskip\parskip
  \message{#1}\centerline{\sectionfont\senum\enspace#1.}
  \nobreak\smallskip}

\def\endsection{\advance\senu by1\penalty-20\smallskip\ssnu=1}

\outer\def\ssection#1\par{\bigskip
  \message{#1}{\noindent\bf\num\ssectionfont\enspace#1.\thinspace}
  \nobreak\normalfont}

\def\endssection{\advance\ssnu by1\smallskip\ifdim\lastskip<\medskipamount
\removelastskip\penalty55\medskip\fi\fonu=1\normalfont\sssn=1}

\outer\def\sssection#1\par{\bigskip
  \message{#1}{\noindent\bf\snum\ssectionfont\enspace#1.\thinspace}
  \nobreak\normalfont}

\def\endsssection{\advance\sssn by1\smallskip\ifdim\lastskip<\medskipamount
\removelastskip\penalty55\medskip\fi\fonu=1\normalfont}

\def\proclaim #1\par{\bigskip
  \message{#1}{\noindent\bf\num\enspace#1.\thinspace}
  \nobreak\claimfont}
\def\sproclaim #1\par{\bigskip
  \message{#1}{\noindent\bf\snum\enspace#1.\thinspace}
  \nobreak\claimfont}

\def\cor{\proclaim Corollary\par}
\def\defi{\proclaim Definition\par}
\def\lemma{\proclaim Lemma\par}
\def\prop{\proclaim Proposition\par}
\def\rmk{\proclaim Remark\par\normalfont}
\def\thm{\proclaim Theorem\par}

\def\endcor{\endssection}
\def\enddefi{\endssection}
\def\endlemma{\endssection}
\def\endprop{\endssection}
\def\endrmk{\endssection}
\def\endthm{\endssection}

\def\scor{\sproclaim Corollary\par}
\def\sdefi{\sproclaim Definition\par}
\def\slemma{\sproclaim Lemma\par}
\def\sprop{\sproclaim Proposition\par}
\def\srmk{\sproclaim Remark\par\normalfont}
\def\sthm{\sproclaim Theorem\par}

\def\endscor{\endsssection}
\def\endsdefi{\endsssection}
\def\endslemma{\endsssection}
\def\endsprop{\endsssection}
\def\endsrmk{\endsssection}
\def\endsthm{\endsssection}

\def\Proof{{\noindent\sl Proof: \/}}


\def\maplefto#1{\ \smash{\mathop{\longleftarrow}\limits^{#1}}\ }
\def\tto{\longrightarrow}
\def\lra{\longrightarrow}
\def\llongrightarrow{\relbar\joinrel\relbar\joinrel\rightarrow}
\def\lllongrightarrow{\hbox to 40pt{\rightarrowfill}}

\def\hooklongrightarrow{\lhook\joinrel\longrightarrow}
\def\twoheadrightarrow{\rightarrow\kern -8pt\rightarrow}

\def\maprighto#1{\smash{\mathop{\longrightarrow}\limits^{#1}}}

\def\mapdownr#1{\Big\downarrow\rlap{$\vcenter{\hbox{$\scriptstyle#1$}}$}}
\def\mapdownl#1{\llap{$\vcenter{\hbox{$\scriptstyle#1$}}$}\Big\downarrow}

\def\llongmaprighto#1{\ \smash{\mathop{\llongrightarrow}\limits^{#1}}\ }

\def\lllongmaprighto#1{\ \smash{\mathop{\lllongrightarrow}\limits^{#1}}\ }

\def\longleftmapsto{\longleftarrow\kern-2pt\mapstochar\;}

\def\llongmapsto{\,\vert\kern-3.2pt\joinrel\longrightarrow\,}
\def\llongmapsto{\,\vert\kern-3.7pt\joinrel\llongrightarrow\,}
\def\lllongmapsto{\,\vert\kern-5.5pt\joinrel\lllongrightarrow\,}

\def\isomarrow{\maprighto{\lower3pt\hbox{$\scriptstyle\sim$}}}
\def\llongisomarrow{\llongmaprighto{\lower3pt\hbox{$\scriptstyle\sim$}}}
\def\lllongisomarrow{\lllongmaprighto{\lower3pt\hbox{$\scriptstyle\sim$}}}

\def\lisomarrow{\maplefto{\lower3pt\hbox{$\scriptstyle\sim$}}}

\font\labprffont=cmtt8
\def\strutdepth{\dp\strutbox}
\def\labtekst#1{\vtop to \strutdepth{\baselineskip\strutdepth\vss\llap{{\labprffont #1}}\null}}
\def\marglabel#1{\strut\vadjust{\kern-\strutdepth\labtekst{#1\ }}}

\def\label #1. #2\par{{\definexref{#1}{\num}{#2}}}
\def\labelf #1\par{{\definexref{#1}{\numfo}{formula}}}
\def\slabel #1. #2\par{{\definexref{#1}{\snum}{#2}}}
\def\slabelf #1\par{{\definexref{#1}{\snumfo}{formula}}}

\def\labels #1. #2\par{{\definexref{#1}{\snum}{#2}}}
\def\labelse #1\par{{\definexref{#1}{\num}{section}}}


\def\fibprod{\mathop\times}
\def\tensor{\mathop\otimes}
\def\dirsum{\mathop\oplus}


\def\AA{{\bf A}}
\def\BB{{\bf B}}
\def\CC {{\bf C}}
\def\EE{{\bf E}}
\def\E{{\EE}}
\def\Etilde{{\widetilde{\E}}}
\def\EEtilde{{\widetilde{\EE}}}
\def\AAtilde{{\widetilde{\AA}}}
\def\FF{{\bf F}}

\def\NN {{\bf N}}
\def\QQ {{\bf Q}}

\def\RR {{\bf R}}
\def\WW {{\bf W}}
\def\ZZ {{\bf Z}}
\def\piO{{\pi_0}}
\def\Sh{{\bf Sh}}

\def\vv{{\bf v}}

\def\bbC {{\Bbb C}}

\def\bbL{{\Bbb L}}
\def\bbK{{\Bbb K}}

\def\gerX{{\euscrfm X}}
\def\gerU{{\euscrfm U}}

\def\hatgerX{\widehat{\gerX}}
\def\hatgerU{\widehat{\gerU}}

\def\gerv{{\widehat{v}}}
\def\germ{{\euscrfm{m}}}
\def\gerbeta{{\widehat{\beta}}}

\def\D{{\rm D}}
\def\Dtilde{\widetilde{\rm D}}
\def\Dtildedagger{\widetilde{\rm D}^\dagger}
\def\G{{\rm G}}
\def\H{{\rm H}}
\def\Q {{\bf Q}}

\def\tildeR{{R}}
\def\tR{{R^{\prime}}}
\def\ttildeR{{\tildeR^{\prime}}}
\def\ttildeRinfty{{\tildeR_\infty^{\prime}}}
\def\tRinfty{{R_\infty^{\prime}}}
\def\tS{{S^{\prime}}}

\def\tSinfty{{S_\infty^{\prime}}}
\def\tTinfty{{T_\infty^{\prime}}}
\def\tAA{\AA^{\prime}}
\def\tAAdagger{\AA^{\prime \dagger}}
\def\tAAtilde{\widetilde{\AA}^{\prime}}
\def\tAAtildedagger{\widetilde{\AA}^{\prime\dagger}}

\def\tGamma{\Gamma^{\prime}}
\def\tgamma{{\widetilde{\gamma}}}

\def\R{{\rm R}}

\def\cA{{\cal A}}
\def\cB{{\cal B}}
\def\cF{{\cal F}}
\def\cR{{\cal R}}
\def\cT{{\cal T}}
\def\ctT{\widetilde{{\cal T}}}
\def\cC{{\cal C}}
\def\cG{{\cal G}}
\def\cH{{\cal H}}
\def\cD{{\euscrfm D}}
\def\cE{{\euscrfm E}}
\def\cDtilde{\widetilde{\euscrfm D}}
\def\cDdagger{{\euscrfm D}^\dagger}
\def\cDtildedagger{\widetilde{\euscrfm D}^\dagger}
\def\cU{{\cal U}}
\def\cV{{\cal V}}
\def\cW{{\cal W}}
\def\cZ{{\cal Z}}

\def\cO{{\cal O}}

\def\barcO{{\overline{\cO}}}
\def\barhatcO{{\overline{\cO}}}
\def\cU{{\cal U}}
\def\cX{{\cal X}}
\def\cY{{\cal Y}}
\def\cS{{\cal S}}

\def\cx{{\widehat{x}}}
\def\barx{{\bar{x}}}
\def\barcx{{\bar{\cx}}}
\def\barpi{{\overline{\pi}}}

\def\arit{{\rm ar}}
\def\rig{{\rm rig}}
\def\form{{\rm form}}
\def\fet{{\rm fet}}
\def\geom{{\rm ge}}
\def\t{{\rm t}}

\def\cont{{\rm cont}}
\def\disc{{\rm disc}}
\def\Kbar{{\overline{K}}}
\def\Vbar{{\overline{V}}}
\def\Rbar{{\overline{R}}}
\def\Rbarhat{{\widehat{\overline{R}}}}
\def\kbar{{\overline{k}}}


\def\et{{\rm et}}
\def\sh{{\rm sh}}
\def\Zar{{\rm Zar}}

\def\EXT {{{\cal E}\kern-.5pt{\it xt} }}

\def\Hom {{\rm Hom}}
\def\HOM {{{\cal H}\kern-.5pt{\it om} }}

\def\Ker {{\rm Ker}}
\def\Coker {{\rm Coker}}
\def\Frac{{\rm Frac}}
\def\Im {{\rm Im}}

\def\min {{\rm min}}

\def\Spec {{\rm Spec}}
\def\Spf {{\rm Spf}}
\def\Spm {{\rm Spm}}
\def\Tr{{\rm Tr}}

\def\Rep{{\rm Rep}}
\def\Mod{{\rm Mod}}
\def\Modet{\Mod^{\rm et}}
\def\trho{\rho^{\prime}}
\def\trhodagger{\rho^{\prime\dagger}}
\def\rhotilde{\widetilde{\rho}}
\def\rhotildedagger{\widetilde{\rho}^\dagger}
\def\trhotilde{\widetilde{\rho}^{\prime}}
\def\trhotildedagger{\widetilde{\rho}^{\prime\dagger}}


\def\phi{{\varphi}}
\def\Inf{{\rm Inf}}

\def\rig{{\rm rig}}

\def\Gal{{\rm Gal}}
\def\e{{\euscrfm e}}
\def\Kinfty{{K_\infty}}
\def\Rinfty{{R_\infty}}
\def\Sinfty{{S_\infty}}

\def\Tinfty{{T_\infty}}
\def\Vinfty{{V_\infty}}

\def\tildeRinfty{{\tildeR_\infty}}

\def\AM{{\rm AM}}
\def\An{{\rm An}}
\def\AnBr{{\rm AB}}
\def\BL{{\rm Be}}
\def\BG{{\rm BG}}
\def\BrinonSen{{\rm Br}}
\def\Ch{{\rm Ch}}
\def\CCINV{{\rm CC1}}
\def\CCJAMS{{\rm CC2}}
\def\Elkik{{\rm El}}
\def\FJAMS{{\rm Fa1}}
\def\FALAN{{\rm Fa2}}
\def\FALAST{{\rm Fa3}}
\def\Fo{{\rm Fo}}
\def\Gabber{{\rm Ga}}
\def\He{{\rm He}}
\def\IoSt{{\rm IS}}
\def\Artin{{\rm Ar}}
\def\Jannsen{{\rm Ja}}
\def\ColmezCMP{{\rm Co}}
\def\Tate{{\rm Ta}}
\def\EGAIV{{\rm EGAIV}}
\def\Sen{{\rm Se}}

\senu=1 \ssnu=1 \fonu=1  \sssn=1

\centerline{\titlefont Global applications of relative
$(\phi,\Gamma)$-modules I}
\bigskip
\centerline{ F{\eightpoint abrizio} A{\eightpoint NDREATTA} \&
A{\eightpoint drian} I{\eightpoint OVITA}}

\bigskip
\bigskip

{\insert\footins{\leftskip\footnotemargin\rightskip\footnotemargin\noindent\eightpoint
$2000$ {\it Mathematics Subject Classification} 11G99,14F20,14F30.
}

\vbox{{\leftskip\abstractmargin \rightskip\abstractmargin
\eightpoint

\noindent A{{\sixrm BSTRACT}} In this paper, given a smooth proper scheme $X$ over a 
$p$-adic dvr and a $p$-power torsion \'etale local system ${\bf L}$ 
on it, we study a family of sheaves associated to the cohomology of local,
relative $(\phi,\Gamma)$-modules of ${\bf L}$ and their cohomology.
As applications we derive descriptions of the \'etale cohomology groups 
on the geometric generic fiber of $X$ with values in ${\bf L}$, as well as of their classical
$(\phi,\Gamma)$-modules, in terms of cohomology of the above mentioned sheaves.

\enspace
}}

\section Introduction\par

\noindent Let $p$ be a prime integer, $K$ a finite extension of
$\Q_p$ and~$V$ its ring of integers. In [\Fo], J.-M.~Fontaine
introduced the notion of $(\phi,\Gamma)$--modules designed to
classify $p$-adic representations of the absolute Galois group
$\G_V$ of~$K$ in terms of semi-linear data. More precisely, if $T$
is a $p$-adic representation of $\G_V$, i.e. $T$ is a finitely
generated $\ZZ_p$-module (respectively a $\QQ_p$-vector space of
finite dimension) with a continuous action of $\G_V$, one
associates to it a $(\phi,\Gamma)$--module, denoted $\D_V(T)$.
This is a finitely generated module over a local ring of dimension
two $\AA_V$ (respectively a finitely generated free module over
$\BB_V:=\AA_V\otimes_{\ZZ_p}\QQ_p$) endowed with a semi-linear
Frobenius endomorphism $\phi$ and a commuting, continuous,
semi-linear action of the group
$\Gamma_V:=\Gal(K(\mu_{p^\infty})/K)$ such that $(\D_V(T), \phi)$
is \'etale. This construction makes the group whose
representations we wish to study simpler with the drawback of
making the coefficients more complicated. It could be seen as a
weak arithmetic analogue of the Riemann--Hilbert correspondence
between representations of the fundamental group of a complex
manifold and vector bundles with integrable connections. The main
point of this construction is that one may recover $T$ with its
$\G_V$-action directly from $\D_V(T)$ and, therefore, all the
invariants which can be constructed from $T$ can be described,
more or less explicitly, in terms of $\D_V(T)$. For example

\

($\ast$) one can express in terms of $\D_V(T)$ the Galois
cohomology groups $\H^i(K, T)=\H^i(\G_V, T)$ of $T$.

\

\noindent More precisely, let us choose a topological generator
$\gamma$ of $\Gamma_V$ and consider the complex
$$
\cC^\bullet(T): \D_V(T)\llongmaprighto{d_0}\D_V(T)\oplus\D_V(T)\llongmaprighto{d_1}\D_V(T)
$$
where $d_0(x)=((1-\phi)(x), (1-\gamma)(x))$ and
$d_1(a,b)=(1-\gamma )(a)-(1-\phi)(b)$. It is proven in [\He] that
for each $i\ge 0$ there is a natural, functorial isomorphism
$$
\H^i(\cC^\bullet(T))\cong \H^i(\G_V, T).
$$
Moreover, for $i=1$ this isomorphism was made explicit in
[\CCJAMS]: let $(x,y)$ be a 1-cocycle for the complex
$\cC^\bullet(T)$ and choose $b\in \AA\otimes_{\ZZ_p}T$ such that
$(\phi-1)(b)=x.$ Define the map $C_{(x,y)}:\G_V\lra
\AA\otimes_{\ZZ_p}T$ by
$$
C_{(x,y)}(\sigma)=(\sigma'-1)/(\gamma-1)y-(\sigma-1)b,
$$
where $\sigma'$ is the image of $\sigma$ in $\Gamma_V$. One can
prove that the image of $C_{(x,y)}$ is in fact contained in $T$,
that $C_{(x,y)}$ is a 1-cocycle whose cohomology class
$[C_{(x,y)}]\in\H^1(\G_V, T)$ only depends on the cohomology class
$[(x,y)]\in\H^1(\cC^\bullet(T))$. Moreover, the isomorphism
$\H^1(\cC^\bullet(T))\cong \H^1(\G_V, T)$ above is then defined by
$[(x,y)]\lra [C_{(x,y)}]$.

As a consequence of ($\ast$) we have explicit descriptions of the
exponential map of Perrin-Riou (or more precisely its ``inverse''
(see [\Fo], [\Ch], [\CCJAMS]) and an explicit relationship with
the ``other world'' of Fontaine's modules: $\D_{\rm dR}(T),
\D_{\rm st}(T), \D_{\rm cris}(T)$ (see [\CCINV], [\BL]).

\noindent Despite being a very useful tool, in fact the only one
which allows the general classification of integral and torsion
$p$-adic representations of $\G_V$, the $(\phi,\Gamma)$--modules
have an unpleasant limitation. Namely, $\D_V(T)$ could not so far
be directly related to geometry when $T$ is the
$\G_V$-representation on a $p$-adic \'etale cohomology group (over
$\Kbar$) of some smooth proper algebraic variety defined over $K$.
Here is a relevant passage from the Introduction to~[\Fo]: ``Il
est claire que ces constructions sont des cas particuliers de
constructions beaucoup plus g\'en\'erales. On doit pouvoir
remplacer les corps que l'on consid\`ere ici par des corps des
fonctions de plusieurs variables ou certaines de leurs
completions. En particulier  (i) la loi de r\'eciprocit\'e
explicite \'enonc\'ee au no. 2.4 doit se g\'en\'eraliser et
\'eclairer d'un jour nouveau les traveaux de Kato sur ce sujet;
(ii) ces constructions doivent se faisceautiser et peut \^etre
donner une approche nouvelle des th\'eor\`emes de comparaison
entre les cohomologies $p$-adiques;''

\bigskip
\noindent The first part of the program sketched above, i.e. the
construction of relative  $(\phi,\Gamma)$--modules was
successfully carried out in~[\An] (at least in the good reduction
case). The main purpose of the present article is to continue
Fontaine's program. In particular  various relative analogues,
local and global, of ($\ast$) are proven.

\noindent Let us first point out that in the relative situation,
over a ``small''-$V$-algebra $R$ (see \S 2)  there are several
variants of $(\phi,\Gamma)$--module functors, denoted $\cD_R( - )$
(arithmetic), $\D_R( - )$ (geometric), $\cDtilde_R( - )$
(tilde--arithmetic), $\Dtilde_R( - )$ (tilde--geometric) and their
overconvergent counterparts  $\cDdagger_R( - )$, $\D^\dagger_R( -
)$, $\cDtildedagger_R( - )$ and $\Dtildedagger_R( - )$. For
simplicity of exposition let us explain our results in terms
of~$\cD_R( - )$ and $\cDtilde_R( - )$.

\bigskip
\noindent I) {\it Local results}. This is carried on in \S 3
together with the appendices \S 7 and \S 8. Let $R$ be a ``small''
$V$-algebra and fix $\eta$ a geometric generic point of $\Spf(R)$.
Let $M$ be a finitely generated $\ZZ_p$-module with continuous
action of $\cG_R:=\pi_1^{alg}\bigl(\Spm(R_K, \eta)\bigr)$ and let
$\D:=\cDtilde_R(M)$. Then $\D$ is a finitely generated
$\AAtilde_R$-module endowed with commuting actions of a
semi-linear Frobenius $\varphi$ and a linear action of the group
$\Gamma_R$ (see \S 2.) As in the classical case, $\Gamma_R$ is a
much smaller group that $\cG_R$. It is the semidirect product of
$\Gamma_V$ and of a group isomorphic to $\ZZ_p^d$ where $d$ is the
relative dimension of~$R$ over~$V$.

Let $\cC^\bullet(\Gamma_R, \D)$ be the standard complex of
continuous cochains computing the continuous $\Gamma_R$-cohomology
of $\D$ and denote $\cT_R^\bullet(\D)$ the mapping cone complex of
the morphism $(\varphi-1)\colon \cC^\bullet(\Gamma_R, \D)\lra
\cC^\bullet(\Gamma_R, \D)$. Then, Theorem~\ref{maintheorem} states
that we have natural isomorphisms, functorial in $R$ and $M$,
$$
\H^i_{cont}(\cG_R, M)\cong
\H^i\bigl(\cT_R^\bullet(\D)\bigr)\quad\hbox{{\rm for all }}i\ge 0.
$$The maps  are defined in \S 3 in an explicit way, following
Colmez's description in the classical case. The input of
Fontaine's construction of the classical $(\phi,\Gamma)$--modules
was to replace modules over perfect, non--noetherain rings with
modules over smaller rings: ``C'est d'ailleurs [...] que j'ai
compris l'int\'er\^et qu'il avait \`a ne pas remplacer
$k(\!(\pi)\!)$ par sa cl\^oture radicielle" Indeed, ``[...]ceci
permet d'introduire des techniques diff\'erentielles". Motivated
by the same needs, in view of applications to comparison
isomorphisms, we show in appendix \S 7 that one can replace the
module $\cDtilde_R(M)$ over the ring~$\AAtilde_R$, which is not
noetherian, with the smaller $(\phi,\Gamma_R)$--module
$\cD_R(M)\subset \cDtilde_R(M)$ over the noetherian, regular
domain~$\AA_R$ of dimension $d+1$. We show that the natural map $$
\H^i_{cont}(\Gamma_R, \cD_R(M)) \llongrightarrow
\H^i_{cont}(\Gamma_R, \cDtilde_R(M)) \quad\hbox{{\rm for all
}}i\ge 0$$is an isomorphism.  The proof follows and slightly
generalizes Tate--Sen's method in [\BrinonSen]. In particular, one
has a natural isomorphism
$$
\H^i_{cont}(\cG_R, M)\cong
\H^i\bigl(\cT_R^\bullet(\cD_R(M))\bigr)\quad\hbox{{\rm for all
}}i\ge 0,
$$
where $ \cT_R^\bullet\bigl(\cD_R(M)\bigr)$ is the mapping cone
complex of the map $(\varphi-1)\colon \cC^\bullet(\Gamma_R,
\cD_R(M))\to \cC^\bullet(\Gamma_R, \cD_R(M))$.

\bigskip
\noindent II) {\it Global results}. This is carried on in \S4,
\S5, \S6. The setting for \S4 and \S5 is the following. Let $X$ be
a smooth, proper, geometrically irreducible scheme of finite type
over $V$  and let $\bbL$ denote a locally constant \'etale sheaf
of $\ZZ/p^s\ZZ$-modules  (for some $s\ge 1$) on the generic fiber
$X_K$ of $X$. Let $\cX$ denote the formal completion of $X$ along
its special fiber and let $X_K^{\rig}$ be the rigid analytic space
attached to $X_K$. Fix a geometric generic point
$\eta=\Spf(\bbC_{\cX})$ and set ${\bf L}$ the fiber of $\bbL$ at
$\eta$.

\noindent To each $\cU\to \cX$ \'etale such that $\cU$ is affine,
$\cU=\Spf(R_{\cU})$, with $R_\cU$ a small $V$-algebra and a choice
of local parameters $(T_1,T_2,\ldots,T_d)$ of $R_\cU$ (as in \S2)
we attach the relative $(\phi,\Gamma)$--module $\cDtilde_\cU({\bf
L}):=\cDtilde_{R_\cU}({\bf L})$. However, the association $\cU\lra
\cDtilde_\cU({\bf L})$ is not functorial because of the dependence
of $\cDtilde_\cU({\bf L})$ on the choice of the local parameters.
In other words the relative $(\phi,\Gamma)$--module construction
does not sheafify.

Nevertheless due to I) above, for every $i\ge 0$, the association
$\cU\lra \H^i(\tilde{\cT}_{R_\cU}^\bullet({\bf L}))$ is functorial
and we denote by $\cH^i({\bf L})$ the sheaf on $\cX^{et}$
associated to it. In \S4 we prove Theorem 4.1: there is a spectral
sequence
$$
E_2^{p,q}=\H^q(\cX^{et}, \cH^p({\bf L})) \Longrightarrow \H^{p+q}((X_K)^{et}, \bbL).
$$
We view this result as a global analogue of ($\ast$): the \'etale
cohomology of $\bbL$ is calculated in terms of local relative
$(\phi,\Gamma)$--modules attached to ${\bf L}$.

\noindent The proof of Theorem~\ref{etalePhiGamma} follows a
roundabout path which was forced on us by lack of enough knowledge
on \'etale cohomology of rigid analytic spaces. More precisely,
for an algebraic, possibly infinite, extension $M$ of $K$
contained in $\Kbar$, Faltings defines in [\FALAST] a Grothendieck
topology $\gerX_M$ on $X$ (see also \S4). The local system $\bbL$
may be thought of as a sheaf on $\gerX_M$ and it follows
from~[\FALAST], see~\ref{compareGaltoXetI}, that there is a
natural isomorphism:
$$
(\ast\ast)\quad \H^i(\gerX_M, \bbL)\cong \H^i((X_M)^{et}, \bbL),
$$
for all $i\ge 0$. The main tool for proving $(\ast\ast)$ is the
result: every point $x\in X_K$ has a neighborhood $W$ which is
$K(\pi,1)$. Such a result, although believed to be true, is yet
unproved in the rigid analytic setting. Therefore the proof of
Theorem~\ref{etalePhiGamma} goes as follows. We define the
analogue Grothendieck topology $\hatgerX_M$ on $\cX$, prove that
there is a spectral sequence with $E_2^{p,q}=\H^q(\cX^{et},
\cH^p({\bf L}^{rig}))$ abutting to $\H^{p+q}(\hatgerX_M,
\bbL^{rig})$, then compare $\H^i(\gerX_M, \bbL)$ to
$\H^i(\hatgerX_M, \bbL^{\rig})$ and $\H^i((X_K)^{et}, \bbL)$ to
$\H^i((X_K^{\rig})^{et}, \bbL^{\rig})$ and in the end use
Faltings's result $(\ast\ast)$.

\bigskip
\noindent In \S5 we introduce a certain family of continuous
sheaves which we call  {\it Fontaine sheaves} and which we denote
by $\barcO_{\_}$, $\cR\bigl(\barcO_{\_}\bigr)$, $A_{\rm
inf}^+(\barcO_{\_})$. There are algebraic and analytic variants of
these: the first are sheaves on $\gerX_M$ and the second on
$\hatgerX_M$. We would like to remark that the local sections of
the Fontaine sheaves are very complicated and they are {\bf not}
relative Fontaine rings. Continuous cohomology of continuous
sheaves  on $\gerX_M$ and $\hatgerX_M$ respectively is developed
in \S 5. As an application a geometric interpretation of
$\cDtilde_V(T_i)$, where $T_i=\H^i((X_{\Kbar})^{et}, \bbL)$, for
$\bbL$ an \'etale local system of $\ZZ/p^s\ZZ$-modules on $X_K$ as
above is given. More precisely, it is proven in \S 5 that there is
a natural isomorphism of classical $(\phi,\Gamma)$--modules:
$$
\H^i\left(\hatgerX_\Kinfty, \bbL^{\rig}\otimes A_{\rm
inf}^+\bigl(\barhatcO_{\hatgerX_\Kinfty}\bigr)\right)\cong
\Dtilde_V\left(\H^i\bigl((X_{\Kbar})^{et}, \bbL\bigr)\right).
$$

\bigskip
\noindent Finally, in \S6 we relax our global assumptions. Now
$\cX$ denotes a formal scheme topologically of finite type over
$V$, smooth and geometrically irreducible, not necessarily
algebrizable, and $X_K^{\rig}$ denotes its generic fiber.

\noindent In \S6 we set up the basic theory for comparison
isomorphisms between the different $p$-adic cohomology theories in
this analytic setting. Our main result is that, if $\bbL^{\rig}$
is a $p$-power torsion local system on $X_K^{\rig}$ and
$\cA^{Font}$ is one of the analytic Fontaine sheaves on
$\hatgerX_M$ listed above, then the cohomology groups
$\displaystyle \H^i(\hatgerX_M, \bbL^{\rig}\otimes \cA^{Font})$
can be calculated as follows. Let us first recall that we fixed a
geometric generic point $\eta=\Spf(\bbC_{\cX})$, where
$\bbC_{\cX}$ is a complete, algebraically closed field which can
be chosen as in \S 5. For each \'etale morphism $\cU\lra \cX$ such
that $\cU$ is affine, $\cU=\Spf(R_\cU)$ with $R_\cU$ a small
$V$-algebra, let $\Rbar_{\cU}$ denote the union of all normal
$R_\cU$-algebras contained in  $\bbC_{\cX}$ which are finite and
\'etale over $R_\cU$ after inverting $p$. Let
$\cA^{Font}(\Rbar_\cU\otimes K)$ denote the Fontaine ring
constructed starting with the pair $(R_\cU,\Rbar_\cU)$ as in~[\Fo]
and denote by ${\bf L}$, as before, the fiber of $\bbL^{\rig}$ at
$\eta$. One can show that the association $\cU\lra
\H^i(\pi_1^{alg}(\cU_M, \eta), {\bf L}\otimes
\cA^{Font}(\Rbar_\cU\otimes K))$ is functorial and denote
$\cH^i_M(\bbL^{\rig}\otimes\cA^{Font})$ the sheaf on $\cX^{et}$
associated to it.

Let us notice that, due to the generalized Tate-Sen's method of \S
7, if $\cA^{\rm Font}=A_{\rm
inf}^+\bigl(\barhatcO_{\hatgerX_M}\bigr)$, the inflation defines
an isomorphism:
$$
\H^i\bigl(\Gamma_{R_\cU}, \cD_{R_\cU}({\bf L})\bigr)\cong
\H^i\bigl(\Gamma_{R_\cU}, \cDtilde_{R_\cU}({\bf L})\bigr)\cong
\H^i\bigl(\pi_1^{alg}(\cU,\eta), {\bf L}\otimes A_{\rm
inf}^+(\Rbar_{\cU}\otimes K)\bigr).
$$
Hence, the sheaf $\cH^i_M\left(\bbL^{rig}\otimes A_{\rm
inf}^+\bigl(\barhatcO_{\hatgerX_M}\bigr)\right)$ is  defined
locally in terms of $\Gamma$-cohomology of relative
$(\phi,\Gamma)$--modules.

\noindent It is proved (Theorem~\ref{cohomologyFontaine})  that
there exists a spectral sequence
$$
(\ast\ast\ast)\quad E_2^{p,q}=\H^q(\cX^{et}, \cH^p_M(\bbL^{\rig}\otimes\cA^{Font}))
\Longrightarrow \H^{p+q}(\hatgerX_M, \bbL^{\rig}\otimes \cA^{Font}).
$$

\bigskip
\noindent At this point we would like to remark that our results
in \S 6 are distinct from those of Faltings in [\FJAMS], [\FALAN],
[\FALAST]. Namely let us consider the following diagram of
categories and functors:
$$
\matrix{\Sh(\hatgerX_M)^{\NN} &
\lllongmaprighto{\gerv_{\cX,M,\ast}^\NN} & \Sh(\cX^{et})^{\NN} \cr
\alpha\downarrow&&\beta\downarrow \cr \Sh(\cX^{et}) &
\lllongmaprighto{\H^0(\cX^{et}, - )}&{\bf AbGr}\cr}
$$
where if $\cC$ is a Grothedieck topology then we denote by
$\Sh(\cC)$ the category of sheaves of abelian groups on $\cC$ and
by $\Sh(\cC)^\NN$ the category of continuous sheaves on $\cC$ (see
\S 5). We also denote $\displaystyle
\alpha=\lim_{\leftarrow}\gerv_{\cX,M,\ast}$ and $\displaystyle
\beta=\lim_{\leftarrow}\H^0(\cX^{et}, -)$.

\noindent We analyze the spectral sequence attached to the
composition of functors: $\displaystyle \H^0(\cX^{et}, - )\circ
\lim_{\leftarrow}(\gerv_{\cX,M,\ast})^\NN$ while it appears,
although very little detail is given, that Faltings considers the
composition of the other two functors in the above diagram (in the
algebraic setting). We believe that our point of view is
appropriate for the applications to relative
$(\phi,\Gamma)$--modules  that we have in mind.

\noindent The analysis in \S 6 and the spectral sequence
$(\ast\ast\ast)$ have already been used in order to construct a
$p$-adic, overconvergent, finite slope Eicher-Shimura isomorphism
and to give a new, cohomological construction of $p$-adic families
of finite slope modular forms in [\IoSt].

\noindent In a sequel paper (``Global arithmetic applications of
relative $(\phi,\Gamma)$-modules, II'') we plan to first extend
the constructions and results in \S 6 of the present paper to
formal schemes over $V$ with semi-stable special fiber and use
them in order to prove comparison isomorphisms between the
different $p$-adic cohomology theories involving Fontaine sheaves
in such analytic settings. We believe that we would be able to
carry on this project for spaces like: the $p$-adic symmetric
domains, their \'etale covers (in the cases where good formal
models exist), the $p$-adic period domains of Rapoport-Zink, etc.

\

{\bf Acknowledgements:} We thank A.~Abbes, V.~Berkovich and W.~Niziol for interesting discussions 
pertaining to the subject of this paper. 
Part of the work on this article was done when the first autor visited the Department of Mathematics 
and Statistics of Concordia University and the second author visited the IHES and il Dipartimento 
di matematica pura ed applicata of the University of Padova. Both we would like to express our 
gratitude  to these institutions for their hospitality.

\

\noindent{\sectionfont Contents.} \spacing \item{{\rm 1.}}
Introduction. \item{{\rm 2.}} Preliminaries. \item{{\rm 3.}}
Galois cohomology and $(\phi,\Gamma)$--modules. \item{{\rm 4.}}
\'Etale cohomology and relative $(\phi,\Gamma)$--modules.
\item{{\rm 5.}} A geometric interpretation of classical
$(\varphi,\Gamma)$--modules. \item{{\rm 6.}} The cohomology of
Fontaine sheaves. \item{{\rm 7.}}  Appendix I: Galois cohomology
via Tate--Sen's method. \item{{\rm 8.}} Appendix II:
Artin--Schreier theory. \item{{ }} References.

\endsection

\

\

\

\centerline{\titlefont I. Local Theory.}

\

\section Preliminaries\par

\label V. section\par \ssection The basic rings\par  Let $V$ be a
complete discrete valuation ring, with perfect residue field~$k$
of characteristic~$p$ and with fraction field~$K=\Frac(V)$ of
characteristic~$0$. Let~$\vv$ be the valuation on~$V$ normalized
so that~$\vv(p)=1$. Let $K \subset \Kbar$ be an algebraic closure
of~$K$ with Galois group~$\Gal(\Kbar/K)=:\G_V$ and denote
by~$\Vbar$ the normalization of\/~$V$ in~$\Kbar$. Define the tower
$$K_0:=K \subset K_1= K(\zeta_p) \subset \cdots \subset
K_n=K(\zeta_{p^n}) \subset \cdots
$$where $\zeta_{p^n}$ is a primitive $p^n$--th root of unity and
$\zeta_{p^{n+1}}^p=\zeta_{p^n}$ for every~$n\in\NN$. Let~$V_n$ be
the normalization of~$V$ in~$K_n$ and define~$\Kinfty:=\cup_n
K_n$. Write $\Gamma_V:=\Gal(\Kinfty/K)$ and
$\H_V:=\Gal(\Kbar/\Kinfty)$ so that~$\Gamma_V=\G_V/\H_V$.

\bigskip

\noindent Let $\tildeR$ be a $V$--algebra such that $k \subset
R\tensor_V k$ is geometrically integral. Assume that~$\tildeR$ is
obtained from~$R^0= V\left\{T_1^{\pm 1},\ldots,T_d^{\pm
1}\right\}$ iterating finitely many times the following
operations:\spacing

\item{{\rm \'et)}} the $p$--adic completion of an \'etale
extension;\spacing

\item{{\rm loc)}} the $p$--adic completion of the localization
with respect to a multiplicative system;\spacing

\item{{\rm comp)}} the completion with respect to an ideal
containing~$p$.\spacing

\noindent Define
$$R_n:= R\tensor_V V_n\left[T_1^{1\over p^n},T_1^{-1\over
p^n},\ldots,T_d^{1\over p^n},T_d^{-1\over p^n}\right],\qquad
\Rinfty:=\cup_n R_n.$$

Let\/~$ \Rbar$ be the direct limit of a maximal chain of normal
$R$-algebras, which are domains and, after inverting~$p$, are
finite and \'etale extensions of\/~$\Rinfty\left[{1\over
p}\right]$.

\spacing Let\/~$m\in\NN$ and let\/~$S$ be a $R_m$--algebra such
that\/~$S$ is finite as $R_m$--module and\/~$R_m\subset S$ is
\'etale after inverting~$p$. Define~$S_n$ as the normalization
of\/~$S\tensor_{R_m} R_n$ in~$S\tensor_{R_m}
R_n\left[p^{-1}\right]$ for every~$n\geq m$.
Let\/~$\Sinfty:=\cup_{n\geq m} S_n$.

Write\/~$\tS_n$ for the normalization of\/~$S_n\tensor_{V_n}
\Vbar$ in~$S_n\tensor_{V_n} \Vbar\bigl[p^{-1}\bigr]$
and\/~$\tSinfty$ for the normalization
of\/~$\Sinfty\tensor_\Vinfty \Vbar$ in~$\Sinfty\tensor_\Vinfty
\Vbar\left[p^{-1}\right]$. We put\/~$\tS:=\tS_m$.

Note that~$\tR=R\tensor_V\Vbar$
and~$\tRinfty=\Rinfty\tensor_{\Vinfty}\Vbar$
\endssection

\label SnSn+1. proposition\par\prop There exist constants~$0 <
\varepsilon < 1$ and~$N=N(S)\in\NN$, depending on~$S$, and there
exists an element~$p^\varepsilon$ of\/~$V_N$ of
valuation~$\varepsilon$ such that\/ $S_{n+1}^p+p^\varepsilon
S_{n+1}\subset S_n+p^\varepsilon S_{n+1}$ (as subrings
of\/~$S_{n+1}$) and\/~$\tS_{n+1}^p+p^\varepsilon \tS_{n+1}\subset
\tS_n+p^\varepsilon \tS_{n+1}$ (as subrings of\/~$\tS_{n+1}$) for
every~$n\geq N$.
\endprop
\Proof The claim concerning~$S_{n+1}$  follows from [\An,
Cor.~3.7].  It follows from~[\An, Prop.~3.6] that there exists a
decreasing sequence of rational numbers~$\{\delta_n(S)\}$ such
that $p^{\delta_n(S)}$ annihilates the trace map $\Tr\colon
\tS_n\to \Hom_{\ttildeR_n}\bigl(\tS_n,\ttildeR_n\bigr)$. This
implies that $p^{\delta_n(S)} \tS_{n+1} \subset
\tS_n\tensor_{\ttildeR_n} \ttildeR_{n+1}$; see loc.~cit. This, and
the fact that the proposition holds for~$\ttildeR$ by direct
check, allows to conclude; see the proof of~[\An, Cor.~3.7] for
details.

\label H. definition\par\defi For every $R$--subalgebra $S\subset
\Rbar$ as in~\ref{V} such that~$\tSinfty$ is an integral domain,
viewed as a subring of\/~$\Rbar$, define
$$\cG_S:=\Gal\left(\Rbar\left[{1\over p}\right]/S\left[{1\over
p}\right]\right), \quad \Gamma_S:=\Gal\left(\Sinfty\left[{1\over
p}\right]/S\left[{1\over p}\right]\right)$$and $$
\cH_S:=\Ker\left(\cG_S\rightarrow\Gamma_S \right).$$Analogously,
let $$\G_S:=\Gal\left(\Rbar\left[{1\over
p}\right]/\tS\left[{1\over p}\right]\right), \quad
\tGamma_S:=\Gal\left(\tSinfty\left[{1\over
p}\right]/\tR\left[{1\over p}\right]\right)$$and $$
\H_S:=\Ker\left(\G_S\rightarrow\tGamma_S
\right)=\Gal\left(\tS_\infty\left[{1\over
p}\right]/\tR\left[{1\over p}\right]\right).$$Note that by
assumption\/~$\cH_S/H_S\cong \H_V$ and that\/~$\Gamma_S$ is
isomorphic to the semidirect product of\/~$\Gamma_V$ and
of\/~$\tGamma_S$. The latter is a finite index subgroup
of\/~$\tGamma_R\cong \ZZ_p^d$. We let $\tgamma_1,\ldots,\tgamma_d$
be topological generators of\/~$\tGamma_R$.

\enddefi

\label en. section\par\ssection {\bf RAE}\par Following
Faltings~[\FJAMS, Def.~2.1] we say that an extension $\Rinfty
\subset \Sinfty$ is almost \'etale if it is  finite and \'etale
after inverting~$p$ and if, for every~$n\in\NN$, the element
$p^{1\over p^n}\e_\infty$ is in the image of $
\Sinfty\tensor_{\Rinfty} \Sinfty$. Here $\e_\infty$ is the
canonical idempotent of the \'etale extension
$\Rinfty\bigl[p^{-1}\bigr] \subset \Sinfty\bigl[p^{-1}\bigr]$.

We say that such extension satisfies (RAE), {\it for refined
almost \'etaleness}, if the following holds. For every~$n\geq m$
let\/~$\e_n$ be the diagonal idempotent associated to the \'etale
extension $R_n\bigl[p^{-1}\bigr] \subset S_n\bigl[p^{-1}\bigr]$.
There exists $\ell\in\NN$, independent of\/~$m$, such that there
exists an element\/~$p^{\ell\over p^n}$ of\/~$V_n$ of
valuation~${\ell\over p^n}$ and $p^{\ell\over p^n} \e_n$ lies in
the image of\/~$S_n \tensor_{R_n} S_n$.

We assume that (RAE) holds for every extension $\Rinfty\subset
\Sinfty$ arising as in~\ref{V}.
\endssection

\label whenRAEholds. remark\par\rmk It is proven in~[\An,
Prop.~5.10 \& Thm.~5.11] that (RAE) holds if $\tildeR$ is of Krull
dimension~$\leq 2$ or if the composite of the
extensions~$V\bigl[T_1^{\pm 1},\cdots, T_d^{\pm 1}\bigr] \to R^0
\to \tildeR$ is flat and has geometrically regular fibers. For
example, this holds if $\tildeR$ is obtained by taking the
completion with respect to an ideal containing~$p$ of the
localization of an \'etale extension of~$V\bigl[T_1^{\pm
1},\cdots, T_d^{\pm 1}\bigr]$; see~[\An, Prop.~5.12].

\endrmk

\label EE. section\par \ssection The rings $\EEtilde_\Sinfty$,
$\EE_S$, $\EEtilde_\tSinfty$ and $\EE_\tS$\par Let\/~$S$ be as
in~\ref{V}. Define
$$\EEtilde_\Sinfty^+:=\lim\bigl(\Sinfty/p^\varepsilon
\Sinfty\bigr),\qquad
\EEtilde_\tSinfty^+:=\lim\bigl(\tSinfty/p^\varepsilon
\tSinfty\bigr)$$where the inverse limit is taken with respect to
Frobenius. Using~\ref{SnSn+1} define the {\it generalized ring of
norms}, $$\EE_S^+\subset \EEtilde_\Sinfty^+,\qquad
\EE_\tS^+\subset \EEtilde_\tSinfty^+$$as the subring consisting of
elements $(a_0,\ldots, a_n,\ldots)$ in~$\EEtilde_\Sinfty^+$
(resp.~in~$\EEtilde_\tSinfty^+$) such that $a_n$ is
in~$S_n/p^\varepsilon S_n$ (resp.~$\tS_n/p^\varepsilon \tS_n$) for
every~$n\geq N(S)$.

By construction~$\Etilde_\Sinfty^+$, $\EE_S^+$,
$\Etilde_\tSinfty^+$ and~$\EE_\tS$ are endowed with a Frobenius
homomorphism~$\varphi$ and a continuous action of~$\Gamma_R$.
Put~$\Etilde_\tSinfty:=\Etilde_\tSinfty^+\bigl[\barpi^{-1}\bigr]$,
$\Etilde_\Sinfty:=\Etilde_\Sinfty^+\bigl[\barpi^{-1}\bigr]$,
$\EE_\tS:=\EE_\tS^+\bigl[\barpi^{-1}\bigr]$ and\/
$\EE_S:=\EE_S^+\bigl[\barpi^{-1}\bigr]$.

Denote by~$\epsilon$ the element
$(1,\zeta_p,\ldots,\zeta_{p^n},\ldots)\in\EE_V^+$ and
by~$\barpi:=\epsilon-1$. By abuse of notation for~$\alpha\in\QQ$,
we denote by~$\piO^\alpha$ a (any) element
$a=(a_0,a_1,\ldots,a_n,\ldots)$ in~$\cup_m \E^+_V(m)$, if it
exists, such that~$\vv(a_i)={\alpha\over p^i}$ for~$i\gg 0$. For
example, $\barpi=\pi_0^{p\over p-1}$; see~[\AnBr, Prop.~4.3(d)].
For every~$i=1,\ldots,d$, let~$x_i:=(T_i,T_i^{1\over
p},T_i^{1\over p^2},\cdots)\in \EE_{R^0}^+$.  The following hold,

\spacing \item{{\it 1)}} the map $\EE_S^+/\pi_0^{p^n\varepsilon}
\EE_S^+\rightarrow S_n/p^\varepsilon S_n$ and the map
$\EEtilde_\Sinfty^+/\pi_0^{p^n\varepsilon}
\EEtilde_\Sinfty^+\rightarrow \Sinfty/p^\varepsilon \Sinfty$,
sending $(a_0,\ldots,a_n,\ldots)\mapsto a_n$, are isomorphisms
for~$n\geq N(S)$; see~[\An, Thm.~5.1];

\spacing \item{{\it 2)}} $\EE_S^+$ is a normal ring, it is finite
as $\EE_{\tildeR}^+$--module and it is an \'etale extension of
$\EE_{\tildeR}^+$, after inverting~$\barpi$, of degree equal to
the generic degree of $\tildeR_m\subset S$; see~[\An, Thm.~4.9 \&
Thm.~5.3];

\spacing \item{{\it 3)}} $\EEtilde_\Sinfty^+$ is normal and
coincides with the $\barpi$--adic completion of the perfect
closure of $\E_S^+$; see~[\An, Cor.~5.4];

\spacing \item{{\it 4)}} we have maps $\barpi^\ell
\EEtilde_{\Sinfty}^+\rightarrow \EE_S^+\tensor_{\EE_R^+}
\EEtilde_\Rinfty^+  \rightarrow \EEtilde_{\Sinfty}^+$. They are
isomorphisms after inverting $\barpi$. In particular,
$\EEtilde_\Sinfty^+\bigl[\barpi^{-1}\bigr]=\EEtilde_\Rinfty^+\tensor_{\EE_R^+}
\EE_S^+\bigl[\barpi^{-1}\bigr]$; see~[\An, Lem.~4.15];

\spacing \item{{\it 5)}} consider the ring
$$\lim_{\infty\leftarrow n} \widehat{\Sinfty}:=\left\{\bigl(
x^{(0)},x^{(1)},\ldots,x^{(m)}, \ldots\bigr)\vert
x^{(m)}\in\widehat{\Sinfty},\enspace
\bigl(x^{(m+1)}\bigr)^p=x^{(m)} \right\},$$where
$\widehat{\Sinfty}$ is the $p$--adic completion of $\Sinfty$, the
transition maps are defined by raising to the $p$--th power, the
multiplicative structure is induced by the one
on~$\widehat{\Sinfty}$ and the additive structure is defined by
$$(\ldots,x^{(m)},\ldots)+(\ldots,y^{(m)},\ldots)=\bigl(\ldots,
\lim_{n\rightarrow\infty}
(x^{(m+n)}+y^{(m+n)})^{p^n},\ldots\bigr).$$The natural map
$\lim_{\infty\leftarrow n} \widehat{\Sinfty}\to
\EEtilde_\Sinfty^+$ is a bijection; see~[\An, Lem.~4.10].

\spacing

\noindent It follows from~(1), see~[\An, Cor.~4.7], that
$$\EE_V^+\cong k_\infty[\![\pi_K]\!]\qquad\hbox{{\rm and}}\qquad \EE_{R^0}^+\cong
\EE_V^+\left\{x_1,\ldots,x_d,{1\over x_1},\ldots,{1\over
x_d}\right\},$$where~$k_\infty$ is the residue field of~$V_\infty$
and~$\pi_K=(\ldots,\tau_n,\tau_{n+1},\ldots)$, with~$\tau_i\in
V_i$ for~$i\gg 0$, is a system of uniformizers
satisfying~$\tau_{i+1}^p\equiv \tau_i$ mod~$p^\varepsilon$. The
convergence in $x_1^{\pm 1},\ldots,x_d^{\pm 1}$ is relative to the
$\barpi$-adic topology on~$\EE_V^+$. Eventually $\EE_{\tildeR}^+$
is obtained from~$\EE_{R^0}^+$ iterating the operations

\item{{\rm \'et)}} the $\barpi$--adic completion of an \'etale
extension;\spacing

\item{{\rm loc)}} the $\barpi$--adic completion of the
localization with respect to a multiplicative system;\spacing

\item{{\rm comp)}} the completion with respect to an ideal
containing a power of~$\barpi$.\spacing

\noindent In particular, $\{\pi_K,x_1,\ldots,x_d\}$ is an absolute
$p$--basis of~$\EE_R^+$.
\endssection

\label tEE. lemma\par\lemma Let\/~$S$ be as in~\ref{V}. The
following hold:

\spacing \item{{\rm 1.}} the maps
$\EEtilde_\tSinfty^+/\pi_0^{p^n\varepsilon}
\EEtilde_\tSinfty^+\rightarrow \tSinfty/p^\varepsilon \tSinfty$
and\/ $\EE_\tS^+/ \pi_0^{p^n\varepsilon} \EEtilde_\tS^+\rightarrow
\tS_n/p^\varepsilon \tS_n $, given by
$(a_0,\ldots,a_n,\ldots)\mapsto a_n$, is an isomorphism for~$n\geq
N(S)$. In particular, $\EE_\ttildeRinfty^+$ coincides with the
$\barpi$--adic completion of\/
$\EEtilde_\tildeRinfty^+\tensor_{\EEtilde_\Vinfty^+}\EEtilde_\Vbar^+$
and $\EE_\ttildeR^+$ coincides with the $\barpi$--adic completion
of\/ $\EE_\tildeR^+\tensor_{\EE_V^+}\EEtilde_\Vbar^+$;

\spacing \item{{\rm 2.}} the extensions $\EE_\tildeR^+\to
\EEtilde_{\tildeRinfty}^+$, $\EE_\tildeR^+\to
\EEtilde_{\ttildeRinfty}^+$ and\/ $\EE_\tildeR^+\to \EE_\tR^+$ are
faithfully flat. For every finitely generated
$\EE_\tildeR^+$--module $M$, the base change via the above
extensions are $\barpi$--adically complete and separated;

\spacing \item{{\rm 3.}} we have maps  $\barpi^\ell
\EEtilde_{\tSinfty}^+\rightarrow \EE_S^+\tensor_{\EE_R^+}
\EEtilde_\tRinfty^+  \rightarrow \EEtilde_{\tSinfty}^+$ and\/
$\barpi^\ell \EE_{\tS}^+\rightarrow
\EE_S^+\tensor_{\EE_R^+}\E_\tR^+  \rightarrow \EE_{\tS}^+$. They
are isomorphisms after inverting~$\barpi$;

\endlemma
\Proof Statements (1) and (3) follow from~\ref{en} arguing as in
the proofs of [\An, Thm.~5.1] and [\An, Lem.~4.15]
respectively.\spacing

(2) By~\ref{EE} we have $\EE_R^+/\pi_0^{\varepsilon p^n} \EE_R^+
\cong R_n/p^\varepsilon R_n$, similarly
$\EE_\Vbar^+/\pi_0^{\varepsilon p^n}\EE_\Vbar^+\cong
\Vbar/p^\varepsilon\Vbar$, $\EE_\tR^+/ \pi_0^{\varepsilon p^n}
\EE_\tR^+\cong \bigl(R_n/p^\varepsilon R_n\bigr)\tensor_{V_n}\Vbar
$ and  $\EEtilde_{\tRinfty}^+/\pi_0^{\varepsilon
p^n}\EEtilde_{\tRinfty}^+\cong \bigl(\Rinfty/p^\varepsilon\Rinfty
\bigr)\tensor_{\Vinfty} \Vbar$. By construction
$\Rinfty\tensor_\Vinfty\Vbar$ is a free
$R_n\tensor_{V_n}\Vbar$--module with basis $\{T_1^{\alpha_1\over
p^m}\cdots T_d^{\alpha_d\over p^m}\} $ with~$m\geq n$ and~$0\leq
\alpha_i< p^{m-n}$. Hence, $\E_\tR^+$
(resp.~$\EEtilde_{\tRinfty}^+$) is the $\barpi$--adic completion
$\EE_R^+\tensor_{\EE_V^+} \E_\Vbar^+$ (resp.~of $\cup_n
\bigl(\EE_R^+\tensor_{\EE_V^+} \E_\Vbar^+\bigr)[x_1^{1\over
p^n},\ldots, x_d^{1\over p^n}]$). Also, $\E_\Vbar^+$ is the
$\pi_0$--adic completion of finite, normal and generically
separable extensions of the dvr $\E_V^+$. Those are free as
$\E_V^+$--module. We may then apply~[\An, Lem.~8.7] to conclude.

\spacing Given an $\Rinfty $--algebra $\Sinfty$, finite and
\'etale over~$\Rinfty\left[{1\over p}\right]$, there
exists~$m\in\NN$ and there exists a $R_m$--algebra $S$, finite and
\'etale over~$R_m\left[{1\over p}\right]$ such that~$\Sinfty$,
defined as in~\ref{V}, is the normalization of $S\tensor_{R_n}
\Rinfty$.

\label equivalenceGalois. theorem\par\thm The functor $\Sinfty \to
\E_S^+$ defines an equivalence of categories from the category
{{\bf $\Rinfty$-AED}} of $\Rinfty$--algebras which are normal
domains, finite and \'etale over~$\Rinfty\left[{1\over p}\right]$
to the category {{\bf $\EE^+_R$-AED}} of $\EE_R^+$--algebras,
which are normal domains, finite and \'etale over~$\EE_R$.\endthm
\Proof See~[\An, Thm.~6.3].

\bigskip

Let $\EE^+_\Rbar$ be $\cup_\Sinfty \EE_S^+$ where the union is
taken over all $\Rinfty$--subalgebras $\Sinfty\subset \Rbar$ such
that $\Sinfty\left[p^{-1}\right]$ is finite \'etale over
$\Rinfty\left[p^{-1}\right]$. Let $\EEtilde^+_\Rbar=\lim
\bigl(\Rbar/p^\varepsilon \Rbar\bigr)$, where the inverse limit is
taken with respect to Frobenius. It coincides with the
$\barpi$--adic completion of~$\cup_\Sinfty \EEtilde_\Sinfty^+$.
Define~$\EE_\Rbar:=\EE_\Rbar^+\left[\barpi^{-1}\right]$
and~$\EEtilde_\Rbar:=\EEtilde_\Rbar^+\left[\barpi^{-1}\right]$.

\label invariantsEE. proposition\par\prop Let\/~$S$ be as
in~\ref{H}. Then, $\widehat{\Rbar}^{\cH_S}= \widehat{\Sinfty}$
and\/~$\widehat{\Rbar}^{\H_S}= \widehat{\tSinfty}$. Furthermore,
$$\bigl(\EE^+_\Rbar\bigr)^{\cH_S}=\EE_S^+,\qquad
\EE_\Rbar^{\cH_S}=\EE_S,\qquad
\bigl(\EEtilde^+_\Rbar\bigr)^{\cH_S}=\EEtilde_\Sinfty^+,\qquad
\EEtilde_\Rbar^{\cH_S}=\EEtilde_\Sinfty$$and
$$\bigl(\EE_\Rbar\bigr)^{\H_S}=\EE_\tS,\qquad
\bigl(\EEtilde_\Rbar^+\bigr)^{\H_S}=\EEtilde_{\tSinfty}^+,\qquad
\bigl(\EEtilde_\Rbar\bigr)^{\H_S}=\EEtilde_{\tSinfty}.$$
Furthermore,   the natural maps
$$\EEtilde_\Rinfty\tensor_{\EE_R}\EE_S\llongrightarrow
\EEtilde_\Sinfty,\quad \EE_\Rbar^{\H_R}\tensor_{\EE_R} \EE_S
\llongrightarrow \EE_\Rbar^{\H_S},\quad
\EE_\Rbar^{\H_R}\tensor_{\EE_R} \EE_S \llongrightarrow
\EEtilde_\tSinfty \llongrightarrow \EEtilde_\Rbar^{\H_S}$$are
isomorphisms.

\endprop
\Proof The fact that $\widehat{\Rbar}^{\cH_S}= \widehat{\Sinfty}$
is proven in~[\An, Lem.~6.13]. The equalities in the first
displayed formula hold due to [\An, Prop.~6.14]. Those in the
second displayed formula follow arguing as in loc.~cit. In fact,
$\EEtilde_{\Rbar}^+$ (resp.~$\EEtilde_\tSinfty^+$) can be written
as in~\ref{EE}(5) as the limit
$\displaystyle\lim_{\infty\leftarrow n} \widehat{\Rbar}$
(resp.~$\displaystyle\lim_{\infty\leftarrow n}
\widehat{\tSinfty}$). To get the last two equalities one is left
to prove that~$\widehat{\Rbar}^{\H_S}= \widehat{\tSinfty}$. This
follows arguing as in [\An, Lem.~6.13]. The fact that the
inclusion~$\EE_\tS\subset \bigl(\EE_\Rbar\bigr)^{\H_S}$ is an
equality can be checked after base change $\EE_\tR^+ \to
\EEtilde_\tRinfty^+$ since the latter is faithfully flat
by~\ref{tEE}(2). But~$\EE_\tS\tensor_{\EE_\tR^+}
\EEtilde_\tRinfty^+\cong \EE_\tSinfty$ by~\ref{tEE}(3) and
$\bigl(\EE_\Rbar^+\bigr)^{\H_S}\tensor_{\EE_\tR^+}
\EEtilde_\tRinfty^+ \subset
\bigl(\EEtilde_\Rbar^+\bigr)^{\H_S}=\EEtilde_{\tSinfty}^+$.

The first equality in the last claim follows from~\ref{EE}(4); the
others follow from the second displayed formula and~\ref{tEE}(3).

\label AA. section\par \label vEEleqN. definition\par\ssection The
rings~$\AAtilde_\Sinfty$, $\AAtilde^\dagger_\Sinfty$, $\AA_S$,
$\AA_S^\dagger$, $\tAAtilde_\Sinfty$, $\tAAtildedagger_\Sinfty$,
$\tAA_S$ and~$\tAAdagger_S$\par
Define~$\AAtilde_\Rbar:=\WW\bigl(\EEtilde_\Rbar\bigr)$. It is
endowed with the following topology, called

{\it The weak topology.} Consider on~$\EEtilde_\Rbar$ the topology
having~$\{\barpi^n\Etilde_\Rbar^+\}_n$ as fundamental system of
neighborhoods of~$0$. On the truncated Witt
vectors~$\WW_m\bigl(\Etilde_\Rbar\bigr)$ we consider the product
topology via the isomorphism~$\WW_m\bigl(\Etilde_\Rbar \bigr)\cong
\bigl(\Etilde_\Rbar\bigr)^m$ given by the phantom components.
Eventually, the weak topology is defined as the projective limit
topology\/~$\displaystyle\WW
\bigl(\Etilde_\Rbar\bigr)=\lim_{\infty\leftarrow m} \WW_m
\bigl(\Etilde_\Rbar\bigr)$.

Alternatively, let~$\pi:=[\varepsilon]-1$ where~$[\varepsilon]$ is
the Teichm\"uller lift of~$\varepsilon$. For every~$n$
and\/~$h\in\NN$ define~$U_{n,h}:=p^n \AAtilde_\Rbar+ \pi^h
\AAtilde_\Rbar^+$. The weak topology on~$\AAtilde_\Rbar$
has~$\{U_{n,h}\}_{n,h\in\NN}$ as fundamental system of
neighborhoods.

Define~$\vv_\EE\colon \EEtilde_\Rbar\to \QQ\cup\{\infty\}$
by~$\vv_\EE(z)=\infty$ if~$z=0$ and~$\vv_\EE(z)=\min\{n\in\QQ\vert
\barpi^n z\in\EEtilde_\Rbar^+\}$.  For~$z=\sum_k [z_k]
p^k\in\AAtilde_\Rbar$ and~$N\in\NN$ we put $$\vv_\EE^{\leq
N}(z):=\inf\{\vv_\EE(z_k)\vert 0\leq k\leq N\}.$$For
every~$N\in\NN$ we have \spacing\item{{\rm (i)}} $\vv_\EE^{\leq
N}(x)=+\infty\Leftrightarrow x=0$; \item{{\rm (ii)}}
$\vv_\EE^{\leq N}(xy)\geq \vv_\EE^{\leq N}(x)+\vv_\EE^{\leq
N}(y)$; \item{{\rm (iii)}} $\vv_\EE^{\leq
N}(x+y)\geq\min(\vv_\EE^{\leq N}(x),\vv_\EE^{\leq N}(y))$ with
equality if~$\vv_\EE^{\leq N}(x)\neq \vv_\EE^{\leq N}(y)$;
\item{{\rm (iv)}} $\vv_\EE^{\leq N}(\barpi)=1$ and $\vv_\EE^{\leq
N}(\barpi x)=\vv_\EE^{\leq N}(\barpi)+\vv_\EE^{\leq N}(x)$;
\item{{\rm (v)}} $\vv_\EE^{\leq N}\bigl(\phi(x)\bigr)=p
\vv_\EE^{\leq N}(x)$; \item{{\rm (vi)}} $\vv_\EE^{\leq
N}\bigl(\gamma(x)\bigr)= \vv_\EE^{\leq N}(x)$ for
every~$\gamma\in\cG_R$.

\spacing \noindent The second claim in~(iii) and~(v) follow
since~$\EEtilde_\Rbar^+$ is by construction the $\barpi$--adic
completion of~$\cup_\Sinfty \EEtilde_\Sinfty^+$ and
each~$\EEtilde_\Sinfty^+$ is normal by~\ref{EE}(3). Note that the
topology on~$\AAtilde_\Rbar/p^{N+1}\AAtilde_\Rbar$ induced from
the weak topology on~$\AAtilde_\Rbar$ coincides with
the~$\vv_\EE^{\leq N}$ topology.\spacing

For every~$S$ as in~\ref{H} define
$$\AAtilde_\Sinfty:=\WW\bigl(\EEtilde_\Sinfty\bigr),\qquad \AAtilde_\Sinfty^+:=\WW\bigl(\EEtilde_\Sinfty^+\bigr),$$
$$\tAAtilde_\Sinfty:=\WW\bigl(\EEtilde_\tSinfty\bigr)\qquad\hbox{{\rm and}}\qquad
\tAAtilde_\Sinfty{}^+:=\WW\bigl(\EEtilde_\tSinfty^+\bigr).$$They
are subrings of~$\WW\bigl(\EEtilde_\Rbar\bigr)$ closed for the
weak topology.

\

{\it Overconvergent coefficients.} For $r\in\QQ_{>0}$ define
$\AAtilde_\Rbar^{(0,r]}$ as the subring of elements
$z=\sum\limits_{k=0}^{\infty}p^{k}[z_{k}]$ of $\AA_\Rbar$ such
that
$$\lim_{k\to\infty}rv_{\EE}(z_{k})+k=+\infty;$$see~[\AnBr, Prop.~4.3]. Write
$\AAtilde^\dagger=\bigcup\limits_{r\in\QQ_{>0}}\AAtilde^{(0,r]}$.
For $z=\sum\limits_{k\in\ZZ}p^{k}[z_{k}]\in\AAtilde^{(0,r]}$, put
$$w_{r}(z)=\cases{\infty & if $z=0$; \cr \inf\limits_{k\in\ZZ}\left( rv_{\EE}(z_{k})+k\right) & otherwise.\cr}$$
Thanks to~[\AnBr, Prop.~4.3] one knows that the map~$w_r\colon
\AAtilde_\Rbar^{(0,r]} \to \RR\cup\{\infty\}$ satisfies \item{{\rm
(i)}} $w_r(x)=+\infty\Leftrightarrow x=0$; \item{{\rm (ii)}}
$w_r(xy)\geq w_r(x)+w_r(y)$; \item{{\rm (iii)}}
$w_r(x+y)\geq\min(w_r(x),w_r(y))$; \item{{\rm (iv)}} $w_r(p)=1$
and $w_r(p x)=w_r(p)+w_r(x)$.

\spacing \noindent For every~$\Sinfty$
define~$\AAtilde_\Sinfty^{(0,r]}:=\AAtilde_\Sinfty\cap
\AAtilde_\Rbar^{(0,r]}$
and~$\AAtilde_\Sinfty^\dagger:=\AAtilde_\Sinfty\cap
\AAtilde_\Rbar^\dagger$. Similarly, define
$\tAAtilde_\Sinfty{}^{(0,r]}:=\tAAtilde_\Sinfty\cap
\AAtilde_\Rbar^{(0,r]}$
and\/~$\tAAtildedagger_\Sinfty:=\tAAtilde_\Sinfty\cap
\AAtilde_\Rbar^\dagger$. By~[\AnBr, Prop.~4.3] they are
$w_r$--adically complete and separated subrings
of~$\AAtilde_\Rbar$.

\

{\it Noetherian coefficients.} In~[\An, Appendix II] a ring~$\AA_S
\subset \WW\bigl(\EEtilde_\Sinfty\bigr)$ has been constructed
whose main properties are: \spacing \item{{\rm i.}} it is complete
and separated for the weak topology. In particular, it is
$p$--adically complete and separated. \spacing \item{{\rm ii.}}
$\AA_S \cap \bigl( p \WW\bigl(\EEtilde_\Rbar\bigr)\bigr)=p\AA_S $;
\spacing\item{{\rm iii.}} $\AA_S/p\AA_S\cong \EE_S$. In
particular, it is noetherian and regular. \spacing \item{{\rm
iv.}} it is endowed with continuous commuting actions
of\/~$\Gamma_R$ and of an operator\/~$\varphi$ lifting those
defined on~$\EE_R$; \spacing \item{{\rm v.}} it contains the
Teichm\"uller lifts of~$\epsilon$, $x_1$, $\ldots$, $x_d$;
\spacing \item{{\rm vi.}} $\AA_S$ is the unique finite and \'etale
$\AA_R$-algebra lifting the finite and \'etale
extension~$\EE_R\subset \EE_S$.\spacing

\noindent In particular, the formation of $\AA_S$ is functorial
in~$\Sinfty$. If one requires the existence of a subring
lifting~$\EE_S^+$, with suitable properties, one can in fact prove
that~$\AA_S$ is unique. We refer to [\AnBr, Prop.~4.49] for
details. Define~$\AA_S^{(0,r]}:=\AA_S\cap \AA_\Rbar^{(0,r]}$ and
$\AA_S^\dagger:=\AA_S\cap \AA_\Rbar^\dagger$.

Let\/~$\AA_\Rbar$ be the completion for the $p$--adic topology
of\/~$\cup_{\Sinfty} \AA_S$, where the union is taken over all
normal $\Rinfty$--subalgebras $\Sinfty\subset \Rbar$ such
that\/~$\Sinfty\left[p^{-1}\right]$ is finite \'etale over
$\Rinfty\left[p^{-1}\right]$. Let~$\tAA_S:=\AA_\Rbar^{\H_S}$
and~$\tAAdagger_S:=\tAA_S\cap \AAtilde_\Rbar^\dagger$.

\endssection

\label reccllASdagger. proposition\par\prop The
extensions~$\AAtilde_\Rinfty^\dagger \subset
\AAtilde_\Sinfty^\dagger$ and\/~$\AA_R^\dagger\subset
\AA_S^\dagger$ are finite and \'etale. Their reduction modulo~$p$
coincide with\/~$\EEtilde_\Rinfty\subset \EEtilde_\Sinfty$
and\/~$\EE_R\subset \EE_S$ respectively.

\endprop
\Proof It is clear that~$\AAtilde^\dagger_\Sinfty$ coincides
with~$\EEtilde_\Sinfty$ modulo~$p$ since it
contains~$\AAtilde_\Sinfty^+$ and $p\AAtilde_\Rbar\cap
\AAtilde_\Sinfty^\dagger= p\AAtilde^\dagger_\Sinfty$. The fact
that~$\AAtilde_\Rinfty^\dagger \subset \AAtilde_\Sinfty^\dagger$
is finite and \'etale is proven in~[\AnBr, Prop.~4.9]. See~[\AnBr,
Prop.~4.34] for the statements regarding~$\AA_R^\dagger\subset
\AA_S^\dagger$.

\label invariantsAA. lemmap\par\lemma The following hold:

\spacing \item{{\rm a)}} $\AA_S=\AA_\Rbar^{\cH_S}$,
$\AAtilde_\Sinfty=\AAtilde_\Rbar^{\cH_S}$,
and\/~$\tAAtilde_\Sinfty=\AAtilde_\Rbar^{\H_S}$. The same
equalities hold considering overconvergent coefficients i.~e.,
$\AA_S^{(0,r]}=\bigl(\AA_\Rbar^{(0,r]}\bigr)^{\cH_S}$,
$\AAtilde^{(0,r]}_\Sinfty=\bigl(\AA_\Rbar^{(0,r]}\bigr)^{\cH_S}$,
$\tAAtilde_\Sinfty{}^{(0,r]}=\bigl(\AA_\Rbar^{(0,r]}\bigr)^{\H_S}$
and $\AA_S^\dagger=\bigl(\AA_\Rbar^\dagger\bigr)^{\cH_S}$,
$\AAtilde^\dagger_\Sinfty=\bigl(\AA_\Rbar^\dagger\bigr)^{\cH_S}$,
$\tAAtildedagger_\Sinfty=\bigl(\AA_\Rbar^\dagger\bigr)^{\H_S}$.
\spacing \item{{\rm b)}} The natural maps
$\AAtilde_\Rinfty\tensor_{\AA_R} \AA_S \to \AAtilde_\Sinfty$,
$\tAA_R\tensor_{\AA_R} \AA_S \to \tAA_S$ and
$\tAAtilde_\Rinfty\tensor_{\AA_R} \AA_S \to \tAAtilde_\Sinfty$ are
isomorphisms.  Similarly, considering overconvergent coefficients,
the maps $\AAtilde^\dagger_\Rinfty\tensor_{\AA^\dagger_R}
\AA_S^\dagger \to \AAtilde^\dagger_\Sinfty$,
$\tAAdagger_R\tensor_{\AA_R^\dagger} \AA_S^\dagger \to
\tAAdagger_S$ and $\tAAtildedagger_\Rinfty\tensor_{\AA^\dagger_R}
\AA_S^\dagger \to \tAAtildedagger_\Sinfty$ are isomorphisms.

\spacing \item{{\rm c)}} We have
$\AA_S^\dagger/p\AA_S^\dagger=\AA_S/p\AA_S=\EE_S$,
$\tAAdagger_S/p\tAAdagger_S=\tAA_S/p\tAA_S=\EE_{\tS}$ and
$\tAAtildedagger_\Sinfty/p\tAAtildedagger_\Sinfty=
\tAAtilde_\Sinfty/p\tAAtilde_\Sinfty=\EEtilde_\tSinfty$.

\spacing \item{{\rm d)}} The maps
$\AAtilde_\Sinfty\tensor_{\AAtilde_\Vinfty}\AAtilde_\Vbar\to
\tAAtilde_\Sinfty$ and\/ $\AA_S\tensor_{\AA_V}\AAtilde_\Vbar\to
\tAAtilde_S$ are injective and have dense image for the weak
topology. The image of
$\AAtilde^{(0,r]}_\Sinfty\tensor_{\AAtilde_\Vinfty^{(0,r]}}
\AAtilde_\Vbar^{(0,r]}\to\tAAtilde_\Sinfty{}^{(0,r]}$ is dense for
the $w_r$--adic topology for every~$r\in\QQ_{>0}$.
\endlemma\Proof (a) We certainly have inclusions
$\AA_S\subset \AA_\Rbar^{\cH_S}$, $\AAtilde_\Sinfty \subset
\AAtilde_\Rbar^{\cH_S}$, $\tAAtilde_\Sinfty\subset
\AAtilde_\Rbar^{\H_S}$ and maps~$\AAtilde_\Rinfty\tensor_{\AA_R}
\AA_S \to \AAtilde_\Sinfty$, $\tAA_R\tensor_{\AA_R} \AA_S \to
\tAA_S$ and $\tAAtilde_\Rinfty\tensor_{\AA_R} \AA_S \to
\tAAtilde_\Sinfty$. The extension $\AA_R \to \AA_S$ is finite and
\'etale and, hence, $\AA_S$ is projective as $\AA_R$--module.
Since~$\AAtilde_\Rinfty$, $\tAA_R$ and~$\tAAtilde_\Rinfty$ are
$p$--adically complete and separated and~$p$ is a not a zero
divisor in these rings, $\AAtilde_\Rinfty\tensor_{\AA_R} \AA_S$,
$\tAA_R\tensor_{\AA_R} \AA_S$ and
$\tAAtilde_\Rinfty\tensor_{\AA_R} \AA_S$ are $p$--adically
complete and separated and~$p$ is a non--zero divisor. The same
holds for~$\AA_S$, $\AA_\Rbar^{\cH_S}$, $\AAtilde_\Sinfty$, $
\AAtilde_\Rbar^{\cH_S}$, $ \AAtilde_\Rbar^{\H_S}$, $\tAA_S$ and
$\tAAtilde_\Sinfty$. To check that all the inclusions and all the
maps above are isomorphisms it then suffices to show it
modulo~$p$.  This follows from~\ref{invariantsEE} if we prove that
$\tAA_S/p\tAA_S=\EE_{\tS}$ and
$\tAAtilde_\Sinfty/p\tAAtilde_\Sinfty=\EEtilde_\tSinfty$. Once
this is established the other statements in~(a) follow.

(c) Since by construction~$\tAAdagger_S\cap p\AAtilde_\Rbar=p
\tAAdagger_S$, $\tAAtildedagger_\Sinfty\cap
p\AAtilde_\Rbar=p\tAAtildedagger_\Sinfty$, $\AA_S^\dagger\cap
p\AAtilde_\Rbar =p\AA_S^\dagger$ and ~$\tAA_S\cap
p\AAtilde_\Rbar=p \tAA_S$, $\tAA_\Sinfty\cap
p\AAtilde_\Rbar=p\tAAtilde_\Sinfty$, $\AA_S\cap p\AAtilde_\Rbar
=p\AA_S$  the maps $\AA_S^\dagger/p\AA_S^\dagger\to
\AA_S/p\AA_S\to\EE_S$,
$\tAAdagger_S/p\tAAdagger_S\to\tAA_S/p\tAA_S\to\EE_{\tS}$ and
$\tAAtildedagger_\Sinfty/p\tAAtildedagger_\Sinfty\to
\tAAtilde_\Sinfty/p\tAAtilde_\Sinfty\to\EEtilde_\tSinfty$ are
injective. It follows from~[\AnBr, \S4.11(e)\&Lem.~4.20] that
there exists $\AA^+_S\subset
\AAtilde^+_\Sinfty\left\{{p^\alpha\over \pi^\beta} \right\}\subset
\AAtilde_\Sinfty^{(0,r]}$, for suitable~$\alpha$, $\beta\in\NN$
and $r\in\QQ_{>0}$, so that $\AA^+_S/p\AA_S^+\cong \EE_S^+$ and
$\AA_S^+$ is complete for the weak topology. Here,
$\AAtilde^+_\Sinfty\left\{{p^\alpha\over \pi^\beta} \right\}$
denotes the completion of $\AAtilde^+_\Sinfty\left[{p^\alpha\over
\pi^\beta} \right]$ with respect to the weak topology. In
particular, since $\pi$--adic convergence in
$\AAtilde^+_\Sinfty\left\{{p^\alpha\over \pi^\beta} \right\}$
implies convergence for the weak topology, $\AA_S^+$ is
$\pi$--adically complete. Note that
$\AA_S^+\tensor_{\AA_V^+}\AAtilde_{\Vbar}^+$ and
$\AAtilde_\Sinfty^+\tensor_{\AAtilde_\Vinfty^+}\AAtilde_\Vbar^+$
map to
$\AAtilde^+_\Sinfty\tensor_{\AA_\Vinfty^+}\AA_\Vbar^+\left\{{p^\alpha\over
\pi^\beta} \right\}$ and that the latter is $\pi$--adically
complete and is contained in $\tAAtilde_\Sinfty{}^{(0,r]}$. We
conclude that $\AA_S^\dagger/p\AA_S^\dagger$,
$\tAAdagger_S/p\tAAdagger_S$ and
$\tAAtildedagger_\Sinfty/p\tAAtildedagger_\Sinfty$ contain the
$\pi$--adic completion of the image of $\AA_S^+$,
$\AA_S^+\tensor_{\AA_V^+}\AAtilde_{\Vbar}^+$ and
$\AAtilde_\Sinfty^+\tensor_{\AAtilde_\Vinfty^+}\AAtilde_\Vbar^+$
respectively. Claim~(c) follows then from~\ref{tEE}.

(b) The fact that $\AAtilde^\dagger_\Rinfty\tensor_{\AA^\dagger_R}
\AA_S^\dagger \cong \AAtilde^\dagger_\Sinfty$ follows
from~\ref{reccllASdagger}.
Since~$\AA_R^\dagger\subset\AA_S^\dagger$ is finite and \'etale
by~\ref{reccllASdagger} there is a unique idempotent $e_{S/R}$
of~$\AA_S^\dagger\tensor_{\AA_R^\dagger}\AA_S^\dagger$ such that
for every~$x\in \AA_S^\dagger\tensor_{\AA_R^\dagger}\AA_S^\dagger$
we have $m(x)=\Tr_{\AA_S^\dagger/\AA_R^\dagger}\bigl(x
e_{S/R}\bigr)$. Here, $m$ is the multiplication map. Write
$e_{S/R}=\sum_{i=1}^u a_i\tensor b_i$ with $a_i$
and~$b_i\in\AA_S^{(0,s]}$ for some~$s\in \QQ_{>0}$. Then,
$e_{S/R}$ is an idempotent of $\tAAdagger_S\tensor_{\tAAdagger_R}
\tAAdagger_S$ and of
$\tAAtildedagger_\Sinfty\tensor_{\tAAtildedagger_\Sinfty}
\tAAtildedagger_\Sinfty$. By the first part of~(b) the extensions
$\tAA_R\subset \tAA_S$ and $\tAAtilde_\Rinfty\subset
\tAAtilde_\Sinfty$ are finite and \'etale and
$m(x)=\Tr_{\tAA_S/\tAA_R}\bigl(x e_{S/R}\bigr)$ and
$m(x)=\Tr_{\tAAtilde_\Sinfty/\tAAtilde_\Rinfty}\bigl(x
e_{S/R}\bigr)$.  We then get that for every~$x\in \tAAdagger_S$
(resp.~$\tAAtildedagger_\Sinfty$) we have $x=m(x\tensor
1)=\sum_{i=1}^u \Tr_{\tAA_S/\tAA_R}(xa_i) b_i$
(resp.~$x=m(x\tensor 1)=\sum_{i=1}^u
\Tr_{\tAAtilde_\Sinfty/\tAAtilde_\Rinfty}(xa_i) b_i$).
Since~$\Tr_{\tAA_S/\tAA_R}$
and~$\Tr_{\tAAtilde_\Sinfty/\tAAtilde_\Rinfty}$ send
overconvergent elements to overconvergent elements, we conclude
that the maps in the second part of~(b) are surjective.

Since the extension~$\AA_R^\dagger\subset \AA_S^\dagger$  is
finite and \'etale, $\AA_S^\dagger$ is projective as
$\AA_R^\dagger$--module. In particular, $p$ is not a zero divisor
in $\tAAdagger_R\tensor_{\AA_R^\dagger} \AA_S^\dagger$ and in
$\tAAtildedagger_\Rinfty\tensor_{\AA^\dagger_R} \AA_S^\dagger$ and
those rings are $p$--adically separated. Thus, to prove that the
maps in the second part of~(b) are injective, and hence are
isomorphisms, it suffices to prove that they are injective
modulo~$p$. This follows from~(c).

(d) Since the extensions $\AA_V\subset \AAtilde_{\Vinfty}\subset
\AAtilde_\Vbar$ are extensions of dvr's, they are flat. Since~$p$
is not a zero divisor in~$\AAtilde_\Sinfty$ and~$\AA_S$, it is not
a zero divisor in
$\AAtilde_\Sinfty\tensor_{\AAtilde_\Vinfty}\AAtilde_\Vbar\to
\tAAtilde_\Sinfty$ and\/ $\AA_S\tensor_{\AA_V}\AAtilde_\Vbar\to
\tAAtilde_S$. Thus, to check the injectivity in~(d) we may reduce
modulo~$p$. The density can be checked modulo~$p^n$ for
every~$n\in\NN$ and, using induction, it suffices in fact to prove
it for~$n=1$. Then, the first claim of~(d) follows
from~\ref{reccllASdagger} and~\ref{tEE}(1).

We prove the second claim of~(d). Suppose that~$r=a/b$ with~$a$
and~$b\in\NN$ and let $A_{\Sinfty,(a,b)}$
(resp.~$A_{\tSinfty,(a,b)}$) denote the $p$--adic completion
 of $\AAtilde_\Sinfty^+\left[p^a\over
[\overline{\pi}]^{\left({p-1\over p}\right)b} \right]$ (resp.~of
$\tAAtilde_\Sinfty{}^+\left[p^a\over
[\overline{\pi}]^{\left({p-1\over p}\right)b} \right]$),
where~$[\barpi]$ is the Teichm\"uller lift of~$\barpi$. Arguing as
in the proof of~[\AnBr, Lem.~4.20] one has
$$A_{\tSinfty,(a,b)} \subseteq
\left\{x\in\tAAtilde_\Sinfty{}^{(0,r]}\vert w_r(x)\geq
0\right\}\subseteq A_{\tSinfty,(a,b)}\left[{1\over
[\barpi]}\right].
$$Since~$w_r([\barpi])>0$, we conclude that $A_{\tSinfty,(a,b)}\left[{1\over
[\barpi]}\right]$ is dense in~$\tAAtilde_\Sinfty{}^{(0,r]}$ for
the $w_r$--adic topology. Since the $[\barpi]$--adic completion
of~$A_{\Sinfty,(a,b)}\tensor_{\AAtilde_\Vinfty^+}\AAtilde_\Vbar^+$
is contained in the $w_r$--adic closure of $
\AAtilde^{(0,r]}_\Sinfty\tensor_{\AAtilde_\Vinfty^{(0,r]}}
\AAtilde_\Vbar^{(0,r]}$ and its image in $A_{\tSinfty,(a,b)}$
contains $\barpi^\ell A_{\tSinfty,(a,b)}$ by~\ref{tEE}, the
conclusion follows.

\label tAsdaggeretale. corollary\par\cor The
extensions~$\tAAtilde_\Rinfty \subset \tAAtilde_\Sinfty$,
$\tAAtildedagger_\Rinfty \subset \tAAtildedagger_\Sinfty$,
$\tAA_R\subset \tAA_S$ and\/~$\tAAdagger_R\subset \tAAdagger_S$
are finite and \'etale. Their reduction modulo~$p$ coincide
with\/~$\EEtilde_\tRinfty\subset \EEtilde_\tSinfty$
and\/~$\EE_\tR\subset \EE_\tS$ respectively.

\endcor

\label phigamma. section\par\ssection $(\phi,\Gamma)$--modules and
Galois representations\par  Let\/~$S$ be as in~\ref{H}. Let
$\Rep\bigl(\cG_S\bigr)$ be the abelian tensor category of finitely
generated~$\ZZ_p$--modules endowed with a continuous action
of~$\cG_S$.

\spacing \noindent Let~$\AA$ be one of the
rings~$\AAtilde_\Sinfty$, $\AAtilde^\dagger_\Sinfty$, $\AA_S$,
$\AA_S^\dagger$, $\tAAtilde_\Sinfty$, $\tAAtildedagger_\Sinfty$,
$\tAA_S$ and~$\tAAdagger_S$ and let~$\Gamma$ be respectively
$\Gamma_S$ or~$\tGamma_S$. Let\/~$(\phi,\Gamma)-\Mod_{\AA}$
(resp.~$(\phi,\Gamma_S)-\Modet_{\AA}$) be the tensor category of
finitely generated~$\AA$-modules $D$ endowed with\spacing
\item{{i.}} a semi-linear action of\/~$\Gamma$;
\spacing\item{{ii.}} a semi-linear homomorphism~$\phi$ commuting
with~$\Gamma$ (resp.~so that~$\phi\tensor 1\colon
D\tensor_{\AA}^\phi \AA \to D $ is an isomorphism of
$\AA$-modules).

\spacing  \noindent  Note that if~$\AA=\AA_S$, then $\AA_S$ is
noetherian and $(\phi,\Gamma_S)-\Mod_{\AA_S}$
(resp.~$(\phi,\Gamma_S)-\Modet_{\AA_S}$) is an abelian category.

\spacing  For any object~$M$ in~$\Rep\bigl(\cG_S\bigr)$, define
$$\cD(M):=\Bigl(\AA_\Rbar\tensor_{\ZZ_p} M \Bigr)^{\cH_S},
\qquad \D(M):=\Bigl(\AA_\Rbar\tensor_{\ZZ_p} M \Bigr)^{\H_S}$$
$$\cDtilde(M):=\Bigl(\AAtilde_\Rbar\tensor_{\ZZ_p} M \Bigr)^{\cH_S},
\qquad \Dtilde(M):=\Bigl(\AAtilde_\Rbar\tensor_{\ZZ_p} M
\Bigr)^{\H_S}.$$Note that $\cD(M)$ (resp.~$\cDtilde(M)$) is an
$\AA_S$-module (resp.~$\AAtilde_\Sinfty$--module) endowed with a
semi-linear action of\/~$\Gamma_S$. Analogously, $\D(M)$
(resp.~$\Dtilde(M)$) is a $\tAA_S$--module (resp.~a
$\tAAtilde_\Sinfty$--module) endowed with a semi-linear action of
$\tGamma_S$. Analogously, define
$$\cD^\dagger(M):=\Bigl(\AA_\Rbar^\dagger\tensor_{\ZZ_p} M \Bigr)^{\cH_S}, \qquad
\D^\dagger(M):=\Bigl(\AA_\Rbar^\dagger\tensor_{\ZZ_p} M
\Bigr)^{\H_S}$$
$$\cDtildedagger(M):=\Bigl(\AAtilde_\Rbar^\dagger\tensor_{\ZZ_p} M \Bigr)^{\cH_S},
\qquad
\Dtildedagger(M):=\Bigl(\AAtilde_\Rbar^\dagger\tensor_{\ZZ_p} M
\Bigr)^{\H_S}.$$Then, $\cD^\dagger(M)$ (resp.~$\cDtildedagger(M)$)
is an $\AA^\dagger_S$-module
(resp.~$\AAtilde_\Sinfty^\dagger$--module) endowed with a
semi-linear action of\/~$\Gamma_S$. Analogously, $\D^\dagger(M)$
(resp.~$\Dtildedagger(M)$) is a $\tAAdagger_S$--module (resp.~a
$\tAAtildedagger_\Sinfty$--module) endowed with a semi-linear
action of $\tGamma_S$.

The homomorphism~$\varphi$ on~$\AAtilde_\Rbar$
and~$\AAtilde^\dagger_\Rbar$ defines a semi-linear action
of~$\varphi$ on all these modules commuting with the action
of~$\Gamma_S$ (resp.~of $\tGamma_S$).
\enddefi

\label MainThm. theorem\par\thm The functor~$\cD$ defines an
equivalence of abelian tensor categories from the
category~$\Rep\bigl(\cG_S\bigr)$ to the
category~$(\phi,\Gamma_S)-\Modet_{\AA_S}$. Let\/~$M$ be a finitely
generated $\ZZ_p$--module endowed with a continuous action
of\/~$\cG_S$. The inverse is defined associating to an \'etale
$(\phi,\Gamma_S)$--module $D$ the $\cG_R$--module
$\cV(D):=\Bigl(\AA_\Rbar\tensor_{\AA_S} D \Bigr)_{(\phi-1)}$.
\endthm

\Proof See~[\An, Thm.~7.11].

\label DMtARAA. lemma\par\lemma Let\/~$M$ be a finitely generated
$\ZZ_p$--module endowed with a continuous action of\/~$\cG_S$.
Then,\spacing

\item{{\rm (i)}} $\cD(M)$ (resp.~$\D(M)$, $\cDtilde(M)$,
$\Dtilde(M)$) is an \'etale $(\phi,\Gamma_S)$--module over~$\AA_S$
(resp.~$\tAA_S$, $\AAtilde_\Sinfty$, $\tAAtilde_\Sinfty$);\spacing

\item{{\rm (i')}} $\cDdagger(M)$ (resp.~$\D^\dagger(M)$,
$\cDtildedagger(M)$, $\Dtildedagger(M)$) is an \'etale
$(\phi,\Gamma_S)$--module over~$\AA^\dagger_S$
(resp.~$\tAAdagger_S$, $\AAtilde^\dagger_\Sinfty$,
$\tAAtildedagger_\Sinfty$);\spacing

\item{{\rm (ii)}} $\displaystyle\cD(M)=\lim_{\infty\leftarrow
n}\cD(M/p^nM)$, $\displaystyle\D(M)=\lim_{\infty\leftarrow
n}\D(M/p^nM)$, $\displaystyle\cDtilde(M)=\lim_{\infty\leftarrow
n}\cDtilde(M/p^nM)$,
$\displaystyle\Dtilde(M)=\lim_{\infty\leftarrow n}\Dtilde(M/p^nM)$
where the limits are inverse limits with respect to~$n\in\NN$.
More precisely, $\cD(M)/p^n\cD(M)\cong\cD(M/p^nM)$,
$\D(M)/p^n\D(M)\cong\D(M/p^nM)$,
$\cDtilde(M)/p^n\cDtilde(M)\cong\cDtilde(M/p^nM)$ and
$\Dtilde(M)/p^n\Dtilde(M)\cong\Dtilde(M/p^nM)$ for
every~$n\in\NN$.

\spacing

\item{{\rm (ii')}} if $M$ is torsion, then $\cDdagger(M)=\cD(M)$,
$\D^\dagger(M)=\D(M)$, $\cDtildedagger(M)=\cDtilde(M)$ and
$\Dtildedagger(M)=\Dtilde(M)$;\spacing

\spacing \item{{\rm (iii)}} the natural maps
$$\cD(M) \tensor_{\AA_S} \tAA_S \llongrightarrow \D(M),
\qquad \cD(M)\tensor_{\AA_S} \AA_\Rbar \llongrightarrow M
\tensor_{\ZZ_p} \AA_\Rbar$$and $$\cD(M) \tensor_{\AA_S}
\AAtilde_\Sinfty \llongrightarrow \cDtilde(M), \qquad
\cD(M)\tensor_{\AA_S} \tAAtilde_\Sinfty \llongrightarrow
\Dtilde(M)
$$are isomorphisms; \spacing

\spacing \item{{\rm (iii')}} the natural maps
$$\cDdagger(M) \tensor_{\AA^\dagger_S} \tAAdagger_S \llongrightarrow \D^\dagger(M),
\qquad \cD^\dagger(M)\tensor_{\AA^\dagger_S} \AA^\dagger_\Rbar
\llongrightarrow M \tensor_{\ZZ_p} \AA^\dagger_\Rbar$$and
$$\cDdagger(M) \tensor_{\AA^\dagger_S} \AAtilde^\dagger_\Sinfty \llongrightarrow
\cDtildedagger(M), \qquad \cD^\dagger(M)\tensor_{\AA^\dagger_S}
\tAAtildedagger_\Sinfty \llongrightarrow \Dtildedagger(M)
$$are isomorphisms; \spacing \item{{\rm (iv)}} the natural maps
$$\cDdagger(M)\tensor_{\AA^\dagger_S}\AA_S\llongrightarrow\cD(M),\qquad
\D^\dagger(M)\tensor_{\tAAdagger_S}\tAA_S\llongrightarrow \D(M)$$
and
$$\cDtildedagger(M)\tensor_{\AAtilde^\dagger_\Sinfty}\AAtilde_\Sinfty\llongrightarrow
\cDtilde(M),\qquad\Dtildedagger(M)\tensor_{\tAAtildedagger_\Sinfty}
\tAAtilde_\Sinfty\llongrightarrow\Dtilde(M)$$are isomorphisms.
\endlemma\Proof We refer the reader to [\An, Thm.~7.11] and~[\AnBr, Thm~4.40] for the
proofs that $\cD(M)$ and~$\cDdagger(M)$ are \'etale
$(\phi,\Gamma_S)$--modules and that
$\cD(M)=\cDdagger(M)\tensor_{\AA_S^\dagger}\AA_S$, that
$\cD(M)\tensor_{\AA_S} \AA_\Rbar\cong M \tensor_{\ZZ_p} \AA_\Rbar$
and $ \cDdagger(M)\tensor_{\AA^\dagger_S} \AA^\dagger_\Rbar\cong M
\tensor_{\ZZ_p} \AA^\dagger_\Rbar$. Claims~(i), (i') and~(iv)
follow from this and the displayed isomorphisms. We prove the
other statements.

Due to~\ref{reccllASdagger}, \ref{invariantsAA}
and~\ref{tAsdaggeretale} to prove (ii'), (iii) and~(iii') one may
pass to an extension $\Sinfty\subset \Tinfty$ in~$\Rbar$ finite,
\'etale and Galois after inverting~$p$. For example,
$\bigl(M\tensor_{\ZZ_p}\AA_\Rbar\bigr)^{\cH_T}=\bigl(M\tensor_{\ZZ_p}\AA_\Rbar\bigr)^{\cH_S}\tensor_{\AA_S}
\AA_T$ and
$\bigl(M\tensor_{\ZZ_p}\AA_\Rbar\bigr)^{\H_T}=\bigl(M\tensor_{\ZZ_p}\AA_\Rbar\bigr)^{\H_S}\tensor_{\tAA_S}
\tAA_T$ by \'etale descent. Hence, if
$\bigl(M\tensor_{\ZZ_p}\AA_\Rbar\bigr)^{\cH_T}\tensor_{\AA_T}\tAA_T\to
\bigl(M\tensor_{\ZZ_p}\AA_\Rbar\bigr)^{\H_T}$ is an isomorphism,
taking the $\cH_S$--invariants, we get the claimed isomorphism
$\bigl(M\tensor_{\ZZ_p}\AA_\Rbar\bigr)^{\cH_S}\tensor_{\AA_S}\tAA_S\to
\bigl(M\tensor_{\ZZ_p}\AA_\Rbar\bigr)^{\H_S}$.

Suppose first that there exists~$N\in\NN$ such that~$p^NM=0$.
Then, there exists an extension~$\Sinfty\subset \Tinfty$ such
that~$\cH_T\subset \cH_S$ acts trivially on~$M$.
Replacing~$\Sinfty$ with~$\Tinfty$, we may then assume
that~$\cH_S$ acts trivially on~$M$. By~\ref{invariantsAA} we have
$\AA^\dagger_S/p^N\AA_S^\dagger=\AA_S/p^N \AA_S $,
$\tAAdagger_S/p^N\tAAdagger_S=\tAA_S/p^N\tAA_S$,
$\AAtilde^\dagger_\Sinfty/p^N
\AAtilde^\dagger_\Sinfty=\AAtilde_\Sinfty/p^N \AAtilde_\Sinfty $
and eventually
$\tAAtildedagger_\Sinfty/p^N\tAAtildedagger_\Sinfty=\tAAtilde_\Sinfty/p^N\tAAtilde_\Sinfty$.
Furthermore, we have in this case $\cD(M)=M\tensor_{\ZZ_p} \AA_S$,
$\D(M)=M\tensor_{\ZZ_p} \tAA_S$, $\cDtilde(M)=M\tensor_{\ZZ_p}
\AAtilde_\Sinfty$, $\Dtilde(M)=M\tensor_{\ZZ_p} \tAAtilde_\Sinfty$
and similarly for the overconvergent $(\phi,\Gamma_S)$--modules.
Then, the claims follow from~\ref{invariantsAA}.

Assume next that~$M$ is free of rank~$n$. It follows from~[\AnBr,
Thm.~4.40] that there exists an extension $\Rinfty\subset \Tinfty$
in~$\Rbar$ finite, \'etale and Galois after inverting~$p$ such
that~$\cD^\dagger(M)\tensor_{\AA_S^\dagger}\AA_T^\dagger $ is a
free $\AA_T^\dagger$--module of rank~$n$.  As we have seen above
we may and will replace~$\Sinfty$ with~$\Tinfty$ so
that~$\cDdagger(M)$ (resp.~$\cD(M)$) is a free
$\AA^\dagger_S$--module (resp.~$\AA_S$--module). Fix a
basis~$\{e_1,\ldots,e_n\}$ of~$\cDdagger(M)$. It is also a
$\AA_R$--basis of~$\cD(M)$. Hence, it is a basis
over~$\AA_\Rbar^\dagger$ (resp.~$\AAtilde_\Rbar^\dagger$, $
\AA_\Rbar$, $\AAtilde_\Rbar$) of $M\tensor_{\ZZ_p}
\AA_\Rbar^\dagger$
(resp.~$M\tensor_{\ZZ_p}\AAtilde_\Rbar^\dagger$, $M\tensor_{\ZZ_p}
\AA_\Rbar$, $M\tensor_{\ZZ_p}\AAtilde_\Rbar$). Since~$\cH_S$
and~$\H_S$ act trivially on~$\{e_1,\ldots,e_s\}$, we get
claims~(iii) and~(iii'). For example,
$\cDtilde(M)=\bigl(M\tensor_{\ZZ_p}\AA_\Rbar\bigr)^{\H_S}=\AAtilde_\Sinfty
e_1\dirsum \cdots \dirsum \AAtilde_\Sinfty
e_n=\cD(M)\tensor_{\AA_S}\AAtilde_\Sinfty$.

We are left to  prove~(ii). We may assume that~$M$ is torsion
free, since the claim for the torsion part is trivial. Note that
$\cD(M)$, $\cDtilde(M)$, $\D(M)$ and~$\Dtilde(M)$ are submodules
of invariants of free modules over $p$--adically complete and
separated rings. For example, $\cD(M)=\bigl(M\tensor_{\ZZ_p}
\AA_\Rbar\bigr)^{\cH_S}\subset M\tensor_{\ZZ_p} \AA_\Rbar$. Hence,
they are themself $p$--adically complete and separated. It
suffices to show that for every~$n\in\NN$ the map from their
quotient modulo~$p^n$ to $\cD\bigl(M/p^n M\bigr)$
(resp.~$\cDtilde\bigl(M/p^n M\bigr)$, $\D\bigl(M/p^n M\bigr)$,
$\Dtilde\bigl(M/p^n M\bigr)$) is an isomorphism. Due to~(iii) it
suffices to show it for~$\cD(M)$ and in this case it follows from
the fact that~$\cD$ is an exact functor by~\ref{MainThm}.

\

We state the following theorem relating the cohomology of the
various $(\phi,\Gamma)$--modules introduced above.

\thm The natural maps
$$\H^n\bigl(\Gamma_S,\cD(M)\bigr)\llongrightarrow
\H^n\bigl(\Gamma_S,\cDtilde(M)\bigr) \llongrightarrow
\H^n\left(\cG_S,M\tensor_{\ZZ_p} \AAtilde_\Rbar\right),$$

$$\H^n\bigl(\tGamma_S,\D(M)\bigr)\llongrightarrow
\H^n\bigl(\tGamma_S,\Dtilde(M)\bigr) \llongrightarrow
\H^n\left(\G_S,M\tensor_{\ZZ_p} \AAtilde_\Rbar\right),$$

$$\H^n\bigl(\Gamma_S,\cDdagger(M)\bigr)\llongrightarrow
\H^n\bigl(\Gamma_S,\cDtildedagger(M)\bigr) \llongrightarrow
\H^n\left(\cG_S,M\tensor_{\ZZ_p}
\AAtilde_\Rbar^\dagger\right)$$and

$$\H^n\bigl(\Gamma_S,\D^\dagger(M)\bigr)\llongrightarrow
\H^n\bigl(\Gamma_S,\Dtilde^\dagger(M)\bigr) \llongrightarrow
\H^n\left(\G_S,M\tensor_{\ZZ_p} \AAtilde_\Rbar^\dagger\right)$$are
all isomorphisms.
\endthm\Proof See~\ref{decompletecohomology}.

\endsection

\section Galois cohomology and $(\phi,\Gamma)$--modules\par In
this section we show how, given a finitely generated
$\ZZ_p$--module $M$ with continuous action of~$\cG_S$, one can
compute the cohomology groups~$\H^n(\cG_S,M)$ and~$\H^n(\G_S,M)$
in terms of the associated $(\phi,\Gamma_S)$--modules~$\cD(M)$,
$\cDtilde(M)$, $\cD^\dagger(M)$, $\cDtildedagger(M)$, $\D(M)$,
$\Dtilde(M)$, $\D^\dagger(M)$ and~$\Dtildedagger(M)$. We start
with the following crucial:

\prop The map $\varphi-1$ on~$\AAtilde_\Rbar$,
$\AAtilde_\Rbar^\dagger$, $\AA_\Rbar$ and\/~$\AA_\Rbar^\dagger$ is
surjective and its kernel is~$\ZZ_p$. Furthermore, the exact
sequence $$ 0\llongrightarrow \ZZ_p \llongrightarrow
\AAtilde_\Rbar \llongmaprighto{\phi-1} \AAtilde_\Rbar
\llongrightarrow 0$$admits a continuous right splitting
$\sigma\colon \AAtilde_\Rbar \rightarrow\AAtilde_\Rbar$ (as sets)
so that $\sigma\bigl(\AA_\Rbar\bigr)\subset \AA_\Rbar$,
$\sigma\bigl(\AAtilde^\dagger_\Rbar\bigr)\subset
\AAtilde^\dagger_\Rbar$
and\/~$\sigma\bigl(\AA^\dagger_\Rbar\bigr)\subset
\AA^\dagger_\Rbar$.\endprop\Proof See~\ref{phiAA}.

\defi Let\/~$\D$ be a $(\phi,\Gamma_S)$-module over~$\AA_S$
or~$\AAtilde_\Sinfty$ or~$\AA_S^\dagger$
or~$\AAtilde_\Sinfty^\dagger$ (resp.~over $\tAA_S$
or~$\tAAtilde_\Sinfty$ or~$\tAAdagger_S$
or~$\tAAtildedagger_\Sinfty$). Define $\cC^\bullet(\Gamma_S,\D)$
(resp.~$\cC^\bullet(\tGamma_S,\D)$) to be the complex of
continuous cochains with values in~$\D$.
\smallskip

Let~$\cT^\bullet(\D)$ (resp.~$\ctT^\bullet(\D)$) be the mapping
cone associated to $\phi-1\colon \cC^\bullet(\Gamma_S,\D)
\rightarrow \cC^\bullet(\Gamma_S,\D)$ (resp.~to $\phi-1\colon
\cC^\bullet(\tGamma_S,\D) \rightarrow \cC^\bullet(\tGamma_S,\D)$).
\enddefi

\label maintheorem. theorem\par\thm There are isomorphisms of
$\delta$-functors from the category of $\Rep(\cG_S)$ to the
category of abelian groups:
$$\rho_i\colon
\H^i\bigl(\cT^\bullet(\cD(\_))\bigr)\isomarrow \H^i(\cG_S,\_),
\qquad \trho_i \colon
\H^i\bigl(\ctT^\bullet(\D(\_))\bigr)\isomarrow \H^i(\G_S,\_),$$
$$\rhotilde_i\colon
\H^i\bigl(\cT^\bullet(\cDtilde(\_))\bigr)\isomarrow
\H^i(\cG_S,\_), \qquad \trhotilde_i \colon
\H^i\bigl(\ctT^\bullet(\Dtilde(\_))\bigr)\isomarrow
\H^i(\G_S,\_),$$
$$\rho_i^\dagger\colon
\H^i\bigl(\cT^\bullet(\cD^\dagger(\_))\bigr)\isomarrow
\H^i(\cG_S,\_), \qquad \trhodagger_i \colon
\H^i\bigl(\ctT^\bullet(\D^\dagger(\_))\bigr)\isomarrow
\H^i(\G_S,\_),$$
$$\rhotildedagger_i\colon
\H^i\bigl(\cT^\bullet(\cDtildedagger(\_))\bigr)\isomarrow
\H^i(\cG_S,\_), \qquad \trhotildedagger_i \colon
\H^i\bigl(\ctT^\bullet(\Dtildedagger(\_))\bigr)\isomarrow
\H^i(\G_S,\_)$$The isomorphisms $\trho_i$, $\trhotilde_i$,
$\trhodagger_i$ and $\trhotildedagger_i$ are $\G_V$--equivariant.

Furthermore, all the maps in the  square
$$\matrix{\H^i\bigl(\cT^\bullet(\cD^\dagger(\_))\bigr) & \llongrightarrow
& \H^i\bigl(\cT^\bullet(\cDtildedagger(\_))\bigr) \cr
\big\downarrow & & \big\downarrow \cr
\H^i\bigl(\cT^\bullet(\cD(\_))\bigr) & \llongrightarrow &
\H^i\bigl(\cT^\bullet(\cDtilde(\_))\bigr),\cr  }
$$induced by the natural inclusions of $(\phi,\Gamma_S)$--modules
$\cD^\dagger(M)\subset \cD(M) \subset \cDtilde(M)$ and
$\cD^\dagger(M)\subset \cD^\dagger(M) \subset \cDtilde(M)$, for
$M\in\Rep(\cG_S)$, are isomorphisms and they are compatible with
the maps~$\rho^\dagger_i$, $\rhotildedagger_i$, $\rho_i$
and\/~$\rhotilde_i$. Similarly, all the maps in the  square
$$\matrix{\H^i\bigl(\cT^\bullet(\D^\dagger(\_))\bigr) &
\llongrightarrow & \H^i\bigl(\cT^\bullet(\Dtildedagger(\_))\bigr)
\cr \big\downarrow & & \big\downarrow \cr
\H^i\bigl(\cT^\bullet(\D(\_))\bigr) & \llongrightarrow &
\H^i\bigl(\cT^\bullet(\Dtilde(\_))\bigr)\cr  }
$$are isomorphisms and are compatible with the maps~$\trhodagger_i$,
$\trhotildedagger_i$, $\trho_i$ and\/~$\trhotilde_i$
\endthm

\Proof First of all we exhibit in~\ref{exhibit} the maps~$\rho_i$
and~$\trho_i$ (with or without~$\widetilde{}$ or~$\dagger$) so
that they are compatible with the displayed squares and~$\trho_i$
they are compatible with the residual action of~$\G_V$ (if one
exists). We then prove that they are isomorphisms
in~\ref{theyareisomorphisms}. Eventually, we show that they are
isomorphisms of $\delta$--functors in~\ref{theyaredeltafunctors}.

\label exhibit. section\par\ssection The maps\par Let $M$ be a
$\ZZ_p$-representation of $\cG_S$. Let~$D(M)$ and~$A$ be \enspace
(1) $\cD(M)$ and~$\AA$, \enspace (2) $\cDtilde$ and~$\AAtilde$,
\enspace (3) $\cDdagger(M)$ and~$\AA^\dagger$ or\enspace (4)
$\cDtildedagger(M)$ and~$\AAtilde^\dagger$. Since in each
case~$D(M)\tensor_{A_S} A_\Rbar\cong M\tensor_{\ZZ_p} A_\Rbar$
by~\ref{DMtARAA} and due to~\ref{phiAA} we have exact sequences of
$\cG_S$--modules \labelf exactMpAA\par
$$0 \llongrightarrow M \llongrightarrow D(M)\tensor_{A_S}
A_\Rbar \llongmaprighto{\phi-1} D(M)\tensor_{A_S} A_\Rbar
\llongrightarrow 0 \eqno{{(\numfo)}}$$\advance\fonu by1\noindent
Similarly, let~$D'(M)$ and~$A'$ be \enspace (1) $\D(M)$
and~$\tAA$, \enspace (2) $\Dtilde$ and~$\tAAtilde$, \enspace (3)
$\D^\dagger(M)$ and~$\tAAdagger$ or\enspace (4) $\Dtildedagger(M)$
and~$\tAAtildedagger$. In each case~$D'(M)\tensor_{A_S'}
A_\Rbar'\cong M\tensor_{\ZZ_p} A_\Rbar'$ by~\ref{DMtARAA}. thanks
to~\ref{phiAA} we get exact sequences of $\cG_S$--modules \labelf
exactMptAA\par$$0 \llongrightarrow M \llongrightarrow
D'(M)\tensor_{A'_S} A_\Rbar' \llongmaprighto{\phi-1}
D'(M)\tensor_{A'_S} A_\Rbar' \llongrightarrow
0\eqno{{(\numfo)}}.$$The maps in both exact sequences are
continuous for the weak topology by~\ref{DMtARAA}.

Let $(\alpha,\beta)$ be an $n$-cochain of~$\cT^\bullet(D(M))$
i.~e., in $\cC^{n-1}(\Gamma_S,D(M))\times \cC^n(\Gamma_S,D(M))$.
Define
$$c_{\alpha,\beta}^n:=\beta+(-1)^n d \bigl(\sigma(\alpha)\bigr) \in
\cC^n(\cG_S,M\tensor_{\ZZ_p} A_\Rbar),$$where  $d$ is the
differential operator on~$\cC^n(\Gamma_S,D(M))$ and~$\sigma$ is
the left inverse of~$\phi-1$ defined in~\ref{phiAA} (for each of
the four possibilities for~$A$).

\

Recall that the derivation on~$\cT^\bullet(D(M))$ is given by~$
d\bigl((\alpha,\beta)\bigr)=\bigl((-1)^n
(\phi-1)(\beta)+d\alpha,d\beta\bigr)$. Thus, $(\alpha,\beta)$ is
an $n$--cocycle if and only if and
satisfies~$(-1)^n(\phi-1)\beta+d \alpha=0$ and~$d \beta=0$. In
this case, $d c_{\alpha,\beta}^n=0$ and $(\phi-1)
c_{\alpha,\beta}^n=(\phi-1)\beta+(-1)^n  d
(\phi-1)\gamma=(\phi-1)\beta+(-1)^n  d \alpha=0$. Thus,
$c_{\alpha,\beta}^n$ is an $n$-cocycle in $\cC^n(\cG_S,M)$
by~(\ref{exactMpAA}).

\

Choose a different left inverse~$\sigma'$ of~$\phi-1$. Then,
$(\phi-1) (\sigma'-\sigma)=0$ so that
$\bigl(\sigma'-\sigma\bigr)(\alpha)$ lies in $\cC^{n-1}(\cG_S,M)$.
Thus, $ \beta+(-1)^n d\bigl(\sigma'(\alpha)\bigr)-
\beta-(-1)^nd\bigl(\sigma(\alpha)\bigr)=(-1)^nd (\sigma' -\sigma)
(\alpha)$. In particular, $c_{\alpha,\beta}^n$ depends on the
choice of~$\sigma$ up to a coboundary with values in~$M$.

\

Let $(\alpha,\beta)=((-1)^{n-1}(\phi-1)b+ da, db)\in
\cC^{n-1}(\Gamma_S,D(M))\times \cC^n(\Gamma_S,D(M))$ be an
$n$-coboundary in~$\cT^\bullet(D(M))$. Then, $c_{\alpha,\beta}^n=
db+(-1)^{2n-1} d (\sigma\circ(\phi-1))(b)+(-1)^n d
(\sigma(d(a)))$. Note that $\bigl(1-(\sigma\circ(\phi-1)\bigr) b$
and~$\sigma\bigl(d(a)\bigr)-d\bigl(\sigma(a)\bigr)$ are
annihilated by~$\phi-1$. Hence, $c_{\alpha,\beta}^n$ is the image
via the differential of $ (1 -
(\sigma\circ(\phi-1))(b)+(-1)^n\left(\sigma\bigl(d(a)\bigr)-d\bigl(\sigma(a)\bigr)\right)
$ lying in~$\cC^{n-1}(\cG_S,M)$. In particular, it is a
coboundary.

\

\noindent We thus get a map
$$r_i^M\colon \H^i\bigl(\cT^\bullet(D(M))\bigr)\llongrightarrow
\H^i(\cG_S,M).$$By construction it is functorial in~$M$. In
case~(1) we get the map~$\rho_i$, in case~(2) we
get~$\rhotilde_i$, in case~(3) we get the map~$\rho^\dagger_i$ and
in case~(4) we get~$\rhotildedagger_i$. By construction they are
compatible with the first commutative displayed square appearing
in~\ref{maintheorem}.

\

Analogously, using~(\ref{exactMptAA}), one gets the claimed
map~$r'_i$. In case~(1) we get the map~$\trho_i$, in case~(2) we
get~$\trhotilde_i$, in case~(3) we get the map~$\trhodagger_i$ and
in case~(4) we get~$\trhotildedagger_i$. They are compatible with
the second commutative displayed square appearing
in~\ref{maintheorem}. Furthermore, we also have actions of~$\G_V$
and we claim that~$r'_i$ is $\G_V$--equivariant.

Indeed, let~$(\alpha,\beta)$ be an $n$--cocycle in
$\cC^{n-1}(\tGamma_S,D'(M))\times \cC^n(\tGamma_S,D'(M))$.
Let~$g\in\G_V$. Then,
$g\bigl((\alpha,\beta)\bigr)=\bigl(g(\alpha),g(\beta)\bigr)$ and
$c^n_{g((\alpha,\beta))}= g(\beta)+(-1)^n
d\bigl(\sigma(g(\alpha))\bigr)$. On the other hand,
$g\bigl(c^n_{\alpha,\beta}\bigr)= g(\beta)+(-1)^n
g\bigl(d(\sigma(a))\bigr)$. Note that $g \circ d=d\circ g$ and
$(\phi-1) (\sigma \circ g- g \circ \sigma )=0$ since~$\phi$ is
$\Gamma_V$--equivariant. Thus,
$c^n_{g((\alpha,\beta))}-g\bigl(c^n_{\alpha,\beta}\bigr)=(-1)^nd
\bigl((\sigma \circ g- g \circ \sigma )(\alpha)\bigr)$ is a
coboundary in~$\cC^n(\cG_S,M)$.

\endssection

\label theyareisomorphisms. proposition\par\prop The maps
$\rho_i$, $\rhotilde_i$, $\rho^\dagger_i$, $\rhotildedagger_i$,
$\trho_i$, $\trhotilde_i$, $\trhodagger_i$ and
$\trhotildedagger_i$ are isomorphisms.\endprop \Proof\advance\ssnu
by-1 We use the notation of~\ref{exhibit}. Since~$\cT^\bullet(D)$
and~$\ctT^\bullet(D')$ are mapping cones, we get exact sequences
\labelf exactcDM\par
$$\H^{n-1}(\Gamma_S,D(M)) \llongmaprighto{\delta'_n} \H^n(\cT^\bullet(D(M))
\longrightarrow \H^n(\Gamma_S,D(M))
\lllongmaprighto{(-1)^n(\phi-1)}
\H^n(\Gamma_S,\cD(M))\eqno{{(\numfo)}}
$$\advance\fonu by1\noindent and\labelf exactDM\par
$$\H^{n-1}(\tGamma_S,D'(M)) \llongmaprighto{\delta'_n} \H^n(\ctT^\bullet(D'(M))
\longrightarrow \H^n(\tGamma_S,D'(M))
\lllongmaprighto{(-1)^n(\phi-1)}
\H^n(\tGamma_S,\D(M)).\eqno{{(\numfo)}}$$\advance\fonu
by1\noindent They are compatible with respect to the natural
inclusions $\cD^\dagger(M)\subset \cD(M) \subset \cDtilde(M)$ and
$\cD^\dagger(M)\subset \cD^\dagger(M) \subset \cDtilde(M)$
(resp.~$\D^\dagger(M)\subset \D(M) \subset \Dtilde(M)$ and
$\D^\dagger(M)\subset \D^\dagger(M) \subset \Dtilde(M)$). Thanks
to~\ref{decompletecohomology} we deduce that the horizontal arrows
in the two displayed squares of~\ref{maintheorem} are
isomorphisms. We are then left to prove that $\rhotilde_i$,
$\rhotildedagger_i$,
 $\trhotilde_i$ and $\trhotildedagger_i$ are isomorphisms.

From the exactness of~(\ref{exactMpAA}) and~(\ref{exactMptAA}) we
get the exact sequences \labelf exactcGRM\par
$$\H^{n-1}(\cG_S, M\tensor_{\ZZ_p}A_\Rbar) \llongmaprighto{\delta_n} \H^n(\cG_S,M)
\llongrightarrow \H^n(\cG_S, M\tensor_{\ZZ_p}A_\Rbar)
\llongmaprighto{\phi-1} \H^n(\cG_S,
M\tensor_{\ZZ_p}A_\Rbar)\eqno{{(\numfo)}}$$\advance\fonu
by1\noindent and\labelf exactGRM\par $$\H^{n-1}(\G_S,
M\tensor_{\ZZ_p}A_\Rbar) \llongmaprighto{\delta_n} \H^n(\G_S,M)
\llongrightarrow \H^n(\G_S, M\tensor_{\ZZ_p}A_\Rbar)
\llongmaprighto{\phi-1} \H^n(\G_S,
M\tensor_{\ZZ_p}A_\Rbar)\eqno{{(\numfo)}}.$$\advance\fonu
by1\noindent Thanks to~\ref{decompletecohomology} the inflation
maps
$$\Inf_n\colon \H^n(\Gamma_S,D(M))\llongrightarrow \H^n(\cG_S,
D(M)\tensor_{A_S}
A_\Rbar)=\H^n(\cG_S,M\tensor_{\ZZ_p}A_\Rbar)$$and
$$ \Inf_n\colon \H^n(\tGamma_S,D'(M))\llongrightarrow \H^n(\G_S,
\D(M)\tensor_{A'_S}
A'_\Rbar)=\H^n(\G_S,M\tensor_{\ZZ_p}A'_\Rbar)$$in cases (2)
and~(4) of~\ref{exhibit} are isomorphisms. Take a cocycle~$\tau$
in $\cC^{n-1}(\Gamma_S, D(M))$. One constructs
$\delta_n\bigl(\Inf_{n-1}(\tau)\bigr)$ as
$d\left(\sigma\bigl(\Inf_{n-1}(\tau)\bigr)\right)$. On the other
hand, $\delta_n'(\tau)=(\tau,0)$ in
$\cC^{n-1}(\Gamma_S,D(M))\times \cC^n(\Gamma_S,D(M))$ and
$c_{\tau,0}^n=(-1)^n d\bigl(\sigma(\tau)\bigr)$. Thus,
$\delta_n\circ (-1)^{n-1} \Inf_{n-1}=\rho_n^M \circ (-1)
\delta'_n$. If~$(\alpha,\beta)$ is an $n$--cocycle
in~$\cT^\bullet(D(M))$, its image in~$\H^n(\tGamma_S,D(M))$ is the
class of~$\beta$. The image of~$c^n_{\alpha,\beta}$ in
$\H^n(\cG_S, M\tensor_{\ZZ_p} A_\Rbar)$ is the class
of~$\beta+(-1)^n d\bigl( \sigma(\alpha)\bigr)$ i.~e., of~$\beta$.
We conclude that the exact sequences~(\ref{exactcDM})
and~(\ref{exactcGRM}) are compatible via~$r_n^M$ and the inflation
maps~$\Inf_n$ and~$\Inf_{n-1}$ i.~e., the following diagram
commutes (the rows continue on the left and on the right):
$$\matrix{ \H^{n-1}(\Gamma_S,D(M)) & \maprighto{-\delta'_n} &
\H^n(\cT^\bullet(D(M))& \longrightarrow & \H^n(\Gamma_S,D(M))&
\maprighto{(-1)^n(\phi-1)}  \cr \mapdownr{(-1)^{n-1}\Inf_{n-1}} &
& \mapdownl{\rho_n^M}& & \mapdownr{\Inf_n} \cr \H^{n-1}(\cG_S,
M\tensor_{\ZZ_p}A_\Rbar) & \maprighto{\delta_n} & \H^n(\cG_S,M)
&\longrightarrow & \H^n(\cG_S, M\tensor_{\ZZ_p}A_\Rbar) &
\maprighto{\phi-1} .\cr}$$An analogous argument implies that the
exact sequences~(\ref{exactDM}) and~(\ref{exactGRM}) are
compatible via~$r'_n$ and the inflation maps~$\Inf_n$
and~$\Inf_{n-1}$. The proposition follows.

\advance\ssnu by1\label theyaredeltafunctors. proposition\par\prop
The functors $\rho_i$, $\rhotilde_i$, $\rho^\dagger_i$,
$\rhotildedagger_i$, $\trho_i$, $\trhotilde_i$, $\trhodagger_i$
and $\trhotildedagger_i$ are $\delta$--functors.
\endprop\Proof We use the notation of~\ref{exhibit}.
We prove the proposition for~$r_i$. The proof for~$r'_i$ is
similar. Let $0 \rightarrow M_1 \rightarrow M_2 \rightarrow M_3
\rightarrow 0$ an exact sequence of $\cG_S$--representations. We
need to prove that the diagram
$$\matrix{ \H^n\bigl(\cT^\bullet(D(M_3))\bigr) &
\llongmaprighto{\delta}&
\H^{n+1}\bigl(\cT^\bullet(D(M_1))\bigr)\cr \mapdownl{r_i^{M_3}}& &
\mapdownr{r_i^{M_1}}\cr \H^i(\cG_S,M_3) & \llongmaprighto{\delta'}
& \H^{i+1}(\cG_S,M_1), \cr}$$where $\delta$ and $\delta'$ are the
connecting homomorphisms, commutes.

Let $(\alpha,\beta)$ be an $n$--cocycle in
$\cT^n(D(M_3))=\cC^{n-1}(\Gamma_S,D(M_3))\times
\cC^n(\Gamma_S,D(M_3))$. Let $(a,b)\in
\cC^{n-1}(\Gamma_S,D(M_2))\times \cC^n(\Gamma_S,D(M_2))$ be a
lifting of~$(\alpha,\beta)$. Then,
$\delta\bigl((\alpha,\beta)\bigr)$ is represented by $d(a,b)=
\bigl((-1)^n(\phi-1)b+d a, d b\bigr)$ and
$$c^{n+1}_{d(a,b)}=db+(-1)^{2n+1} d \bigl(\sigma\circ
(\phi-1)(b)\bigr)+(-1)^{n+1}d \bigl(\sigma(d a)\bigr).$$

On the other hand, $c^n_{(\alpha,\beta)}= \beta+(-1)^n d
\bigl(\sigma(\alpha)\bigr)$. Consider $c^n_{(a,b)}= b+(-1)^n d
\bigl(\sigma(a)\bigr)$
in~$\cC^n(\cG_S,M_2\tensor_{\ZZ_p}A_\Rbar)$. Then,
$\gamma:=c^n_{(a,b)}-\sigma\bigl((\phi-1)\bigl(c^n_{(a,b)}\bigr)\bigr)$
lies in~$\cC^n(\cG_S,M_2)$. Furthermore, it
lifts~$c^n_{(\alpha,\beta)}$ since~$ \sigma\bigl( (\phi-1)
c^n_{(\alpha,\beta)}\bigr)=0$ because~$(\alpha,\beta)$ is a
cocycle. Then, the class
of~$\delta'\bigl(c^n_{(\alpha,\beta)}\bigr)$ is~$d \gamma$. To
compute it we may take the differential
in~$\cC^n(\cG_S,M_2\tensor_{\ZZ_p}A_\Rbar)$ i.~e.,
$$dc^n_{(a,b)}-d
\sigma\bigl((\phi-1)\bigl(c^n_{(a,b)}\bigr)\bigr)= d b- d
\bigl(\sigma\circ(\phi-1)(b)\bigr)+(-1)^{n+1} d \bigl(\sigma(d
a)\bigr).$$Here, we used
$d\sigma\bigl((\phi-1)\bigl(d(\sigma(a))\bigr)\bigr)= d
\bigl(\sigma(d a)\bigr)$. The conclusion follows.

\endsection

\

\

\centerline{\titlefont Global Theory.}

\

\section \'Etale cohomology and relative $(\phi,\Gamma)$--modules\par
\noindent As in the Introduction, let $X$ denote a smooth,
geometrically irreducible and proper scheme over $\Spec(V)$. Fix a
field extension $K\subset M \subset \Kbar$. In this section we
review a Grothendieck topology on~$X$, introduced by Faltings
in~[\FALAST] and denoted $\gerX_M$, and its relation to \'etale
cohomology; see \ref{etaletheorem}. We also define the analogue
Grothendieck topology, $\hatgerX_M$ on the formal completion $\cX$
of $X$ along its special fiber. In this section we study $p$-power
torsion sheaves on these Grothendieck topologies and compare their
cohomology theories. The main result of this section is the
following. Let $\bbL$ be an \'etale local system of
$\ZZ/p^s\ZZ$-modules on $X_K$, for some $s\ge 1$. Fix a geometric
generic point $\eta$ of $X_K^\rig$ and denote by ${\bf L}$ the
fiber of $\bbL^\rig$ (the corresponding \'etale local system on
the rigid space $X_K^\rig$ attached to $X_K$) at $\eta$. For every
\'etale morphism $\cU\tto \cX$ such that $\cU$ is affine,
$\cU=\Spf(R_{\cU})$, with $R_{\cU}$ a small $V$-algebra (see \S2)
we let $\cD_{\cU}({\bf L})$, $\D_{\cU}({\bf L})$, $\cDtilde
_{\cU}({\bf L})$, $\Dtilde_\cU({\bf L})$ denote the respective
$(\phi,\Gamma)$--modules over $R_{\cU}$. For each $i\ge 0$ the
associations $\cU\tto \H^i(\cT^\bullet(\cD_{\cU}({\bf L})))$,
$\cU\tto \H^i(\cT^\bullet(\D_{\cU}({\bf L})))$, $\cU\tto
\H^i(\cT^\bullet(\cDtilde_{\cU}({\bf L})))$, $\cU\tto
\H^i(\cT^\bullet(\Dtilde_{\cU}({\bf L})))$ are functorial and we
denote by $\cH^{i,\arit}(\bbL)$, $\cH^{i,\geom}(\bbL)$,
$\cH^{i,\t,\arit}(\bbL)$ and $\cH^{i,\t,\geom}(\bbL)$ respectively
the associated sheaves on $\cX^\et$.  We have

\label etalePhiGamma. theorem\par\thm There are spectral sequences
$$
i)\quad E_2^{p,q}=\H^q\Bigl(\cX^\et,
\cH^{p,\ast,\geom}\bigl(\bbL\bigr)\Bigr)\Longrightarrow
\H^{p+q}\Bigl(\bigl(X_{\Kbar}\bigr)^\et, \bbL\Bigr).
$$
$$
ii)\quad E_2^{p,q}=\H^q\Bigl(\cX^\et,
\cH^{p,\ast,\arit}\bigl(\bbL\bigr)\Bigr)\Longrightarrow
\H^{p+q}\Bigl(\bigl(X_{K}\bigr)^\et, \bbL\Bigr).
$$
where $\ast$ stands for $\emptyset$ or $\t$. Moreover, the
spectral sequence i) is compatible with the residual $\G_V$-action
on all of its terms.
\endthm

\noindent
The proof of theorem \ref{etalePhiGamma} will take the rest of this section.

\label toposes. ssection\par\ssection Some Grothendieck topologies
and associated sheaves\par  Following {\rm [\FALAST, \S 4, p.~214]} we define
the following site:

\smallskip

We denote by $X_M^\et$ the small \'etale site of\/~$X_M$ and
by~$\Sh(X_M^\et)$ the category of sheaves of abelian groups
on~$X_M^\et$.

\smallskip

\noindent {\it The site $\gerX_M$.}\enspace The objects consists of
pairs $\bigl(U,W\bigr)$ where\spacing

\item{{\rm (i)}} $U\to X$ is \'etale; \spacing \item{{\rm (ii)}}
$W \to U\tensor_V M$ is a finite \'etale cover.\spacing

\noindent The maps are compatible maps of pairs and the coverings of
a pair~$(U,W)$ are finite families~$\{(U_\alpha,W_\alpha)\}_\alpha$
over~$(U,W)$ such that $\amalg_\alpha U_\alpha \to U$ and
$\amalg_\alpha W_\alpha \to W$ are surjective. It is easily checked
that we get a noetherian Grothendieck topology, in the sense of
[\Artin, I.0.1 \& II.5.1]. Note that one has a final object,
namely~$(X,X_M)$. Let $\Sh(\gerX_M)$ be the category of sheaves of
abelian groups in~$\gerX_M$.

\smallskip

Let $\cX$ denote the formal scheme associated to~$X$ i.~e., the
formal completion of $X$ along its special fiber. Denote by
$\cX^\et$ the small \'etale site on~$\cX$ and by~$\Sh(\cX^\et)$ the
category of sheaves of abelian groups on~$\cX^\et$.

\

\noindent {\it The sites $\cU^{M,\fet}$ and\/
$\cU_M^\fet$.}\enspace Let $\cU\to \cX$ be an \'etale map of
formal schemes. Define~$\cU^{M,\fet}$ to be the category whose
objects are pairs $\bigl(\cW,L\bigr)$ where

\spacing \item{{\rm (i)}} $L$ is a finite extension of\/~$K$
contained in~$M$;\spacing \item{{\rm (ii)}} $\cW \to
\cU^\rig\tensor_K L$ is a finite \'etale cover of $L$--rigid
analytic spaces; here $\cU^\rig$ denotes the $K$--rigid analytic
space associated to~$\cU$.\spacing

\noindent Define $\Hom_{\cU^{M,\fet}}\bigl((\cW',L'),
(\cW,L)\bigr)$ to be empty if~$L \not\subset L'$ and to be the set
of morphisms $g\colon \cW'\to \cW\tensor_L L'$ of $L'$--rigid
analytic spaces if~$L\subset L'$. The coverings of a
pair~$(\cW,L)$ in~$\cU^{M,\fet}$ are finite families of
pairs~$\{(\cW_\alpha,L_\alpha)\}_\alpha$ over~$(\cW,L)$ such that
$\amalg_\alpha \cW_\alpha\to \cW$ is surjective. Define the fiber
product of two pairs $(\cW',L')$ and\/~$(\cW'',L'')$ over a pair
$(\cW,L)$ to be~$(\cW'\fibprod_{\cW} \cW'', L''') $ with~$L'''$
equal to the composite of~$L'$ and~$L''$. It  is the fiber product
in the category~$\cU^{M,\fet}$

Let~$\cU_2\to\cU_1$ be a map of formal schemes over~$\cX$. Assume
that they are \'etle over~$\cX$. We then have a morphism of
Grothendieck topologies $\rho_{\cU_1,\cU_2}\colon \cU_1^{M,\fet}
\to \cU_2^{M,\fet}$ given on objects by sending~$(\cW,L) \mapsto
\left(\cW\fibprod_{\cU_1^\rig} \cU_2^\rig,L\right)$. It is clear
how to define such a map for morphisms and that it sends covering
families to covering families.\smallskip

Let~${\cal S}_\cU$ be the system of morphisms in~$\cU^{M,\fet}$ of
pairs $(\cW',L')\to (\cW,L)$ such that $g\colon \cW'\to
\cW\tensor_L L'$ is an isomorphism. Then,

\slabel multsys. lemma\par\slemma The following hold:\spacing

\item{{\rm i)}}  the composite of two composable elements
of\/~${\cal S}_\cU$ is in~${\cal S}_\cU$;\spacing

\item{{\rm ii)}} given a map $\cU_2\to \cU_1$ of formal schemes
\'etale over~$\cX$, we have $\rho_{\cU_1,\cU_2}\bigl({\cal
S}_{\cU_1}\bigr)\subset {\cal S}_{\cU_2} $;

\item{{\rm iii)}} the base change of an element of\/~${\cal
S}_\cU$ via a morphism in~$\cU^{M,\fet}$ is again an element
of\/~${\cal S}_\cU$;

\item{{\rm iii)}} if $f\colon (\cW_1,L_1) \to (\cW,L)$ and
$g\colon (\cW_2,L_2) \to(\cW,L)$ are morphisms lying in~${\cal
S}_\cU$ and if\/~$h\colon (\cW_1,L_1)\to (\cW_2,L_2)$ is a
morphism in~$\cU^{M,\fet}$ such that~$f=g\circ h$, then $h$ is
in~${\cal S}_\cU$;

\endslemma

\Proof Left to the reader.

\smallskip

Thanks to~\ref{multsys} the category $\cU^{M,\fet}$ localized with
respect to~${\cal S}_\cU$ exists and we denote it by~$\cU_M^\fet$.
Note that the fiber product of two pairs over a given one exists
in~$\cU_M^\fet$ and it coincides with the fiber product
in~$\cU^{M,\fet}$. The coverings of a pair~$(\cW,L)$
in~$\cU_M^\fet$ are still finite families of
pairs~$\{(\cW_\alpha,L_\alpha)\}_\alpha$ over~$(\cW,L)$ such that
$\amalg_\alpha \cW_\alpha\to \cW$ is surjective. By~\ref{multsys}
the category $\cU_M^\fet$ and the given families of covering
define a noetherian Grothendieck topology. By abuse of notation we
will simply write~$\cW$ for an object~$(\cW,L)$ of~$\cU_M^\fet$.

\smallskip

We recall that, given pairs $(\cW_1,L_1)$ and\/~$(\cW_2,L_2)$
in~$\cU_M^\fet$, one defines the set of homomorphisms
$$\Hom_{\cU_M^\fet}\bigl((\cW_1,L_1),(\cW_2,L_2)\bigr):=\lim_{(\cW',L')\to (\cW_1,L_1)}
\Hom_{\cU^{M,\fet}}\bigl((\cW',L'),(\cW_2,L_2)\bigr),$$where the
direct limit is taken over all morphisms $(\cW',L')\to
(\cW_1,L_1)$ in~${\cal S}_\cU$. Equivalently, due
to~\ref{multsys}, it is the set of classes of morphisms
$(\cW_1,L_1) \leftarrow (\cW',L') \rightarrow (\cW_2,L_2)$, where
$(\cW',L')\to (\cW_1,L_1) $ is in~${\cal S}_\cU$, and two such
diagrams $(\cW_1,L_1) \leftarrow (\cW',L') \rightarrow
(\cW_2,L_2)$ and $(\cW_1,L_1) \leftarrow (\cW'',L'') \rightarrow
(\cW_2,L_2)$ are equivalent if and only if there is a third one
$(\cW_1,L_1) \leftarrow (\cW''',L''') \rightarrow (\cW_2,L_2)$
mapping to the two. If~$(\cW_1,L_1) \leftarrow (\cW',L')
\rightarrow (\cW_2,L_2)$ and $(\cW_2,L_2) \leftarrow (\cW'',L'')
\rightarrow (\cW_3,L_3)$ are two homomorphisms, the composite
$(\cW_1,L_1) \leftarrow (\cW''',L''') \rightarrow (\cW_3,L_3)$ is
defined by taking $(\cW''',L''')$ to be the fiber product
of~$(\cW',L')$ and $(\cW'',L'')$ over~$(\cW_2,L_2)$.
\smallskip

If~$\cU_2\to\cU_1$ is a map of formal schemes over~$\cX$ and they
are \'etale over~$\cX$, due to~\ref{multsys} the map
$\rho_{\cU_1,\cU_2}\colon \cU_1^{M,\fet} \to \cU_2^{M,\fet}$
extends to the localized categories and defines a map of
Grothendieck topologies $\cU_{2,M}^\fet \to \cU_{1,M}^\fet$ which,
by abuse of notation, we write $\cW \to \cW\fibprod_{\cU_1^\rig}
\cU_2^\rig$.

\

\noindent {\it The site $\hatgerX_M$.}\enspace Define $\hatgerX_M$
to be the category of pairs~$(\cU,\cW)$ where $\cU\to \cX$ is an
\'etale map of formal schemes and $\cW$ is an object
of~$\cU_M^\fet$. A morphism of pairs~$(\cU',\cW')\to (\cU,\cW)$ is
defined to be a morphism $\cU'\to \cU$ as schemes over~$\cX$ and a
map $\cW' \to \cW\fibprod_{\cU^\rig} \cU^{'\rig}$
in~$\cU_M^{'\fet}$. A morphism is a covering if $\cU'\to \cU$ is
\'etale and surjective and $\cW' \to \cW\fibprod_{\cU^\rig}
\cU^{'\rig}$ is a covering in~$\cU_M^{'\fet}$.

Remark that the fiber product~$(\cU''',\cW''')$ of two pairs
$(\cU',\cW')$ and~$(\cU'',\cW'')$ over a pair~$(\cU,\cW)$ exists
putting~$\cU''':=\cU'\fibprod_{\cU}\cU''$ and~$\cW'''$ to be the
fiber product in~$\cU_M^{'''\fet}$
of~$\cW'\fibprod_{\cU^{'\rig}}\cU^{'''\rig}$
and~$\cW''\fibprod_{\cU^{''\rig}}\cU^{'''\rig}$
over~$\cW\fibprod_{\cU^{\rig}}\cU^{'''\rig}$. The
pair~$(\cX,(\cX^\rig,K))$ is a final object in~$\hatgerX_M$. We
let\/ $\Sh\bigl(\hatgerX_M\bigr)$ be the category of sheaves of
abelian groups on~$\hatgerX_M$.

\

\noindent We remark that in all the categories~$\Sh(\_)$ introduced
above AB3$^\ast$ and AB5 hold and the representable objects provide
families of generators. In particular, one has enough injectives;
see [\Artin, Thm.~II.1.6 \& \S~II.1.8].

\endssection

\ssection Morphisms of Grothendieck topologies\par One has natural
functors:\spacing

\item{{\rm I}} $u_{X,M}\colon {\gerX}_M \longrightarrow
\left(X\tensor_V M \right)^\et  $ with $u_{X,M}(U,W):=W$;\spacing

\item{{\rm II.a}} $v_{X,M}\colon X^\et\longrightarrow {\gerX}_M$
given by $v_{X,M}(U):=\bigl(U,U\tensor_V M\bigr)$;

\spacing

\item{{\rm II.b}} $\gerv_{\cX,M}\colon \cX^\et \longrightarrow
{\hatgerX}_M $ given by
$\gerv_{\cX,M}(\cU):=(\cU,(\cU^\rig,K))$;\spacing

\item{{\rm III}} $\mu_{X,M}\colon {\gerX}_M \longrightarrow
{\hatgerX}_M$ given by $\mu_{X,M}\bigl(U,W\bigr):=(\cU,(\cW,L))$
where $\cU$ is the formal scheme associated to~$U$ and if the
cover $W\to U\tensor_V M $ is defined over a finite extension
$K\subset L$, contained in~$M$, then $\cW\to \cU^\rig_L$ is the
pull--back via $\cU^\rig_L\to U^\rig_L$ of the associated finite
and \'etale cover of rigid analytic spaces~$W^\rig\to U^\rig_L$.

\spacing

\item{{\rm IV}} $\nu_X\colon X^\et \to \cX^\et$ given by
$\nu_{X}(U)=\cU$ where $\cU$ is the formal scheme associated to~$U$.

\spacing \noindent Let~$K\subset M_1 \subset M_2\subset \Kbar$ be
field extensions. Define\spacing

\item{{\rm V.a}} $\beta_{M_1,M_2}\colon {\gerX}_{M_1}\to \gerX_{M_2}$
by $\beta_{M_1,M_2}(U,W)=\bigl(U,W\tensor_{M_1}M_2\bigr)$.

\spacing \item{{\rm V.b}} $\gerbeta_{M_1,M_2}\colon
{\hatgerX}_{M_1}\to \hatgerX_{M_2}$ by
$\gerbeta_{M_1,M_2}(\cU,\cW)=\bigl(\cU,\cW\bigr)$.

\spacing

Due to the definition of~$\hatgerX_M$, the functors~$\mu_{X,M}$
and~$\gerbeta_{M_1,M_2}$ are well defined. More precisely,
given~$(U,W)$, the image $\mu_{X,M}\bigl(U,W\bigr)$ does not
depend on the subfield~$L\subset M$ to which~$W$ descends.
Analogously, given $\cU\in\cX^\et$, then $\gerbeta_{M_1,M_2}$
sends the multiplicative system~$\cS_\cU$, used to
define~$\hatgerX_{M_1}$, to the multiplicative system~$\cS_\cU$
used to define~$\hatgerX_{M_2}$.

It is clear that the above functors send covering families to
covering families and commute with fiber products.  In particular,
they are morphisms of topologies, see~{\rm [\Artin, Def.~II.4.5]}.
Given any such functor, denote it by~$g$, we let\/~$g_\ast$
and\/~$g^\ast$ be the induced morphisms of the associated category
of sheaves for the given topologies; see~{\rm [\Artin,
p.~41--42]}. Note that the functors above preserve final objects.
Then, the induced functor $g^\ast$ on topoi of sheaves is exact
by~{\rm [\Artin, Thm.~II.4.14]}.

We work out an example. If $\cF$ is a sheaf on~$\hatgerX_M$,
then~$\mu_{X,M,\ast}(\cF)$ is the sheaf on~$\gerX_M$ defined
by~$(U,W)\mapsto \cF\bigl(\mu_{X,M}\bigl(U,W\bigr)\bigr)$. If
$\cF$ is a sheaf on~$\gerX_M$, then~$\mu_{X,M}^\ast(\cF)$ is the
sheaf associated to the separated presheaf defined
by~$(\cU,\cW)\mapsto \lim_{(U',W')}\cF(U',W')$ where the limit is
the direct limit taken over all pairs~$(U',W')$ in~$\gerX_M$ and
all maps~$(\cU,\cW) \to \mu_{X,M}\bigl(U',W'\bigr) $
in~$\hatgerX_M$.

\

{\bf Notation:} If $\cF$ is a sheaf on~$X^\et$ or is
in~$\Sh(X^\et)^\NN$ (see section 5 for the definition),
we write $\cF^\form$ for~$\nu_X^\ast(\cF)$,
respectively for~$\nu_X^{\ast,\NN}(\cF)$.

\spacing

If~$\bbL$ is a locally constant sheaf on~$X_M^\et$, we still
denote~$\bbL$ its push forward $u_{X,M,\ast}(\bbL)\in
\Sh(\gerX_M)$. It is a locally constant sheaf on~$\gerX_M$.

\spacing

If~$\cF$ is a  sheaf on~$\gerX_M$ or is in~$\Sh(\gerX_M)^\NN$, we
denote by~$\cF^\rig$ the pull--back~$\mu_{X,M}^\ast(\cF)$. Note
that if~$\cF\in\Sh(\gerX_M)$ is locally constant, then~$\cF^\rig$
is also a locally constant sheaf on~$\hatgerX_M$.

\endssection

\label definegeopoints. ssection\par\ssection Stalks {\rm
[\FALAST, p.~214]}\par Let\/~$K_x$ be a finite field extension
of\/~$K$ contained in~$\Kbar$ and denote by~$V_x$ its valuation
ring.

Fix a map\/~$x\colon \Spec(V_x)\to X$   of $V$--schemes and denote
by~$\barx\colon \Spec(\Vbar)\to X$ the composite of\/~$x$ with the
natural map~$\Spec(\Vbar)\to \Spec(V_x)$. Taking $p$--adic
completions, $x$ (resp.~$\barx$) defines a unique morphism of formal
$V$--schemes $\cx\colon \Spf(V_x)\to \cX$ (resp.~$\barcx\colon
\Spf(\widehat{\Vbar})\to \cX$).

\spacing

Let\/~$\cO_{X,x}^\sh$ be the the direct limit~$\lim_i R_i$ taken
over all pairs~$\{(R_i, f_i)\}_i$ where~$\Spec(R_i)$ is \'etale
over~$X$ and $f_i\colon R_i\to \Vbar$ defines a point
over~$\barx$. Let\/~$\cF$ be a sheaf in $\Sh\bigl(X^\et\bigr)$.
The stalk\/~$\cF_{x}$ of\/~$\cF$ at~$x$ is defined as
$\cF_{x}=\cF\bigl( \cO_{X,x}^\sh\bigr)$ by which we mean the
direct limit~$\lim_i \cF\bigl(\Spec(R_i)\bigr)$.

Define~$\barcO_{X,x,M}$ as the direct limit~$\lim_{i,j} R_{i,j}'$
over the pairs~$\{(R_{i,j}', R_{i,j}'\to \Vbar)\}_{i,j}$ where (1)
$R_{i,j}'$ is an integral  $R_i$--algebra and is normal as a ring,
(2) $R_{i,j}'\tensor_V K$ is a finite and \'etale extension
of\/~$R_i\tensor_V M$, (3) the composite~$R_i\tensor_V M \to
R_{i,j}'\tensor_V K \to \Kbar$ is~$r\tensor \ell \mapsto f_i(r)
\cdot \ell$. If\/~$\cF$ is a sheaf in $\Sh\bigl(\gerX_M\bigr)$, we
then write $\cF_x$ or equivalently~$\cF\bigl(\barcO_{X,x,M}\tensor_V
K\bigr)$ for the direct
limit~$\lim_{i,j}\cF\bigl(\Spec(R_i),\Spec(R_{i,j}'\tensor_V
K)\bigr)$. We call it the stalk of\/~$\cF$ at~$x$.

Let\/~$G_{x,M}$ be the Galois group of\/~$\barcO_{X,x,M}\tensor_V
K$ over~$\cO_{X,x}^\sh\tensor_V M$. Then, $\cF_x$ is endowed with
an action of\/~$G_{x,M}$.

\spacing

Analogously, let\/~$ \cO_{\cX,\cx}^\sh$ be the direct
limit~$\lim_i S_i$ over all pairs~$\{(S_i, g_i)\}_{i\in I}$ such
that\/~$S_i$ is $p$--adically complete and separated $V$--algebra,
$\Spf(S_i)\to \cX$ is an \'etale map of formal schemes and
$g_i\colon S_i\to \widehat{\Vbar}$ defines a formal point
over~$\barcx$. If\/~$\cF$ is a sheaf in $\Sh\bigl(\cX^\et\bigr)$,
the stalk\/~$\cF_{\cx}$ of\/~$\cF$ at\/~$\cx$ is defined to be the
direct limit~$\cF\bigl( \cO_{\cX,\cx}^\sh\bigr):=\lim_{i\in
I}\cF\bigl(\Spf(R_i)\bigr)$.

Write~$\barcO_{\cX,\cx,M}$ for the direct limit~$\lim_{i,j}
S_{i,j}'$ over all triples~$\{(S_{i,j}',S_{i,j}'\to
\widehat{\Vbar}, L_{i,j})\}_{i,j}$ where (1)~$L_{i,j}$ is a finite
extension of\/~$K$ contained in~$M$, (2) $S_{i,j}'$ is an integral
extension of\/~$S_i$ and is normal as a ring, (3)
$S_{i,j}'\tensor_V K$ is a finite and \'etale $S_i\tensor_V
L_{i,j}$--algebra, (4) the composite $S_i\tensor_V L_{i,j} \to
S_{i,j}'\tensor_V K \to \widehat{\Kbar}$ is $a\tensor \ell \mapsto
g_i(a)\cdot \ell$. Given a sheaf\/~$\cF$ in
$\Sh\bigl(\hatgerX_M\bigr)$, denote by\/ $\cF_\cx$, or
equivalently~$\cF\bigl( \barcO_{\cX,\cx,M}\tensor_V K\bigr)$, the
direct limit~$\lim_{i,j}\cF\bigl(\Spf(S_i),
(\Spm(S_{i,j}'\tensor_V K),L_{i,j})\bigr)$. We call it the stalk
of\/~$\cF$ at~$\cx$.

Denote by\/~$G_{\cx,M}$ the Galois group
of\/~$\barcO_{\cX,\cx}\tensor_V K$
over~$\cO_{\cX,\cx}^\sh\tensor_V M$. Then, $\cF_\cx$ is endowed
with an action of\/~$G_{\cx,M}$.

\slabel sh. lemma\par\slemma Let\/~$k(x)$ (resp.~$\kbar$) be the
residue field of\/~$V_x$ (resp.~$\Vbar$) and denote
by\/~$x_k\colon \Spec(k_x)\to X_k$ (resp.~$\barx_k\colon
\Spec(\kbar)\to X_k$) the points induced by~$x$ (resp.~$\barx$).
Then, \spacing

\item{{\rm i.}} $\cO_{X,x}^\sh$ coincides with the strict
henselization of\/~$\cO_{X,x_k}$ and\/~$\cO_{\cX,\cx}^\sh$
coincides with the strict formal henselization
of\/~$\cO_{\cX,x_k}$; \spacing

\item{{\rm ii.}} $\bigl(\cO_{X,x}^\sh,(p)\bigr)$
and\/~$\bigl(\cO_{\cX,\cx}^\sh,(p)\bigr)$ are noetherian henselian
pairs and the natural map $\cO_{X,x}^\sh\to \cO_{\cX,\cx}^\sh$ is
an isomorphism after taking $p$--adic completions;

\item{{\rm iii.}} the base change functor from the category of
finite extensions of\/~$\cO_{X,x}^\sh$, \'etale after
inverting~$p$, to the category of finite extensions
of\/~$\cO_{\cX,\cx}^\sh$, \'etale after inverting~$p$, is an
equivalence of categories;

\spacing

\item{{\rm iv.}} the maps  $\cO_{X,x}^\sh/p \cO_{X,x}^\sh \to
\cO_{\cX,\cx}^\sh/p \cO_{\cX,\cx}^\sh$ and\/ $\barcO_{X,x,M}/p
\barcO_{X,x,M}\to \barcO_{\cX,\cx,M}/p \barcO_{\cX,\cx,M}$ are
isomorphisms.

\spacing

\item{{\rm v.}} Frobenius on $ \barcO_{\cX,\cx,M}/p
\barcO_{\cX,\cx,M}$ is surjective with kernel $p^{1\over p}
\barcO_{\cX,\cx,M}/p \barcO_{\cX,\cx,M}$.

\endslemma
\Proof (i) The strict henselization of~$\cO_{X,x_k}$ is defined as
the direct limit~$\lim_j T_j$ over all pairs~$\{(T_j,t_j)\}_j$ where
$\Spec(T_j)$ is \'etale over~$X$ and~$t_j\colon T_j\to \kbar$ is a
point over~$\barx_k$. In particular, we get a map
$\cO_{X,x}^\sh=\lim_i R_i \to  \cO_{X,x_k}^\sh=\lim_j T_j$ by
associating to a pair~$\bigl(R_i,f_i\colon R_i\to \Vbar\bigr)$ the
pair~$\bigl(R_i, R_i\to \Vbar\to \kbar\bigr)$. To conclude that
such a map is an isomorphism it suffices to show that for any
pair~$(T_j,t_j)$ there is a unique map of $V$--algebras $T_j\to
\Vbar$ lifting~$t_j$ and inducing the point~$\barx$. The base
change of~$T_j$ via~$\barx$ defines an \'etale $\Vbar$--algebra
$A_j$ and $t_j$ induces a map of $\Vbar$--algebras $A_j \to
\kbar$. By \'etalness of~$A_j$ the latter lifts uniquely to a map
of $\Vbar$--algebras $A_i \to \widehat{\Vbar}$ which, since~$T_j$
is of finite type over~$V$, factors via~$\Vbar$.

The strict formal henselization of~$\cO_{\cX,X_k}$ is defined as the
direct limit~$\lim_j Q_j$ over all pairs~$\{Q_j,q_j\}_j$ where $Q_j$
is a $p$-adically complete and separated $V$--algebra, $\Spf(Q_j)\to
\cX$ is an \'etale map of formal schemes and~$q_j\colon Q_j\to
\kbar$ is a point over~$\barx_k$. The proof that~$\cO_{\cX,\cx}^\sh$
is the strict formal henselization of~$\cO_{\cX,X_k}$ is similar to
the first part of the proof and left to the reader.

(ii) It follows from~(i) that~$\cO_{X,x}^\sh$
(resp.~$\cO_{\cX,\cx}^\sh$) is a local ring with residue
field~$\kbar$ and maximal ideal~$\germ_x$ (resp.~$\germ_\cx$)
generated by the maximal ideal of~$\cO_{X,x_k}$. In particular,
the graded rings ${\rm gr}_{\germ_x}\cO_{X,x}^\sh$ and~${\rm
gr}_{\germ_\cx}\cO_{\cX,\cx}^\sh$ are noetherian so
that~$\cO_{X,x}^\sh$ and~$\cO_{\cX,\cx}^\sh$ are noetherian.

We claim that~$(\cO_{\cX,\cx}^\sh,\germ_\cx)$ is a henselian
pair; see~[\Elkik, \S0.1]. This amounts to prove that any \'etale
map~$\cO_{\cX,\cx}^\sh\to B$, such
that~$\kbar=\cO_{\cX,\cx}^\sh/\germ_\cx \cO_{\cX,\cx}^\sh\to
B/\germ_x B$ is an isomorphism, admits a section. Note that there
exists~$i$ and an \'etale extension $S_i\to A$ such that~$B$ is
obtained by base change of~$A$ via $S_i\to \cO_{\cX,\cx}^\sh$.
Via~$a\colon A/\germ_x A\to B/\germ_x B\cong \kbar$ the
pair~$(\widehat{A},a)$, where~$\widehat{A}$ is the $p$--adic
completion of~$A$, appears in the inductive system used to define
the strict formal henselization of~$\cO_{\cX,x_k}$ so that, thanks
to~(i), we get a natural map $A \to \cO_{\cX,\cx}^\sh$ and, hence
base--changing, a map of $\cO_{\cX,\cx}^\sh$--algebras $B \to
\cO_{\cX,\cx}^\sh$. Analogously, one proves
that~$(\cO_{X,x}^\sh,\germ_x)$ is a henselian pair.

Note that~$p$ is contained in~$\germ_x$, so
that~$\bigl(\cO_{X,x}^\sh,(p)\bigr)$
and\/~$\bigl(\cO_{\cX,\cx}^\sh,(p)\bigr)$ are henselian pairs.
Let~$\cO_{X\tensor \FF_p,x_k}^\sh$ be the strict henselization of
the local ring of~$X\tensor_V V/pV$ at~$x_k$. By construction we
have natural injective maps $\cO_{X,x}^\sh/p \cO_{X,x}^\sh\to
\cO_{\cX,\cx}^\sh/p \cO_{\cX,\cx}^\sh \to \cO_{X\tensor
\FF_p,x_k}^\sh$. We claim that such maps are isomorphisms. It
suffices to show that the composite is surjective. Using~(i) this
is equivalent to prove that the map from the strict henselization
of\/~$\cO_{X,x_k}$ to the strict henselization of\/~$\cO_{X\tensor
\FF_p,x_k}$ is surjective. This amounts to show that given an
\'etale map $f\colon \Spec(R) \to \cO_{X\tensor \FF_p,x_k}$, there
exists an \'etale map $g\colon \Spec(S)\to \cO_{X,x_k}$ reducing
to~$f$ modulo~$p$. By the Jacobian criterion of \'etalness we
have~$R=\cO_{X\tensor \FF_p,x_k}[T_1,\ldots,T_d]/(h_1,\ldots,h_d)$
with~$\det\left(\partial h_i/\partial T_j\right)_{i,j=1}^d$
invertible in~$R$. Then,
$S:=\cO_{X,x_k}[T_1,\ldots,T_d]/(q_1,\ldots,q_d)\left[\det\left(\partial
q_i/\partial T_j\right)^{-1}\right]$, with~$q_i$ lifting~$h_i$, is
an \'etale $\cO_{X,x_k}$--algebra and lifts~$R$  as wanted.
Since~$p$ is not a zero divisor in~$\cO_{X,x}^\sh$ and
in~$\cO_{\cX,\cx}^\sh$, we conclude that the graded rings ${\rm
gr}_p \cO_{X,x}^\sh$ and in~${\rm gr}_p\cO_{\cX,\cx}^\sh$ are
isomorphic, concluding the proof of~(ii).

(iii) Let~$\widehat{\cO_{X,x}^\sh}$
(resp.~$\widehat{\cO_{\cX,\cx}^\sh}$) be the $p$--adic completion
of~$\cO_{X,x}^\sh$ (resp.~$\cO_{\cX,\cx}^\sh$). Thanks to~[\Elkik,
Thm.~5] one knows that the category of finite extensions
of\/~$\cO_{X,x}^\sh$, \'etale after inverting~$p$ (resp.~the
category of finite extensions of\/~$\cO_{\cX,\cx}^\sh$, \'etale
after inverting~$p$), is equivalent to the category of finite
extensions of $\widehat{\cO_{X,x}^\sh}
=\widehat{\cO_{\cX,\cx}^\sh}$, \'etale after inverting~$p$. The
claim follows from~(ii).

(iv) The first claim follows from~(ii). The second follows from
the first and~(iii).

(v) Note that $p^\alpha\in \barcO_{\cX,\cx,M}$ for
every~$\alpha\in \QQ_{>0}$. It follows from~[\FALAST, \S 3, Lemma
5] that Frobenius is surjective on $\barcO_{\cX,\cx,M}/p^\alpha
\barcO_{\cX,\cx,M}$ for every~$0<\alpha<1$. Let~$a\in
\barcO_{\cX,\cx,M}$. Write $a=b^p+ p^{1\over p} c$ with~$b$
and~$c\in \barcO_{\cX,\cx,M}$. Write $c=d^p+ p^{1-{1\over p}} e$
with~$e\in \barcO_{\cX,\cx,M}$. Then, $a\equiv (b+p^{1\over p^2}
d)^p $ modulo~$p \barcO_{\cX,\cx,M}$.

Let~$a\in \barcO_{\cX,\cx,M}$ be such that~$a^p\in
p\barcO_{\cX,\cx,M}$. Then, ${a^p\over p}=\left({a\over p^{1\over
p}}\right)^p$ lies in~$\barcO_{\cX,\cx,M}$. Since the latter is a
normal ring, this implies that~${a\over p^{1\over p}}\in
\barcO_{\cX,\cx,M}$ as claimed.

\slabel stalks. sprop\par\slabel whyGalois. proposition\par\sprop
The notation is as above;\spacing

\item{{\rm 1)}} a sequence of sheaves $\cF\to \cG \to \cH$
on~$X^\et$ (resp.~$\gerX_M$, resp.~$\cX^\et$,  resp.~$\hatgerX_M$)
is exact if and only if for every point~$x$ of\/~$X$, defined over
a finite extension of~$K$, the induced sequence of stalks
$\cF_x\to \cG_x \to \cH_x$ (resp.~$\cF_\cx\to \cG_\cx \to
\cH_\cx$) is exact;\spacing

\item{{\rm 2)}} let $\cF$ be in $\Sh(X^\et)$. Then, $\nu_X^\ast(\cF)_\cx\cong
\cF_x$;\spacing

\item{{\rm 3)}} if\/~$\cF$ is in~$\Sh({\gerX}_M)$, then,
$\mu_{X,M}^\ast(\cF)_\cx\cong \cF_x$;

\spacing

\item{{\rm 4)}} fix field extensions $K\subset M_1 \subset M_2
\subset \Kbar$. Then, $\beta_{M_1,M_2}^\ast$
(resp.~$\gerbeta_{M_1,M_2}^\ast$) of a flasque sheaf is flasque.
Furthermore, if\/~$\cF$ is in~$\Sh\bigl(\gerX_{M_1}\bigr)$
(resp.~$\cF$ is in $\Sh\bigl(\hatgerX_{M_1}\bigr)$), then one has
natural identifications: \itemitem{{\rm a)}}
$\beta_{M_1,M_2}^\ast(\cF)_x=\cF_x$
(resp.~$\gerbeta_{M_1,M_2}^\ast(\cF)_\cx=\cF_\cx$); \itemitem{{\rm
b)}} if~$M_1\subset M_2$ is Galois with group~$G$, then
$\H^0(\gerX_{M_1},\cF)=\H^0\left(\gerX_{M_2},
\beta_{M_1,M_2}^\ast(\cF)\right)^{G}$
(resp.~$\H^0(\hatgerX_{M_1},\cF)=\H^0\left(\hatgerX_{M_2},
\gerbeta_{M_1,M_2}^\ast(\cF)\right)^G$).

\spacing \noindent Assume that\/~$K_x$ is contained in~$M$.
Then,\spacing

\item{{\rm 5)}} we have a natural isomorphisms $G_{\cx,M}\cong
G_{x,M}$ and, if\/~$\cF$ is in~${\gerX}_M$, the isomorphism
$\mu_{X,M}^\ast(\cF)_\cx\cong \cF_x$ is compatible with the actions
of\/~$G_{\cx,M}$ and\/~$G_{x,M}$

\spacing

\item{{\rm 6)}} let\/~$\cF$ be a sheaf in~$\gerX_M$. Then, $\left(\R^q v_{X,M,\ast}(\cF)\right)_x\cong
\H^q\bigl(G_{x,M},\cF_x\bigr)$;

\spacing

\item{{\rm 7)}} let\/~$\cF$ be a sheaf in~$\hatgerX_M$. Then, $\left(\R^q
\gerv_{\cX,M,\ast}(\cF)\right)_\cx\cong \H^q\bigl(G_{\cx,M},
\cF_\cx\bigr)$.

\endsprop
\Proof (1) In each case it suffices to prove that a sheaf is
trivial if and only if all its stalks are.

We give a proof for a sheaf on~$\gerX_M$ and leave the other cases
to the reader. Let~$\cF\in\Sh(\gerX_M)$ such that for every
point~$x$ of\/~$X$, defined over a finite extension of~$K$, we
have~$\cF_x=0$. Let~$(U,W)\in\gerX_M$ and let~$\alpha\in\cF(U,W)$.
Then, for every~$x\colon \Spec(V_x)\to U$ and every point~$y\colon
\Spec(K_y)\to W$ over~$x\tensor_V K$, which exists since~$W\to
U_M$ is finite, there exists~$(U_x,W_y)\in \gerX_M$ and a map
$\bigl(U_x,W_y\bigr)\to (U,W)$ such that (1) $x\tensor_V \Vbar$
factors via~$U_x$, (2) $y\tensor_K \Kbar$ factors via~$W_y$ and
(3) the image of~$\alpha$ in~$\cF\bigl(U_x,W_y\bigr)$ is~$0$.

Choose finitely many points~$x_1,\cdots,x_n$ and~$y_1,\ldots,y_n$
such that $\amalg_{i=1}^n(U_{x_i},W_{y_i})\to (U,W)$ is a covering
of~$(U,W)$ in~$\gerX_M$. Since~$\cF$ is a sheaf, the homomorphism
$\cF(U,W)\to \prod_{i=1}^n \cF\bigl(U_{x_i},W_{y_i}\bigr)$ is
injective. Hence, $\alpha=0$ to start with.

(2) Since any sheaf is the direct limit of representable sheaves
and direct limits commute with~$\nu_X^\ast$ and with taking
stalks, we may assume that~$\cF$ is represented by an \'etale
$X$--scheme $Y\to X$. In particular, $\nu_X^\ast(\cF)$ is
represented by the formal scheme~$\cY$ associated to~$Y$.
Let~$Y_x$ (resp.~$\cY_\cx$) be the pull back of~$Y$ (resp.~$\cY$)
to~$\Spec(\cO_{X,x})$ (resp.~$\Spec(\cO_{\cX,\cx})$). We then have
the following diagram
$$\matrix{ \cF_x & & \nu_X^\ast(\cF)_\cx\cr \Vert & & \Vert \cr
\Hom_{\cO_{X,x}}\left(\cO_{Y_x},\cO_{X,x}^\sh\right) &
\longrightarrow & \Hom_{\cO_{\cX,\cx}}
\left(\cO_{\cY_\cx},\cO_{\cX,\cx}^\sh\right) & \longrightarrow &
\Hom_k \bigl(\cO_{Y_x}\tensor_V k,\kbar\bigr).\cr}$$By~\ref{sh}(i)
these maps are bijective as claimed.

(4) If~$\cF$ is in~$\Sh(\gerX_{M_1})$,
then~$\beta_{M_1,M_2}^\ast(\cF)$ is the sheaf in~$\gerX_{M_2}$
associated to the presheaf $\beta_{M_1,M_2}^{-1}(\cF)$ defined
by~$(U,W)\mapsto\lim \cF(U',W')$ where the limit is taken over all
the pairs~$(U',W')$ in~$\gerX_{M_1}$ and all the maps~$(U,W)\to
(U',W'\tensor_{M_1} M_2)$. This is equivalent to take the direct
limit over all pairs~$(U,W')$ in~$\gerX_{M_1}$ and over all
map~$(U,W)\to (U,W')$ as $U_{M_1}$--schemes. If~$M_1 \subset M_2$ is
finite, there exists an initial pair, namely~$(U,W)$ itself, viewed
in~$\gerX_{M_1}$ via the finite and \'etale map $W\to U\tensor_V M_2
\to U\tensor_V M_1$, so that
$\beta_{M_1,M_2}^{-1}(\cF)(U,W)=\cF\bigl(U,W\bigr)$. In general,
there exists a finite extension~$M_1\subset L$ contained in~$M_2$
and a pair~$(U,W_L)$ in~$\gerX_L$ such that~$W=W_L\tensor_L M_2$.
Since any morphism of pairs in~$\gerX_{M_2}$ descends to a finite
extension of~$M_1$, we conclude
that~$\beta_{M_1,M_2}^{-1}(\cF)(U,W)=\cF\bigl(U,W_L\tensor_LM_2\bigr)$,
defined as the direct limit~$\lim_{L'}\cF\bigl(U,W_L\tensor_L
L'\bigr)$ taken over all finite extensions~$L\subset L'$ contained
in~$M_2$, considering~$\bigl(U,W_L\tensor_L L'\bigr)$
in~$\gerX_{M_1}$ via the finite and \'etale map $W_L\tensor_L L' \to
U\tensor_V L \to U\tensor_V K$. In any case, we conclude
that~$\beta_{M_1,M_2}^{-1}(\cF)$ is already a sheaf i.~e.,
$\beta_{M_1,M_2}^{-1}(\cF)=\beta_{M_1,M_2}^{\ast}(\cF)$.
Furthermore, $\beta_{M_1,M_2}^{\ast}$ preserves flasque objects and
satisfies (a).

For~(b), recall that~$\gerX_{M_1}$ and~$\gerX_{M_2}$ have final
objects so that global sections can be computed using the final
objects. Since $X_{M_2}\to X_{M_1}$ is a limit of finite and
\'etale covers with Galois group~$G$ and~$\cF$ is a sheaf
on~$\gerX_{M_1}$, one has~$\cF(X,X_{M_1})=\cF(X,X_{M_2})^{G}$.
Then,
$\H^0(\gerX_{M_1},\cF)=\cF(X,X_{M_1})=\cF(X,X_{M_2})^{G}=\H^0\left(\gerX_{M_2},
\beta_{M_1,M_2}^\ast(\cF)\right)^{G}$ and~(b) follows.

A similar argument works for~$\gerbeta_{M_1,M_2}^\ast$. Details are
left to the reader.

(3)\&(5)  The first claim in~(5) follows from~\ref{sh}(iii). To
get the second claim and~(3), one argues as in~(2) reducing to the
case of a sheaf represented by a pair~$(U,W)$, so
that~$\mu_{X,M}^\ast(U,W)=(\cU,\cW,L)$, and using~\ref{sh}(iii).

(6) Consider the functor $\Sh(\gerX_M)\to (G_{x,M}-\hbox{\rm
Modules})$, associating to a sheaf~$\cF$ its stalk~$\cF_x$. It is an
exact functor. Recall
that~$\cF_x=\lim_{i,j}\cF\bigl(\Spec(R_i),\Spec(R_{i,j}'\tensor_V
K)\bigr)$. Thus, the continuous Galois cohomology
$\H^\ast\bigl(G_{x,M}, \cF_x\bigr)$ is the direct limit over~$i$
and~$j$   of the Ch\v{e}ch cohomology of~$\cF$ relative to the
covering $\bigl(\Spec(R_i),\Spec(R_{i,j}'\tensor_V K)\bigr)$. In
particular, if~$\cF$ is injective, it is flasque and
$\H^q\bigl(G_{x,M}, \cF_x\bigr)=0$ for~$q\geq 1$.

Both~$\left\{\left(\R^q v_{\cX,M,\ast}(\cF)\right)_x\right\}_q$
and~$ \left\{\H^q\bigl(G_{x,M}, {\cF}_x\bigr)\right\}_q$ are
$\delta$--functors from~$\Sh(\gerX_M)$ to the category of abelian
groups. Also~$\left(\R^q v_{\cX,M,\ast}(\cF)\right)_x$ is zero
for~$q\geq 1$ and~$\cF$ injective. For~$q=0$ we have $$\left(\R^0
v_{\cX,M,\ast}(\cF)\right)_x=\lim_i
\cF\bigl(\Spec(R_i),\Spec(R_i\tensor_V
K)\bigr)=\cF\bigl(\barcO_{X,x,M}\tensor_V K\bigr)^{G_{x,M}}.$$This
proves the claim.

(7) The proof  is similar to the proof of~(6) and left to the
reader.

\slabel fincompareformaltoalgebraic. lemma\par\slemma We have the
following equivalences of $\delta$--functors  : \spacing
\item{{\rm i.}} $\R^q\bigl(\nu_X^\ast \circ v_{X,M,\ast}\bigr)=
\nu_X^\ast \circ \R^q v_{X,M,\ast}$ and $\R^q\bigl(
\gerv_{\cX,M,\ast}\circ\mu_{X,M}^\ast\bigr)=
\bigl(\R^q\gerv_{\cX,M,\ast}\bigr) \circ \mu_{X,M}^\ast$;

\spacing \item{{\rm ii.}} $\nu_X^\ast \circ\R^q v_{X,M,\ast}
\isomarrow \bigl(\R^q \gerv_{\cX,M,\ast}\bigr) \circ
\mu_{X,M}^\ast$.

\spacing

\endslemma
\Proof (i) Since~$\nu_X^\ast$ and~$\mu_{X,M}^\ast$ are exact
and~$v_{X,M,\ast}$ and~$\gerv_{\cX,M,\ast}$ are left exact, the
derived functors of~$\nu_X^\ast\circ v_{X,M,\ast}$
and~$\gerv_{\cX,M,\ast}\circ \mu_{X,M}^\ast$ exist. By
\ref{whyGalois} we have
$$\nu_X^\ast\left(\R^q
v_{X,M,\ast}(\cF)\right)_\cx\cong \left(\R^q
v_{X,M,\ast}\left(u_{X,M,\ast}(\cF)\right)\right)_x\cong
\H^q\left(G_{x,M},\cF_x\right)$$ and $$\R^q
\gerv_{\cX,M,\ast}\left(\mu_{X,M}^\ast(\cF)\right)_\cx\cong
\H^q\left(G_{\cx,M},\mu_{X,M}^\ast(\cF)_x\right).$$This implies
that if~$\cF$ is injective, $\nu_X^\ast\left(\R^q
v_{X,M,\ast}(\cF)\right)=0$ and $\R^q
\gerv_{\cX,M,\ast}\left(\mu_{X,M}^\ast(\cF)\right)=0$ for~$q\geq
1$. Hence, $\R^q\bigl(\nu_X^\ast \circ v_{X,M,\ast}\bigr)=
\nu_X^\ast \circ \R^p v_{X,M,\ast}$ (resp.~$\R^q\bigl(
\gerv_{\cX,M,\ast}\circ\mu_{X,M}^\ast\bigr)=
\R^q\gerv_{\cX,M,\ast} \circ \mu_{X,M}^\ast$). Indeed,  they are
both $\delta$--functors since~$\nu_X^\ast$
(resp.~$\mu_{X,M}^\ast$) is exact, they are both erasable and they
coincide for~$q=0$.

(ii) We construct a map $\gamma_\cF\colon\nu_X^\ast\left(
v_{X,M,\ast}(\cF)\right)\longrightarrow
\gerv_{\cX,M,\ast}\left(\mu_{X,M}^\ast(\cF)\right)$ functorial
in~$\cF$.  The sheaf~$\nu_X^\ast\left( v_{X,M,\ast}(\cF)\right)$
is the sheafification  of the presheaf $\cF_1$ which associates to
an object~$\cU$ in~$\cX^\et$ the direct limit $\lim \cF(U,U_K)$
taken over all~$ U\in X^\et$ and all maps from~$\cU$ to the formal
scheme associated to~$U$. On the other hand, the {\it presheaf}
$\cF_2:=\gerv_{\cX,M,\ast}\left(\mu_{X,M}^\ast(\cF)\right)$
($\mu_{X,M}^\ast(\cF)$ is taken as presheaf) associates
to~$\cU\in\cX^\et$ the direct limit $\lim \cF(U,W)$ over
all~$(U,W)$ in~$\gerX_M$ and all maps from~$\cU$ to the formal
scheme associated to~$U$ and from~$\cU^\rig$ to the rigid analytic
space defined by~$W^\rig\fibprod_{U^\rig} \cU^\rig$. We thus get a
morphism at the level of presheaves $\cF_1\to \cF_2$. Passing to
the associated sheaves we get the claimed map.

The map~$\gamma_\cF$ induces $\R^q \gamma_\cF\colon
\R^q\bigl(\nu_X^\ast \circ v_{X,M,\ast}\bigr)(\cF)\longrightarrow
\R^q\bigl( \gerv_{\cX,M,\ast}\circ\mu_{X,M}^\ast\bigr)(\cF)$.
Using (i), we get a natural transformation of $\delta$--functors as
claimed in~(ii). We are left to prove that it is an isomorphism.
This can be checked on stalks and, as explained in the proof
of~(i), it amounts to prove that for any sheaf~$\cF$ one has
$\H^q\left(G_{x,M},\cF_x\right)\cong
\H^q\left(G_{\cx,M},\mu_{X,M}^\ast(\cF)_x\right)$. The conclusion
follows since~$\mu_{X,M}^\ast(\cF)_\cx\cong \cF_x$
and~$G_{x,M}\cong G_{\cx,M}$ thanks to~\ref{whyGalois}.

\endssection

We next show that the sites introduced above are very useful in
order to compute \'etale cohomology:

\label compareGaltoXetI. proposition\par\prop {\rm
(Faltings)}\enspace Let $\bbL$ be a finite locally constant
\'etale sheaf on $X_M$ annihilated by~$p^s$. For every~$i$
 the map
$\H^i\left(\gerX_M,\bbL\right) \longrightarrow
\H^i\left(X_M^\et,\bbL\right)$, induced by $u_{X,M}$, is an
isomorphism. \endprop

\Proof {\rm [\FALAST, Rmk.~p.~242]}\enspace Put
$G_M:=\Gal(\Kbar/M)$. We have a spectral sequence
$$\H^p\bigl(G_M,\H^q(X_\Kbar^\et,\bbL)\Longrightarrow
\H^{p+q}\bigl(X_M^\et,\bbL\bigr)
$$and, thanks to~\ref{stalks}(4.b), a spectral sequence
$$\H^p\bigl(G_M,\H^q(\gerX_\Kbar,\bbL)\Longrightarrow
\H^{p+q}\bigl(\gerX_M,\bbL\bigr).$$Hence, it suffices to prove the
proposition for~$M=\Kbar$. Let $z_{X,\Kbar}\colon X^\et \to
X_\Kbar^\et$ be the map $U\to U\tensor_V \Kbar$. One knows from
[\FJAMS, Cor.~II.2.2] that $\Spec\left(\cO_{X,x}^\sh\tensor_V
\Kbar\right)$ is $K(\pi,1)$.  Hence, the stalk~$\bigl(\R^q
z_{X,\Kbar,\ast}(\bbL)\bigr)_x$  is
$\H^q\bigl(G_{x,\Kbar},\bbL_x\bigr)$. By~\ref{stalks} also the
stalk $\bigl(\R^q v_{X,\Kbar,\ast}(\bbL)\bigr)_x$ coincides
with~$\H^q\bigl(G_{x,\Kbar},\bbL_x\bigr)$. Hence, $\R^q
z_{X,\Kbar,\ast}(\bbL)\cong \R^q v_{X,\Kbar,\ast}(\bbL)$. Using
the spectral sequences $$\H^p\bigl(X^\et, \R^q
z_{X,\Kbar,\ast}(\bbL)\bigr)\Longrightarrow
\H^{p+q}\left(X_\Kbar^\et,\bbL\right),\quad \H^p\bigl(X^\et, \R^q
v_{X,\Kbar,\ast}(\bbL)\bigr)\Longrightarrow
\H^{p+q}\left(\gerX_\Kbar^\et,\bbL\right)$$the proposition
follows.

\label compforandalg. section\par\ssection Comparison between
algebraic cohomology and formal cohomology\par Since~$\nu_X^\ast$
is an exact functor, given an injective resolution~$I^\bullet$
of~$\cF$, then~$0\to \nu_X^\ast(\cF)\to \nu_X^\ast(I^\bullet)$ is
exact so that given an injective resolution~$J^\bullet$
of~$\cF^\form=\nu_X^\ast(\cF)$ we can extend the identity map
on~$\cF$ to a morphism of complexes~$\nu_X^\ast(I^\bullet)\to
J^\bullet$. Since~$\nu_X$ sends the final object~$X$ of~$X^\et$ to
the final object~$\cX$ of~$\cX^\et$, one has a natural
map~$I^\bullet(X)\to \nu_X^\ast(I^\bullet)(\cX)$. Then,

\label fincomparisonformaalgebraiccoho. definition\par\sdefi One has
natural maps of $\delta$--functors $$\rho_{X,\cX}^q(\cF)\colon
\H^q\bigl(X^\et, \cF \bigr)\to \H^q\bigl(\cX^\et,\cF^\form \bigr),$$

$$\rho_{\gerX_M,\hatgerX_M}^q(\cF)\colon \H^q\bigl(\gerX_M, \cF
\bigr)\to \H^q\bigl(\hatgerX_M,\cF^\rig \bigr)$$
\endsdefi

\noindent Note that one has spectral sequences \labelf
specHqv\par$$\H^q\bigl(X^\et, \R^p
v_{X,M,\ast}(\cF)\bigr)\Longrightarrow
\H^{p+q}(\gerX_M,\cF),\eqno{{(\numfo)}}$$\advance\fonu
by1\noindent and\labelf specHqnu\par
$$\H^q\bigl(\cX^\et,\nu_X^\ast \R^p
v_{\cX,M,\ast}(\cF)\bigr)=  \H^q\bigl(\cX^\et, \R^p
\gerv_{\cX,M,\ast}\bigl(\cF^\rig\bigr)\bigr)\Longrightarrow
\H^{p+q}\bigl(\hatgerX_M,\cF^\rig\bigr).\eqno{{(\numfo)}}$$\advance\fonu
by1\noindent where the equality on the left hand side is due
to~\ref{fincompareformaltoalgebraic}.

\slabel finGabber. proposition\par\sprop  The following hold:\spacing

\item{{\rm a.}} If\/~$\cF$ in~$\Sh(X^\et)$ is torsion, the map
$\rho_{X,\cX}^q(\cF)$ is an isomorphism. \spacing\item{{\rm b.}}
the spectral sequences (\ref{specHqv}) and\/~(\ref{specHqnu})  are
compatible via $\rho_{X,\cX}^q$
and~$\rho_{\gerX_M,\hatgerX_M}^{p+q}$; \spacing \item{{\rm c.}}
if\/~$\cF$ is a torsion sheaf on~$\gerX_M$, the map
$\rho_{\gerX_M,\hatgerX_M}^q(\cF)$ is an isomorphism.

\endsprop\Proof (a) Let~$X_k:=X\tensor_V k$ and denote
by~$\iota\colon X^\et \to X_k^\et$ and~$\widehat{\iota}\colon
\cX^\et \to X_k^\et$ the morphisms of topologies induced by the
closed immersions $X_k\subset X$ and~$X_k\subset \cX$
respectively. In fact, $\widehat{\iota}$ is an equivalence, since
the \'etale sites of~$\cX$ and of~$X_k$ coincide, and
$\widehat{\iota} \circ \nu_X=\iota$. For any sheaf~$\cF$
on~$X^\et$ denote~$\cF_k:=\iota^\ast(\cF)$ or, equivalently,
$\widehat{\iota}^\ast\bigl(\cF^\form\bigr)$. We then have
$\H^q\bigl(X^\et, \cF \bigr)\to \H^q\bigl(\cX^\et,\cF^\form \bigr)
\cong \H^q\bigl(X_k^\et,\cF_k \bigr)$ where the first map
is~$\rho_\cF^q$. The composite is defined by restriction and is an
isomorphism if~$\cF$ is a torsion sheaf due to~[\Gabber, Cor.~1]
and the fact that~$X$ is proper over~$V$. The conclusion follows.

(b) left to the reader.

(c) The left hand side of the spectral sequences~(\ref{specHqv})
and\/~(\ref{specHqnu}) are isomorphic by~(a) since~$\R^p
v_{X,M,\ast}$ sends a torsion sheaf to a torsion sheaf. The
conclusion follows from~(b).

\slabel comparegerXKbartohatgerXKabrI. corollary\par\scor Let
$\bbL$ be a locally constant  sheaf on~$\gerX_M$ annihilated
by~$p^s$. Then, the two sides of the Leray spectral sequences
$$ \H^j\Bigl(\cX^\et, \R^i
\gerv_{\cX,M,\ast}\bigl(\bbL^\rig\bigr) \Bigr)
 \Longrightarrow \H^{i+j}\left(\hatgerX_M,
\bbL^\rig\right)$$and
$$\H^j\Bigl( X^\et, \R^i
v_{X,M,\ast}\bigl(\bbL \bigr)\Bigr) \Longrightarrow
\H^{i+j}\left(\gerX_M,\bbL\right),$$obtained from the morphisms of
topoi $\gerv_{\cX,M}\colon \hatgerX_M \to \cX^\et$ and
$v_{X,M}\colon \gerX_M\to X^\et$, are naturally isomorphic.
\endscor \Proof The statements follow from~\ref{finGabber}.

\endssection

\label proof1. section\par\ssection The proof of Theorem
\ref{etalePhiGamma}\par The proofs of i) and ii) are very similar
therefore let us prove i). It follows from~\ref{maintheorem} that
$\cH^{i,\ast,\geom}_{\Kbar}(\bbL)$ is the sheaf associated to the
functor: $\cU=\Spf(R_{\cU})\lra
\H^i\bigl(\pi_1^{alg}(\cU^\rig_{\Kbar},\eta), {\bf L}\bigr)$
where~$\cU\to \cX$ run through the \'etale maps such
that~$\cU=\Spf(R_\cU)$ is affine and~$R_\cU$ satisfies the
assumptions of~\ref{V}. Therefore its stalk at a point~$\cx$
of~$\cX$ can be described as follows:
$$
\bigl(\cH^{i,\ast,\geom}_{\Kbar}(\bbL)\bigr)_{\cx}\cong
\H^i(\G_{\cx,\Kbar}, \bigl(\bbL^\rig\bigr)_\cx).
$$
On the other hand for every small $\cU=\Spf(R_\cU)$ we have a natural, functorial map
$$
\alpha^i_\cU\colon\H^i(\pi_1^{alg}(\cU^\rig_{\Kbar},\eta), {\bf
L})\lra \R^i\gerv_{\cX,\Kbar,\ast}(\bbL^\rig)(\cU),
$$
which induces a morphism of sheaves on $\cX^\et$
$$
\alpha^i\colon\cH^{i,\ast,\geom}_{\Kbar}(\bbL)\lra
\R^i\gerv_{\cX,\Kbar,\ast}(\bbL^\rig).
$$
For every point $\cx$ of $\cX$ the map $\alpha^i_\cx$ induced by $\alpha_i$ on stalks
is the canonical morphism
$$
\alpha^i_\cx:\bigl(\cH^{i,\geom}_\Kbar\bigr)_\cx\cong \H^i(\G_{\cx,\Kbar},
\bigl(\bbL^\rig\bigr)_\cx)\lra
\bigl(\R^i\gerv_{\cX,\Kbar,\ast}(\bbL^\rig)\bigr)_\cx,
$$
which by proposition \ref{stalks}(7) is an isomorphism. Therefore,
$\alpha^i$ induces an isomorphism of sheaves on $\cX^\et$:
$\cH^{i,\ast,\geom}_\Kbar(\bbL)\cong
\R^i\gerv_{\cX,\Kbar,\ast}(\bbL^\rig)$ and, thus, the Leray
spectral sequence produces the spectral sequence
$$
E_2^{p,q}=\H^q(\cX^\et,
\cH^{p,\ast,\geom}_\Kbar(\bbL)\Longrightarrow
\H^{p+q}(\hatgerX_\Kbar,\bbL^\rig).
$$
By~\ref{comparegerXKbartohatgerXKabrI}
$$
\H^{p+q}(\hatgerX_\Kbar,\bbL^\rig)\cong \H^{p+q}(\gerX_\Kbar, \bbL)
$$
and by~\ref{compareGaltoXetI}
$$
\H^{p+q}(\hatgerX_\Kbar,
\bbL)\cong \H^{p+q}(\bigl(X_\Kbar\bigr)^\et, \bbL).
$$
All these isomorphisms are equivariant for the residual action of $\G_V$.
This proves the claim.

\endssection

\endsection

\section A geometric interpretation of classical $(\varphi,\Gamma)$--modules\par

\noindent Let the notations be as in the previous section and fix
as before $M$ an algebraic extension of $K$ contained in $\Kbar$.
In this section we work with continuous sheaves on all our
topologies (see \S 4). We define  families of continuous sheaves
denoted $\barcO_{\gerX_M}$, $\cR\bigl(\barcO_{\gerX_M}\bigr)$,
$A_{\rm inf}^+\bigl(\barcO_{\gerX_M}\bigr)$ and call them {\sl
algebraic Fontaine sheaves} on $\gerX_M$ (respectively
$\barhatcO_{\hatgerX_M}, \cR\left(\barhatcO_{\hatgerX_M}\right)$,
$A_{\rm inf}^+\left(\barhatcO_{\hatgerX_M}\right)$ called {\sl
analytic Fontaine sheaves} on $\hatgerX_M$) and study their
properties. In this section we compare the cohomology on $\gerX_M$
of an \'etale local system $\bbL$ of $\ZZ/p^s\ZZ$-modules on $X_K$
tensored by one of the algebraic Fontaine sheaves with the
cohomology on $\hatgerX_M$ of its analytic analogue. As a
consequence we derive the following result.

\noindent Let us fix $M=K_\infty=K(\mu_{p^\infty})$ and consider
the following sheaf  $\cF_\infty:=\bbL^\rig\otimes A_{\rm
inf}^+\bigl(\barhatcO_{\hatgerX_\Kinfty}\bigr)$ on
$\hatgerX_{K\infty}$.

\label geometricPhiGamma. theorem\par\thm We have natural
isomorphisms of classical $(\phi,\Gamma)$--modules
$$
\H^i(\hatgerX_{K_\infty}, \cF_\infty)\cong \cDtilde_V\Bigl(
\H^i(\bigl(X_\Kbar\bigr)^\et, \bbL)\Bigr),
$$
for all $i\ge 0$.
\endthm

\noindent
The proof of theorem \ref{geometricPhiGamma} will take the rest of
this section.

\label topologicalsheaves. section\par\ssection Categories of
inverse systems\par We review some of the results of~[\Jannsen]
which will be needed in the sequel. Let~$\cA$ be an abelian
category. Denote by $\cA^\NN$ the category of inverse systems
indexed by the set of natural numbers. Objects are inverse systems
$\{A_n\}_n:=\ldots \to A_{n+1}\to A_n \ldots A_2 \to A_1$, where
the $A_i$'s are objects of~$\cA$ and the arrows denote morphisms
in~$\cA$. The morphisms in $\cA^\NN$ are commutative diagrams
$$\matrix{\ldots & \to & A_{n+1} & \to & A_n & \ldots &
A_2  &\to & A_1\cr  & & \downarrow & & \downarrow & & \downarrow &
& \downarrow & \cr \ldots & \to & B_{n+1} & \to & B_n & \ldots &
B_2 &\to & B_1,}$$where the vertical arrows are morphisms
in~$\cA$. Then, $\cA^\NN$ is an abelian category with kernels and
cokernels taken componentwise and if~$\cA$ has enough injectives,
then~$\cA^\NN$ also has enough injectives; see~[\Jannsen, Prop.~1.1].
Furthermore, there is a fully faithful and exact functor $\cA\to
\cA^\NN$ sending an object~$A$ of~$\cA$ to the inverse system
$\{A\}_n:=\ldots \to A\to A \ldots  \to A$ with transition maps
given by the identity and a morphism $f\colon A\to B$ of~$\cA$ to
the map of inverse systems $\{A\}_n\to \{B\}_n$ defined by~$f$ on
each component. By~[\Jannsen, Prop.~1.1] such map preserves
injective objects.

Let~$h\colon \cA \to \cB$ be a left exact functor of abelian
categories. It induces a left exact functor $h^\NN\colon \cA^\NN\to
\cB^\NN$ which, by abuse of notation and if no confusion is
possible, we denote again by~$h$. If~$\cA$ has enough injectives,
one can derive the functor~$h^\NN$. It is proven in~[\Jannsen,
Prop.~1.2] that~$\R^i \bigl(h^\NN\bigr)=\bigl(\R^i h\bigr)^\NN$.

If inverse limits over~$\NN$ exist in~$\cB$, define the left exact
functor $\displaystyle\lim_\leftarrow h\colon \cA^\NN\to \cB$ as
the composite of~$h^\NN$ and the inverse limit functor
$\displaystyle\lim_\leftarrow\colon \cB^\NN \to \cB$. Assume
that~$\cA$ and~$\cB$ have enough injectives. For
every~$A=\{A_n\}_n\in \cA^\NN$ one then has a spectral sequence
$${\lim_\leftarrow}^{(p)} \R^q h (A_n) \Longrightarrow \R^{p+q}
\bigl(\lim_{\leftarrow} h(A)\bigr),$$where
$\displaystyle{{\lim_\leftarrow}}^{(p)}$ is the $p$--th derived
functor of~$\displaystyle\lim_\leftarrow$ in~$\cB$. If in~$\cB$
infinite products exist and are exact functors, then
$\displaystyle{{\lim_\leftarrow}}^{(p)}=0$ for~$p\geq 1$ and the
above spectral sequence reduces to the simpler exact sequence
\labelf inversecohcohinverse\par$$0 \longrightarrow
{\lim_\leftarrow}^{(1)} \R^{i-1} h(A_n)\longrightarrow
\R^i\bigl(\lim_\leftarrow h) (A) \to \lim_\leftarrow \R^i h(A_n)
\longrightarrow 0;\eqno{{(\numfo)}}$$\advance\fonu by1\noindent
see~[\Jannsen, Prop.~1.6].

Note that via the map $\cA\to \cA^\NN$ given above, if~$A\in\cA$
then $\R^i h^\NN \bigl(\{A\}_n\bigr)=\{\R^ih (A)\}_n$
and~$\displaystyle\R^i\lim_\leftarrow h \left(\{A\}_n\right)=\R^i
h(A)$.
\endssection

\label groupscontcoho. section\par\ssection Example {\rm
[\Jannsen, \S 2]}\par Let~$G$ be a profinite group. Let~$\cA$ be
the category of discrete modules with continuous action of~$G$ and
let~$\cB$ be the category of abelian groups. For every~$j$ let $
\H^j\bigl(G,\_\bigr)\colon \cA^\NN\to \cB$ be the $j$--th derived
functor of~$\displaystyle\lim_\leftarrow \H^0(G,\_)$ on~$\cA^\NN$.
By loc.~cit.~for every inverse system~$T=\{T_n\}_n\in\cA^\NN$
we have an exact sequence \labelf lim1\par$$ 0 \longrightarrow
{\lim_\leftarrow}^{(1)} \H^{j-1}\bigl(G,T_n\bigr) \longrightarrow
\H^j\bigl(G,T\bigr) \longrightarrow \lim_\leftarrow
\H^j\bigl(G,T_n\bigr)\longrightarrow 0.\eqno{{(\numfo)}}
$$\advance\fonu by1\noindent Moreover given~$\{(N_n,d_n)\}_n\in \cB^\NN$,
one computes $\displaystyle{\lim_\leftarrow}^{(1)} N_n$ as the
cokernel of the map \labelf lim1ascokernel\par$$\prod_n({\rm
Id}-d_n)\colon \prod_n N_n \llongrightarrow \prod_n
N_n.\eqno{{(\numfo)}}$$\advance\fonu by1\noindent For later use we
remark the following. Assume that each $N_n$ is a module over a
ring~$C$ and that $d_n\colon N_{n+1}\to N_n$ is a homomorphism of
$C$--modules. Suppose that for every~$n$ there exists an
element~$c_n\in C$ annihilating the cokernel of~$d_n$. One then
proves by induction on~$m\in\NN$ that the cokernel of
$\prod_n({\rm Id}-d_n)\colon \prod_{n\leq m} N_n \to \prod_{n\leq
m} N_n$ is annihilated by~$c_1\cdots c_m$. In particular, if~$C$
is a complete local domain  and $c_n=c^{1\over p^n}\in C$ for every~$n$
with~$c$ in the maximal ideal of~$C$, then $c^{p\over p-1}=\prod_m
c_m$ annihilates $\displaystyle{\lim_\leftarrow}^{(1)} N_n$.

For every~$\{(T_n,d_n)\}\in \cA^\NN$ one
defines~$\displaystyle\H^j_\cont\left(G,\lim_{\infty\leftarrow n}
T_n\right)$ as the continuous cohomology defined by continuous
cochains modulo continuous coboundaries with values in
$\displaystyle\lim_{\infty\leftarrow n} T_n$ endowed with the
inverse limit topology considering on each~$T_n$ the discrete
topology. As explained in~[\Jannsen, Pf.~of Thm.~2.2] there exists
a canonical complex~$D^\bullet(G,T_n)$ whose $G$--invariants
define the continuous cochains~$C^\bullet(G,T_n)$ of~$G$ with
values in~$T_n$ and such that each~$D^i(G,T_n)$ is $G$--acyclic.
This resolution is functorial so that we get a resolution
$$(T_n,d_n) \subset (D^1(G,T_n),d^1_n) \to (D^2(T,T_n),d^2_n) \to
\cdots.$$The continuous cohomology
$\displaystyle\H^i_\cont\bigl(G,\lim_{\infty\leftarrow n}
T_n\bigr)$ is obtained by
applying~$\displaystyle\lim_{\infty\leftarrow n} \H^0(G,\_)$ to
this resolution and taking homology. Due to~(\ref{lim1}),
since~$D^i(G,T_n)$ is $G$--acyclic, we have \labelf
Hicont\par$$\H^i_\cont\bigl(G,(D^j(G,T_n),d^1_n) \bigr)=\cases{ 0
& if~$i\geq 2$\cr \lim_\leftarrow^{(1)} C^j(G,T_n) & for~$i=1$\cr
\lim_{\infty \leftarrow n} C^j(G,T_n)& if
$i=0$.\cr}\eqno{{(\numfo)}}$$\advance\fonu by1\noindent In
particular, if the system $\{T_n\}_n$ is Mittag--Leffler, then
$(D^j(G,T_n),d^1_n)$ is acyclic for every~$j$ and we obtain
$$\H^i_\cont\bigl(G,\lim_{\infty\leftarrow n} T_n\bigr)\isomarrow
\H^i(G,T).$$Next, assume as before that there exists a complete
local domain~$C$ such that~$T_n$ is a $C$--module and $d_n$ is a
homomorphism of $C$--modules. Suppose also that there is~$c$ in
the maximal ideal of~$C$ such that~$c^{1\over p^n}\in C$
and~$c^{1\over p^n}$ annihilates the cokernel of~$d_n$. Then,
$c^{1\over p^n}$ annihilates also the cokernel of
$C^i(G,T_{n+1})\to C^i(G,T_n)$ so that~$c^{p\over p-1}$
annihilates $\H^1_\cont\bigl(G,(D^j(G,T_n),d^1_n) \bigr)$. This
implies that if we invert~$c$ we have an isomorphism \labelf
Hicontiso\par $$ \H^i_\cont\bigl(G,\lim_{\infty\leftarrow n}
T_n\bigr)\bigl[c^{-1}\bigr]\isomarrow
\H^i(G,T)\bigl[c^{-1}\bigr].\eqno{{(\numfo)}}$$

\endssection

\ssection Fontaine sheaves on $\gerX_M$ and $\hatgerX_M$\par  We now come
to the definition of a family of sheaves on $\gerX_M$ and $\hatgerX_M$ which
will play a crucial role in the sequel. See \ref{etaletheorem}.

\slabel Obar. definition\par \sdefi {\rm [\FALAST,
p.~219--221]}\enspace Let\/  $\barcO_{\gerX_M}$ be the sheaf of
rings on ${\gerX}_M$ defined requiring that for every object
$(U,W)$ in~$\gerX_M$ the ring $\barcO_{\gerX_M}(U,W)$ consists of
the normalization of $\Gamma\bigl(U,\cO_U\bigr)$
in~$\Gamma\bigl(W,\cO_W\bigr)$.

Denote by\/~$\cR\bigl(\barcO_{\gerX_M}\bigr)$ the sheaf of rings
in\/~$\Sh(\gerX_M)^\NN$ given by the inverse system $
\left\{\barcO_{\gerX_M}/p\barcO_{\gerX_M}\right\}$, where the
transition maps are given by Frobenius.

For every~$s\in\NN$ define the sheaf of rings $A_{\rm
inf,s}^+\left(\barcO_{\gerX_M}\right)$ in\/~$\Sh(\gerX_M)^\NN$ as
the inverse
system~$\left\{\WW_s\bigl(\barcO_{\gerX_M}/p\barcO_{\gerX_M}\bigr)\right\}$.
Here, $\WW_s\bigl(\barcO_{\gerX_M}/p\barcO_{\gerX_M}\bigr)$ is the
sheaf~$\bigl(\barcO_{\gerX_M}/p\barcO_{\gerX_M}\big)^s$ with ring
operations defined by Witt polynomials and the transition maps in
the inverse system are defined by Frobenius. Define $A_{\rm
inf}^+\left(\barcO_{\gerX_M}\right)$ to be the inverse system of
sheaves
$\left\{\WW_n\bigl(\barcO_{\gerX_M}/p\barcO_{\gerX_M}\bigr)\right\}_n$
where the transition maps are  defined as the composite of the
projection
$\WW_{n+1}\bigl(\barcO_{\gerX_M}/p\barcO_{\gerX_M}\bigr) \to
\WW_n\bigl(\barcO_{\gerX_M}/p\barcO_{\gerX_M}\bigr)$ and Frobenius
on $\WW_n\bigl(\barcO_{\gerX_M}/p\barcO_{\gerX_M}\bigr)$.

\spacing Similarly, $\barhatcO_{\hatgerX_M}$ is the sheaf of rings
on ${\hatgerX}_M$ associating to an object $(\cU,\cW,L)$
in~$\hatgerX_M$ the ring $\barhatcO(\cU,\cW)$ defined as the
normalization of\/ $\Gamma\bigl(\cU,\cO_\cU\bigr)$
in~$\Gamma\bigl(\cW,\cO_\cW\bigr)\tensor_L M$.

Let\/~$\cR\left(\barhatcO_{\hatgerX_M}\right)$ be the sheaf of
rings in~$\Sh(\hatgerX_M)^\NN$ given by
$\left\{\bigl(\barhatcO_{\hatgerX_M}/p\barhatcO_{\hatgerX_M}\bigr)\right\}$,
where the inverse system is taken using Frobenius as transition
map.

For~$s\in\NN$ define the sheaf of rings $A_{\rm
inf,s}^+\left(\barhatcO_{\hatgerX_M}\right)$
in~$\Sh(\hatgerX_M)^\NN$ as the inverse system
$\left\{\WW_s\bigl(\barhatcO_{\hatgerX_M}/p\barhatcO_{\hatgerX_M}\bigr)\right\}$
with transition maps given by Frobenius. Eventually, let $A_{\rm
inf}^+\left(\barhatcO_{\hatgerX_M}\right)$
in~$\Sh(\hatgerX_M)^\NN$ be the sheaf
$\left\{\WW_n\bigl(\barhatcO_{\hatgerX_M}/p\barhatcO_{\hatgerX_M}\bigr)\right\}$
where the transition maps are  defined as the composite
$$\WW_{n+1}\bigl(\barhatcO_{\hatgerX_M}/p\barhatcO_{\hatgerX_M}\bigr)
\longrightarrow
\WW_n\bigl(\barhatcO_{\hatgerX_M}/p\barhatcO_{\hatgerX_M}\bigr)\longrightarrow
\WW_n\bigl(\barhatcO_{\hatgerX_M}/p\barhatcO_{\hatgerX_M}\bigr),$$
where the first map is the natural projection and the second is
Forbenius.

\spacing We denote by~$\varphi$ the Frobenius operator acting on
the sheaves, or inverse systems of sheaves, introduced above.

\endsdefi

\label AinRObarKbar. remark\par\srmk Note that if~$X=V$
and~$M=\Kbar$, one has
$\H^0_\cont\bigl((V,\Kbar),\cR(\barcO_\Kbar)\bigr)=\EEtilde_\Vbar^+$,
$\H^0_\cont\left((V,\Kbar),A_{\rm
inf,s}^+\bigl(\barcO_\Kbar\bigr)\right)=\WW_s\bigl(\EEtilde_\Vbar^+\bigr)$
and $\H^0_\cont\left((V,\Kbar),A_{\rm
inf}^+\bigl(\barcO_\Kbar\bigr)\right)=\AAtilde_\Vbar^+$.

For later use, we recall that we  denote by~$\pi$ the
element~$[\varepsilon]-1$ of\/~$\AAtilde_\Vbar^+$
where~$\varepsilon$ is the
element~$(1,\zeta_p,\zeta_{p^2},\cdots)\in \EE_V^+$
and~$[\varepsilon]$ is its Teichm\"uller lift.
\endsrmk

\sssection Notation\par If~$\cF$ is in~$\Sh(\gerX_M)^\NN$
(resp.~$\Sh(\hatgerX_M)$) write~$\H^i_\cont(\gerX_M,\cF)$
(respectively $\H^i_\cont(\hatgerX_M,\cF)$) for the $i$--th
derived functor $\displaystyle \lim_\leftarrow \H^0(\gerX_M,\_)$
(resp.~$\displaystyle \lim_\leftarrow \H^0(\hatgerX_M,\_) $)
applied to~$\cF$. Note that if~$\cF=\{\cG\}_n$ with~$\cG\in
\Sh(\gerX_M)$ (resp.~in~$\Sh(\hatgerX_M)$), then
$\H^i_\cont(\gerX_M,\cF)=\H^i(\gerX_M,\cG)$
(resp.~$\H^i_\cont(\hatgerX_M,\cF)=\H^i(\hatgerX_M,\cG)$).

\endsssection

\slabel compareAinfhatoverlineandoverline. lemma\par\slemma One
has $A_{\rm inf,\ast}^+\bigl(
\barcO_{\gerX_M}\bigr)^\rig\isomarrow A_{\rm
inf,\ast}^+\bigl(\barhatcO_{\hatgerX_M}\bigr)$ where
$\ast=s\in\NN$ or $\ast=\emptyset$.
\endslemma
\Proof Consider a  pair $(U,W)$ in~$\gerX_M$, with~$W$ defined
over some finite extension $K\subset L$ contained in~$M$. Recall
from section 4 that~$\mu_{X,M}(U,W):=(\cU,\cW,L)$. We have a
natural map $\barcO_{\gerX_M}(U,W) \to
\mu_{X,M,\ast}\left(\barhatcO_{\hatgerX_M}\right)(U,W)$ from the
normalization of~$\Gamma(U,\cO_U)$ in~$\Gamma(W,\cO_W)\tensor_L M$
to the normalization of~$\Gamma(\cU,\cO_\cU)$
in~$\Gamma(\cW,\cO_\cW)\tensor_L M$. This induces a natural map
$\barcO_{\gerX_M} \to \mu_{X,M,\ast}\left(\barhatcO\right)$ and,
hence, a map
$\mu_{X,M}^\ast\left(\mu_{X,M,\ast}\left(\barhatcO_{\hatgerX_M}
\right)\right)\longrightarrow \barhatcO_{\hatgerX_M}$, coming from
adjunction of~$\mu_{X,M,\ast}$ and~$\mu_{X,M}^\ast$. We then get a
homomorphism
$$\mu_{X,M}^{\ast,\NN}\left(A_{\rm inf,\ast}^+\bigl(
\barcO_{\gerX_M}\bigr)\right)\to A_{\rm
inf,\ast}^+\bigl(\barhatcO_{\hatgerX_M}\bigr).$$We claim that
these maps are isomorphisms. It suffices to prove it componentwise
and by devissage it is enough to show that
$\mu_{X,M}^\ast\left(\mu_{X,M,\ast}\left(\barhatcO_{\hatgerX_M}/p
\barhatcO_{\hatgerX_M} \right)\right)\longrightarrow
\barhatcO_{\hatgerX_M}/p\barhatcO_{\hatgerX_M} $ is an
isomorphism. Due to~\ref{stalks}(3) this amounts to prove that,
for every point~$x$ of~$X$ as in~\ref{definegeopoints}, the
natural map $\barcO_{X,x,M}/p \barcO_{X,x,M}\to
\barhatcO_{\cX,\cx,M}/p \barhatcO_{\cX,\cx,M}$ is an isomorphism.
This follows from~\ref{sh}(iv).

\slabel compareformaltoalgebraic. lemma\par\slemma We have the
following equivalences of $\delta$--functors  : \spacing
\item{{\rm i.}} $\R^q\bigl(\nu_X^{\ast,\NN} \circ v_{X,M,\ast}^\NN\bigr)=
\nu_X^{\ast,\NN} \circ \R^p v_{X,M,\ast}^\NN$ and $\R^q\bigl(
\gerv_{\cX,M,\ast}^\NN\circ\mu_{X,M}^{\ast,\NN}\bigr)=
\bigl(\R^q\gerv_{\cX,M,\ast}^\NN\bigr) \circ \mu_{X,M}^{\ast,\NN}$;

\spacing \item{{\rm ii.}} $\nu_X^{\ast,\NN} \circ\R^q v_{X,M,\ast}^\NN
\isomarrow \bigl(\R^q \gerv_{\cX,M,\ast}^\NN\bigr) \circ
\mu_{X,M}^{\ast,\NN}$.

\spacing

\Proof The result follows for lemma
\ref{fincompareformaltoalgebraic} and~\ref{topologicalsheaves}.

\endssection

\label contcompforandalg. section\par\ssection Comparison between
algebraic and formal cohomology of continuous sheaves\par Since
$\nu_X^{\ast,\NN}$ is an exact functor, as in section
\ref{compforandalg}, given an injective resolution~$I^\bullet$ of
a continuous sheaf~$\cF$, then~$0\to \nu_X^\ast(\cF)\to
\nu_X^{\ast,\NN}(I^\bullet)$ is exact so that given an injective
resolution~$J^\bullet$ of~$\cF^\form=\nu_X^{\ast,\NN}(\cF)$ we can
extend the identity map on~$\cF$ to a morphism of
complexes~$\nu_X^{\ast,\NN}(I^\bullet)\to J^\bullet$.
Since~$\nu_X$ sends the final object~$X$ of~$X^\et$ to the final
object~$\cX$ of~$\cX^\et$, one has a natural map~$I^\bullet(X)\to
\nu_X^{\ast,\NN}(I^\bullet)(\cX)$. Then,

\label comparisonformaalgebraiccoho. definition\par\sdefi One has
natural maps of $\delta$--functors
$$\rho_{X,\cX}^{\cont,q}(\cF)\colon \H^q_\cont\bigl(X^\et, \cF
\bigr)\to \H^q_\cont\bigl(\cX^\et,\cF^\form \bigr),$$ and $$
\rho_{\gerX_M,\hatgerX_M}^{\cont,q}(\cF)\colon
\H^q_\cont\bigl(\gerX_M, \cF\bigr)\to
\H^q_\cont\bigl(\hatgerX_M,\cF^\rig \bigr).$$
\endsdefi

\noindent Note that one has spectral sequences \labelf
specHqvNN\par$$ \H^q_\cont\bigl(X^\et,\R^p
v_{X,M,\ast}^\NN(\cF)\bigr)\Longrightarrow
\H^{p+q}_\cont(\gerX_M,\cF),\eqno{{(\numfo)}}$$\advance\fonu
by1\noindent and\labelf specHqnuNN\par
$$\H^q_\cont\bigl(\cX^\et,\nu_X^{\ast} \R^p
v_{\cX,M,\ast}^\NN(\cF)\bigr)=  \H^q_\cont\bigl(\cX^\et, \R^p
\gerv_{\cX,M,\ast}^\NN \bigl(\cF^\rig\bigr)\bigr) \Longrightarrow
\H^{p+q}_\cont\bigl(\hatgerX_M,\cF^\rig\bigr), \eqno{{(\numfo)}}$$
\noindent where the equality on the left hand side is due
to~\ref{compareformaltoalgebraic}.

\slabel Gabber. proposition\par\sprop  The following hold:\spacing
\item{{\rm a.}} If\/~$\cF$ is a torsion sheaf on~$\Sh(X^\et)^\NN$,
then~$\rho_{X,\cX}^{\cont,q}(\cF)$ is an isomorphism.
\spacing\item{{\rm b.}} the spectral sequences (\ref{specHqvNN})
and\/~(\ref{specHqnuNN}) are compatible via
$\rho_{X,\cX}^{\cont,q}$
and~$\rho_{\gerX_M,\hatgerX_M}^{\cont,p+q}$; \spacing \item{{\rm
c.}} if\/~$\cF$ is a torsion sheaf in~$\Sh(\gerX_M)^\NN$, the map
$\rho_{\gerX_M,\hatgerX_M}^{\cont,q}(\cF)$ is an isomorphism.

\endsprop\Proof (a) follows from~\ref{finGabber} (a) and the exact
sequence~(\ref{inversecohcohinverse}) noting that the inverse
limit of a torsion inverse system of sheaves is itself torsion;
(b) is left to the reader; (c) is proven similarly
to~\ref{finGabber} (c).

\slabel comparegerXKbartohatgerXKabrII. corollary\par\scor Let
$\bbL$ be a locally constant  sheaf on~$\gerX_M$ annihilated
by~$p^s$. Then, the two sides of the Leray spectral sequences
$$\eqalign{ \H^j_\cont\Bigl(\cX^\et, \R^i
\gerv_{\cX,M,\ast}^\NN\Bigl(\bbL^\rig  \tensor  A_{\rm
inf,s}^+\bigl(& \barhatcO_{\hatgerX_M}\bigr)\Bigr)\Bigr)
 \Longrightarrow \cr & \H^{i+j}_\cont\left(\hatgerX_M,
\bbL^\rig \tensor A_{\rm
inf,s}^+\bigl(\barhatcO_{\hatgerX_M}\bigr)\right)\cr }$$and
$$ \H^j_\cont\Bigl( X^\et, \R^i
v_{X,M,\ast}^\NN\Bigl(\bbL \tensor  A_{\rm
inf,s}^+\bigl(\barcO_{\gerX_M}\bigr) \Bigr)\Bigr) \Longrightarrow
\H^{i+j}_\cont\left(\gerX_M, \bbL \tensor A_{\rm
inf,s}^+\bigl(\barcO_{\gerX_M}\bigr)\right)$$are isomorphic.
\endscor \Proof The statements follow from~\ref{Gabber}
and~\ref{compareAinfhatoverlineandoverline}.

\endssection

\label compareGaltoXetII. proposition\par\prop {\rm
(Faltings)}\enspace Let $\bbL$ be a finite locally constant
\'etale sheaf on $X_\Kbar$ annihilated by~$p^s$. For every~$i$ the
kernel and the cokernel of the induced map of
$\WW_s\bigl(\EEtilde_\Vbar^+\bigr)$--modules
$$\H^i\left(\gerX_\Kbar,\bbL\right)\tensor
\WW_s\bigl(\EEtilde_\Vbar^+\bigr) \llongrightarrow
\H^i_\cont\left(\gerX_\Kbar,\bbL\tensor A_{\rm
inf,s}^+\bigl(\barcO_{\gerX_\Kbar}\bigr)\right) $$are annihilated
by the Teichm\"uller lift of any element in the maximal ideal
of~$\EEtilde_\Vbar^+$.
\endprop
\Proof By devissage one reduces to the case~$s=1$. The statement
follows then from~[\FALAST, \S3, Thm.~3.8].

\label formalfaltings. proposition\par\prop We have a commutative
square
$$\matrix{ \H^i\left(\hatgerX_\Kbar,\bbL^\rig\right)\tensor
\WW_s\bigl(\EEtilde_\Vbar^+\bigr) & \llongrightarrow &
\H^i_\cont\left(\hatgerX_\Kbar,\bbL^\rig\tensor A_{\rm
inf,s}^+\bigl(\barhatcO_{\hatgerX_\Kbar}\bigr)\right) \cr
\mapdownl{\Vert} & &\Big\downarrow\cr \lim_\leftarrow
\H^i\left(\hatgerX_\Kbar,\bbL^\rig\right)\tensor
\WW_s\bigl(\Vbar/p\Vbar\bigr) &\llongrightarrow  & \lim_\leftarrow
\H^i\left(\hatgerX_\Kbar,\bbL^\rig\tensor
\WW_s\left(\barhatcO_{\hatgerX_\Kbar}/p
\barhatcO_{\hatgerX_\Kbar}\right)\right), \cr}$$where the inverse
limits are taken with respect to Frobenius. The kernel and the
cokernel of any two maps in the square are annihilated by the
Teichm\"uller lift of any element in the maximal ideal
of~$\EEtilde_\Vbar^+$. Furthermore, each map \labelf
kernelandcokernelannihilatedbypi\par$$
\H^i\left(\hatgerX_\Kbar,\bbL^\rig\right)\tensor
\WW_s\bigl(\Vbar/p\Vbar\bigr) \llongrightarrow
\H^i\left(\hatgerX_\Kbar,\bbL^\rig\tensor
\WW_s\left(\barcO_{\hatgerX_\Kbar}/p
\barcO_{\hatgerX_\Kbar}\right)\right),\eqno{{(\numfo)}}$$appearing
in the inverse limits in the displayed square, has kernel and
cokernel annihilated by the Teichm\"uller lift of any element in
the maximal ideal of~$\EEtilde_\Vbar^+$.

\endprop \Proof We first of all construct the maps in the square.
The top horizontal  map is defined by the natural map
$\bbL^\rig\to \bbL^\rig\tensor A_{\rm
inf,s}^+\left(\barhatcO_{\hatgerX_\Kbar}\right)$. Similarly, the
lower horizontal arrow is induced by the map $\bbL^\rig \to
\bbL^\rig\tensor \WW_s\left(\barhatcO_{\hatgerX_\Kbar}/p
\barhatcO_{\hatgerX_\Kbar}\right) $. Note
that~$\H^0_\cont(\hatgerX_\Kbar,\_)$ is the composite of the
functors $\displaystyle\lim_\leftarrow
\H^0_\NN(\hatgerX_\Kbar,\_)$. This gives a spectral sequence in
which the derived
functors~$\displaystyle{{\lim_\leftarrow}}^{(i)}$
of~$\displaystyle\lim_\leftarrow$ on the category of abelian
groups appear. Since~$\displaystyle{{\lim_\leftarrow}}^{(i)}=0$
for~$i\geq 2$, see~[\Jannsen, \S1], we get an exact sequence
$$\eqalign{ 0  \to  {\lim_\leftarrow}^{(1)}
\H^{i-1}_\NN\Bigl(\hatgerX_\Kbar, \bbL^\rig \tensor \WW_s &
\left(\barhatcO_{\hatgerX_\Kbar}/
p\barhatcO_{\hatgerX_\Kbar}\right)\Bigr)
\Bigl(\barhatcO_{\hatgerX_\Kbar}\Bigr)\Bigr) \longrightarrow \cr &
\longrightarrow \H^i_\cont\Bigl(\hatgerX_\Kbar, \bbL^\rig \tensor
A_{\rm inf,s}^+\Bigl(\barhatcO_{\hatgerX_\Kbar}\Bigr)\Bigr)
\longrightarrow \cr & \longrightarrow \lim_\leftarrow
\H^i_\NN\Bigl(\hatgerX_\Kbar, \bbL^\rig \tensor
\WW_s\left(\barhatcO_{\hatgerX_\Kbar}/
p\WW_s\barhatcO_{\hatgerX_\Kbar}\right)\Bigr)\longrightarrow
0.\cr}
$$This provides the right vertical map in the square. Clearly the
square commutes. The fact that the top horizontal arrow has kernel
and  cokernel annihilated by the Teichm\"uller lift of any element
in the maximal ideal of~$\EEtilde_\Vbar^+$ follows
by~\ref{comparegerXKbartohatgerXKabrII}
and~\ref{compareGaltoXetII}. The equality on the left hand side
follows since $\H^i\left(\hatgerX_\Kbar,\bbL^\rig\right)$ is a
finite group being isomorphic
to~$\H^i\left(X_\Kbar^\et,\bbL\right)$
by~\ref{comparegerXKbartohatgerXKabrI} and~\ref{compareGaltoXetI}.

\noindent To conclude the proof, it suffices to show that the
kernel and cokernel of~(\ref{kernelandcokernelannihilatedbypi})
are annihilated by the Teichm\"uller lift of any element in the
maximal ideal of~$\EEtilde_\Vbar^+$.  We may reduce to the
case~$s=1$. For any integer~$m\geq 1$
let~$\left(\barhatcO_{\hatgerX_\Kbar}/p\barhatcO_{\hatgerX_\Kbar}\right)^{\geq
m}$ be the inverse system $
\left\{\barhatcO_{\hatgerX_\Kbar}/p\barhatcO_{\hatgerX_\Kbar}\right\}$
where the transition maps are the identity in degree~$\geq m$ and
are Frobenius in degree~$<m$. Let $\beta_m\colon
\cR\left(\barhatcO_{\hatgerX_\Kbar}\right) \longrightarrow
\left(\barhatcO_{\hatgerX_\Kbar}/p\barhatcO_{\hatgerX_\Kbar}\right)^{\geq
m}$ be the map of inverse systems whose $n$--th component
is~$\varphi^{n-m}\colon
\barhatcO_{\hatgerX_\Kbar}/p\barhatcO_{\hatgerX_\Kbar} \to
\barcO_\Kbar/p\barcO_\Kbar$ for~$n> m$ and is the identity for~$n<
m$. We claim that~$\beta_m$ is surjective. It suffices to check it
componentwise and, for each component, to check surjectivity
of~$\varphi^n\colon
\barhatcO_{\hatgerX_\Kbar}/p\barhatcO_{\hatgerX_\Kbar}\to
\barhatcO_{\hatgerX_\Kbar}/p\barhatcO_{\hatgerX_\Kbar}$ on stalks.
This follows from~\ref{sh}(v). Consider~$\pi_0^{p^m}
\cR\left(\barhatcO_{\hatgerX_\Kbar}\right)$ with~$\pi_0:=(p,p^{1\over p},p^{1\over
p^2},\cdots)$. Then, $\pi_0^{p^m}
\cR\left(\barhatcO_{\hatgerX_\Kbar}\right) $ is the inverse system
$\{p^{1\over
p^{n-m}}\barhatcO_{\hatgerX_\Kbar}/p\barhatcO_{\hatgerX_\Kbar}\}_n
$ with transition map given by Frobenius. We claim that
$\Ker(\beta)=\pi_0^{p^m}
\cR\left(\barhatcO_{\hatgerX_\Kbar}\right)$. This also can be
checked component--wise, for each component it can be checked on
stalks and it follows from~\ref{sh}(v).
Note that
$$\H^i_\cont\left(\hatgerX_M,
\left(\barhatcO_{\hatgerX_\Kbar}/p\barhatcO_{\hatgerX_\Kbar}\right)^{\geq
m}\right)\cong \H^i\left(\hatgerX_M,
\barhatcO_{\hatgerX_\Kbar}/p\barhatcO_{\hatgerX_\Kbar}\right).$$Indeed,
by~[\Jannsen, Prop.~1.1] an injective resolution of
$(\barhatcO_{\hatgerX_\Kbar}/p\barhatcO_{\hatgerX_\Kbar})^{\geq
m}$ is given by an injective resolution of each component of this
inverse system which is constant in degree~$n\geq m$. Take the
long exact sequence of the groups
$\H^i_\cont\left(\hatgerX_\Kbar,\_\right)$ associated to the short
exact sequence $$0\longrightarrow
\bbL^\rig\tensor\cR\left(\barhatcO_{\hatgerX_\Kbar}\right)
\lllongmaprighto{1\tensor\pi_0^{p^m}}
\bbL^\rig\tensor\cR\left(\barhatcO_{\hatgerX_\Kbar}\right)
\lllongmaprighto{1\tensor\beta_m} \left(\bbL^\rig\tensor
\barhatcO_{\hatgerX_\Kbar}\right)^{\geq m} \longrightarrow 0.$$We
get the exact sequence$$\eqalign{ \H^i_\cont &
\left(\hatgerX_\Kbar,\bbL^\rig\tensor
\cR\left(\barhatcO_{\hatgerX_\Kbar}\right)\right)
\llongmaprighto{\pi_0^{p^m}}
\H^i_\cont\left(\hatgerX_\Kbar,\bbL^\rig\tensor
\cR\left(\barhatcO_{\hatgerX_\Kbar}\right)\right) \llongrightarrow
\cr & \H^i\left(\hatgerX_\Kbar,\bbL^\rig\tensor
\barhatcO_{\hatgerX_\Kbar}/p\barhatcO_{\hatgerX_\Kbar}\right)
\llongmaprighto{\delta_i}
\H^{i+1}_\cont\left(\hatgerX_\Kbar,\bbL^\rig\tensor
\cR\left(\barhatcO_{\hatgerX_\Kbar}\right)\right) \cr &
 \llongmaprighto{\pi_0^{p^m}}
\H^{i+1}_\cont\left(\hatgerX_\Kbar,\bbL^\rig\tensor
\cR\left(\barhatcO_{\hatgerX_\Kbar}\right)\right)\cr}$$which we will compare
with the exact sequence
$$\eqalign{
\H^i\left(\hatgerX_\Kbar,\bbL^\rig\right) \tensor \EEtilde_\Vbar^+
\maprighto{\pi_0^{p^m}} &
\H^i\left(\hatgerX_\Kbar,\bbL^\rig\right)\tensor \EEtilde_\Vbar^+
\longrightarrow \H^i\left(\hatgerX_\Kbar,\bbL^\rig\right)\tensor
\bigl(\Vbar/p\Vbar\bigr)\maprighto{0} \cr & \llongrightarrow
\H^{i+1}\left(\hatgerX_\Kbar,\bbL^\rig\right)\tensor
\EEtilde_\Vbar^+ \maprighto{\pi_0^{p^m}}
\H^{i+1}\left(\hatgerX_\Kbar,\bbL^\rig\right)\tensor
\EEtilde_\Vbar^+\cr}$$via the maps
$f_j:\H^j\left(\hatgerX_\Kbar,\bbL^\rig\right) \tensor
\EEtilde_\Vbar^+\to
\H^j_\cont\left(\hatgerX_\Kbar,\bbL^\rig\tensor
\cR\left(\barhatcO_{\hatgerX_\Kbar}\right)\right) $ defined
in~\ref{compareGaltoXetII} for~$j=i$ or~$j=i+1$.

\noindent
Set $\delta_{-1}=0$ and let us denote for the rest of the proof
$\cF:=\bbL^\rig\tensor \barhatcO_{\hatgerX_\Kbar}/p\barhatcO_{\hatgerX_\Kbar}$,
$\cG:=\bbL^\rig\tensor
\cR(\barhatcO_{\hatgerX_\Kbar})$ and $\EE:=\EEtilde_\Vbar^+$. Fix $m\ge 1$ and
$i\ge 0$ and consider the diagram
$$\matrix{ \H^i\left(\hatgerX_\Kbar,\bbL^\rig\right)\tensor
\Vbar/p\Vbar & \maprighto{0} &
\H^{i+1}\left(\hatgerX_\Kbar,\bbL^\rig\right)\tensor \EE&
\maprighto{\pi_0^{p^m}}&\H^{i+1}\left(\hatgerX_\Kbar, \bbL^\rig\right)\tensor \EE\cr
f_i\downarrow & &f_{i+1}\downarrow& &f_{i+1}\downarrow\cr
\H^i_\cont\left(\hatgerX_\Kbar,\cF\right)
& \maprighto{\delta_i} & \H^{i+1}_\cont\left(\hatgerX_\Kbar,\cG\right)&\maprighto{\pi_0^{p^m}}& \H^{i+1}_\cont
\left(\hatgerX_\Kbar,\cG\right)\cr}$$
Let us remark that the right square of the diagram is commutaive and that the rows are exact.
We claim that the image of $\delta_i$ is annihilated by every element of the maximal ideal of
$\EE$, i.e. that $\delta_i$ is ``almost zero''.
For every
$\epsilon\in \QQ$ with $\epsilon >0$
let us denote by $\pi_0^\epsilon$ any element $r$ of $\EE$
such that $v_{\EE}(r)=\epsilon$. Let us fix any such $\epsilon$ and let
$x\in \H^i_\cont\left(\hatgerX_\Kbar,\cF\right)$. Denote by $y=\delta_i(x)\in \Ker(\pi_0^{p^m})$.
As the cokernel of $f_{i+1}$ is annihilated by any element of the maximal ideal of $\EE$,
$\pi_0^{\epsilon/2}y=f_{i+1}(t)$ for some $t\in \H^{i+1}\left(\hatgerX_\Kbar,\bbL^\rig\right)\tensor \EE$
and therefore $0=\pi_0^{p^m}(\pi_0^{\epsilon/2}y)=\pi_0^{p^m}f_{i+1}(t)=f_{i+1}(\pi_0^{p^m}t)$.
As the kernel of $f_{i+1}$ is also annihilated by every element of the maximal ideal of $\EE$
we have $0=\pi_0^{\epsilon/2}(\pi_0^{p^m}t)=\pi_0^{p^m}(\pi_0^{\epsilon/2}t)$ and because
multiplication by $\pi_0^{p^m}$ is injective on the top row of the diagram, we deduce
$\pi_0^{\epsilon/2}t=0$. Thus $\pi_0^\epsilon\delta_i(x)=\pi_0^{\epsilon/2}(f_{i+1}(t))=
f_{i+1}(\pi_0^{\epsilon/2}t)=0$, which proves the claim.

\noindent
Now we consider the diagram.
$$\matrix{ 0 \rightarrow \H^i\left(\hatgerX_\Kbar,\bbL^\rig\right)\tensor
\EE & \maprighto{\pi_0^{p^m}} &
\H^{i}\left(\hatgerX_\Kbar,\bbL^\rig\right)\tensor \EE&
\rightarrow &\H^{i}\left(\hatgerX_\Kbar, \bbL^\rig\right)\tensor
\Vbar/p\Vbar \rightarrow 0\cr \overline{f}_i\downarrow &
&f_{i}\downarrow& &g_i\downarrow\cr 0\rightarrow
\H^i_\cont\left(\hatgerX_\Kbar,\cG\right)/M_{i-1} &
\maprighto{\pi_0^{p^m}} &
\H^{i}_\cont\left(\hatgerX_\Kbar,\cG\right)&\rightarrow&
\H^{i}_\cont
\left(\hatgerX_\Kbar,\cF\right)\maprighto{\delta_i}M_i\cr}$$ where
for every $i\ge 0$ we denoted by $M_i$ the image of $\delta_i$ in
$\H^{i+1}_\cont\left(\hatgerX_\Kbar, \cG\right)$ and
$\overline{f}_i$ is the composition of $f_i$ with the natural
projection. It is clear that the diagram is commutaive and the
rows are exact. Moreover, the snake lemma and the fact that
$\delta_i\circ g_i=0$ give the following exact sequence of
$\EE$-modules.
$$
\Ker(f_i)\rightarrow\Ker(g_i)\rightarrow \Coker(\overline{f}_i)\rightarrow\Coker(f_i)\rightarrow
\Coker(g_i)\rightarrow M_i.
$$
As $\Coker(\overline{f}_i)$ is a quotient of $\Coker(f_i)$, we deduce that
the modules $\Ker(f_i)$, $\Coker(\overline{f}_i)$, $\Coker(f_i)$ and $M_i$ are annihilated
by every element of the maximal ideal of $\EE$, and therefore the same holds for
$\Ker(g_i)$ and $\Coker(g_i)$. This finishes the proof of Proposition \ref{formalfaltings}.


\label etaletheorem. theorem\par\thm Let\/~$\bbL$ be a locally
constant sheaf on $X_M$ annihilated by $p^s$. We have a first
quadrant spectral sequence:

$$\H^j\Bigl(\cX^\et, \R^i \gerv_{\cX,M,\ast}^\cont\Bigl(\bbL^\rig
\tensor A_{\rm inf,s}^+
\Bigl(\barhatcO_{\hatgerX_M}\Bigr)\Bigr)\Bigr)
 \Longrightarrow  \H^{i+j}_\cont\left(\gerX_M, \bbL \tensor A_{\rm
inf,s}^+\Bigl(\barcO_{\gerX_M}\Bigr)\right).$$If~$M=\Kbar$, there
is a map of\/ $\WW_s\bigl(\EEtilde_\Vbar^+\bigr)$--modules
$$\H^n\bigl(X^\et_\Kbar, \bbL\bigr) \tensor \WW\bigl(\EEtilde_\Vbar^+\bigr)
\llongrightarrow\H^n_\cont \left(\gerX_\Kbar, \bbL \tensor A_{\rm
inf,s}^+\bigl(\barcO_{\gerX_\Kbar}\bigr)\right),$$which is an
isomorphism after inverting~$\pi$.

\endthm

\Proof The first spectral sequence abuts
to~$\H^{i+j}_\cont\left(\hatgerX_M, \bbL^\rig \tensor A_{\rm
inf,s}^+\Bigl(\barhatcO_{\hatgerX_M}\Bigr)\right)$. The first
statement follows then from~\ref{comparegerXKbartohatgerXKabrII}.
The second one is the content of~\ref{compareGaltoXetII}.

\label proof2. section\par\ssection Proof of theorem
\ref{geometricPhiGamma}\par The groups
$\H^n_\cont\left(\hatgerX_\Kinfty, \bbL^\rig \tensor A_{\rm
inf,s}^+\bigl(\barhatcO_{\hatgerX_\Kinfty}\bigr)\right)\bigl[\pi^{-1}\bigr]$,
are  modules
over~$\WW\bigl(\EEtilde_\Vbar\bigr)^{\cH_V}=\WW\bigl(\EEtilde_\Vinfty\bigr)$
and have residual action of $\Gamma_V$ and $\phi$. By~\ref{stalks}
the functor $\beta_{\Kinfty,\Kbar}^{\ast,\NN}\colon
\Sh\bigl(\hatgerX_\Kinfty\bigr)^\NN \to
\Sh\bigl(\hatgerX_\Kbar\bigr)^\NN$ is exact, sends flasque objects
to flasque objects and~$\H^0_\cont(\hatgerX_\Kinfty,\cF)$ is equal
to $\displaystyle\lim_\leftarrow \H^0\left(\cH_V,
\H^0_\NN\bigl(\hatgerX_\Kbar,\beta_{\Kinfty,\Kbar}^{\ast,\NN}(\cF)\bigr)\right)$
for every~$\cF$ in $\Sh\bigl(\hatgerX_\Kinfty\bigr)^\NN$. Here,
$\H^0_\NN(\hatgerX_\Kbar,\_\bigr)$ is the functor from
$\Sh\bigl(\hatgerX_\Kbar\bigr)^\NN $ to the category of inverse
systems of $\cH_V$--modules mapping~$\{\cG_n\}_n\mapsto
\{\H^0\bigl(\hatgerX_\Kbar,\cG_n\bigr)\}$. We then get a spectral
sequence \labelf HJHINNspectral\par$$\eqalign{ \H^j\Bigl(\cH_V,
\H^i_\NN\Bigl(\hatgerX_\Kbar, \bbL^\rig \tensor & A_{\rm
inf,s}^+\Bigl( \barhatcO_{\hatgerX_\Kbar}\Bigr)\Bigr)\Bigr) \cr &
\Longrightarrow \H^{i+j}_\cont\left(\hatgerX_\Kinfty, \bbL^\rig
\tensor A_{\rm
inf,s}^+\Bigl(\barhatcO_{\hatgerX_\Kinfty}\Bigr)\right).\cr}\eqno{{(\numfo)}}$$Here,
$ \H^j\bigl(\cH_V,\_\bigr)$ stands for the $j$--th derived functor
of~$\displaystyle\lim_\leftarrow \H^0(\cH_V,\_)$ on the category
of inverse systems of $\cH_V$--modules.

Put~$M:=\H^i_\NN\Bigl(\hatgerX_\Kbar, \bbL^\rig \tensor A_{\rm
inf,s}^+\Bigl(\barhatcO_{\hatgerX_\Kbar}\Bigr)\Bigr)$. Then, $M$
is the inverse system~$\{M_n\}_n$
with~$M_n:=\H^i\Bigl(\hatgerX_\Kbar, \bbL^\rig \tensor
\WW_s\bigl(\barhatcO_{\hatgerX_\Kbar}/p\barhatcO_{\hatgerX_\Kbar}\bigr)\Bigr)$
and transition maps $d_n\colon M_{n+1}\to M_n$ given by Frobenius.
By~\ref{formalfaltings} each~$d_n$ has cokernel annihilated the
Teichm\"uller lift of any element in the maximal ideal
of~$\EEtilde_\Vbar^+$ for every~$n\in\NN$.
Let~$C^\bullet(\cH_V,M_n)$ be the complex of continuous cochains
with values in~$M_n$. For every~$i\in\NN$ the transition maps in
$\{C^i(\cH_V,M_n)\}_n$ are given by Frobenius and their cokernels
are also annihilated the Teichm\"uller lift of any element in the
maximal ideal of~$\EEtilde_\Vbar^+$ for every~$n\in\NN$. We deduce
from~(\ref{Hicontiso}) and the following discussion that we have a
canonical isomorphism
$$\H^i_\cont\bigl(\cH_V, \lim_{\infty\leftarrow n} M_n\bigr)\bigl[\pi^{-1}\bigr]\isomarrow
\H^i(\cH_V,M)\bigl[\pi^{-1}\bigr],$$where
$\displaystyle\H^i_\cont\bigl(\cH_V, \lim_{\infty\leftarrow n}
M_n\bigr)$ is continuous cohomology. Eventually, we conclude
from~\ref{etaletheorem} that
$$\H^j(\cH_V,M)\bigl[\pi^{-1}\bigr]\cong \H^j\left(\cH_V,\H^i(X^\et_\Kbar, \bbL) \tensor
\WW\bigl(\EEtilde_\Vbar\bigr)\right).$$By~\ref{induced} the latter
is zero for~$j\geq 1$ and is equal to the invariants under~$\H_V$
for~$j=0$. In particular, the spectral
sequence~(\ref{HJHINNspectral}) degenerates if we invert~$\pi$.
If~$\bbL$ is defined on~$X_K$ the isomorphism one gets is
compatible with respect to the residual action of~$\Gamma_V$ and
the action of Frobenius. The $\cH_V$--invariants of
$\H^n(X^\et_\Kbar, \bbL) \tensor \WW\bigl(\EEtilde_\Vbar\bigr)$
coincide by definition with $\cDtilde_V\bigl(\H^n(X^\et_\Kbar,
\bbL) \bigr)$.

\endssection

\endsection

\section The cohomology of Fontaine sheaves\par

In this section $\cX$ denotes a formal scheme topologically of
finite type, smooth and geometrically irreducible over~$V$ and let
$X_K^{rig}$ be its generic fiber. We study the cohomology on
$\hatgerX_M$ of continuous sheaves satisfying certain assumptions
(see \ref{assumption}). For example, it follows
from~\ref{forsomesheavesitholds} that these sheaves $\cF$ can be
taken to be of the following form:\spacing

\item{{1)}} If $\bbL$ is a $p$-power torsion \'etale local system
on $X_K^{rig}$ we set $\cF:=\bbL\otimes A_{\rm
inf}^+\left(\barhatcO_{\hatgerX_M}\right)$.

\spacing

\item{{2)}} If $\bbL$ is an  \'etale sheaf on $X_K^{rig}$ such
that $\displaystyle\bbL=\lim_\leftarrow \bbL_n$, with each
$\bbL_n$ a locally constant $\ZZ/p^n\ZZ$-module and we set
$\cF:=\bbL\hat{\tensor}\barhatcO_{\hatgerX_M}$.

\spacing   \noindent Then the cohomology groups
$\H^i\bigl(\hatgerX_M, \cF\bigr)\bigl[\pi^{-1}\bigr]$ can be
calculated as follows (here $\pi$ is $\bigl[\epsilon\bigr]-1\in
A_{\rm inf}^+(\Vbar)$ if $\cF$ is of the first type and $\pi$ is
$p$ if $\cF$ is of the second).

Let us fix a geometric generic point $\eta=\Spm(\CC_{\cX})$ as in
\S5 and for each small formal scheme $\cU=\Spf(R_\cU)$
(see~\ref{small}) with a map $\cU\lra \cX$ which is \'etale,
define $\Rbar_\cU$ to be the union of all finite, normal
$R_\cU$-algebras contained in $\CC_{\cX}$, which are \'etale after
inverting $p$. Denote by $\cF(\Rbar_\cU\otimes K)$ the inductive
limit of the sections $\cF(\cU,\cW)$, where $\cW$ runs over all
objects of $\cU_K^\fet$. Then $\cF(\Rbar\otimes_VK)$ is a
continuous representation of $\pi_1^{alg}(\cU_K,\eta)$. Moreover
(see~\ref{barhatcOMRbar=Rbar})
$\barhatcO_{\hatgerX_M}(\Rbar_\cU\otimes_VK)\bigl[p^{-1}\bigr] \cong
\widehat{\Rbar}_\cU\bigl[p^{-1}\bigr]$  and $A_{\rm
inf}^+\bigl(\barhatcO_{\hatgerX_M}\bigr)(\Rbar_\cU\otimes_VK)\bigl[\pi^{-1}\bigr]$
is isomorphic to the relative Fontaine ring $A_{\rm inf}^+$ (in which $\pi$ was inverted)
constructed using the pair $(R_\cU,\Rbar_\cU)$. We
make the following local assumption on~$\cX$.

\spacing

\noindent {\it $\cX$ admits a covering by small
objects~$\{\cU_i\}_i$ for which\/~Assumption~(ii)
of~\ref{assumption} holds.}

\spacing

\noindent Assumption~(ii) in~\ref{assumption} is equivalent to
requiring the existence of a basis of~$\cX$ by small objects
satisfying a technical compatibility condition. For any
such~$\cU=\Spf(R_\cU)$, the association $\cU\lra
\H^i\bigl(\pi_1^{alg}(\cU_K,\eta),
\cF(\Rbar_\cU\otimes_VK)\bigr)\bigl[\pi^{-1}\bigr]$ is functorial
and we denote by $\cH^i_{\Gal_M}(\cF)$ the sheaf on $\cX^\et$
associated to it. Then the main result of this section is:

\label cohomologyFontaine. theorem\par\thm Assume that the above assumption
holds. Then, there exists a spectral sequence
$$
E_2^{p,q}=\H^q\bigl(\cX^{et}, \cH^p_{\Gal_M}(\cF)\bigr)
\Longrightarrow \H^{p+q}(\hatgerX_M, \cF).
$$
\endthm

\noindent As mentioned in the Introduction, theorem
\ref{cohomologyFontaine} is the main technical tool needed to
prove comparison isomorphisms relating different $p$-adic
cohomology theories on $X_K^{rig}$. The proof of theorem
\ref{cohomologyFontaine} will take the rest of the section.

\label Zarikisites. section\par\ssection Zariski sites\par Denote
by $\cX^\Zar$ the Zariski topology on~$\cX$.\spacing

{\it The site $\hatgerX_M^\Zar$.} Let the underlying category of $\hatgerX_M^\Zar$ be the
full subcategory of the category of $\hatgerX_M$ defined
in~\ref{toposes} whose objects are pairs~$(\cU,\cW)$
with~$(\cU,\cW)\in \hatgerX_M$ and~$\cU\to \cX$ is a Zariski open
formal subscheme. We define a family of maps in $\hatgerX_M^\Zar$ to be a covering
family if it is a covering family when considered in $\hatgerX_M$.
We let $$\iota\colon \hatgerX_M^\Zar
\llongrightarrow \hatgerX_M$$be the natural functor. We also
denote by
$$\gerv_{\cX,M}\colon \cX^\Zar \llongrightarrow \hatgerX_M^\Zar $$the
map of Grothendieck topologies given by
$\gerv_{\cX,M}(\cU):=(\cU,(\cU^\rig,K))$. Since $\iota$ sends
covering families to covering families, it is clear that
$\iota_\ast\colon \Sh\bigl(\hatgerX_M^\Zar\bigr) \to
\Sh\bigl(\cX_M^\Zar\bigr)$ and $\iota_\ast^\NN\colon
\Sh\bigl(\hatgerX_M^\Zar\bigr)^\NN \to
\Sh\bigl(\cX_M^\Zar\bigr)^\NN$ send flasque objects to
flasque objects.\spacing

{\it Stalks.} Let~$\cx\colon \Spf(V_\cx) \to \cX$ be a closed
immersion of formal schemes with~$V\subset V_\cx(\subset \Kbar)$ a
finite extension of discrete valuation rings. Let~$\cO_{\cX,\cx}$ be the local ring
of~$\cO_\cX$ at~$\cx$. Define~$\barcO_{\cX,\cx,M}^\Zar$ to be the
limit $\lim_{i,j} S_{i,j}$ over all quadruples
$(R_i,S_{i,j},S_{i,j}\to \widehat{\Vbar},L_{i,j})$ where (1)
$\Spf(R_i)\subset \cX$ is a Zariski open neighborhood of~$\cx$,
(2) $L_{i,j}$ is a finite extension of\/~$K$ contained in~$M$, (3)
$R_i\subset S_{i,j}$ is an integral extension with~$S_{i,j}$
normal, (4) $S_{i,j}\tensor_V K$ is a finite and \'etale
$R_i\tensor_V L_{i,j}$--algebra, (5) the composite $R_i\tensor_V
L_{i,j}\to S_{i,j}\tensor_V K\to \widehat{\Kbar}$ is $a\tensor
\ell \mapsto \cx^\ast(a) \cdot \ell$. If~$\cF$ is a sheaf
on~$\hatgerX_M^\Zar$, define the stalk of~$\cF$ at~$\cx$ to be
$$\cF_\cx=\cF(\barcO_{\cX,\cx,M}^\Zar):=\lim_{i,j}\cF\bigl(\Spf(R_i),(\Spm(S_{i,j}\tensor_V
K),L_{i,j})\bigr).$$A sequence of sheaves on~$\hatgerX_M^\Zar$ is
exacts if and only if the induced sequence of stalks is exact for
every closed immersion~$\cx\colon \Spf(V_\cx) \to \cX$ as above.
As in~\ref{stalks} one proves that $\left(\R^q
\gerv_{\cX,M,\ast}(\cF)\right)_\cx\cong \H^q\bigl(G_{\cx,M},
\cF_\cx\bigr)$ where
$G_{\cx,M}:=\Gal\bigl(\barcO_{\cX,\cx,M}^\Zar/\cO_{\cX,\cx}\tensor_V
K\bigr)$.\endssection

\label pointedsites. section\par\ssection Pointed \'etale
sites\par Let\/~$\bbK$ be an algebraic closure of the field of
fractions of\/~$\cX\tensor_V k$. Let\/~$\WW(\bbK)$ be the Witt
vectors of\/~$\bbK$ and let\/~$\bbC_\cX$ be the $p$--adic
completion of an algebraic closure of the fraction field
of\/~$\WW(\bbK)$ containing~$\Kbar$. \spacing

{\it The site $\cX^{\et\,\bullet}$.} Denote by
$\cX^{\et\,\bullet}$ the following Grothendieck topology. As a
category it consists of pairs~$(\cU,s)$ where~$\cU\to \cX$ is an
\'etale morphism of formal schemes and $s$ is a $V$--morphism
$\Gamma(\cU,\cO_\cU)\to \bbC_\cX$. A map of pairs~$(\cU,s)\to
(\cU',s')$ is a map of $\cX$--schemes $\cU\to \cU'$ such that the
composite of\/ $\Gamma(\cU',\cO_{\cU'})\to \Gamma(\cU,\cO_\cU)$
with~$s$ is~$s'$. A covering $\amalg_{i\in I} (\cU_i,s_i)\to
(\cU,s)$ is defined to a map of pairs~$(\cU_i,s_i)\to (\cU,s)$ for
every~$i$ such that~$I$ is finite and~$\amalg_i \cU_i\to \cU$ is
\'etale surjective.\spacing

Fix~$\cx\colon \Spf(V_x)\to \cX$ as in~\ref{definegeopoints} and
{\it choose} a homomorphism $\eta_\cx\colon \cO_{\cX,\cx}^\sh \to
\bbC_\cX$ of $\cO_{\cX,\cx}$--algebras. Given a sheaf~$\cF$
on~$\cX^{\et\bullet}$ define~$\cF_\cx$ to
be~$\cF\bigl(\cO_{\cX,\cx}^\sh\bigr)$ as in~\ref{definegeopoints}.
One then proves that a sequence of sheaves on~$\cX^{\et\bullet}$
is exact if and only if the associated sequence of stalks is exact
for every~$\cx\colon \Spf(V_x)\to \cX$.

\spacing {\it The site $\hatgerX_M^\bullet$.} Define
$\hatgerX_M^\bullet$ to be the following Grothendieck topology.
Its objects are the pairs $\bigl((\cU,s),\cW,L\bigr)$
where~$(\cU,\cW,L)$ is an object of\/~$\hatgerX_M$ and\/~$(\cU,s)$
is an object of\/~$\cX^{\et\,\bullet}$. A morphism
$\bigl((\cU,s),\cW,L\bigr) \to \bigl((\cU',s'),\cW',L'\bigr)$
in~$\hatgerX_M^\bullet$ is a morphism $(\cU,\cW,L)\to
(\cU',\cW',L')$ in $\hatgerX_M^\bullet$ such that the induced map
$\cU\to \cU'$ arises from a map~$(\cU,s)\to (\cU',s')$
in~$\cX^{\et\,\bullet}$. A covering $\amalg_{i\in I}
\bigl((\cU_i,s_i),\cW_i,L_i\bigr) \to
\bigl((\cU',s'),\cW',L'\bigr)$ is the datum of morphisms
$\bigl((\cU_i,s_i),\cW_i,L_i\bigr) \to
\bigl((\cU',s'),\cW',L'\bigr)$ for every~$i\in I$ such that~$I$ is
finite and  the map $\amalg_i (\cU_i,\cW_i,L_i)\to (\cU'\cW',L') $
is a covering in~$\hatgerX_M$.\spacing

Fix~$\cx\colon \Spf(V_x)\to \cX$ as in~\ref{definegeopoints} and
{\it choose} a homomorphism $\overline{\eta}_\cx\colon
\barcO_{\cX,\cx,M} \to \bbC_\cX$ of $\cO_{\cX,\cx}^\sh$--algebras.
Given a sheaf~$\cF$ on~$\hatgerX_M^\bullet$ let~$\cF_\cx$
be~$\cF\bigl(\barcO_{\cX,\cx,M}\tensor_V K\bigr)$, defined as
in~\ref{definegeopoints}. Then, a sequence of sheaves
on~$\hatgerX_M^\bullet$ is exact if and only if the associated
sequence of stalks is exact for every~$\cx\colon \Spf(V_\cx)\to
\cX$. \spacing

\noindent We have  functors

\item{{\rm i)}} $a\colon \cX^{\et\bullet}\to \cX^\et$ given by
$a\bigl(\cU,s\bigr)=\cU$;

\item{{\rm (ii)}} $b\colon\hatgerX_M^\bullet \to \hatgerX_M$ given
by $b\bigl((\cU,s),\cW,L\bigr)=(\cU,\cW,L)$;

\item{{\rm (iii)}} $\gerv_{\cX,M}\colon \cX^{\et\bullet}\to
\hatgerX_M^\bullet$ given by
$\gerv_{\cX,M}(\cU,s)=\bigl((\cU,s),\cU^\rig,K\bigr) $.

\spacing \noindent As in~\ref{stalks}(7) one proves that for every
point~$\cx\colon \Spf(V_x)\to \cX$, $$\bigl(\R^i
\gerv_{\cX,M,\ast}(\cF)\bigr)_\cx\cong \H^i\left(G_{\cx,M},
\cF_\cx\right),\qquad
G_{\cx,M}:=\Gal\bigl(\barcO_{\cX,\cx}\tensor_V
K/\cO_{\cX,\cx}^\sh\tensor_V M\bigr).$$
Then:

\slabel okforpointed. lemma\par\slemma Let\/~$\cF$ be a sheaf
on~$\cX^\et$ (resp.~$\hatgerX_M$). We have a natural isomorphism
of $\delta$--functors
$\H^i\bigl(\cX^{\et\bullet},a_\ast(\cF)\bigr) \cong
\H^i(\cX^\et,\cF)$ (resp.~
$\H^i\bigl(\hatgerX_M^\bullet,b_\ast(\cF)\bigr) \cong
\H^i(\hatgerX_M,\cF)$).

Analogously, if~$\cF\in \Sh(\cX^\et)^\NN$
(resp.~$\Sh(\hatgerX_M)^\NN$). Then, we have a natural isomorphism
of $\delta$--functors
$\H^i_\cont\bigl(\cX^{\et\bullet},a_\ast^\NN(\cF)\bigr) \cong
\H^i_\cont(\cX^\et,\cF)$ (resp.~
$\H^i_\cont\bigl(\hatgerX_M^\bullet,b_\ast^\NN(\cF)\bigr) \cong
\H^i_\cont(\hatgerX_M,\cF)$).

\endslemma

\Proof We  have functors
$$a_\ast\colon \Sh\bigl(\cX^\et \bigr)\llongrightarrow
\Sh\bigl(\cX^{\et\bullet} \bigr),\qquad b_\ast\colon\Sh\bigl(
\hatgerX_M\bigr) \llongrightarrow \Sh\bigl(\hatgerX_M^\bullet
\bigr),$$which send flasque objects to flasque objects. Since~$a$
and~$b$ are surjective, $a_\ast$ and $b_\ast$ are also exact.
The lemma follows.

\

This allows us to work with pointed sites, better suited for
Galois cohomology as we will see.

\endssection

\label defUbarrepresentation. ssection\par\ssection The site
$\cU_M^\fet$\par Let $\cU\subset \cX$ be a Zariski open formal
subscheme or an object of~$\cX^{\et\bullet}$. Let~$\cU_M^\fet$ be
the Grothendieck topology~$\cU_M^\fet$ introduced
in~\ref{toposes}. It is a full subcategory of $\hatgerX_M^\Zar$
(resp.~$\hatgerX_M$). If~$\cU'\to \cU$ is a morphism in~$\cX^\Zar$
(resp.~$\cX^{\et\bullet}$), we have a map of Grothendieck
topologies
$$\rho_{\cU,\cU'}\colon \cU_M^\fet  \llongrightarrow
\cU_M^{'\fet}$$letting $\rho_{\cU,\cU'}\bigl(\cU,\cW\bigr)$ be the
pair~$(\cU',\cW')$ where~$\cW':=\cW\fibprod_{\cU^\rig}
\cU^{'\rig}$; see~\ref{toposes}.

\

Assume that\/~$\cU=\Spf(R_\cU)$ is affine. We have an
inclusion~$R_\cU\subset \bbC_\cX$ (this way we work
with~$\cX^{\et\bullet}$ instead of~$\cX^\et$). Let\/~$R_\cU\subset
\Rbar_\cU$ be the union of all finite and normal
$R_\cU$--subalgebras of\/~$\bbC_\cX$, which are \'etale after
inverting~$p$. We then have an inclusion~$R_\cU\subset \bbC_\cX$.
If~$\cU=\amalg_i \cU_i$, with~$\cU_i$ of the type above for
every~$i$, define $\Rbar_{\cU}:=\prod_i \Rbar_{\cU_i}$.\spacing

Define~$\pi_1(\cU_M)$ to be $\Gal\left(\Rbar_\cU\tensor_V
K/R_\cU\tensor_V M\right)$ and
let\/~$\Rep_\disc\bigl(\pi_1(\cU_M)\bigr)$ be the category of
abelian groups, with the discrete topology, endowed with a
continuous action of\/~$\pi_1(\cU_M)$.

\slabel hatgerXMbullet. lemma\par\slemma The category~$\cU_M^\fet$
is equivalent, as Grothendieck topology, to the category of finite
sets with continuous action
of\/~$\pi_1(\cU_M):=\Gal\left(\Rbar_\cU\left[{1\over
p}\right]/R_\cU\tensor_V M\right)$. In particular,\spacing

\item{{\rm 1)}} the functor
$$\Sh\left(\cU_M^\fet\right) \llongrightarrow
\Rep_\disc\bigl(\pi_1(\cU_M)\bigr), \qquad \cF\mapsto
\cF(\Rbar_\cU\tensor_V K),
$$with $\cF(\Rbar_\cU\tensor_V K):=\lim_{(\cU,\cW)}
\cF(\cU,\cW)$ where the direct limit is over all elements
of\/~$\cU_M^\fet$, defines an equivalence of categories;\spacing

\item{{\rm 2)}} for~$\cF\in\Sh\left(\cU_M^\fet\right)$ we have
$\H^i\bigl(\cU_M^\fet,\cF\bigr)=\H^i\bigl(\pi_1(\cU_M),\cF(\Rbar_\cU\tensor_V
K) \bigr)$, where the latter is the derived functor
of\/~$\Rep_\disc\bigl(\pi_1(\cU_M)\bigr) \ni A\mapsto
A^{\pi_1(\cU_M)}$(the Galois invariants of\/~$A$);\spacing

\item{{\rm 3)}} the functor
$$\Sh\left((\cU_M)^\fet\right)^\NN\llongrightarrow
\Rep_\disc\bigl(\pi_1(\cU_M)\bigr)^\NN,\qquad \{\cF_n\}\to
\{\cF_n(\Rbar_\cU\tensor_V K)\}$$is an equivalence of
categories;\spacing

\item{{\rm 4)}} for every~$\cF\in\Sh\left((\cU_M)^\fet\right)^\NN$
we have
$$\H^i_\cont\bigl(\cU_M^\fet,\cF\bigr)=\H^i\bigl(\pi_1(\cU_M),\cF(\Rbar_\cU\tensor_V K)
\bigr),$$where the latter is the $i$--th derived functor
of\/~$\Rep_\disc\bigl(\pi_1(\cU_M)\bigr)^\NN\to {\rm AbGr}$ given
by $\displaystyle \{A_n\}\mapsto \lim_{\infty\leftarrow n}
A_n^{\pi_1(\cU_M)}$.

\endslemma

\Proof  The first claim follows noting that~$\cU_M^\fet$ is the
category of finite and \'etale covers of~$R_\cU\tensor_V M$. By
Grothendieck's reformulation of Galois theory the latter is
equivalent to the category of finite sets with continuous action
of\/~$\pi_1(\cU_M)$.

Claims~(1) and~(3) follow from this. For example for claim~(1), an
inverse of the functor given in~(1) is given as follows. Let~$G\in
\Rep_\disc\bigl(\pi_1(\cU_M)\bigr)$.
Let~$(\cU,\amalg_i(\cW_i,L_i))\in \hatgerX_M^\Zar$
(resp.~$\hatgerX_M^\bullet$) with~$\cW_i=\Spm(S_i)$
and~$S_i\tensor_{L_i} M$ a domain and fix an embedding $f_i\colon
S_i\tensor_{L_i} M\to \Rbar_\cU\tensor_V K$.
Let~$H_i:=\Gal(\Rbar_\cU\tensor_V K/S_i\tensor_L M)\subset
\pi_1(\cU_M) $ which is independent of~$f_i$. Then, define
$\cG\bigl(\cU,\amalg_i(\cW_i,L_i)\bigr)=\dirsum_i G^{H_i}$. One
verifies that~$\cG$ is a sheaf and that the two functors are the
inverse one of the other.

For claims~(2) and~(4) we note that the cohomology groups
appearing are universal $\delta$--functors coinciding for~$i=0$.

\

\slabel FRbarK. definition\par\sdefi Let\/~$\cF$ be
in~$\Sh\bigl(\hatgerX_M^\Zar\bigr)$ ( or
$\Sh(\hatgerX_M^\bullet)$, or
$\Sh\bigl(\hatgerX_M^\Zar\bigr)^\NN$, or
$\Sh\bigl(\hatgerX_M^\bullet\bigr)^\NN$). We
define\/~$\cF(\Rbar_\cU\tensor_V K)$ as the image of\/~$\cF$
in~$\Rep_\disc\bigl(\pi_1(\cU_M)\bigr)$ (or
in~$\Rep_\disc\bigl(\pi_1(\cU_M)\bigr)^\NN$) of\/~$\cF$ via the
pull--back maps $\Sh(\hatgerX_M^\Zar)\to \Sh(\cU_M^\fet) \cong
\Rep_\disc\bigl(\pi_1(\cU_M)\bigr)$ (respectively via the
pull--back ~$\Sh(\hatgerX_M^\bullet)\to
\Rep_\disc\bigl(\pi_1(\cU_M)\bigr)$ , etc.).
\endsdefi

\endssection

{\bf Convention:} From now on we simply write~$\hatgerX_M^\ast$
for~$\hatgerX_M^\Zar$ or~$\hatgerX_M^\bullet$ and~$\cX^\ast$
for~$\cX^\Zar$ or, respectively,~$\cX^{\et\bullet}$.

\label CfRbar. ssection\par\ssection The sheaf
$\cH_{\Gal_M}^i$\par Let\/~$\cU'\to \cU$ be a map in $\cX^\ast$
with~$\cU'$ and~$\cU$ affine. We then get an induced map
$$\H^i\bigl(\pi_1(\cU_M), \cF(\Rbar_\cU\tensor_V
K)\bigr)\llongrightarrow
\H^i\bigl(\pi_1(\cU_M'),\cF(\Rbar_{\cU'}\tensor_V K)\bigr).$$In
particular, $\cU\to \H^i\bigl(\pi_1(\cU_M), \cF(\Rbar_\cU\tensor_V
K)\bigr)$  is a controvariant functor on the category of affine
objects of~$\cX^\ast$.

\slabel defHiGalM. definition\par\sdefi
For~$\cF\in\Sh(\hatgerX_M^\ast)$ define $\cH_{\Gal_M}^i(\cF)$ to
be the sheaf on~$\cX^\ast$ associated to the controvariant functor
given by $\cU\to \H^i\bigl(\pi_1(\cU_M), \cF(\Rbar_\cU\tensor_V
K)\bigr)$ for~$\cU$ affine.\endsdefi

\slabel standardres. ssection\par\sssection The standard
resolution\par Let~$\cG$ be a presheaf on~$\hatgerX_M^\ast$.
For~$i\in\NN$ and for~$\cU=\Spf(R_\cU)$ an affine object
of~$\cX^\ast$,  define
$$E^i(\cG)_\cU:=\Hom_\ZZ\left(
\ZZ\left[\pi_1(\cU_M)^{i+1}\right],\cG(\Rbar_\cU\tensor_V
K)\right).$$It is naturally endowed with an action
of~$\pi_1(\cU_M)$ defining for every~$\gamma,g_0,\ldots,g_i\in
\pi_1(\cU_M)$ and every~$f\in E^i(\cG)_\cU$ the action~$\gamma
\cdot f(g_0,\ldots,g_i)=\gamma^{-1}\bigl(f(\gamma
g_0,\ldots,\gamma g_i)\bigr)$.  Denote by~$C^i(\cG)_\cU\subset
E^i(\cG)_\cU$ the subgroup of invariants for the action
of~$\pi_1(\cU_M)$. Consider the map $$d_i\colon
\ZZ\left[\pi_1(\cU_M)^{i+1}\right]\to
\ZZ\left[\pi_1(\cU_M)^i\right], \quad (g_0,\ldots,g_i) \mapsto
\sum_{j=0}^{i}(-1)^j
(g_0,\ldots,g_{j-1},g_{j+1},\ldots,g_i)$$for~$i\geq 1$ and given
by~$g_0\mapsto 1$ for~$i=0$. We then get an exact sequence of
$\pi_1(\cU_M)$--modules
$$ \cdots \longrightarrow \ZZ\left[\pi_1(\cU_M)^{2}\right]\longrightarrow
\ZZ\left[\pi_1(\cU_M)^{1}\right] \longrightarrow \ZZ
\longrightarrow 0.$$Taking $\Hom\bigl(\_,\cG(\Rbar_\cU\tensor_V
K)\bigr)$ we get an exact sequence of $\pi_1(\cU_M)$--modules
\labelf exactcgRbar\par
$$0 \longrightarrow \cG(\Rbar_\cU\tensor_V K) \longrightarrow E^0(\cG)_\cU
\longrightarrow E^1(\cG)_\cU\longrightarrow
\cdots\eqno{{(\numfo)}}
$$which provides a resolution of~$\cG(\Rbar_\cU\tensor_V K)$ by
acyclic $\pi_1(\cU_M)$--modules. Using~\ref{hatgerXMbullet} we
define the sheaf~$\cW\mapsto E^i(\cG)(\cU,\cW)$ on the
category~$\cU^\fet$ associated to~$E^i(\cG)_\cU$. Furthermore,
$(\cU,\cW) \mapsto E^i(\cG)(\cU,\cW)$ is a controvariant functor
defined on the subcategory of~$\hatgerX_M^\ast$ of
pairs~$(\cU,\cW)$ with~$\cU$ affine.

\endsssection
\slabel sheafDicG. definition\par\sdefi
Let\/~$\cF\in\Sh(\hatgerX_M^\ast)$. For every~$i\in\NN$
define~$\cE^i(\cF)$ to be the sheaf on~$\hatgerX_M^\ast$
associated to the contravariant functor $(\cU,\cW)\to
E^i(\cF)(\cU,\cW)$ for $(\cU,\cW)$ affines. Define~$\cC^i(\cF)$ to
be the sheaf on~$\cX^\ast$ associated to the contravariant functor
associating to an affine~$\cU$ the continuous $i$--th cochains
of\/~$\pi_1(\cU_M)$ with values in~$\cF(\Rbar_\cU\tensor_V K)$
i.~e., $\cC^i(\cF)(\cU)=E^i(\cF)_\cU^{\pi_1(\cU_M)}$.\endsdefi

\slabel propsheafDicG. proposition\par\sprop The following hold:
\spacing \item{{\rm i)}} the differentials $d_i$ of \ref{standardres}
define an exact sequence of sheaves on $\hatgerX_M^\ast$
$$0\llongrightarrow \cF \llongrightarrow \cE^0(\cF) \llongrightarrow
\cE^1(\cF)\llongrightarrow \cE^2(\cF) \llongrightarrow \cdots ;$$

\spacing \item{{\rm ii)}} for every~$j\geq 1$ and every~$i$ one
has $\R^j \gerv_{\cX,M,\ast} \cE^i(\cF)=0$;

\spacing \item{{\rm iii)}} for every~$i$ one has
$\gerv_{\cX,M,\ast} \cE^i(\cF)=\cC^i(\cF)$.
\endsprop
\Proof (i) let $(\cU,(\cW,L))\in \hatgerX_M^\ast$ with~$\cU$
affine. Suppose that~$\cW=\Spm(S)$  with~$S\tensor_L M$  an
integral domain. Write~$\Gal_M(\cW):=\Gal\bigl(\Rbar_\cU\tensor_V
K/S\tensor_L M\bigr)$. Then, $E^i(\cF)(\cU,\cW)$ is
$E^i(\cF)_\cU^{\Gal_M(\cW)}$. In particular,
using~(\ref{exactcgRbar}), it follows that the kernel of
$E^0(\cF)(\cU,\cW)\to E^1(\cF)(\cU,\cW)$ is
$\cF\bigl(\Rbar_\cU\tensor_V K\bigr)^{\Gal_M(\cW)}$. This
coincides with~$\cF(\cU,\cW)$ since~$\cF$ is a sheaf thanks
to~\ref{hatgerXMbullet}. In particular, the kernel of $\cE^0(\cF)
\to \cE^1(\cF)$ is~$\cF$. To check the exactness of the sequence
in~(i) it is enough to pass to the stalks. Given~$\cx\colon
\Spf(V_\cx) \to \cX$ as in~\ref{Zarikisites}
or~\ref{pointedsites}, the stalk~$\cE^i(\cF)_\cx$ is the direct
limit~$\lim E^i(\cF)(\cU,\cW)$ over all~$(\cU,\cW)$ with~$\cU$ an
affine neighborhood of~$\cx$ and~$\cW=\Spm(S)$ with~$S\tensor_L
M\subset \Rbar_\cU\tensor_V K$. Hence, $\cE^i(\cF)_\cx=\lim
\cE^i(\cF)_\cU$ where the limit is now taken over all affine open
neighborhoods~$\cU$ of~$\cx$. Since for any
such~(\ref{exactcgRbar}) is exact, we conclude that the stalk
at~$\cx$ of the sequence in~(i) is exact as well.

(ii) The claim can be checked on stalks. As explained
in~\ref{Zarikisites} or~\ref{pointedsites}, given~$\cx\colon
\Spf(V_\cx) \to \cX$ as before, one has $\left(\R^q
\gerv_{\cX,M,\ast}(\cE^i(\cF))\right)_\cx\cong
\H^q\bigl(G_{\cx,M}, \cE^i(\cF)_\cx\bigr)$. But $\cE^i(\cF)_\cx$
coincides with the direct limit~$\lim \cE^i(\cF)_\cU$ taken over
all affine neighborhoods~$\cU$ of~$\cx$. Hence,

$$\eqalign{ \cE^i(\cF)_\cx & =\lim_{\cx\in\cU} \cE^i(\cF)_\cU =
\lim_{\cx\in\cU}
\Hom\left(\ZZ\left[\pi_1(\cU_M)^{i+1}\right],\cF(\Rbar_\cU\tensor_V
K)\right) =\cr & =\Hom\left(\lim_{\leftarrow}
\ZZ\left[\pi_1(\cU_M)^{i+1}\right],\lim_{\rightarrow}\cF(\Rbar_\cU\tensor_V
K)\right)
=\Hom\left(\ZZ\left[\bigl(G_\cx\bigr)^{i+1}\right],\cF_\cx
\right),\cr}$$where $G_\cx$ is~$G_{\cx,M}^\Zar$ or~$G_{\cx,M}$
depending whether~$\hatgerX_M^\ast$ is~$\hatgerX_M^\Zar$
or~$\hatgerX_M^\bullet$. In particular, $\left(\R^q
\gerv_{\cX,M,\ast}(\cE^i(\cF))\right)_\cx=0$ if~$q\geq 1$.
Claim~(ii) follows.

(iii) For every affine open~$\cU\subset \cX$ there exists a map
from the group of $i$--th cochains
$C^i\bigl(\pi_1(\cU_M),\cF(\Rbar_\cU\tensor_V
K)\bigr)=\bigl(E^i_\cU\bigr)^{\pi_1(\cU_M)}$
to~$\gerv_{\cX,M,\ast} \cE^i(\cF)(\cU)$. This provides a natural
map $\cC^i(\cF)\to \gerv_{\cX,M,\ast} \cE^i(\cF)$. On the other
hand, it follows from the discussion above that
$\left(\gerv_{\cX,M,\ast}(\cE^i(\cF))\right)_\cx$ is equal to the
group of $i$--th cochains $C^i\bigl(G_{\cx,M}^\Zar, \cF_\cx\bigr)$
i.~e., the stalk of~$\cC^i(\cF)$. The claim follows.

\slabel provethatlimHiUFisRivF. corollary\par\scor If  $\cF\in
\Sh(\hatgerX_M^\ast)$, then $\R^i\gerv_{\cX,M,\ast} (\cF) \cong
\cH_{\Gal_M}^i(\cF)$ functorially in~$\cF$.\endscor

\endssection

\ssection The sheaf $\cH_{\Gal_M,\cont}^i$\par We wish to prove an
analogue of~\ref{provethatlimHiUFisRivF} in the case of a
continuous sheaf~$\cF=\{\cF_n\}_n\in \Sh(\hatgerX_M^\ast)^\NN$. We need some
assumptions.

\slabel small. definition\par\sdefi An object\/~$\cU$
of\/~$\cX^\ast$ is called {\it small}
if~$\cU:=\Spf\bigl(R_\cU\bigr)$ is affine and\/~$R_\cU$ satisfies
the assumptions of~\ref{V} and~{\rm (RAE)} (see~\ref{en}).

Define~$R_{\cU,M,\infty}$ to be the normalization
of\/~$R_{\cU,\infty}$ in the subring of~$\Rbar_\cU\tensor_V K$
generated by\/~$M$ and\/~$R_{\cU,\infty}$, where~$R_{\cU,\infty}$
is defined as in~\ref{V}. Denote by~$\Gamma_{\cU,M}$ the group
$\Gal\bigl(R_{\cU,M,\infty}\tensor_V K/R_\cU\tensor_V M\bigr)$.
Let\/~$\cH_{\cU,M}$ be the kernel of the map~$\pi_1(\cU_M\bigr)\to
\Gamma_{\cU,M}$. Let us remark that the definitions of
$R_{\cU,\infty},R_{\cU,M,\infty}, \cH_{\cU,M},\Gamma_{\cU,M}$
depend on a choice of local parameters of $R_\cU$ and so are not
canonical.

\endsdefi

\slabel inftysite. ssection\par\sssection The site
$\hatgerU_M^\ast(\infty)$ \par Let~$\cU$ be a small object
of~$\cX^\ast$. For every map $\cU'\to \cU$
with~$\cU':=\Spf(R_{\cU'})$ affine and~$R_{\cU'}\tensor_V K$ an
integral domain, we let~$\cH_{\cU',M}$ be the kernel
of~$\pi_1(\cU'_M)\to \Gamma_\cU$. Note that such a map is
surjective.

Let~$\hatgerU_M^\ast(\infty)$ be the following full subcategory
of~$\hatgerU_M^\ast$. Let~$(\cU',\cW')\in \hatgerU_M^\ast$ and
assume that~$\cU':=\cup_i \cU'_i$ with~$\cU'_i$ connected. Then,
$\cW'$ lies in~$\cU^{'\fet}$ which, via the equivalence
of~\ref{hatgerXMbullet}, is equivalent to  the category of finite
sets with continuous action of~$\pi_1(\cU'_M)=\prod_i
\pi_1(\cU'_{i,M})$. We then say that~$(\cU',\cW')$ lies
in~$\hatgerU_M^\ast(\infty)$ if and only if~$\cW'$ lies in the
subcategory of finite sets with continuous action
of~$\prod_i\Gamma_\cU$ (viewed as a quotient of~$\pi_1(\cU'_M)$).
We then have natural maps of Grothendieck topologies
$\cU^\ast\maprighto{\alpha}
\hatgerU_M^\ast(\infty)\maprighto{\beta} \hatgerU_M^\ast$ giving
rise to maps on the category of sheaves
$$\Sh(\hatgerU_M^\ast) \llongmaprighto{\beta_\ast}
\Sh(\hatgerU_M^\ast(\infty)) \llongmaprighto{ \alpha_\ast}
\Sh(\cU^\ast)
$$whose composite is $\gerv_{\cU,M,\ast}$.
As in~\ref{Zarikisites} or~\ref{pointedsites} one has a notion of
stalks in~$\Sh(\hatgerU_M^\ast(\infty))$. For~$\cx\colon
\Spf(V_\cx)\to \cX$ a point as in~\ref{definegeopoints},
let~$H_{\cx,M}$ be the kernel of the map~$G_{\cx,M}\to \Gamma$.
If~$\cF\in \Sh(\hatgerU_M)$ and~$\cF_\cx$ is its stalk, one proves
as in~\ref{stalks} that $$R^q\beta_\ast(\cF)_\cx\cong
\H^q\bigl(H_{\cx,M},\cF_\cx\bigr).$${\it Caveat:} The
site~$\hatgerU_M(\infty)$ depends on the choice of an extension $R_\cU
\subset R_{\cU,\infty}$. In particular, if $\{\cU_i\}_i$ is a
covering of~$\cX$ by small objects, the sites~$\hatgerU_{i,M}(\infty)$ do
not necessarily glue  so that the site $\hatgerX_M(\infty)$ is not
defined in general.
\endsssection

\slabel assumption. ssection\par\sssection {\bf Assumption}\par We
suppose that  \spacing \item{{\rm i)}} $\{\cF_n\}_{n\in\NN}$ is a
sheaf of $A_{\rm inf}^+\left(\Vinfty\right)$--modules
(resp.~of~$\{\Vinfty/p^n\Vinfty\}_n$--modules) on
$\hatgerX_M^\ast$; \spacing \item{{\rm ii)}} $\cX$ admits
\itemitem{{a)}} a covering~$\cS:=\{\cW_i\}_{i}$ in~$\cX^\ast$ by
small objects, \itemitem{{b)}} a choice $R_{\cU_i}\subset
R_{\cU_i,\infty}$ as in~\ref{V}, \itemitem{{c)}} for every~$i$ a
basis $\cT_i:=\{\cU_{i,j}\}_{j}$ of~$\cW_i$ by small objects such
that, putting~$R_{\cU_{i,j},\infty}$ to be the normalization of
$R_{\cU_{i,j}}\tensor_{R_{\cW_i}} R_{\cW_i,\infty}$,
condition~(RAE) holds for~$R_{\cU_{i,j},\infty}$.

\spacing \noindent Furthermore,  for every~$i$, $j$ and~$n\in\NN$,
putting~$\cU:=\cU_{i,j}$, the following hold:\spacing

\spacing \item{{\rm iii)}}  the cokernel of
$\cF_{n+1}(\Rbar_\cU\tensor_V K)\to \cF_n(\Rbar_\cU\tensor_V K)$
is annihilated by any element of the maximal ideal
of~$\WW\bigl(\Vbar/p\Vbar)$ (resp.~$\Vbar$);

\spacing \item{{\rm iv)}} for every~$q\geq 1$ the group
$\H^q\left(\cH_{\cU,M},\cF_n\bigl(\Rbar_{\cU}\tensor_V
K\bigr)\right)$ is annihilated by any element of the maximal ideal
of $\WW\bigl(\Vbar/p\Vbar)$ (resp.~$\Vbar$);

\spacing \item{{\rm v)}} the cokernel of the transition maps
$\cF_{n+1}(R_{\cU,M,\infty}\tensor_V K)\to
\cF_n(R_{\cU,M,\infty}\tensor_V K)$ is annihilated by any element
of the maximal ideal of~$\WW\bigl(\Vinfty/p\Vinfty)$
(resp.~$\Vinfty$);

\spacing \item{{\rm vi)}} for every covering~$\cZ\to \cU$ by small
obiects in~$\cX^\ast$ and every~$q\geq 1$ the Chech cohomology
group $\H^q\bigl(\cZ\to \cU,\cF_n(R_{\cZ,\infty}\tensor_V
K)\bigr)$ is annihilated  by any element of the maximal ideal
of~$\WW\bigl(\Vinfty/p\Vinfty)$ (resp.~$\Vinfty$).
\endsssection

\noindent As usual we write~$\pi$ for the element~$[\epsilon]-1$
in~$\AAtilde_\Vinfty^+$. We put~$\pi=p$ if $\{\cF_n\}_{n\in\NN}$
is a sheaf of~$\{\Vinfty/p^n\Vinfty\}_n$--modules. It follows
from~(iii) and~\ref{groupscontcoho} that we have an isomorphism
$$\H^i\bigl(\pi_1(\cU_M),\cF(\Rbar_\cU\tensor_V
K)\bigr)\bigl[\pi^{-1}\bigr] \cong
\H^i_\cont\bigl(\pi_1(\cU_M),\lim_{\infty\leftarrow n}
\cF_n(\Rbar_\cU\tensor_V
K)\bigr)\bigl[\pi^{-1}\bigr].$$If\/~$\cU'\to \cU$ is a map in
$\cX^\ast$ with~$\cU'$ and~$\cU$ small objects in~$\cT_i$, we then
get an induced map
$$\H^i\bigl(\pi_1(\cU_M), \cF(\Rbar_\cU\tensor_V
K)\bigr)\bigl[\pi^{-1}\bigr]\llongrightarrow
\H^i\bigl(\pi_1(\cU_M'),\cF(\Rbar_{\cU'}\tensor_V
K)\bigr)\bigl[\pi^{-1}\bigr].$$As in~\ref{defHiGalM}, we define

\slabel defHiGalmcont. definition\par\sdefi Assume that\/~$\cF$
satisfies the assumption above. Let\/~$\cH_{\Gal_M,\cont}^i(\cF)$
be the sheaf on~$\cX^\ast$ associated to the controvariant functor
sending an object~$\cU$ of $\cX^\ast$, with~$\cU\in\cup_i \cT_i$,
to $ \H^i\bigl(\pi_1(\cU_M), \cF(\Rbar_\cU\tensor_V
K)\bigr)\bigl[\pi^{-1}\bigr]$.\endsdefi

We want to prove the following:

\slabel provethatlimHiUFisRivFcont. theorem\par\sthm Let\/ $\cF\in
\Sh(\hatgerX_M^\ast)^\NN$ be such that  the conditions
of~\ref{assumption} are fulfilled. Then, $\R^i\gerv_{\cX,M,\ast}
(\cF)\bigl[\pi^{-1}\bigr] \cong \cH_{\Gal_M}^i(\cF)$. The
isomorphism is functorial in~$\cF$.\endsthm \Proof It suffices to
prove that for every small object~$\cW_i\in \cS$, we have an
isomorphism $\R^i\gerv_{\cX,M,\ast}
(\cF)\bigl[\pi^{-1}\bigr]\vert_{\cW_i} \cong
\cH_{\Gal_M}^i(\cF)\vert_{\cW_i}$ functorially in~$\cW_i$
and~$\cF$. We construct the isomorphism and leave it to the reader
to check the functoriality in~$\cW_i$ and~$\cF$.

We may and will, {\it till the end of this section}, assume
that~$\cX=\cW_i$ {\it is small}. We put~$\cT:=\cT_i$ and we
write~$\Gamma$ for~$\Gamma_{\cW_i}$. Consider the maps on the
category of sheaves
$$\Sh(\hatgerX_M^\ast)^\NN \llongmaprighto{\beta_\ast^\NN}
\Sh(\hatgerX_M^\ast(\infty))^\NN \llongmaprighto{\lim_\leftarrow
\alpha_\ast} \Sh(\cX^\ast),
$$introduced in~\ref{inftysite}. The
composite is $\displaystyle\lim_\leftarrow \gerv_{\cX,M,\ast}$.
Since~$\alpha_\ast$ and~$\beta_\ast$ are left exact
and~$\beta_\ast$ sends injective to injective, we have a spectral
sequence \labelf spectralphabeta\par$$\R^p \lim_\leftarrow
\alpha_\ast \bigl( \R^q \beta_\ast^\NN(\cF)\bigr) \Longrightarrow
\R^{p+q}\bigl(\lim_\leftarrow
\gerv_{\cX,M,\ast}\bigr)(\cF).\eqno{{(\numfo)}}$$

\slabel Rqalpha=0. slemma\par\slemma For every~$q\geq 1$ the group
$\R^q \beta_\ast(\cF)$ is annihilated by any element of the
maximal ideal of\/~$\WW\bigl(\Vinfty/p\Vinfty)$ (resp.~$\Vinfty$).
\endslemma
\Proof Since~$\R^q\beta_\ast^\NN=\bigl(\R^q\beta_\ast\bigr)^\NN$
as remarked in~\ref{topologicalsheaves}, it suffices to prove that
for every~$n\in\NN$ and every~$q\geq 1$ the sheaf
$R^q\beta_\ast(\cF_n)$ is annihilated by any element of the
maximal ideal of~$\WW\bigl(\Vinfty/p\Vinfty)$ (resp.~$\Vinfty$).
It suffices to prove the vanishing on stalks. But for~$\cx\colon
\Spf(\cx)\to \cX$ a point as in~\ref{definegeopoints}, we have
$R^q\beta_\ast(\cF_n)_\cx\cong \H^q\bigl(H_{\cx,M},\cF_\cx\bigr)$
as explained in~\ref{inftysite}.  The latter coincides with the
direct limit $\lim
\H^q\left(\cH_{\cU,M},\cF_n\bigl(\Rbar_{\cU}\tensor_V
K\bigr)\right)$ taken over the small objects belonging to a basis
of~$\cX^\ast$ containing~$\cx$. The claim then follows
from~\ref{assumption}(i)\&(iv).

\

\noindent Using~\ref{Rqalpha=0} and~(\ref{spectralphabeta}) we
conclude that
$$\R^p \lim_\leftarrow \alpha_\ast \bigl(\beta_\ast^\NN(\cF)\bigr)\bigl[\pi^{-1}\bigr] \cong
\R^p\bigl(\lim_\leftarrow
\gerv_{\cX,M,\ast}\bigr)(\cF)\bigl[\pi^{-1}\bigr].$$We are left to
compute~$\displaystyle\R^p \lim_\leftarrow \alpha_\ast$. For this
we use the analogue of~\ref{standardres}
on~$\hatgerX_M^\ast(\infty)$.

Given~$\cU$ in~$\cT$, write $R_{\cU,M,\infty}$ as the union
$\cup_n R_{\cU,M,n} $ of finite $R_\cU$--algebras such that
$R_\cU\tensor_V K \subset \R_{\cU,M,n}\tensor_V K$ is finite and
\'etale. Then, for every covering $\cU'\to \cU$ with~$\cU'\in
\cT$, we have $R_{\cU',M,\infty}\tensor_V K\cong \cup_n R_{\cU'}
\tensor_{R_\cU} R_{\cU,M,n}  \tensor_V K$ by construction.
Let~$\cU''\to \cU'\fibprod_{\cU} \cU'$ be a covering with~$\cU''$
in~$\cT$. Then, we also have $R_{\cU'',M,\infty}\tensor_V K\cong
\cup_n R_{\cU''} \tensor_{R_\cU} R_{\cU,M,n} \tensor_V K$.
Since~$\cF_n$ is a sheaf, we conclude that the sequence
$$0\lllongrightarrow \cF_n(R_{\cU,M,\infty}\tensor_V K)
\llongrightarrow \cF_n(R_{\cU',M,\infty}\tensor_V K)
\llongrightarrow \cF_n(R_{\cU'',M,\infty}\tensor_V K)
$$is exact i.~e., $\cU \to \cF_n(R_{\cU,M,\infty}\tensor_V K)$ satisfies the sheaf
property with respect to coverings $\cU'\to \cU$ with~$\cU'$
and~$\cU$ small and lying in~$\cT$. Then, the following makes
sense:

\slabel Diccfisaresolution. definition\par\sdefi For every small
object~$\cU\to \cX$ lying in~$\cT$ and every~$i$, $n\in\NN$ define
$E^i(\Gamma,\cF_n)_\cU$ to be $\Hom_\ZZ\left(
\ZZ\left[\Gamma^{i+1}\right],\cF_n(R_{\cU,M,\infty}\tensor_V
K)\right)$. Define\/~$\cE^i(\Gamma,\cF_n)$ to be the sheaf
on~$\hatgerX_M^\ast(\infty)$ characterized by the property that,
for every small object~$\cU\in \cT$, its restriction to
$\cU_M^\fet$ (see~\ref{defUbarrepresentation}) is~$
E^i(\Gamma,\cF_n)_\cU$ as representation of~$\Gamma_\cU$.
Let\/~$\cE^i(\Gamma,\cF):=\bigl\{\cE^i(\Gamma,\cF_n)\bigr\}_n $.

Let\/ $C^i(\Gamma,\cF_n)_\cU\subset E^i(\Gamma,\cF_n)_\cU$ be the
subgroup of invariants for the action of\/~$\Gamma$ i.~e., the
group of $i$--th cochains of\/~$\Gamma$ with values
in~$\cF_n(R_{\cU,M,\infty}\tensor_V K)$. Denote
by~$\cC^i(\Gamma,\cF_n)$ the unique sheaf on~$\cX^\ast$ whose
value for every small object~$\cU$ is~$C^i(\Gamma,\cF_n)_\cU$.
Eventually, let~$\cC^i(\Gamma,\cF):=\bigl\{
\cC^i(\Gamma,\cF_n)\bigr\}_n$.
\endsdefi

\slabel standardisacyclic. propostion\par\sprop Assume
that\/~$\cF$ satisfies~\ref{assumption}. Then:\spacing \item{{\rm
i)}} we have an exact sequence
in~$\Sh(\hatgerX_M^\ast(\infty))^\NN$
$$0 \longrightarrow \beta_\ast^\NN(\cF) \longrightarrow \cE^0(\Gamma,\cF)
\longrightarrow \cE^1(\Gamma,\cF) \longrightarrow
\cdots;$$\spacing \item{{\rm ii)}} $\displaystyle\R^q
\lim_\leftarrow \alpha_\ast\bigl(
\cE^i(\Gamma,\cF)\bigr)\bigl[\pi^{-1}\bigr]=0$ for every~$q\geq 1$
and every~$i$;\spacing

\item{{\rm iii)}} $\displaystyle\lim_\leftarrow
\alpha_\ast\bigl(\cE^i(\Gamma,\cF)\bigr)\bigl[\pi^{-1}\bigr]$ is
the sheaf associated to the controvariant functor sending a small
object\/~$\cU$ to  $\displaystyle\lim_{\infty\leftarrow n}
\cC^i\left(\Gamma,\cF_n\bigl(R_{\cU,M,\infty}\tensor_V
K)\right)\bigl[\pi^{-1}\bigr]$.

\spacing\noindent In particular, $\displaystyle\R^q
\lim_\leftarrow \alpha_\ast
(\beta_\ast^\NN(\cF))\bigl[\pi^{-1}\bigr]$ is the $q$--th
cohomology of the complex $$\lim_{\infty\leftarrow n}
\cC^0(\Gamma,\cF_n)\bigl[\pi^{-1}\bigr] \longrightarrow
\lim_{\infty\leftarrow n} \cC^1(\Gamma,\cF_n)\bigl[\pi^{-1}\bigr]
\longrightarrow \lim_{\infty\leftarrow n}
\cC^2(\Gamma,\cF_n)\bigl[\pi^{-1}\bigr] \longrightarrow \cdots
,$$proving~\ref{provethatlimHiUFisRivFcont}.

\endsprop
\Proof Claim~(i) can be checked componentwise and then it follows
as in the proof of~\ref{propsheafDicG}(i).

(ii)--(iii) We use the spectral sequence $${\lim}^{(p)}
\bigl(\R^q\alpha_\ast(\cE^i(\Gamma,\cF_n))\bigr) \Longrightarrow
\R^{p+q}\lim_\leftarrow \alpha_\ast(\cE^i(\Gamma,\cF))$$given
in~\ref{topologicalsheaves}. Since each~$\cF_n$ is a sheaf we have
$\cE^i(\Gamma,\cF_n)(R_{\cU,M,\infty}\tensor_V
K)=E^i(\Gamma,\cF_n)_\cU$. Hence,
$\H^q\bigl(\Gamma,\cE^i(\Gamma,\cF_n)(R_{\cU,M,\infty}\tensor_V
K)\bigr)$ is~$0$ for every~$q\geq 1$ and it coincides with the
cochains~$\cC^i(\Gamma,\cF_n)(R_{\cU,M,\infty}\tensor_V
K)=C^i(\Gamma,\cF_n)_\cU$ for~$q=0$.

Arguing as in~\ref{propsheafDicG}(ii) we conclude that
$\R^q\alpha_\ast\bigl(\beta_\ast(\cE^i(\Gamma,\cF_n))\bigr)=0$
for~$q\geq 1$. We are left to compute $\lim^{(p)}
\alpha_\ast\bigl(\cE^i(\Gamma,\cF_n)\bigr)=\lim^{(p)}
\cC^i(\Gamma,\cF_n) $.

Due to~\ref{assumption}(vi), for every small object~$\cU\in \cT$
and every~$n$ the Chech cohomology group $\H^q\left(\cZ\to \cU,
\cF_n(R_{\cZ,\infty}\tensor_V K)\right)$, relative to every
covering~$\cZ\to \cU$ by small objects lying in~$\cT$, is
annihilated by any element of the maximal ideal
of~$\WW\bigl(\Vinfty/p\Vinfty)$ (resp.~$\Vinfty$). But we have
\advance\fonu by1\labelf FC^i(F)\par$$\cC^i(\Gamma,
\cF_n(R_{\_,\infty}\otimes K))=\lim_{m\rightarrow\infty}
\cC^i(\Gamma/p^m\Gamma, \cF_n(R_{\_,\infty}\otimes
K))=\lim_{m\rightarrow\infty}
\bigl(\prod_{\Gamma/p^m\Gamma}\cF_n(R_{\_,\infty}\otimes
K)\bigr).\eqno{{(\numfo)}}$$ As both inductive limit and finite
products are exact functors we deduce that the Chech cohomology
group $\H^q\left(\cZ\to \cU, \cC^i\bigl(\Gamma,
\cF_n(R_{\cZ,\infty}\tensor_V K)\bigr)\right)$ relative to every
covering~$\cZ\to \cU$ with~$\cZ$ and~$\cU\in\cT$ is annihilated by
any element of the maximal ideal of~$\WW\bigl(\Vinfty/p\Vinfty)$
(resp.~$\Vinfty$). Hence, the restriction of
$\cC^i\bigl(\Gamma,\cF_n(R_{\_,\infty}\tensor_V K)$ to~$\cU$ is
flasque, see~[\Artin, II.4.2], up to multiplication by any element
of the maximal ideal of~$\WW\bigl(\Vinfty/p\Vinfty)$
(resp.~$\Vinfty$). In particular,
$\H^q\bigl(\cU,\cC^i(\Gamma,\cF_n)\bigr)$ is almost zero for
every~$q\geq 1$; see~[\Artin, II.4.4]. Due to~\ref{assumption}(v)
the projective system
$\left\{\cF_n(R_{\cU,M,\infty}\otimes_VK)\right\}_n$ is almost
Mittag-Lefler and using once again~(\ref{FC^i(F)}) we also have
that the projective system
$\left\{\cC^i\bigl(\Gamma,\cF_n)\right\}_n$ is almost
Mittag--Leffler. Hence, $\lim^{(1)}\cC^i(\Gamma, \cF_n)$ is almost
zero.

By~[\Jannsen, Lemma 3.12] the sheaf~$\lim^{(q)}
\cC^i(\Gamma,\cF_n)$ is the sheaf associated to the presheaf
$\cU\mapsto
\H^q\bigl(\cU,\bigl(\cC^i(\Gamma,\cF_n)\bigr)_n\bigr)$. We have,
for each $q\ge 1$, exact sequences $$0\to {\lim}^{(1)}
\H^{q-1}\bigl(\cU,\cC^i(\Gamma,\cF_n)\bigr) \lra\H^{q}\bigl(\cU,
\bigl(\cC^i(\Gamma,\cF_n)\bigr)_n\bigr)\lra \lim_{\infty
\leftarrow n}\H^{q}\bigl(\cU, \cC^i(\Gamma,\cF_n)\bigr)\to 0.$$
For $q\ge 2$, using the fact proved above that
$\H^{s}\bigl(\cU,\cC^i(\Gamma,\cF_n)\bigr)$ is almost zero for
$s\ge 1$, we deduce that $\H^{q}\bigl(\cU,
\bigl(\cC^i(\Gamma,\cF_n)\bigr)_n\bigr)$ is almost zero. We
conclude that $ \lim^{(q)}\cC^i(\Gamma, \cF_n)$ is annihilated by
every element of the maximal ideal~$\WW\bigl(\Vinfty/p\Vinfty)$
(resp.~$\Vinfty$) for~$q\geq 2$.

\noindent
Thus,
$$
R^p\lim_\leftarrow\alpha_\ast\bigl(\cE^j(\Gamma,\cF)\bigr)\bigl[\pi^{-1}\bigr]=0
$$
for $p\ge 1$ and all $j\ge 0$, and
$$
\lim_\leftarrow\alpha_\ast\cE^\bullet(\Gamma,
\cF)\bigl[\pi^{-1}\bigr]\cong \lim_{\infty \leftarrow
n}\cC^\bullet(\Gamma,\cF_n)\bigl[\pi^{-1}\bigr].
$$
The conclusion follows.

\endssection

\

\label forsomesheavesitholds. theorem\par\thm Let\/~$\cX$ be
formally smooth, topologically of finite type and geometrically
irreducible over~$V$ for which assumption ii) of \ref{assumption}
is satisfied. Let\/ $\bbL=(\bbL_n)_n$ be a projective
system of sheaves such that\/~$\bbL_n\cong
\bbL_{n+1}/p^n\bbL_{n+1}$ for every~$n$. Let
$\cF\in\Sh(\cX^\ast)^\NN$ is a sheaf of one of
the following types:\spacing

\item{{\rm A)}} $\cF$ is~$\bbL^\rig \tensor A_{\rm inf}^+
\Bigl(\barhatcO_{\hatgerX_M}\Bigr):=\left(\bbL_n\tensor
\WW_n\bigl(\barhatcO_{\hatgerX_M}/p\barhatcO_{\hatgerX_M}
\bigr)\right)_n$;\spacing

\item{{\rm B)}} $\cF:=\left(\bbL_n\tensor
\bigl(\barhatcO_{\hatgerX_M}/p^n
\barhatcO_{\hatgerX_M}\bigr)\right)_n$ where
$\barhatcO_{\hatgerX_M}/p^{n+1} \barhatcO_{\hatgerX_M}\to
\barhatcO_{\hatgerX_M}/p^n \barhatcO_{\hatgerX_M}$ is the natural
projection for each~$n\in\NN$. \spacing

\noindent Then, the assumptions in~\ref{assumption} hold.\endthm

\label relatetophigamma. remark\par\rmk Concerning
assumption ii) of \ref{assumption}, we note that the existence of a basis by affine
subschemes satisfying the assumptions of~\ref{V} follows from the
fact that~$\cX$ is formally smooth and topologically of
finite type over~$V$. The content of ii) is~(RAE).
See~\ref{whenRAEholds} for
examples when it holds.

\

In case~(A), assume further that~$\bbL$ is a $p$--power torsion
i.~e., annihilated by~$p^s$ for some~$s$. Then, one can compute
the sheaf~$\cH_{\Gal_M,\cont}^i(\cF)$ introduced
in~\ref{defHiGalmcont} via relative $(\phi,\Gamma)$--modules.
Indeed, assume that~$\cU=\Spf(R_\cU)$ is small and that (RAE)
holds for~$R_{\cU,\infty}$.\spacing

For~$M=K$ we have~$\pi_1(\cU_M)=\cG_{R_\cU}$ and
by~\ref{decompletecohomology} the inflation
$$\H^i(\Gamma_{R_\cU}, \cD({\bf
L})\bigr)\llongrightarrow \H^i\bigl(\pi_1(\cU_K),{\bf
L}\tensor_{\ZZ_p}\AAtilde_{\Rbar_\cU}\bigr)\cong
\H^i\bigl(\pi_1(\cU_K), \cF(\Rbar_\cU\tensor_V
K)\bigr)\bigl[\pi^{-1}\bigr]$$is an isomorphism.

\spacing

Analogously, for~$M=\Kbar$ we have~$\pi_1(\cU_M)=\G_{R_\cU}$ so
that
$$\H^i(\tGamma_{R_\cU}, \D({\bf L})\bigr)\isomarrow \H^i\bigl(\pi_1(\cU_\Kbar),{\bf
L}\tensor_{\ZZ_p}\AAtilde_{\Rbar_\cU}\bigr)\cong
\H^i\bigl(\pi_1(\cU_\Kbar), \cF(\Rbar_\cU\tensor_V
K)\bigr)\bigl[\pi^{-1}\bigr].$$

\spacing

For~$M=\Kinfty$, the group~$\pi_1(\cU_M)$ is the subgroup
of~$\cG_{R_\cU}$ generated by~$\G_{R_\cU}$ and~$\H_V$. It follows
from~\ref{induced} that~$\H^i\bigl(\cH_{R_\cU},{\bf
L}\tensor_{\ZZ_p}\AAtilde_{\Rbar_\cU}\bigr)=0$ for~$i\geq 1$.
Hence, we conclude that $\H^i\bigl(\pi_1(\cU_M),{\bf
L}\tensor_{\ZZ_p}\AAtilde_{\Rbar_\cU}\bigr)=\H^i(\tGamma_R,
\cDtilde({\bf L})\bigr)$. Using~\ref{compareDicDi} we get that the
latter coincides with~$\H^i(\tGamma_R, \cD({\bf L})\bigr)$ so that
$$\H^i(\tGamma_{R_\cU}, \cD({\bf L})\bigr)\isomarrow
\H^i\bigl(\pi_1(\cU_\Kinfty),{\bf
L}\tensor_{\ZZ_p}\AAtilde_{\Rbar_\cU}\bigr)\cong
\H^i\bigl(\pi_1(\cU_\Kinfty), \cF(\Rbar_\cU\tensor_V
K)\bigr)\bigl[\pi^{-1}\bigr].$$
\endrmk

\ssection Proof of\/ \ref{forsomesheavesitholds}\par We start with
some preliminary results.

\slabel technicalneededinS6. lemma\par\slemma Let\/~$R$ be as
in~\ref{V}. Let\/~$\Sinfty \subset \Tinfty$ be integral extensions
of~$\Rinfty$ such that $\Sinfty\tensor_V K= \Tinfty\tensor_V K$
and\/ $\Rinfty\subset \Sinfty$ is almost \'etale (see~\ref{en}).
Then, the cokernel of\/ $\Sinfty \subset \Tinfty$ is annihilated
by any element of the maximal ideal of\/~$\Vinfty$.\endslemma
\Proof Let $\e_\infty$ be the canonical idempotent of the \'etale
extension $\Rinfty\bigl[p^{-1}\bigr] \subset
\Sinfty\bigl[p^{-1}\bigr]=\Tinfty\bigl[p^{-1}\bigr]$. Since
$\Rinfty\subset \Sinfty$ is assume to be almost \'etale, for
every~$\alpha\in\QQ_{>0}$ we may write $p^\alpha \e_\infty$ as a
finite sum~$\sum_i a_i\tensor b_i$ with~$a_i$ and~$b_i$
in~$\Sinfty$. Let $m\colon \Sinfty\bigl[p^{-1}\bigr]
\tensor_\Rinfty \Sinfty\bigl[p^{-1}\bigr] \to
\Sinfty\bigl[p^{-1}\bigr]$ be the multiplication map and let
$\Tr\colon \Sinfty\bigl[p^{-1}\bigr]\to \Rinfty\bigl[p^{-1}\bigr]$
be the trace map. Then, $\e_\infty$ is characterized by the
property that $m(x\tensor y)=(\Tr\tensor {\rm Id})\bigl((x\tensor
y)\cdot \e_\infty\bigr)$. In particular, for every~$x\in \Tinfty$
we have $p^\alpha x= m\bigl(p^\alpha x\tensor 1\bigr)= \sum_i
\Tr(a_i x) b_i$. But~$\Tr(a_i x)\in\Rinfty$ since~$x$ and~$a_i$
are integral over~$\Rinfty$ and~$\Rinfty$ is integrally closed.
Hence, $p^\alpha x= \sum_i \Tr(a_i x) b_i$ lies in~$\Sinfty$ as
claimed.

\slabel RcU'inftyisbasechenageofRcU. lemma\par\slemma
Let\/~$\cU=\Spf(R_\cU)$ be an affine small object of\/~$\cX^\ast$
and let\/~$\cU'\to \cU$ be a covering with\/~$\cU'$ affine. Then,
$R_{\cU',M, \infty} \cong R_{\cU,M,\infty}\tensor_{R_\cU}
R_{\cU'}$.
\endslemma \Proof Write the composite of~$M$
and~$\Kinfty$ (in~$\Kbar$) as the union~$\cup_n M_n$ where~$M_0=K
\subset \cdots \subset  M_n \subset \cdots$ and~$K\subset M_n$ is
a finite extension for every~$n$. Let~$W_n$ be the ring of
integers of~$M_n$ and let~$\FF_n$ be its residue filed.
Let~$T_1,\ldots,T_d\in R_\cU$ be parameters as in~\ref{V}.
Since~$R_\cU\tensor_V k$ is a smooth $k$--algebra, then
$R_{\cU}\bigl[T_1^{1\over p^n},\ldots, T_d^{1\over
p^n}\bigr]\tensor_V \FF_n$ is a smooth $k$--algebra. Hence,
$R_{\cU}\tensor_V W_n\bigl[T_1^{1\over p^n},\ldots, T_d^{1\over
p^n}\bigr] $ is a regular ring modulo the maximal ideal of~$W_n$
and, hence, it is a regular ring itself. In particular, it is
normal. This implies that~$R_{\cU,M,\infty}\cong \cup_n
R_{\cU}\tensor_V W_n\bigl[T_1^{1\over p^n},\ldots, T_d^{1\over
p^n}\bigr]$.

Since~$\cU'\to \cU$ is formally \'etale, then $R_{\cU'}\tensor_V
k$ is a smooth $k$--algebra. Reasoning as above we conclude
that~$R_{\cU',M,\infty}\cong \cup_n R_{\cU}\tensor_V
W_n\bigl[T_1^{1\over p^n},\ldots, T_d^{1\over p^n}\bigr]$. The
lemma follows.

\

Let\/~$\cU=\Spf(R_\cU)$ be an affine small object of\/~$\cX^\ast$.
Let\/~$A$ be the union of some collection of almost \'etale,
integral $R_{\cU,M,\infty}$--subalgebras of\/~$ \Rbar_\cU$. Write
$$\bigl(\barhatcO_{\hatgerX_M}/p\barhatcO_{\hatgerX_M}\bigr)(A\tensor_V K):=\lim
\bigl(\barhatcO_{\hatgerX_M}/p\barhatcO_{\hatgerX_M}\bigr)(\cU,(\cW,L))$$where
the direct limit is taken  over all~$(\cU,(\cW,L))\in\hatgerX_M$
with\/~$\cW=\Spm(S_\cW)$ such that\/~$S_\cW\tensor_L M \subset
A\tensor_V K$.

\slabel barhatcOMRbar=Rbar. proposition\par\sprop Assume
that\/~$R_{\cU}$ is small over $V$.
Then, the natural map
$$A/pA\llongrightarrow
\bigl(\barhatcO_{\hatgerX_M}/p\barhatcO_{\hatgerX_M}\bigr)(A\tensor_V
K)$$has kernel and cokernel annihilated by any element of the
maximal ideal of\/~$\Vinfty$.
\endsprop
\Proof The presheaf $\barhatcO_{\hatgerX_M}/p
\barhatcO_{\hatgerX_M}$ is separated i.~e., if $(\cU',\cW',L')\to
(\cU,\cW,L)$ is a covering map, the natural map
$$\barhatcO_{\hatgerX_M}(\cU,\cW,L)/p \barhatcO_{\hatgerX_M}(\cU,\cW,L)
\lra \barhatcO_{\hatgerX_M}(\cU',\cW',L')/p
\barhatcO_{\hatgerX_M}(\cU',\cW',L')$$is injective. This implies
that we have an injective map
$$A/pA= \barhatcO_{\hatgerX_M}(A\tensor_V K)
/p\barhatcO_{\hatgerX_M}(A\tensor_V K) \hooklongrightarrow
\bigl(\barhatcO_{\hatgerX_M}/p\barhatcO_{\hatgerX_M}\bigr)(A\tensor_V
K).$$We also get that the sheaf associated to the presheaf
associating to a triple $(\cU,\cW, L)$ the ring $
\barhatcO_{\hatgerX_M}(\cU,\cW,L)/p
\barhatcO_{\hatgerX_M}(\cU,\cW,L)$ is defined by taking
$\bigl(\barhatcO_{\hatgerX_M}/p\barhatcO_{\hatgerX_M}\bigr)(\cU,\cW,L)$
to be the direct limit, over all coverings $(\cU',\cW',L')$
of~$(\cU,\cW,L)$ with~$\cU'$ affine, of the elements $b$ in the
group $
\barhatcO_{\hatgerX_M}(\cU',\cW',L')/p\barhatcO_{\hatgerX_M}(\cU',\cW',L')$
such that the image of~$b$
in~$\barhatcO_{\hatgerX_M}(\cU'',\cW'',L'')/p
\barhatcO_{\hatgerX_M}(\cU'',\cW'',L'')$ is~$0$
where~$(\cU'',\cW'',L'') $ is the fiber product
of~$(\cU',\cW',L')$ with itself over~$(\cU,\cW,L)$. Hence,
$$\bigl(\barhatcO_{\hatgerX_M}/p\barhatcO_{\hatgerX_M}\bigr)(A\tensor_V
K)=\lim_{S,T}\Ker_{S,T}$$where the notation is as follows. The
direct limit is taken over all normal
$R_{\cU,M,\infty}$--subalgebras~$S$ of~$A$, finite and \'etale
after inverting~$p$ over~$R_{\cU,M,\infty}$, all covers~$\cU'\to
\cU$ and all normal extensions~$R_{\cU',M,\infty} \tensor_{R_\cU}
S \to T$, finite, \'etale and Galois after inverting~$p$.
Eventually, we put $\cU'':=\Spf(R_{\cU''})$ to be the fiber
product of~$\cU'$ with itself over~$\cU$ i.~e.,
$R_{\cU''}:=\widehat{R_{\cU'}\tensor_{R_\cU}R_{\cU'}}$. We let
$$\Ker_{S,T}:=\Ker\left(T/p T \maprighto{\llongrightarrow} \widetilde{T}_S/p
\widetilde{T}_S \right),$$where~$\widetilde{T}_S$ is the
normalization of the base change to~$R_{\cU''}$ of
$T\tensor_{(R_{\cU',M,\infty} \tensor_{R_{\cU,M,\infty}} S)} T$.
For every~$S$ and~$T$ as above,
write~$G_{S,T}:=\Gal\bigl(T\tensor_V K/
S\tensor_{R_{\cU,M,\infty}}R_{\cU',M,\infty}\tensor_V K\bigr)$.
Then, $\widetilde{T}_S$  is the product~$\prod_{g\in G_{S,T}}
\widetilde{T\tensor_{R_{\cU'}}R_{\cU''}}$ where tilde stands for
the normalization (of~$T\tensor_{R_{\cU'}} R_{\cU''}$) and we
view~$R_{\cU''}$ as $R_{\cU'}$--algebra choosing the left action.
Hence,
$$\Ker_{S,T}=\Ker\left(T/pT
\maprighto{\llongrightarrow} \prod_{g\in G_{S,T}}
\widetilde{T\tensor_{R_{\cU'}} R_{\cU''}\over p
\widetilde{T\tensor_{R_{\cU'}} R_{\cU''}}}\right).$$The two maps
in the display are $a\mapsto (a,\cdots,a)$ and $a\mapsto
\bigl(g(a)\bigr)_{g\in G_{S,T}}$.

For the rest of this proof we make the following notations: if $B$
is a $R_{\cU,M,\infty}$-algebra we denote by
$B':=B\otimes_{R_{\cU,M,\infty}}R_{\cU',M,\infty}=B\otimes_{R_\cU}R_{\cU'}$,
by
$B'':=B\otimes_{R_{\cU',M,\infty}}R_{\cU'',M,\infty}=B\otimes_{R_{\cU'}}R_{\cU''}$
(the second equalities above follow form
\ref{RcU'inftyisbasechenageofRcU}) and by $\widetilde{B}$ the
normalization of $B$ in $B\bigl[p^{-1}\bigr]$. We then get a
commutative diagram \labelf commdiagS'T'\par
$$\matrix{0 & \to & S/pS & \llongrightarrow &
S'/pS' & \maprighto{\llongrightarrow} & S''/pS'' \cr & &
\big\downarrow & & \mapdownl{\alpha} & &  \mapdownr{\beta} \cr 0 &
\to & \Ker_{S,T} & \llongrightarrow & T/pT &
\maprighto{\llongrightarrow}
&\widetilde{T}_S/p\widetilde{T}_S=\prod_{g\in G_{S,T}}
\bigl({\widetilde{T}''/ p \widetilde{T}''}\bigr).
\cr}\eqno{{(\numfo)}}$$The top row is exact by \'etale descent
and the bottom row is exact by construction. We claim that the
kernel and cokernel of the map $S/pS\lra \Ker_{S,T}$ are
annihilated by any element of the maximal ideal of~$\Vinfty$. To
do this we analyze the maps~$\alpha$ and~$\beta$.\spacing

{\it Analysis of\/ $\Ker(\alpha)$ and\/ $\Ker(\beta)$.}\enspace
Note that he extension~$R_{\cU,M,\infty} \subset S$ is integral
and almost \'etale by~\ref{whenRAEholds}. Hence, the
extensions~$R_{\cU',M,\infty} \subset S'$ and~$R_{\cU'',M,\infty}
\subset S''$ are integral and almost \'etale as well. Since the
extension $R_{\cU}\to R_{\cU'}$ (resp.~$R_\cU\to R_{\cU''}$) is
faithfully flat, the rings $S'$ and~$S''$ have no non--trivial
$p$--torsion. In particular, $S'$ (resp.~$S''$) injects into its
normalization $\widetilde{S}'$ (resp.~$\widetilde{S}''$) which
is~$T^{G_{S,T}}$. Thanks to Lemma~\ref{technicalneededinS6} the
cokernel of $S'\to \widetilde{S}'= T^{G_{S,T}}$ (resp.~$S''\to
\widetilde{S}''$) is annihilated by any element of the maximal
ideal of~$\Vinfty$.

Consider the following. If $0\lra B\lra C\lra D\lra 0$ is an exact
sequence of abelian groups then the kernel of the induced map
$B/pB\to C/pC$ is the image in $B/pB$ of the group of $p$-torsion
elements of $D$. In particular if $B,C,D$ are $\Vinfty$-modules
and $D$ is annihilated by an element $a\in \Vinfty$
then~$\Ker(B/pB\to C/pC)$ is also annihilated by $a$. It follows
from this obvious  fact that the kernel of the map $S'/pS'\to
\widetilde{S}'/p\widetilde{S}'$ and the kernel of the map
$S''/pS''\to \widetilde{S}''/p\widetilde{S}''$ are annihilated by
any element of the maximal ideal of $\Vinfty$.

The map $\widetilde{S}'/ p\widetilde{S}'\to T/ pT$
(resp.~$\widetilde{S}''/ p\widetilde{S}''\to \widetilde{T}_S/
p\widetilde{T}_S $) is injective since $\widetilde{S}'\to T$
(resp.~$\widetilde{S}''\to \widetilde{T}_S$) is an integral
extension of normal rings. Hence, the kernel of~$\alpha$ and the
kernel of~$\beta$ are annihilated by any element of the maximal
ideal of~$\Vinfty$.

\spacing {\it Analysis of the image of\/  $\Coker(S/pS\lra
\Ker_{S,T})$ in $\Coker(\alpha)$.}\enspace Define $Z$ as
$Z:=\Coker\bigl(S'/pS'\to (T/pT)^{\G_{S,T}}\bigr)\subset
\Coker(\alpha)$. Since $\Ker_{S,T}$ is $\G_{S,T}$-invariant (by
definition), the image of\/  $\Coker(S/pS\lra \Ker_{S,T})$ in
$\Coker(\alpha)$ is contained in~$Z$. Put
$Y:=\Coker\bigl(S'/pS'\to \widetilde{S}'/p\widetilde{S}'\bigr)$.
Let us remark that we have an exact sequence of groups:
$$
0\lra Y\lra Z\lra \Coker\bigl(\widetilde{S}'/p\widetilde{S}'
\lra (T/pT)^{\G_{S,T}}\bigr)\lra 0.
$$
We know that $Y$ is annihilated by any element of the maximal
ideal of $V_\infty$, so let us examine the last term of the
sequence. This is the same as
$\Coker\bigl(T^{\G_{S,T}}/pT^{\G_{S,T}}\lra
(T/pT)^{\G_{S,T}}\bigr)$. Consider the exact sequence
$$ 0\llongrightarrow T^{G_{S,T}}/p
T^{G_{S,T}} \llongrightarrow \bigl(T/pT\bigr)^{G_{S,T}}
\llongrightarrow \H^1\bigl(G_{S,T},T\bigr).$$Since
$R_{\cU',M,\infty}\to T$ is almost \'etale, the group
$\H^1\bigl(G_{S,T},T\bigr)$ is annihilated by any element of the
maximal ideal of~$\Vinfty$; see~[\FJAMS, Thm.~I.2.4(ii)]. Hence,
the cokernel of $T^{G_{S,T}}/p T^{G_{S,T}} \tto
\bigl(T/pT\bigr)^{G_{S,T}}$ is annihilated by any element of the
maximal ideal of~$\Vinfty$. We deduce that the same is true for
the module $Z$ above. \spacing

Now using the snake lemma applied to the commutaive
diagram~(\ref{commdiagS'T'}), we get that the kernel and cokernel
of the map $S/pS\lra \Ker_{S,T}$ are annihilated by any element of
the maximal ideal of $V_\infty$ as claimed.

\spacing This concludes the proof in the case that~$A$ is the union of
almost \'etale, integral and normal
$R_{\cU,M,\infty}$--subalgebras of\/~$ \Rbar_\cU$. In the general
case, assume that~$Q$ is an almost \'etale, integral
$R_{\cU,M,\infty}$--subalgebra of~$A$ and let~$S$ be its
normalization. Then, the cokernel of $Q\to S$ annihilated by any
element of the maximal ideal of~$\Vinfty$ by
Lemma~\ref{technicalneededinS6}. The same then applies to the
kernel and the cokernel of $Q/pQ\to S/pS$. The conclusion follows.

\sssection End of proof of\/~\ref{forsomesheavesitholds}\par
Assumption~(i)  clearly holds. We let~$\{\cW_i\}_i=\cS$ be a
covering of~$\cX$  and let~$\cT_i:=\{\cU_{ij}\}_j$ be a basis
of~$\cW_i$ as in~\ref{assumption}(ii). Let~$\cU\in\cT_i$ for
some~$i$.

\spacing

(iii) The group~$\bbL_n(\Rbar_\cU\tensor_V K)$ is constant on the
connected components of~$\cU$ and does not depend on~$\cU$ itself.
It then suffices to verify assumption~(iii) for~$\bbL_n$ the
constant sheaf i.~e, $\bbL_n=\ZZ/p^s\ZZ$ for some~$s$ in case~(A)
or~$\bbL_n=(\ZZ/p^n\ZZ)$ in case~(B). In this case (iii) follows
from~\ref{barhatcOMRbar=Rbar} with~$A=\Rbar_\cU$.

\spacing

(iv) Due to~\ref{barhatcOMRbar=Rbar} it suffices to prove that
$\H^q\bigl(\cH_\cU, \bbL_n(\Rbar_\cU\tensor_V K)\tensor
\WW_n(\Rbar_\cU/p\Rbar_\cU)\bigr)$ (resp.~$\H^q\bigl(\cH_\cU,
\bbL_n(\Rbar_\cU\tensor_V K)\tensor (\barhatcO_{\hatgerX_M}/p^n
\barhatcO_{\hatgerX_M})\bigr)$ is annihilated by any power of the
maximal ideal of~$\WW_n\bigl(\Vinfty/p\Vinfty)$ (resp.~$\Vinfty$)
for every~$q\geq 1$. In both cases, one reduces by devissage to
the case~$n=1$. The claim then follows from~\ref{induced}
and~\ref{variantcomputecoh}.

\spacing

(v) Given~$n\in\NN$, let~$R_\cU\tensor_V K \subset B$ be a finite
and \'etale extension such that~$\bbL_{n+1}$ and~$\bbL_n$ are
constant on the \'etale site of~$B_M$. We may assume that~$B$ is
defined over a finite extension~$K\subset L$ contained in~$M$ and
that~$R_\cU\tensor_V M \subset B\tensor_L M$ is a Galois extension
of integral domains. Define~$A_{\cU}$ as the normalization
of~$R_{\cU}$ in the subring of~$\Rbar_{\cU}\tensor_V M$ generated
by~$R_{\cU,M,\infty}$ and~$B$. Let~$G$ be the Galois group
of~$A_\cU\tensor_V K $ over~$R_{\cU,M,\infty}\tensor_V K $.

Then, assumption~(v), with~$A_\cU$ in place of~$R_{\cU,M,\infty}$,
holds due to~\ref{barhatcOMRbar=Rbar} since we may reduce to the
that~$\bbL$ is trivial. Let~$D_n$ be the kernel
of~$\cF_{n+1}(A_\cU\tensor_V K) \to \cF_n(A_\cU\tensor_V K)$.
Let~$E_n$ be the tensor product of the kernel of~$\bbL_{n+1}(B)\to
\bbL_n(B)$ with
$\WW_{n+1}\bigl(\barhatcO_{\hatgerX_M}/p\barhatcO_{\hatgerX_M}\bigr)(A_{\cU}\tensor_V
K)$ (resp.~$\barhatcO_{\hatgerX_M}/p^n
\barhatcO_{\hatgerX_M})(A_{\cU}\tensor_V K)$). Let~$F_n$ be the
tensor product of~$L_{n+1}(B)$ with the kernel of Frobenius on $
\WW_n\left(\barhatcO_{\hatgerX_M}/p \barhatcO_{\hatgerX_M}
\right)(A_{\cU}\tensor_V K)$ i.~e.,
$\WW_n\left(\barhatcO_{\hatgerX_M}/p^{1\over
p}\barhatcO_{\hatgerX_M}\right) (A_{\cU}\tensor_V K)$, in
case~(A). In case~(B) define~$F_n$ as the tensor product
of~$L_{n+1}(B)$ with the kernel of the  projection from
$(\barhatcO_{\hatgerX_M}/p^{n+1}\barhatcO_{\hatgerX_M})(A_{\cU}\tensor_V
K)$ to $
(\barhatcO_{\hatgerX_M}/p^n\barhatcO_{\hatgerX_M})(A_{\cU}\tensor_V
K)$ i.~e.,
$(\barhatcO_{\hatgerX_M}/p\barhatcO_{\hatgerX_M})(A_{\cU}\tensor_V
K)$.

Then, $D_n$ is generated by the images of~$E_n$ and~$F_n$
in~$\cF_{n+1}(A_\cU\tensor_V K)$.  It follows from~\ref{induced}
and~\ref{variantcomputecoh}  that~$\H^q(\cH_\cU,E_n)$
and~$\H^q(\cH_\cU,F_n)$ are annihilated by any element of the
maximal ideal of~$\WW\bigl(\Vinfty/p\Vinfty)$ (resp.~$\Vinfty$)
for~$q\geq 1$. Thus, the same applies to~$\H^q(\cH_\cU,D_n)$ and,
hence, to the cokernel of the map from $\cF_{n+1}(A_\cU\tensor_V
K)^{\cH_\cU}=\cF_{n+1}(R_{\cU,M,\infty}\tensor_V K)$ to
$\cF_n(A_\cU\tensor_V K)^{\cH_\cU}=\cF_n(R_{\cU,M,\infty}\tensor_V
K)$. This concludes the proof of~(v).

\spacing

(vi) For every covering~$\cZ\to \cU$ in~$\cX^\ast$ with~$\cZ\in
\cT_i$ define~$\H^q_{\cF_n}(\cZ\to \cU)$ as the Chech cohomology
group
$$\H^q_{\cF_n}(\cZ\to \cU):=\H^q\bigl(\cZ\to \cU,\bbL_n(B)\tensor
\WW_n(A_{\cU}\tensor_{R_\cU} R_{\cZ}/ p R_{\cZ}) \bigr)$$
respectively $$\H^q_{\cF_n}(\cZ\to \cU):=\H^q\bigl(\cZ\to
\cU,\bbL_n(B)\tensor (A_{\cU}\tensor_{R_\cU} R_{\cZ}/ p^n R_{\cZ})
\bigr).$$See the proof of~(v) for the notation. For every~$q\geq
1$ the group~$\H^q_{\cF_n}(\cZ\to \cU)$ is~$0$ since the sheaves
considered are quasi--coherent.

Due to~\ref{barhatcOMRbar=Rbar} we conclude that assumption~(vi)
holds using~$A_{\cU}\tensor_{R_\cU} R_{\cZ}\tensor_V K$ instead
of~$R_{\cZ,M,\infty}$. Consider the spectral sequence
$$\H^p\bigl(G,\H^q\bigl(\cZ\to \cU, \cF_n(A_{\cU}\tensor_{R_\cU}
(R_{\cZ}\tensor_V K)\bigr)\Longrightarrow \H^{p+q}\bigl(\cZ\to
\cU,\cF_n(R_{\cZ,\infty}\tensor_V K)\bigr).$$Then, up to
multiplication by any element of the maximal ideal
of~$\WW\bigl(\Vinfty/p\Vinfty)$ (resp.~of $\Vinfty$), the group
$\H^p\bigl(\cZ\to \cU,\cF_n(R_{\cZ,\infty}\tensor_V K)\bigr)$ is
isomorphic to the group $\H^p\bigl(G,\H^0\bigl(\cZ\to
\cU,\cF_n(A_{\cU}\tensor_{R_\cU} (R_{\cZ}\tensor_V K)\bigr)\bigr)$
i.~e., $\H^p\bigl(G,\cF_n(A_\cU\tensor_V K)\bigr)$. Let~$C$ be the
kernel of the surjective map $\cH_\cU\to G$. Consider the spectral
sequence
$$\H^p\bigl(G,\H^q(C,\bbL_n(B)\tensor \Rbar_\cU/p
\Rbar_\cU\bigr)\Longrightarrow \H^{p+q}\bigl(\cH_U,
\bbL_n(B)\tensor \Rbar_\cU/p \Rbar_\cU\bigr).$$Note that
$\H^q(C,\bbL_n(B)\tensor \Rbar_\cU/p \Rbar_\cU\bigr)$ and
$\H^q(\cH_\cU,\bbL_n(B)\tensor \Rbar_\cU/p \Rbar_\cU\bigr)$ are
annihilated by multiplication by any element of the maximal ideal
of~$\WW\bigl(\Vinfty/p\Vinfty)$ (resp.~$\Vinfty$) for~$q\geq 1$
due to~\ref{induced} and~\ref{variantcomputecoh}. Hence, the same
must hold for $\H^q\bigl(G,\bbL_n(B)\tensor A_{\cU}/p
A_{\cU}\bigr)$. By devissage and~\ref{barhatcOMRbar=Rbar} one
concludes that the same must hold
for~$\H^q\bigl(G,\cF_n(A_\cU\tensor_V K)\bigr)$. Thus, (vi) holds.

\endsssection

\endssection

\label proof3. section\par\ssection Proof of theorem
\ref{cohomologyFontaine}\par By theorem
\ref{provethatlimHiUFisRivFcont} if $\cF$ is a sheaf of
$\Sh(\hatgerX_M^\ast)^\NN$ such that the assumptions
\ref{assumption} are satisfied then $\R^i\gerv_{\cX,M,\ast}
(\cF)\bigl[\pi^{-1}\bigr] \cong \cH_{\Gal_M}^i(\cF)$. Using this
isomorphism  the Leray spectral sequence for the composition of
functors $\H^0(\cX^\ast, - )\circ\gerv_{\cX,M,\ast}$ becomes
$$
E_2^{p,q}=\H^q\bigl(\cX^\ast, \cH^p_{\Gal_M}(\cF)\bigr)\Longrightarrow \H^{p+q}(\hatgerX_M^\ast, \cF).
$$
In particular, we obtain a spectral sequence for $\ast=\bullet$.
Now theorem \ref{cohomologyFontaine} follows as the functors
$\H^i(\cX^{et,\bullet}, - )$ and $\H^i(\cX^{et}, - )$ are
canonically isomorphic, same as the functors
$\H^i(\hatgerX_M^\bullet, - )$ and $\H^i(\hatgerX_M, - )$;
see~\ref{okforpointed}.

\endssection
\endsection

\section Appendix I: Galois cohomology via Tate--Sen's method\par
The goal of this section is to prove Proposition~\ref{induced}
stating that, if\/~$M$ is a $\ZZ_p$-representation of\/~$\cG_S$,
then the groups $\H^i(\cH_S,\cD(M)\tensor_{\AA_S}\AA_\Rbar)$,
$\H^i(\cH_S,\cDtilde(M)\tensor_{\AAtilde_S}\AAtilde_\Rbar)$,
$\H^i(\H_S,\Dtilde(M)\tensor_{\tAA_S}\AA_\Rbar)$ and
$\H^i(\H_S,\Dtilde(M)\tensor_{\tAAtilde_S}\tAAtilde_\Rbar)$ are
trivial for~$i\geq 1$.  This will be the key tool to compute the
Galois cohomology of~$M$ in terms of the associated
$(\phi,\Gamma_S)$--modules.

To treat all the cases above, we follow the axiomatic approach
started by Colmez in  [\ColmezCMP, \S3.2\&3.3] and developed
in~[\AnBr, \S2].

\label axioms. section\par\ssection The Axioms\par Let~$\cG$ be a
profinite group and~$\cH$ a closed normal subgroup of~$\cG$ such
that~$\Gamma=\cG/\cH$ is endowed with a continuous homomorphism
$\chi\colon\Gamma\to\ZZ_{p}^{\times}$ with open image and kernel
isomorphic to $\ZZ_{p}^{d}$. Suppose that $\gamma
g\gamma^{-1}=g^{\chi(\gamma)}$ for every~$g\in\Ker(\chi)$ and
every~$\gamma\in\Gamma$. Let\/~$\G\subset \cG$ be a closed normal
subgroup, put\/~$\H:=\G\cap \cH$ and assume that~$\tGamma:=\G/\H$
is~$\Ker(\chi)$.

Let $\widetilde{\Lambda}$ be $\ZZ_{p}$-algebra which is an
integral domain and is endowed with a map
$v\colon\widetilde{\Lambda}\to\RR\cup\{ +\infty\}$ such that:

\item{{\rm (i)}} $v(x)=+\infty\Leftrightarrow x=0$; \item{{\rm
(ii)}} $v(xy)\geq v(x)+v(y)$; \item{{\rm (iii)}}
$v(x+y)\geq\min(v(x),v(y))$; \item{{\rm (iv)}} $v(p)>0$ and
$v(px)=v(p)+v(x)$.

We endow $\widetilde{\Lambda}$ with the (separated) topology
induced by~$v$. We assume that $\widetilde{\Lambda}$ is complete
for this topology and that it is endowed with a continuous action
of~$\cG$ such that $v(g(x))=v(x)$ for $x\in\widetilde{\Lambda}$
and for $g\in \cG$.

\spacing

We introduce the following axioms  \'a la {\sl Tate-Sen}:

\spacing

\item{{\rm (TS1)}} there exists $c_{1}\in\RR_{>0}$ such that for
every open normal subgroups $H_{1}\subseteq H_{2}$ of $\cH$
(resp.~of $\H$), there exists
$\alpha\in\widetilde{\Lambda}^{H_{1}}$ such that
$v(\alpha)>-c_{1}$ and $\sum\limits_{\tau\in
H_{2}/H_{1}}\tau(\alpha)=1$;\spacing

\noindent Write $\Ker(\chi)$ as a subgroup of
$\dirsum_{i=1}^{d}\ZZ_{p}\gamma_{i}$ so that there is~$m_0\in\NN$
with $p^{m_0} \dirsum_{i=1}^{d}\ZZ_{p}\gamma_{i} \subset
\Ker(\chi)$ and let~$\gamma_0\in\Gamma$ be such
that~$\Im(\chi)=\ZZ_p\chi(\gamma_0)\dirsum F$ with~$F$ a finite
group. Let~$\cH' \subset \cH$ be an open normal subgroup. Assume
that there exists an integer $m_{0,\cH'}\geq m_0$ and a lifting
$p^{m_{0,\cH'}-m_0}\Gamma \subset \cG/\cH'$ as a normal subgroup
centralizing~$\cH/\cH'$ and that for every~$i\in\{ 0,\ldots,d\}$
one has an increasing sequence $\left(
\Lambda_{m,\cH'}^{(i)}\right) _{m\geq m_{0,\cH'}}$ of {\sl closed}
subrings of $\widetilde{\Lambda}^{\cH'}$ stable under~$\cG/\cH'$
and  maps
$\left(\tau^{(i)}_{m,\cH'}\colon\widetilde{\Lambda}^{\cH'}\to\Lambda_{m,\cH'}^{(i)}\right)_{m\geq
m_{0,\cH'}}$ such that: \spacing

\item{{\rm (TS2)}} for every $i$ and~$j\in\{ 0,\ldots,d\}$ and
every $m\geq m_{0,\cH'}$:

\itemitem{{\rm (a)}} $\tau^{(i)}_{m,\cH'}$ is
$\Lambda_{m,\cH'}^{(i)}$-linear and $\tau^{(i)}_{m,\cH'}(x)=x$ if
$x\in\Lambda_{m,\cH'}^{(i)}$; \itemitem{{\rm (b)}} there exists
$c_{2,\cH'}\in\RR_{>0}$ such that for every
$x\in\widetilde{\Lambda}^{\cH'}$, one has $v\left(
\tau^{(i)}_{m,\cH'}(x)\right) \geq v(x)-c_{2,\cH'}$ and
$\displaystyle\lim_{m\to\infty}\tau^{(i)}_{m,\cH'}(x)=x$ ;
\itemitem{{\rm (c)}} $\tau^{(i)}_{m,\cH'}$ commutes with
$\tau^{(j)}_{m^{\prime},\cH'}$ and with the action of $\cG/\cH'$;
\itemitem{{\rm (d)}} for every open normal subgroup $\cH''\subset
\cH'$ we have $\Lambda_{m,\cH'}^{(i)} \subset
\Lambda_{m,\cH''}^{(i)} $, as subrings of~$\widetilde{\Lambda}$,
and the following diagram commutes

$$\matrix{
\widetilde{\Lambda}^{\cH'}& \llongmaprighto{\tau^{(i)}_{m,\cH'}} &
\Lambda_{m,\cH'}^{(i)}\cr \big\downarrow & & \big\downarrow \cr
\widetilde{\Lambda}^{\cH''}& \llongmaprighto{\tau^{(i)}_{m,\cH''}}
& \Lambda_{m,\cH''}^{(i)}};$$

\spacing

\item{{\rm (TS3)}} let $X_{m,\cH'}^{(i)}=\left(
1-\tau^{(i)}_{m,\cH'}\right) \left(
\widetilde{\Lambda}^{\cH'}\right)$. Then, $1-\gamma_{i}^{p^{m}}$
is invertible on $X_{m,\cH'}^{(i)}$ and there exists
$c_{3,\cH'}\in\RR_{>0}$ such that for every $x\in
X_{m,\cH'}^{(i)}$, one has $v\left( \left(
1-\gamma_{i}^{p^{m}}\right) ^{-1}(x)\right) \geq v(x)-c_{3,\cH'}$.
Furthermore, there exists $c_{4,\cH'}\in\RR_{>0}$ such that for
every $i\in\{ 1,\ldots,d\}$ and every
$x\in\Lambda^{(i)}_{m,\cH'}$, one has $v\left( \left(
\gamma_{i}^{p^{m}}-1\right) (x)\right) \geq v(x)+c_{4,\cH'}$.

\item{{\rm (TS4)}} Let~$\H' \subset \H$ be an open normal
subgroup. Assume that there exists an integer $m_{0,\H'}\geq m_0$
and a lifting $p^{m_{0,\H'}-m_0}\Gamma' \subset \G/\H'$ as a
normal subgroup centralizing~$\H/\H'$ and that for every~$i\in\{
1,\ldots,d\}$ one has an increasing sequence $\left(
\Lambda_{m,\H'}^{(i)}\right) _{m\geq m_{0,\H'}}$ of {\sl closed}
subrings of $\widetilde{\Lambda}^{\H'}$ stable under~$\G/\H'$ and
maps
$\left(t^{(i)}_{m,\H'}\colon\widetilde{\Lambda}^{\H'}\to\Lambda_{m,\H'}^{(i)}\right)_{m\geq
m_{0,\cH'}}$ such  that the analogues of~(TS2) and~(TS3) hold.

\endssection

\noindent We followed closely~[\AnBr, \S2] with the differences
that we added~(TS4) and that we require~(TS2) and~(TS3) to hold
not only for~$\cH$ as in loc.~cit., but for every open normal
subgroup of~$\cH$ as well.

\label notDLambdatilde. section\par\ssection Notation\par Let~$W$
be a free $\widetilde{\Lambda}$--module of finite rank~$a$ i.~e.,
$W:=\widetilde{\Lambda}^a$. We consider it as a topological module
with respect to the (separated) topology defined as the product
topology considering on~$\widetilde{\Lambda}$ the $v$--adic
topology. Note that such topology is independent of the choice
of~$\widetilde{\Lambda}$--basis of~$W$. For every positive~$n\in
\QQ$ write~$\widetilde{\Lambda}_{\geq n}$ for the subgroup
of~$\widetilde{\Lambda}$ consisting of elements~$x$ such
that~$v(x)\geq n$. They are a fundamental system of neighborhoods
for the topology on~$\widetilde{\Lambda}$ for~$n\to \infty$.
Let~$W_{\geq n}$ be the image of $\widetilde{\Lambda}_{\geq n}^a$
in~$W$; they form a fundamental system of neighborhoods for the
given topology on~$W$. Assume that~$W$ is endowed with a
continuous action of~$\cG$. We consider continuous cohomology of a
closed subgroup~$H'$ of~$\cG$ with values in~$W$. If~$f\in
C^r\bigl(H',W\bigr)$ is a continuous cochain, with~$r\geq 0$ and
with the profinite topology on~$H'$, write~$v(f):=\min\bigl\{n\in
\NN\vert f(g_1,\ldots,g_r)\in W_{\geq n} \forall g_1,\ldots,g_r\in
H'\bigr\}$. We write $\partial \colon C^r\bigl(H',W\bigr) \to
C^{r+1}\bigl(H',W\bigr)$ for the boundary map.

\endssection

\label Hhasapproximatedtricialcoho. lemma\par \lemma {\rm [\Tate,
\S3.2]}\enspace  Let\/~$H_0$ be an open subgroup of\/~$H$
(resp.~of\/~$\cH$) and let\/~$f$ be an $r$--cochain of\/~$H_0$
with values in~$W$ for~$r\geq 1$.

\item{{\rm (1)}} Assume that there exists an open normal
subgroup~$H_1\subset H_0$ such that\/~$f$ factors via an
$r$--cochain of\/~$H_0/H_1$. Then, there exists an $r-1$--cochain
$h$ of\/~$H_0/H_1$ with values in~$W$ such
that\/~$v\bigl(f-\partial h\bigr) > v(\partial f)-c_1$ and\/~$v(h)
> v(f)-c_1$.

\item{{\rm (2)}} There exists a sequence of open normal
subgroups~$H_n\subset H_0$ and continuous cochains $f_n\in
C^r\bigl(H_0/H_n,W\bigr)$ for~$n\in \NN$ such that $f\equiv f_n$
modulo~$W_{\geq n}$ for~$n\to \infty$.
\endlemma
\Proof We work out the case of $H_0\subset H$. For $H_0\subset
\cH$ the argument is analogous and the details are left to the
reader.

(1) Let~$\alpha\in \widetilde{\Lambda}^{H_1}$ be an element
satisfying~(TS1).  Define the $r-1$--cochain $\alpha \cup f$
of~$H_0/H_1$ with values in~$W$ by
$$\bigl(\alpha \cup f\bigr) \bigl(g_1,\ldots,g_{r-1}\bigr):=(-1)^r \sum_{t\in
H_0/H_1} g_1 \cdots g_{r-1} t(\alpha)\cdot
f\bigl(g_1,\ldots,g_{r-1},t\bigr).$$One computes that
$$\eqalign{\partial(\alpha\cup f)(g_1,\ldots,g_r) & = g_1 \bigl((\alpha\cup
f)(g_2,\ldots,g_r)\bigr) +\sum_{j=1}^{r-1} (-1)^j (\alpha\cup
f)(\ldots,g_j g_{j+1},\ldots)+ \cr & + (-1)^r (\alpha\cup
f)(g_1,\ldots,g_{r-1})=\cr & = (-1)^r \sum_{t\in H_0/H_1} g_1
\cdots g_r t(\alpha)\cdot g_1 f\bigl(g_2,\ldots,g_r,t\bigr)+ \cr &
+ (-1)^r \sum_{j=1}^{r-1}  \sum_{t\in H_0/H_1} (-1)^j g_1 \cdots
g_r t(\alpha)\cdot f\bigl(g_1,\ldots,g_jg_{j+1},g_r,t\bigr)+ \cr &
+ \sum_{t\in H_0/H_1} g_1 \cdots g_{r-1} t(\alpha)\cdot
f\bigl(g_1,\ldots,g_{r-1},t\bigr)\cr}$$and
$$\eqalign{ \bigl(\alpha \cup\partial f\bigr) (g_1,\ldots,g_r) & = (-1)^{r+1} \sum_{t\in
H_0/H_1} g_1 \cdots g_r t(\alpha)\cdot
\partial f\bigl(g_1,\ldots,g_r,t\bigr)=\cr & = (-1)^{r+1} \sum_{t\in
H_0/H_1} g_1 \cdots g_r t(\alpha)\cdot g_1 f(g_2,\ldots,g_r, t) +
\cr & + (-1)^{r+1} \sum_{t\in H_0/H_1} g_1 \cdots g_r
t(\alpha)\cdot \sum_{j=1}^{r-1}(-1)^j
f\bigl(g_1,\ldots,g_jg_{j+1},g_r,t\bigr)+\cr & - \sum_{t\in
H_0/H_1} g_1 \cdots g_r t(\alpha)\cdot f\bigl(g_1,\ldots, g_r
t\bigr) +\cr &+  \sum_{t\in H_0/H_1} g_1 \cdots g_r t(\alpha)\cdot
f(g_1,\ldots,g_r) \cr}$$Since $\sum_{t\in H_0/H_1} t(\alpha)=1$,
we have $\bigl(\alpha \cup\partial
f\bigr)=f-\partial\bigl(\alpha\cup f\bigr)$.  Put~$h=\alpha\cup
f$. Then, $v(h) > v(f)-c_1$ and $v\bigl(f-\partial h\bigr)=v\bigl(
\alpha \cup\partial f\bigr)\geq v(\alpha)+ v(\partial f)$.

(2)  Since~$f$ is continuous  there exists an open normal
subgroup~$H_n$ such that the composite $\bar{f}_n\colon H_0^r \to
W \to W/W_{\geq n}$ factors via~$\bigl(H_0/H_n\bigr)^r$. Let~$f_n$
be the composite of $\bar{f}_n$ with a splitting $W/W_{\geq n}\to
W$ (as sets). Then, $f_n$ is a continuous cochain and
$v(f-f_n)\geq n$.

\label computecoh. proposition\par \prop {\rm [\Tate,
Prop.~10]}\enspace We have $\H^r\bigl(\H,W\bigr)=0$ for~$r\geq 1$
and\/ $\H^r\bigl(\cH,W\bigr)=0$ for~$r\geq 1$. In particular,
$\H^r\bigl(G,W\bigr)=\H^r\bigl(\Gamma,W^H\bigr)$ and\/
$\H^r\bigl(\cG,W\bigr)=\H^r\bigl(\tGamma,W^\cH\bigr)$.
\endprop
\Proof The last statement follows from the first one and from the
spectral sequences $\H^j(\Gamma,\H^j(\cH,\_))\Longrightarrow
\H^{i+j}(\cG,\_)$ and $\H^j(\tGamma,\H^j(\H,\_))\Longrightarrow
\H^{i+j}(\G,\_)$. Let~$H_0:=H$ or~$\cH$. Let~$f$ be an $r$--th
cochain of~$H_0$, for~$r\geq 1$, with values in~$W$.
Let~$\{H_n,f_n\}_n$ be as in~\ref{Hhasapproximatedtricialcoho}(2)
and, for each~$n$, write~$h_n$ for the continuous $r-1$-cochain
satisfying~\ref{Hhasapproximatedtricialcoho}(1) i.~e.,
$v\bigl(f_n-\partial h_n\bigr) > v(\partial f_n)-c_1$
and\/~$v(h_n)
> v(f_n)-c_1$. Then, $\{h_n\}$ is Cauchy so that
it converges to a continuous $r-1$--cochain~$h$. Furthermore,
$v(f_n-\partial h_n) \geq n-c_1$ for every~$n$ so that~$\partial
h_n \to f$ for~$n\to \infty$. We conclude that~$f=\partial h$ as
claimed.

\

Let~$\widetilde{\Lambda}_{\geq n}$ be the subset
of~$\widetilde{\Lambda}$ consisting of elements~$b$ such
that~$v(b)\geq n$. Then, $\widetilde{\Lambda}_{\geq 0}$ is a ring
and $\widetilde{\Lambda}_{\geq n}$ is an ideal for every~$n\geq 0$
due to the properties of~$v$. We write $\overline{\Lambda}_n$ for
 the quotient $\widetilde{\Lambda}_{\geq 0}/
\widetilde{\Lambda}_{\geq n}$. Assume that the following
strengthening of (TS1) holds:\spacing

\item{{\rm (TS1')}} for every $c \in\RR_{>0}$ and for every open
normal subgroups $H_{1}\subseteq H_{2}$ of $\cH$ (resp.~of $\H$),
there exists $\alpha\in\widetilde{\Lambda}_{\geq 0}^{H_{1}}$ such
that $v\left(\sum\limits_{\tau\in
H_{2}/H_{1}}\tau(\alpha)\right)\leq c$;\spacing

\noindent One then has the following variant of~\ref{computecoh}:

\label variantcomputecoh. proposition\par\prop Let\/~$W$ be a free
$\overline{\Lambda}_n$--module of finite rank~$a$ endowed with a
continuous action of\/~$\cG$. Then, for every~$c\in \RR_{\geq 0}$
and every integer~$r\geq 1$ there exists an element~$\gamma_c\in
\widetilde{\Lambda}_{\geq 0}^\cH$ of valuation~$v(\gamma_c)<c$
such that $\gamma_c\cdot \H^r\bigl(\H,W\bigr)=0$ and\/ $\gamma_c
\cdot \H^r\bigl(\cH,W\bigr)=0$.
\endprop

\label  defDicDi. section\par \ssection Decompletion\par The
notation is as in~\ref{notDLambdatilde}. Write~$\cD(W):=W^\cH$
and~$\D(W):=W^\H$. They are closed subgroups of~$W$ endowed with
the topology induced from~$W$.

It is proven in~[\AnBr, Cor.~2.3] that~(TS1) implies that there
exists an open normal subgroup~$\cH_W\subset \cH$ and a
$\widetilde{\Lambda}$--basis $e_1,\ldots,e_a$ of~$W$ such
that~$W^{\cH_W}\cong \widetilde{\Lambda}^{\cH_W} e_1\dirsum
\cdots\dirsum \widetilde{\Lambda}^{\cH_W} e_a$. For
every~$i=0,\ldots,d$ and every~$m\geq m_W=m_{0,\cH_W}$ define the
map $\tau_{m,\cH_W}^{(i)}\colon W^{\cH_W} \rightarrow W^{\cH_W}$
by~$\sum_{i=1}^a \beta_i e_i \mapsto \sum_{i=1}^a
\tau_{m,\cH_W}^{(i)}(\beta_i) e_i$. Due to (TS2)(c)\&(d) such map
is independent of~$\cH_W$ and the basis~$e_1,\ldots,e_a$ and it
descends to a map on~$\cD(W)=W^\cH$. Due to~(TS2)(b) it is
continuous for the topology on~$ W^{\cH_W}$ induced from~$W$. We
then drop the index~$\cH_W$ and we write simply
$$\tau_m^{(i)}\colon \cD(W) \llongrightarrow \cD(W)\quad \hbox{{\rm for }} i=0,\ldots,d, \quad m\geq m_W.$$
Using~(TS4) and repeating the construction above, we get similarly
continuous maps
$$ t_m^{(i)}\colon \D(W) \llongrightarrow \D(W) \quad \hbox{{\rm for }}
i=1,\ldots,d, \quad m\geq m_W.$$For every~$m\geq m_M$ due to~(TS2)
we have a decomposition
$$\cD(W):=\cD_m(W)\dirsum \cD^{(0)}_m(W)\dirsum \cdots
\dirsum  \cD^{(d)}_m(W),$$where
$\cD^{(d)}_m(W):=\bigl(1-\tau^{(d)}_m\bigr)\bigl(\cD(W)\bigr)$,
$\cD^{(d-1)}_m(W):=\bigl(1-\tau^{(d-1)}_m\bigr)\left(\cD(W)_{\tau^{(d)}_m=1}\right)$,
$\ldots$, $\cD^{(0)}_m(W):=\bigl(1-\tau^{(0)}_m\bigr)
\left(\cD(W)_{\tau^{(d)}_m=1,\ldots,\tau^{(1)}_m=1}\right)$
and\/~$\cD_m(W)=\cD(W)_{\tau^{(d)}_m=1,\ldots,\tau^{(0)}_m=1}$.
They are closed $\Gamma$--submodules of~$\cD(W)$. We endowed them
with the induced topology. By~(TS2) the decomposition above is an
isomorphism of topological $\Gamma$--modules.  Similarly, we have
an isomorphism of topological $\tGamma$--modules
$$\D(W):=\D_m(W)\dirsum \D^{(1)}_m(W)\dirsum \cdots \dirsum
\D^{(d)}_m(W),$$where
$\D^{(d)}_m(W):=\bigl(1-t^{(d)}_m\bigr)\bigl(\cD(W)\bigr)$,
$\ldots$, $\D^{(1)}_m(W):=\bigl(1-t^{(1)}_m\bigr)
\left(\D(W)_{\tau^{(d)}_m=1,\ldots,\tau^{(2)}_m=1}\right)$
and\/~$\D_m(W)=\D(W)_{t^{(d)}_m=1,\ldots,t^{(1)}_m=1}$ are closed
$\tGamma$--submodules of~$\D(W)$.

\endssection

\label decompletionofcoho. proposition\par\prop There exists an
integer~$N\geq m_W$ such that for every~$n\geq N$,
if~$\gamma_i^{p^n}\in \Gamma$ the map $\gamma_i^{p^n}-1$ is
bijective with continuous inverse on $\cD^{(i)}_m(W)$
for~$i=0,\ldots,d$ (resp.~$\D^{(i)}_m(W)$ for~$i=1,\ldots,d$).

Then, the maps of continuous cohomology groups
$\H^j\bigl(\Gamma,\cD_m(W)\bigr)\to \H^j\bigl(\Gamma,
\cD(W)\bigr)$ and\/~$\H^j\bigl(\tGamma,\D_m(W)\bigr)\to
\H^j\bigl(\tGamma, \D(W)\bigr)$ are isomorphisms.
\endprop
\Proof We deduce from the first statement that $\H^j\left(\ZZ_p
\gamma_i^{p^n},\cD^{(i)}_m(W)\right)=0$ and that $\H^j\left(\ZZ_p
\gamma_i^{p^n},\D^{(i)}_m(W)\right)=0$ for every~$j\geq 0$. We get
from  the Hochschild---Serre spectral sequence that
$\H^j\left(\Gamma,\cD^{(i)}_m(W)\right)=0$ and
$\H^j\left(\tGamma,\D^{(i)}_m(W)\right)=0$ for~$i\geq 1$.
For~$i=0$ the first statement implies that $\gamma_0^{p^n}-1$ is
bijective with continuous inverse on the group
$\H^j(\tGamma,\cD^{(0)}_m(W)$. By Hochschild--Serre
$\H^j\left(\Gamma,\cD^{(i)}_m(W)\right)=0$ for~$i$ as well. The
second statement follows.

Since~$\gamma_i^{p^t}-1=(\gamma_i^{p^s}-1)
\bigl(\sum_{j=0}^{p^{t-s}-1} \gamma_i^{p^sj}\bigr)$ for~$t\geq s$,
if~$\gamma_i^{p^t}-1$ is bijective with continuous inverse
on~$\cD^{(i)}_m(W)$ (resp.~$\D^{(i)}_m(W)$)
also~$\gamma_i^{p^s}-1$ is. Hence, it suffices to prove that
$\gamma_i^{p^m}-1$ is invertible with continuous inverse.

We prove the statement for~$\cD^{(i)}_m(W)$. The proof
for~$\D^{(i)}_m(W)$ is similar and the details are left to the
reader. Write $W^{\cH_W}\cong \widetilde{\Lambda}^{\cH_W}
e_1\dirsum \cdots\dirsum \widetilde{\Lambda}^{\cH_W} e_a$ as
in~\ref{defDicDi} and write
$W^{\cH_W,(i)}_m:=\bigl(1-\tau^{(i)}_m\bigr)
\left(W^{\cH_W}_{\tau^{(i-1)}_m=1,\ldots,\tau^{(1)}_m=1}\right)$.
Due to the assumptions in~\ref{axioms} we have a lifting
$p^{m_{0,\cH_W}-m_0}\Gamma\subset \cG/\cH_W$ so
that~$p^{m_{0,\cH_W}-m_0}\Gamma$ and~$\cH/\cH_W$ commute.
Since~$\cD^{(i)}_m(W)=\bigl(W^{\cH_W,(i)}_m\bigr)^{\cH}$
by~(TS2)(c) it then suffices to prove that $\gamma_i^{p^m}-1$ is
invertible with continuous inverse on~$W^{\cH_W,(i)}_m$.

Extend $v$ on~$\widetilde{\Lambda}e_1\dirsum \cdots
\widetilde{\Lambda}e_a$ by~$v\bigl(\sum_{j=1}^a z_j
e_j\bigr):=\inf\{v(z_j)\vert j=1,\ldots,a\}$. It defines the weak
topology on~$\widetilde{\Lambda}e_1\dirsum \cdots
\widetilde{\Lambda}e_a$. Since the action of~$\cG/\cH_W$
on~$W^{\cH_W}$ is continuous, there exists an integer~$N\geq
m_{0,\cH_W}$ such that~$\gamma_i^{p^N}$ acts trivially
on~$W^{\cH_W}/W_{\geq c_{3,\cH_W}+1}^{\cH_W}\cong \dirsum_{j=1}^a
\widetilde{\Lambda}^{\cH_W}/ \widetilde{\Lambda}^{\cH_W}_{\geq
c_{3,\cH_W}+1} e_j$. Take~$m\geq N$. Following~[\CCINV,
Prop.~II.6.4] define
$$f_i\colon W^{\cH_W,(i)}_m \llongrightarrow W^{\cH_W,(i)}_m, \quad
f_i\bigl(\sum_{j=1}^a z_j e_j\bigr):=\sum_{j=1}^a
\bigl(1-\gamma_i^{p^m}\bigr)^{-1}(z_j) e_j.$$It is well defined,
continuous, bijective and with continuous inverse (for the weak
topology) due to~(TS3). Then,
$z-f_i\bigl((1-\gamma_i^{p^m})(z)\bigr)=-f_i\left(\sum_{j=1}^a
\gamma_i^{p^m}(z_j) \bigl(1-\gamma_i^{p^m}\bigr)(e_j) \right)$.
Write $$g_{i,z}\colon W^{\cH_W,(i)}_m \llongrightarrow
W^{\cH_W,(i)}_m, \quad
g_{i,z}(y):=y-f_i\left(\bigl(1-\gamma_i^{p^m}\bigr)(y)-z\right).
$$  Then,
$v\bigl(f_i(z)\bigr)\geq v(z)-c_{3,\cH_W}$ by~(TS3) and
$v\bigl(g_{i,0}(y)\bigr)\geq v(y)+
\inf\{v\bigl((1-\gamma_i^{p^m})(e_i)\bigr)\vert
i=1,\ldots,a\}-c_{3,\cH_W}\geq 1$. Hence,
$v\bigl(g_{i,z}(y_1)-g_{i,z}(y_2)\bigr)=v\bigl(g_{i,0}(y_1-y_2)
\geq v(y_1-y_2) +1$. This implies that~$g_{i,z}$ is a contracting
operator for the $v$--adic topology so that there exists a unique
fixed point~$y_z$. Since~$f_i$ is bijective, we get that~$y_z$ is
the only solution of~$\bigl(1-\gamma_i^{p^m}\bigr)(y)=z$. We
deduce that~$1-\gamma_i^{p^m}$ is bijective on~$W^{\cH_W,(i)}_m$.
Furthermore, since~$y_z$ is the limit of the sequence
$g_{i,z}^n(z)$ and~$g_{i,z}(z)-z=f_i(\gamma_i^{p^m}(z)\bigr)$, we
have~$v(y_z-z)\geq v\bigl(g_{i,z}(z)-z\bigr)\geq
v\bigl(\gamma_i^{p^m}(z)\bigr)-c_{3,\cH_W}$. Hence,
$\bigl(1-\gamma_i^{p^m}\bigr)^{-1}$ is continuous
on~$W^{\cH_W,(i)}_m$.

\

We are ready to apply the considerations above in the cases of
interest to us. Let\/~$S$ be as in~\ref{H}. Let\/~$M$ be a
$\ZZ_p$--representation of\/~$\cG_S$. Let~$M\cong \ZZ_p^a
\dirsum_{i=1}^b \ZZ_p/p^{c_i}\ZZ_p$. For $A=\Rbarhat\tensor_V K$,
$\AA_\Rbar$, $\AA_\Rbar^\dagger$, $\AAtilde_\Rbar$
or~$\AAtilde_\Rbar^\dagger$, then $M\tensor_{\ZZ_p}
A=A^a\dirsum_{i=1}^b \bigl(A/p^{c_i}A\bigr)$. We
consider~$M\tensor_{\ZZ_p}A$ as topological module for the product
topology considering on~$A$ the topology induced from the
$p$--adic topology on~$\Rbarhat\tensor_V K$ or form the weak
topology on~$\AA_\Rbar$ and considering on each~$A/p^{c_i}A$ the
quotient topology.

\label induced. proposition\par \prop We have: \spacing

\item{{\rm 1)}}  the ring~$\widetilde{\Lambda}:=\Rbarhat\tensor_V
K$ with~$v(b):=\min\left\{\alpha\in\QQ\vert {b\over p^\alpha}\in
\Rbarhat\right\}$ satisfies (TS1). Furthermore, the following
holds

\itemitem{{\rm (TS1') }} for every $c \in\RR_{>0}$ and every open
normal subgroups $H_{1}\subseteq H_{2}$ of\/ $\cH$ (resp.~of\/
$\H$), there exists $\alpha\in  \Rbarhat^{H_1}$ such that
$\sum\limits_{\tau\in H_{2}/H_{1}}\tau(\alpha)$ is an element
of\/~$\Vbar$ of valuation~$\leq c$;

\spacing

\item{{\rm 2)}} for every~$r\in\QQ_{>0}$ the
ring~$\widetilde{\Lambda}:=\AAtilde_\Rbar^{(0,r]}$ with~$v=w_r$
satisfies (TS1);

\spacing

\item{{\rm 3)}} for every~$N\in\NN$ the
ring~$\widetilde{\Lambda}:=\AAtilde_\Rbar/p^{N+1}\AAtilde_\Rbar$
with~$v=\vv_\EE^{\leq N}$ satisfies (TS1);\spacing

\item{{\rm 4)}}
$\H^i\bigl(\cH_S,M\tensor_{\ZZ_p}(\Rbarhat\tensor_V K)\bigr)=0$
for every $i\geq 1$;\spacing

\item{{\rm 5)}} {{\rm (a)}}
$\H^i\bigl(\cH_S,M\tensor_{\ZZ_p}\AAtilde_\Rbar\bigr)=0$, {{\rm
(b)}} $ \H^i\Bigl(\cH_S,M\tensor_{\ZZ_p}
\AAtilde_\Rbar^\dagger\Bigr)=0$ for every $i\geq 1$;

\spacing

\item{{\rm 6)}}{{\rm (a)}}
$\H^i\bigl(\H_S,M\tensor_{\ZZ_p}\AAtilde_\Rbar\bigr)=0$, {{\rm
(b)}} $\H^i\Bigl(\H_S,M\tensor_{\ZZ_p}
\AAtilde_\Rbar^\dagger\Bigr)=0$ for every $i\geq 1$.

\spacing\noindent
\endprop \Proof For open normal
subgroups~$\cH_1\subset \cH_2$ of~$\cH$ claim~(1) follows
from~[\AnBr, Prop~3.4 \& Rmk.~3.5] and claim~(2) follows
from~[\AnBr, Prop.~4.6].

Let~$\H_1\subset \H_2$ be normal subgroups of~$\H$. They
correspond to extensions~$\tRinfty\subset \tSinfty^2 \subset
\tSinfty^1$ which are finite and Galois
over~$\tRinfty\bigl[p^{-1}\bigr]$ of degree~$d_1$ and~$d_2$
respectively. In particular, there exists an extension
$\Vinfty\subset \Vinfty^\natural$, finite and Galois after
inverting~$p$, such that they arise by taking the normalization of
the base change of extensions of~$\Rinfty\tensor_{\Vinfty}
{\Vinfty^\natural}$ finite and Galois after inverting~$p$ of
degree~$d_1$ and~$d_2$ respectively. This is equivalent to require
that there exist open normal subgroups of~$\cH_1\subset \cH_2$
of~$\cH$ such that~$\cH_1\cap \H=\H_1$, $\cH_2\cap \H=\H_2$
and~$\H_1/\H_2\cong \cH_1/\cH_2$. Then, (TS1) (resp.~(TS1'))
for~$\cH_1\subset \cH_2$ implies (TS1) (resp.~(TS1'))
for~$\H_1\subset \H_2$. Hence, (1) and~(2) follow.

(3) Let~$H_1 \subset H_2$ be open normal subgroups of~$\cH$
(resp.~$\H$) and let $\alpha_r$ be an element
of~$\AAtilde_\Rbar^{(0,r]}$ satisfying~(TS1). If we write
$\alpha_r:=\sum_k p^k [z_k]$ with~$z_k\in \EEtilde_\Rbar$,
since~$w_r(\alpha_r)=\inf\{r \vv_\EE(z_k)+r\}\geq -c_1$, we get
that $\vv_\EE(z_k)\geq {-c_1-k \over r}$. Hence, for
every~$N\in\NN$ we have $\vv_\EE^{\leq N}(\alpha_r)\geq
{-c_1-N\over r}$. In particular, the claim follows.

(4)--(6) Claim~(4) follows from~\ref{computecoh}.
Since~$\AAtilde_\Rbar^\dagger/p^{N+1}
\AAtilde_\Rbar^\dagger=\AAtilde_\Rbar/p^{N+1} \AAtilde_\Rbar$,
claims~(5) and~(6)  follow from~\ref{computecoh} for~$M$ a
$\cG_S$--representation which is free as
$\ZZ_p/p^{N+1}\ZZ_p$--module. In particular, (5) and~(6) hold for
torsion representations.

We may then assume that~$M$ is torsion free.
Let~$A:=\AAtilde_\Rbar$ or~$\AAtilde_\Rbar^\dagger$ and
let~$H=\cH_S$ or~$\H_S$. Let~$f$ be an $i$--th cocycle of~$H$ with
values in~$M\tensor_{\ZZ_p} A$, continuous for the weak topology.
If~$A=\AAtilde_\Rbar^\dagger$, then~$\AAtilde_\Rbar^\dagger=\cup_r
\AAtilde_\Rbar^{(0,r]}$  and~$\AAtilde_\Rbar^{(0,r]}$ is open
in~$\AAtilde_\Rbar$ for the weak topology since it
contains~$\AAtilde_\Rbar^+$. In particular, since~$H$ is compact
in this case~$f$ takes values in~$M\tensor_{\ZZ_p}
\AAtilde_\Rbar^{(0,r]}$ for some~$r$.

Since~$f$ is continuous, for every~$n\in\NN$ there exists an open
normal subgroup~$H_n$ of~$H$ such that the composite
$\bar{f}_n\colon H^r \to M\tensor_{\ZZ_p} A \to M\tensor_{\ZZ_p}
\bigl(A/ (U_{n+1, [{c_1+ n\over r}]+n } \cap A)\bigr)$ factors
via~$\bigl(H_0/H_n\bigr)^r$; see~\ref{AA} for the
notation~$U_{n,h}$. Here, for~$u\in\QQ$ we write $[u]$ for the
smallest positive integer bigger or equal to~$u$. Let~$f_n$ be the
composite of $\bar{f}_n$ with a splitting $M\tensor_{\ZZ_p}
\bigl(A/ (U_{n+1, [{c_1+ n\over r}]+n }\cap A)\bigr) \to
M\tensor_{\ZZ_p} A$ (as sets). Then, $f_n$ is a continuous cochain
and we also have $f_n\equiv f$ in $M\tensor_{\ZZ_p} \bigl(A/
(U_{n+1,  [{c_1+ n\over r}]+n }\cap A) \bigr)$. For every~$n$
let~$h_n:=\alpha\cup f_n$ be the continuous $i-1$--cochain defined
as in the proof of~\ref{Hhasapproximatedtricialcoho}(1). The
computations in loc.~cit.~show that $f_n-\partial h_n\equiv \alpha
\cup \partial f_n \equiv 0$ and $h_{n+1}\equiv h_n $ in
$M\tensor_{\ZZ_p} \bigl(A/ (U_{n+1, n }\cap A) \bigr)$. Then,
$\{h_n\}$ is Cauchy for the weak topology and $\{\partial h_n\}_n$
converges to~$f$ for the weak topology. In particular, $h_n$
converges to a continuous $i-1$--cochain~$h$ with values
in~$M\tensor_{\ZZ_p}\AAtilde_\Rbar$ and~$\partial h=f$.

If~$A=\AAtilde_R$, this concludes the proof.
If~$A=\AAtilde_\Rbar^\dagger$, since $p^n \AAtilde_\Rbar\cap
\AAtilde_\Rbar^{(0,r]}=p^n \AAtilde_\Rbar^{(0,r]}$ and
since~$w_r(p)=1$ and~$w_r(\pi)\geq {p r\over p-1}$ by~[\AnBr,
Prop.~4.3(d)], we conclude that $\{h_n\}$ is Cauchy for the
$w_r$--adic topology as well. Since~$\AAtilde_\Rbar^{(0,r]}$ is
complete and separated for the $w_r$--adic topology by~[\AnBr,
Prop.~4.3(c)], we conclude that~$h$ in fact takes values
in~$M\tensor_{\ZZ_p}\AAtilde_\Rbar^{(0,r]}$. The conclusion
follows.

\

\ssection Sen's theory for $\Rbarhat\bigl[p^{-1}\bigr]$\par Before
passing to the $(\phi,\Gamma)$--modules,  we first show that our
theory applies in the case of
$\Rbarhat\bigl[p^{-1}\bigr]$--representations. These results are
due to Sen [\Sen], in the classical case of a dvr with perfect
residue field, and are due to [\BrinonSen] for a dvr with
imperfect residue field. The key point is of course to show
that~\ref{decompletionofcoho} applies. This follows essentially
from results proven in~[\AnBr]. We review some of the basic
definitions and properties from loc.~cit.

Let~$S$ be a $R$--algebra as in~\ref{V}. Fix $m_{0,S}\in\NN$ such
that
$p^{m_{0,S}}\bigoplus\limits_{i=1}^{d}\ZZ_{p}\gamma_{i}\subseteq\Gamma_S$.
Then,  for every $m\geq m_{0,S}$, the ring $S_{m+1}[p^{-1}]$ is a
free $S_m[p^{-1}]$-module of rank $p^{d+1}$ (resp.~$S_{m+1}.
\Vbar[p^{-1}] $ is a free $S_m.\Vbar[p^{-1}]$-module of rank
$p^d$). For every $i\in\{ 1,\ldots,d\}$ and every $n\in\NN$,
define
$$S_{n,\ast}^{(i)}=S\big[ T_1^{{1\over
p^n}},\ldots,T_{i-1}^{{1\over p^n}},T_{i+1}^{{1\over
p^n}},\ldots,T_{d}^{{1\over p^n}}\big].V_{n}$$and
$$S_{n,\ast}^{'(i)}=\tS\big[ T_1^{{1\over
p^n}},\ldots,T_{i-1}^{{1\over p^n}},T_{i+1}^{{1\over
p^n}},\ldots,T_{d}^{{1\over p^n}}\big].$$For~$i=0$, one puts
$S_{n,\ast}^{(i)}=S\big[ T_1^{{1\over p^n}},\ldots,T_{d}^{{1\over
p^n}}\big]$. Eventually, let $S_{\infty,\ast}^{
(i)}=\bigcup\limits_{n\in\NN}R_{n,\ast}^{ (i)}$ and
$S_{\infty,\ast}^{' (i)}=\bigcup\limits_{n\in\NN}S_{n,\ast}^{'
(i)}$. For every $i\in\{ 0,\ldots,d\}$ and $m\in\NN$, one defines
$$\widehat{S}_{\infty,m,K}^{(i)}=\cases{
\Bigl( \widehat{S_{\infty,\ast}^{(0)}}.V_{m}\Bigr) [p^{-1}] & if
$i=0$,\cr \Bigl( \widehat{S_{\infty,\ast}^{(i)}}.S_{m}\Bigr)
[p^{-1}] & if $i\in\{ 1,\ldots,d\}$,\cr}$$where the hat stands for
$p$--adic completion. Similarly, for $i\in\{ 1,\ldots,d\}$ and
$m\in\NN$, put
$$\widehat{S}_{\infty,m,K}^{'(i)}=\Bigl( \widehat{S_{\infty,\ast}^{'(i)}}.S_{m}\Bigr)
[p^{-1}].$$Note that
$\widehat{S}_{\infty,m,K}^{(i)}\subset\widehat{S}_{\infty}[p^{-1}]$
for every $i\in\{ 0,\ldots,d\}$ and $m\in\NN$ and that
$\widehat{S}_{\infty,m,K}^{'(i)}\subset
\widehat{\tSinfty}[p^{-1}]$. For $n\geq m\geq m_0$ and $x\in
S_{n}[p^{-1}]$, one puts
$$\tau^{(i)}_{m}(x)=\cases{
{1\over p^{n-m}}\Tr_{S_{n}/S_{n,\ast}^{\prime}.V_{m}}(x) & if
$i=0$,\cr {1\over p^{n-m}}\Tr_{S_{n}/S_{n,\ast}^{(i)}.S_{m}}(x) &
if $i\in\{ 1,\ldots,d\}$.\cr}$$For $n\geq m\geq m_0$ and $x\in
\tS_{n}[p^{-1}]$ and every~$i=1,\ldots,d$ define $$t^{(i)}_m(x)=
{1\over p^{n-m}}\Tr_{\tS_{n}/S_{n,\ast}^{' (i)}.S_{m}}(x).$$Such
maps do not depend on~$n$ for~$n\gg 0$ so that they are defined on
$\Sinfty[p^{-1}]$ (resp.~$\tSinfty[p^{-1}]$).

\sprop For every~$i=0,\ldots,d$ and every~$m\geq m_0$ the map
$\tau^{(i)}_{m}$ is continuous for the $p$--adic topology so that
it extends to a unique $\widehat{S}_{\infty,m,K}^{(i)}$--linear
map
$$\tau^{(i)}_{m}\colon \widehat{\Sinfty}[p^{-1}]\llongrightarrow
\widehat{S}_{\infty,m,K}^{(i)}.$$Analogously, for
every~$i=1,\ldots,d$ and every~$m\geq m_0$ the map $t^{(i)}_{m}$
is continuous for the $p$--adic topology so that it extends to a
unique $\widehat{S}_{\infty,m,K}^{'(i)}$--linear map
$$\tau^{(i)}_{m}\colon \widehat{\tSinfty}[p^{-1}]\llongrightarrow
\widehat{S}_{\infty,m,K}^{'(i)}.$$
\endsprop
\Proof The claim for~$\tau_{m}^{(i)}$ follows from~[\AnBr,
Lem.~3.8]. Since~$t_m^{(i)}$ is obtained from~$\tau_m^{(i)}$ by
base--change from~$\Vinfty$ to~$\Vbar$, the claim for~$t_m^{(i)}$
follows as well.

\

Note that~$\bigl(\widehat{\Rbar}\bigl[p^{-1}\bigr]\bigr)^{\cH_S}=
\widehat{\Sinfty}\bigl[p^{-1}\bigr]$
and\/~$\bigl(\widehat{\Rbar}\bigl[p^{-1}\bigr]\bigr)^{\H_S}=
\widehat{\tSinfty}\bigl[p^{-1}\bigr]$ due to~\ref{invariantsEE}.
Furthermore,

\sprop The rings~$\widehat{S}_{\infty,m,K}^{(i)}$ and the
applications~$\tau^{(i)}_{m}$ satisfy~(TS2) and\/ (TS3). The
rings~$\widehat{S}_{\infty,m,K}^{'(i)}$ with the
applications~$t^{(i)}_{m}$ satisfy~(TS4).
\endsprop
\Proof The fact that (TS2) and~(TS3) hold is proven in~[\AnBr,
Prop.~3.9] and in~[\AnBr, Prop.~3.10]. Axiom~(TS4) follows from
this since~$t_m^{(i)}$ is obtained from~$\tau_m^{(i)}$ by
base--change from~$\widehat{\Vinfty}$ to~$\widehat{\Vbar}$.

\slemma We have $\bigcup_m \left( \bigcap_i
\widehat{S}_{\infty,m,K}^{(i)} \right)=\Sinfty\bigl[p^{-1}\bigr]$.
Analogously, we also have\/ $\bigcup_m \left( \bigcap_i
\widehat{S}_{\infty,m,K}^{'(i)} \right)=\Sinfty\tensor_V
\widehat{\Vbar}\bigl[p^{-1}\bigr]$.\endslemma\Proof The first
claim follows from~[\AnBr, Lem.~3.11]. The second is proven as
in~loc.~cit.

\

Let\/~$M$ be a $\ZZ_p$--representation of\/~$\cG_S$ and
let~$W:=M\tensor_{\ZZ_p}\Rbarhat\bigl[p^{-1}\bigr]$. Due
to~\ref{induced} we know that the natural maps
$$
\H^n\bigl(\Gamma_S,W^{\cH_S}\bigr) \longrightarrow
\H^n\left(\cG_S,W\right)\quad\hbox{{\rm
and}}\quad\H^n\bigl(\tGamma_S,W^{\H_S}\bigr) \llongrightarrow
\H^n\left(\G_S,W\right)$$are isomorphisms. Furthermore,

\slabel classicalSen. theorem\par\sthm There exists a finitely
generated, projective $\Sinfty\bigl[p^{-1}\bigr]$--submodule
$N\subset W^{\cH_S}$, stable under~$\Gamma_S$, such that
$N\tensor_{\Sinfty} \widehat{\Sinfty}\cong W^{\cH_S}$ and the
natural map
$$\H^n(\Gamma_S,N) \llongrightarrow \H^n\bigl(\Gamma_S,W^{\cH_S}\bigr)
$$is an isomorphism. Furthermore, if \/ $N':=N\tensor_\Vinfty \widehat{\Vbar}$, then
$N'\tensor_{\tSinfty} \widehat{\tSinfty}\cong W^{\H_S}$ and the
natural map
$$\H^n\bigl(\tGamma_S,N'\bigr) \llongrightarrow \H^n\bigl(\tGamma_S,W^{\H_S}\bigr)
$$is an isomorphism.
\endsthm
\Proof Put~$N$ to be the base change of~$\cD_m(W)$, as defined
in~\ref{defDicDi}, via the natural map $\bigcap_i
\widehat{S}_{\infty,m,K}^{(i)}\to \Sinfty\bigl[p^{-1}\bigr]$.
Similarly, put $N'$ to be the base change of~$\D_m(W)$ via
$\bigcap_i \widehat{S}_{\infty,m,K}^{('i)}\to
\tSinfty\bigl[p^{-1}\bigr]$.

Due to~[\AnBr, Thm.~3.1] there exists an
$\Rbarhat\bigl[p^{-1}\bigr]$--basis~$e_1,\ldots,e_a$ of
$M\tensor_{\ZZ_p} \Rbarhat\bigl[p^{-1}\bigr]$ stable under an open
subgroup~$\cH_W$ of~$\cH_S$, normal in~$\cG_S$.
Let~$\Sinfty[p^{-1}]\subset \Tinfty[p^{-1}]$ be the corresponding
Galois extension. Then~$\cD_m(W)$ (resp.~$\D_m(W)$) is by
construction the set of $\cH_S/\cH_W$--invariants
(resp.~$\H_S/\H_W$--invariants) of the free $\bigcap_i
\widehat{T}_{\infty,m,K}^{(i)}$--module (resp.~$\bigcap_i
\widehat{T}_{\infty,m,K}^{'(i)}$--module) with basis
$e_1,\ldots,e_a$. By~[\An, Cor.~3.11] we have
$\widehat{\Tinfty}[p^{-1}]\cong
\widehat{\Sinfty}[p^{-1}]\tensor_{\Sinfty}\Tinfty$ so that the
extension $\widehat{\Sinfty}[p^{-1}]\subset
\widehat{\Tinfty}[p^{-1}] $ is finite, \'etale and Galois with
group~$\cH_S/\cH_W$. Then, the claims that $N'=N\tensor_\Vinfty
\widehat{\Vbar}$ and that $N$ and $N'$ satisfy the requirements of
the theorem follow from~\ref{decompletionofcoho} and  \'etale
descent.
\endssection

\ssection Sen's theory for $\AAtilde_\Rbar$ and
$\AAtilde_\Rbar^\dagger$\par We  recall some facts proven
in~[\AnBr] needed in order to prove that (TS2) and~(TS3) hold also
for the rings $\AAtilde_\Rbar$ and $\AAtilde_\Rbar^\dagger$.
Let~$S$ be a $R$--algebra as in~\ref{V}. For every~$i=0,\ldots,d$
let
$$\AA_S^{(i)}(\infty):=\cup_n \AA_S\left[\bigl[x_0\bigr]^{1\over
p^n}, \ldots, \bigl[x_{i-1}\bigr]^{1\over
p^n},\bigl[x_{i+1}\bigr]^{1\over
p^n},\ldots,\bigl[x_d\bigr]^{1\over p^n}\right]$$and let
$\overline{\AA_S^{(i)}(\infty)}$ be the closure of
$\AA_S^{(i)}(\infty)$ in $\AAtilde_\Sinfty$ for the weak topology.
Here, we write~$x_0$ for the element~$\epsilon$ and we write
$\bigl[x_i\bigr]$ for the Teichm\"uller lift of~$x_i$. Then,
\endssection

\label recallAnBr. proposition\par\prop For every~$m\geq 0$ and
every~$i=0,\ldots,d$ there exists a homomorphism
$$ \tau^{(i)}_m=\tau^{(i)}_{S,m}\colon\AAtilde_\Sinfty \longrightarrow
\overline{\AA_S^{(i)}(\infty)}\left[ [x_i]^{1\over
p^m}\right],$$called the generalized trace \'a la Tate, such that
\spacing \item{{\rm (i)}} it is
$\overline{\AA_S^{(i)}(\infty)}\left[ [x_i]^{1\over
p^m}\right]$-linear and it is the identity on
$\overline{\AA_S^{(i)}(\infty)}\left[ [x_i]^{1\over p^m}\right]$;
\spacing \item{{\rm (ii)}} it is continuous for the weak topology;
\spacing \item{{\rm (iii)}} it commutes with the action
of\/~$\Gal\left(\Sinfty/R\right)$ and\/~$\tau_m^{(i)}\circ
\tau_n^{(j)}=\tau_n^{(j)}\circ \tau_m^{(i)}$ for~$m$, $n\in\NN$
and\/~$i$, $j\in\{0,\ldots,d\}$; \spacing \item{{\rm (iv)}} for
every $n\in\ZZ$ such that~$m+n\geq 0$ we have $\varphi^n\circ
\tau_{m+n}^{(i)}=\tau_m^{(i)}\circ \varphi^n$; \spacing \item{{\rm
(v)}} it is compatible for varying~$S$ i.~e., given a map of
$R$--algebras $S\to T$ as in~\ref{V} we have that
$\tau^{(i)}_{m,T}$ restricted to~$\AAtilde_\Sinfty$ coincides
with~$\tau^{(i)}_{m,S}$;  \spacing \item{{\rm (vi)}} there
exists~$r_S\in \QQ_{>0}$ such that\/~(TS2) and\/~(TS3) hold for
every~$0<r<r_S$ with\/~$\widetilde{\Lambda}:=\AA_\Rbar^{(0,r]}$
and\/~$v=w_r$, taking\/~$\widetilde{\Lambda}_m^{(i)}$ for~$m\geq
0$ to be the closure of $\AAtilde_\Sinfty^{(0,r]}\cap
\AA_S^{(i)}(\infty)\left[ [x_i]^{1\over p^m}\right]$
in~$\AAtilde_\Rbar^{(0,r]}$ and taking\/~$\tau^{(i)}_{S,m}$
for~$m\geq 0$ to be the restriction of the maps defined in~(i);

\spacing \item{{\rm (vii)}} for every~$N\in\NN$ (TS2) and\/~(TS3)
hold for\/~$\widetilde{\Lambda}:=\AA_\Rbar/p^{N+1} \AA_\Rbar$
and\/~$v=\vv_\EE^{\leq N}$, taking\/~$\widetilde{\Lambda}_m^{(i)}$
for~$m\geq 0$ to be
$\overline{\AA_S^{(i)}(\infty)}/p^{N+1}\overline{\AA_S^{(i)}(\infty)}\left[
[x_i]^{1\over p^m}\right]$  and taking\/~$\tau^{(i)}_{S,m}$
for~$m\geq 0$ to be the reduction modulo~$p^{N+1}$ of the maps
defined in~(i);

\spacing \item{{\rm (viii)}} there exists $m_S\in\NN$ such that
for~$m\geq m_S$ the map $\gamma_i^{p^m}-1$ is an isomorphism
on~$\bigl(1-\tau^{(i)}_m\bigr)\bigl(\AAtilde_\Sinfty\bigr)$ with
continuous inverse (for the weak topology);
\endprop

\Proof Claims~(i)--(v) follow from~[\AnBr, Prop.~4.15]. The
verification of~(TS2) (resp.~of (TS3)) in~(vi) follows
from~[\AnBr, Prop.~4.24] (resp.~[\AnBr, Prop.~4.30\&Prop.~4.32]).
The fact that~(TS2) holds in~(vii) follows from~(ii) and the fact
that the weak topology on~$\AA_\Rbar/p^{N+1} \AA_\Rbar$ is the
$\vv_\EE^{\leq N}$--adic topology.

Since~$\bigl(1-\tau^{(i)}_m\bigr)\bigl(\AAtilde_\Sinfty\bigr)$ is
$p$--adically complete and separated, the fact that
$\gamma_i^{p^m}-1$ is bijective can be verified modulo~$p$ and
follows from~[\AnBr, Prop.~4.30\&Prop.~4.32]. In particular,
$(\gamma_i^{p^m}-1)^{-1}$ is bijective on
$p^{N+1}\bigl(1-\tau^{(i)}_m\bigr)\bigl(\AAtilde_\Sinfty\bigr)$
and, consequently, on
$\bigl(1-\tau^{(i)}_m\bigr)\bigl(\AAtilde_\Sinfty/p^{N+1}\AAtilde_\Sinfty\bigr)$
for every~$N\in\NN$. Note that for every~$h\in\NN$ the
group~$\pi^h\AAtilde_\Sinfty^+$ is contained in the subgroup of
elements~$x\in\AAtilde_\Sinfty^{(0,r]}$ such that~$w_r(x)\geq h
w_r(\pi)$. By~(TS3) for~$\AAtilde_\Rbar^{(0,r]}$ there exist
constants~$c_{3,S}$ and~$c_{4,S}$ such that for every element $z$
in $ \pi^h(1-\tau_m^{(i)})\AAtilde_\Sinfty^+$
(resp.~$\overline{\AA_S^{(i)}(\infty)}\left[ [x_i]^{1\over
p^m}\right]\cap \pi^h\AAtilde_\Sinfty^+$) one has
$$r\vv_\EE(z_k)^{\leq N}\bigl((1-\gamma_i^{p^m})^{-1}(z)\bigr) \geq
w_r\bigl((1-\gamma_i^{p^m})^{-1}(z)\bigr)-N \geq
w_r(z)-c_{3,S}-N$$and, respectively,
$$r\vv_\EE(z_k)^{\leq N}\bigl((1-\gamma_i^{p^m})(z)\bigr)
\geq w_r\bigl((1-\gamma_i^{p^m})(z)\bigr)-N\geq
w_r(z)+c_{4,S}-N.$$Since~$w_r(z)\geq r \vv_\EE^{\leq N}(z)$, we
conclude that $\vv_\EE(z_k)^{\leq
N}\bigl((1-\gamma_i^{p^m})^{-1}(z)\bigr)\geq \vv_\EE^{\leq
N}(z)-{c_{3,S}+N\over r} $ and $\vv_\EE(z_k)^{\leq
N}\bigl((1-\gamma_i^{p^m})(z)\bigr)\geq \vv_\EE^{\leq
N}(z)+{c_{4,S}-N\over r}$. Hence, (vii) and~(viii) follow.

\

Similarly, given  a $R$--algebra~$S$ as in~\ref{V}, for
every~$i=1,\ldots,d$ let ${\tAA_S}^{(i)}(\infty):=\cup_n
\tAA_S\left[\bigl[x_1\bigr]^{1\over p^n}, \ldots,
\bigl[x_{i-1}\bigr]^{1\over p^n},\bigl[x_{i+1}\bigr]^{1\over
p^n},\ldots,\bigl[x_d\bigr]^{1\over p^n}\right]$ and let
$\overline{{\tAA_S}^{(i)}(\infty)}$ be the closure of
${\tAA_S}^{(i)}(\infty)$ in $\tAAtilde_\Sinfty$ for the weak
topology. Note that since~$i\geq 1$ the
ring~$\overline{\AA_S^{(i)}(\infty)}$ contains the closure for the
weak topology of~$\cup_n \AA_V\left[\bigl[x_0\bigr]^{1\over
p^n}\right]$ which is~$\AAtilde_\Vinfty$ by~[\AnBr, Cor.~4.13].
Then,
$\overline{\AA_S^{(i)}(\infty)}\tensor_{\AAtilde_\Vinfty}\AAtilde_\Vbar$
maps to~$\overline{{\tAA_S}^{(i)}(\infty)}$ and the image is dense
for the weak topology. Recall that
$\AAtilde_\Sinfty\tensor_{\AAtilde_\Vinfty}\AAtilde_\Vbar$ injects
and is dense in~$\tAAtilde_\Sinfty$
and~$\AA_S\tensor_{\AA_V}\AA_\Vbar$ injects and is dense
in~$\tAA_S$ by~\ref{invariantsAA}. Hence, for every~$i=1,\ldots,d$
we may base--change $\tau^{(i)}_{S,m}$
via~$\tensor_{\AAtilde_\Vinfty} \AAtilde_{\Vbar}$ and complete
with respect to the weak topology. We obtain a map $$
t^{(i)}_m=t^{(i)}_{S,m}\colon\tAAtilde_\Sinfty \longrightarrow
\overline{{\tAA_S}^{(i)}(\infty)}\left[ [x_i]^{1\over
p^m}\right].$$

\label recallAnBr2. proposition\par\prop The analogues of the
statements~(i)--(viii) of\/~\ref{recallAnBr} hold
for~$\tAAtilde_\Sinfty$, the
rings~$\overline{{\tAA_S}^{(i)}(\infty)}\left[ [x_i]^{1\over
p^m}\right]$ and the maps~$t^{(i)}_m$.
\endprop
\Proof The proposition follows from~\ref{recallAnBr}, from the
construction of~$t^{(i)}_m$ and density arguments. For~(vi) note
that~$\AAtilde_\Sinfty^{(0,r]}\tensor_{\AAtilde_\Vinfty^{(0,r]}}
\AAtilde_\Vbar^{(0,r]}$ maps to ${\tAAtilde_\Sinfty}{ }^{(0,r]}$
and has dense image for the $w_r$--adic topology
by~\ref{invariantsAA}(d).

\

\label AAtildetau=1. lemma\par\lemma We have
$\AAtilde_\Sinfty^{\tau^{(0)}_0=1,\ldots,\tau^{(d)}_0=1}=\AA_S$
and
$\bigl(\tAAtilde_\Sinfty\bigr)^{t^{(1)}_0=1,\ldots,t^{(d)}_0=1}=\tAA_S$.
\endlemma
\Proof By~[\AnBr, Cor.~4.13] the monomials
$\{\bigl[x_0\bigr]^{\alpha_1\over
p^n}\cdots,\bigl[x_d\bigr]^{\alpha_d\over p^n}\}_{0\leq \alpha_i<
p^n}$ form an $\AA_S$--basis of
$\varphi^{-n}(\AA_S)=\AA_S\left[\bigl[x_0\bigr]^{1\over p^n},
\ldots,\bigl[x_d\bigr]^{1\over p^n}\right]$ and $\cup_n
\varphi^{-n}(\AA_S)$ is dense in~$\AAtilde_\Sinfty$ for the weak
topology. In particular, $\cup_n
\varphi^{-n}(\AA_S)^{\tau^{(0)}_0=1,\ldots,\tau^{(d)}_0=1}$ is
dense in $\AAtilde_\Sinfty^{\tau^{(1)}_0=1,\ldots,\tau^{(d)}_0=1}$
and $\cup_n
\varphi^{-n}(\AA_S\tensor_{\AA_V}\AAtilde_\Vbar)^{t^{(0)}_0=1,\ldots,t^{(d)}_0=1}$
is dense
in~$\bigl(\tAAtilde_\Sinfty\bigr)^{t^{(1)}_0=1,\ldots,t^{(d)}_0=1}$
respectively. From the fact
that~$\tau_0^{(i)}\bigl([x_i]^{\alpha\over p^n}\bigr)=0$
for~$0<\alpha<p^n$, we get that $\cap_{i=0}^d \AA_S^{(i)}(\infty)=
\AA_S$ is dense
in~$\AAtilde_\Sinfty^{\tau^{(0)}_0=1,\ldots,\tau^{(d)}_0=1}$ and
$\AA_S\tensor_{\AA_V}\AAtilde_\Vbar$ is dense in
$\bigl(\tAAtilde_\Sinfty\bigr)^{t^{(1)}_0=1,\ldots,t^{(d)}_0=1}$
respectively. The conclusion follows from~\ref{invariantsAA}.

\

\label deftauandtonphiGammamodules. section\par\ssection The
operators~$\tau^{(i)}_m$ and\/~$t^{(i)}_m$ on
$(\varphi,\Gamma)$--modules\par Let\/~$S$ be as in~\ref{H}.
Let\/~$M$ be a $\ZZ_p$--representation of\/~$\cG_S$. If~$M\cong
\ZZ_p^a\dirsum_{i=1}^p \ZZ_p/p^{c_i}\ZZ_p$, then
$M\tensor_{\ZZ_p}\AAtilde_\Rbar\cong
\AAtilde_\Rbar^a\dirsum_{i=1}^b\AAtilde_\Rbar/p^{c_i}\AAtilde_\Rbar$
and on the latter we have the product topology considering
on~$\AAtilde_\Rbar$ the weak topology. Recall that we have defined
$\cDtilde(M)=\bigl(M\tensor_{\ZZ_p}\AAtilde_\Rbar\bigr)^{\cH_S}$
and
$\Dtilde(M):=\bigl(M\tensor_{\ZZ_p}\AAtilde_\Rbar\bigr)^{\H_S}$.
We then define the weak topology on~$\cDtilde(M)$ and $\Dtilde(M)$
to be the topology induced from the inclusions~$\cDtilde(M)\subset
\Dtilde(M) \subset M\tensor_{\ZZ_p}\AAtilde_\Rbar$.

Assume first that~$M\cong \bigl( \ZZ/p^{N+1}\ZZ\bigr)^a$.
Since~$\widetilde{\Lambda}:=\AAtilde_\Rbar/p^{N+1}\AAtilde_\Rbar$
satisfies~(TS1)--(TS4) due
to~\ref{recallAnBr}\&\ref{recallAnBr2}(vii), we may
apply~\ref{defDicDi} and define the operators~$\tau^{(i)}_m$
(resp.~$t^{(i)}_m$) on~$\cDtilde(M)$ and on~$\Dtilde(M)$ and we
get decompositions $$\cDtilde(M):=\cDtilde_m(M)\dirsum
\cDtilde^{(0)}_m(M)\dirsum \cdots \dirsum \cDtilde^{(d)}_m(M)$$and
$$\Dtilde(M):=\Dtilde_m(M)\dirsum \Dtilde^{(1)}_m(M)\dirsum \cdots
\dirsum  \Dtilde^{(d)}_m(M).$$

By devissage we get the operators~$\tau^{(i)}_m$
(resp.~$t^{(i)}_m$) on~$\cDtilde(M)$ and on~$\Dtilde(M)$ and the
decomposition above for any torsion $\cG_S$--representation~$M$.

If~$M$ is torsion free,
$\displaystyle\cDtilde(M):=\lim_{\infty\leftarrow
n}\cDtilde(M/p^nM)$ and
$\displaystyle\Dtilde(M):=\lim_{\infty\leftarrow n}
\Dtilde(M/p^nM)$ by~\ref{DMtARAA}. Using the construction for the
torsion case and passing to the limit, we get the
operators~$\tau^{(i)}_m$ (resp.~$t^{(i)}_m$) on~$\cDtilde(M)$ and
on~$\Dtilde(M)$ and the decomposition above.
\endssection

\label compareDicDi. proposition\par\prop Let\/~$S$ be as
in~\ref{H} and let\/~$M$ be a $\ZZ_p$--representation
of\/~$\cG_S$. Then,

\spacing\item{{\rm 1)}}  $\cDtilde_0(M)=\cD(M)$ and
$\Dtilde_0(M)=\D(M)$;

\spacing\item{{\rm 2)}} the operators~$\tau^{(i)}_m$
(resp.~$t^{(i)}_m$) on~$\cDtilde(M)$ (resp.~$\Dtilde(M)$) are
continuous for the weak topology;

\spacing\item{{\rm 3)}} the operators~$\tau^{(i)}_m$
(resp.~$t^{(i)}_m$) preserve~$\cDtildedagger(M)$
(resp.~$\Dtildedagger(M)$). In particular, we have
$$\cDtildedagger(M):=\cDtildedagger_m(M)\dirsum \cDtilde^{\dagger,(0)}_m(M)\dirsum
\cdots \dirsum \cDtilde^{\dagger,(d)}_m(M)$$and
$$\Dtildedagger(M):=\Dtildedagger_m(M)\dirsum \Dtilde^{\dagger,(1)}_m(M)\dirsum \cdots
\dirsum  \Dtilde^{\dagger,(d)}_m(M),$$where the modules on the
right hand side are defined as in~\ref{defDicDi};

\spacing\item{{\rm 4)}} if~$\gamma_i^{p^n}\in \Gamma_S$
(resp.~$\tGamma_S$) then $\gamma_i^{p^n}-1$ is bijective on
$\cDtilde^{(i)}_m(M)$ (resp.~$\Dtilde^{(i)}_m(M)$) with continuous
inverse (for the weak topology);

\spacing\item{{\rm 5)}} if~$\gamma_i^{p^n}\in \Gamma_S$
(resp.~$\tGamma_S$) then $\gamma_i^{p^n}-1$ is bijective on
$\cDtilde^{\dagger,(i)}_m(M)$ (resp.~$\Dtilde^{\dagger,(i)}_m(M)$)
with continuous inverse (for the weak topology);

\spacing\item{{\rm 6)}} $\cDtildedagger_m(M)= \cDtilde_m(M)\cap
\cDtildedagger(M)$,
$\cDtilde^{\dagger,(i)}_m(M)=\cDtilde^{(i)}_m(M)\cap
\cDtildedagger(M)$,
$\Dtilde^{\dagger,(i)}_m(M)=\Dtilde^{(i)}_m(M)\cap\Dtildedagger(M)$
and $\Dtildedagger_m(M)= \Dtilde_m(M)\cap \Dtildedagger(M)$. In
particular, $\cDtildedagger_0(M)=\cDdagger(M)$ and
$\Dtildedagger_0(M)=\D^\dagger(M)$.
\endprop
\Proof Since~$\gamma_i^{p^t}-1=(\gamma_i^{p^s}-1)
\bigl(\sum_{j=0}^{p^{t-s}-1} \gamma_i^{p^sj}\bigr)$ for~$t\geq s$,
it suffices to prove the bijectivity and the existence of a
continuous inverse in~(4) and~(5) for~$n\gg 0$. Assuming~(3), we
have $\cDtilde^{\dagger,(i)}_m(M)=\cDtilde^{(i)}_m(M)\cap
\cDtildedagger(M)$ and
$\Dtilde^{\dagger,(i)}_m(M)=\Dtilde^{(i)}_m(M)\cap\Dtildedagger(M)$.
Then, (4) and the bijectivity in~(5) imply the existence of a
continuous inverse in~(5). Claim~(6) follows from the others. For
every~$m$ and~$n\in\NN$ the maps
$$\phi^n\tensor 1\colon
\cDtilde(M)\tensor_{\AA_S}^{\phi^n} \AA_S\isomarrow
\cDtilde(M),\qquad \phi^n\tensor 1\colon
\cDtildedagger(M)\tensor_{\AA_S^\dagger}^{\phi^n}
\AA_S^\dagger\isomarrow \cDtildedagger(M)$$and $$\phi^n\tensor
1\colon \Dtilde(M)\tensor_{\tAA_S}^{\phi^n} \tAA_S\isomarrow
\Dtilde(M),\qquad \phi^n\tensor 1\colon
\Dtildedagger(M)\tensor_{\tAAdagger_S}^{\phi^n}
\tAAdagger_S\isomarrow \Dtildedagger(M)$$are isomorphisms
by~\ref{DMtARAA}(i). It follows from~\ref{recallAnBr}
and~\ref{recallAnBr2}(iv) that  $\bigl(\varphi^n\tensor
1\bigr)\circ \tau_{m+n}^{(i)}=\tau_m^{(i)}\circ
\bigl(\varphi^n\tensor 1\bigr) $ and $\bigl(\varphi^n\tensor
1\bigr)\circ \tau_{m+n}^{(i)}=t_m^{(i)}\circ
\bigl(\varphi^n\tensor 1\bigr)$ and that $\varphi^n\tensor 1$
defines an isomorphism from
$\cDtilde_{m+n}(M)\tensor_{\AA_S}^{\phi^n} \AA_S$ (respectively
from $\cDtilde^{(i)}_{m+n}(M)\tensor_{\AA_S}^{\phi^n} \AA_S$,
respectively from $\Dtilde_{m+n}(M)\tensor_{\tAA_S}^{\phi^n}
\tAA_S$, respectively from
$\Dtilde^{(i)}_{m+n}(M)\tensor_{\AA_S}^{\phi^n} \AA_S$) to
$\cDtilde_m(M)$ (respectively $\cDtilde^{(i)}_m(M)$,
$\Dtilde_m(M)$, $\Dtilde^{(i)}_m(M)$). Hence, it suffices to prove
claims~(2), (4) and~(5) for~$m\gg 0$ to deduce it for
every~$m\in\NN$.

Since $\displaystyle\cDtilde(M):=\lim_{\infty \leftarrow
n}\cDtilde(M/p^nM)$ and
$\displaystyle\Dtilde(M):=\lim_{\infty\leftarrow n}
\Dtilde(M/p^nM)$ by~\ref{DMtARAA} and the operators~$\tau^{(i)}_m$
(resp.~$t^{(i)}_m$) are constructed on each $\cDtilde(M/p^nM)$
(resp.~$ \Dtilde(M/p^nM)$) passing to the limit, to prove~(1), (2)
and~(4) one may assume that~$M$ is a torsion representation. By
devissage one may also assume that~$M$ is a free
$\ZZ/p^{N+1}\ZZ$--module for some~$N\in\NN$. Note
that~$\tau^{(i)}_m$ and~$t^{(i)}_m$ commute with the Galois action
and are compatible with extensions~$\Sinfty\subset\Tinfty$
and~$\tSinfty\subset \tTinfty$ by~\ref{recallAnBr}
and~\ref{recallAnBr2}.  Due to~\ref{reccllASdagger},
\ref{invariantsAA} and~\ref{tAsdaggeretale} and \'etale descent,
it then suffices to prove~(1), (2) and~(4) passing to an extension
$\Sinfty\subset \Tinfty$ in~$\Rbar$ finite, \'etale and Galois
after inverting~$p$ i.~e., for
$\bigl(M\tensor_{\ZZ_p}\AAtilde_\Rbar\bigr)^{\cH_T}$ instead
of~$\cDtilde(M)$ and
$\bigl(M\tensor_{\ZZ_p}\AAtilde_\Rbar\bigr)^{\H_T}$ instead
of~$\Dtilde(M)$. We may then assume that~$\cH_T$, and
hence~$\H_T$, act trivially on~$M$. Claim~(1) follows then
from~\ref{AAtildetau=1}. Claim~(2) follows
from~\ref{recallAnBr}(ii) and~\ref{recallAnBr2}(ii). Claim~(4)
for~$m\gg 0$ follows from~\ref{decompletionofcoho}
since~$\widetilde{\Lambda}:=\AAtilde_\Rbar/p^{N+1}\AAtilde_\Rbar$
satisfies~(TS1)--(TS4). This concludes the proof of~(1), (2)
and~(4).

If~$M$ is a torsion representation, then~(3) and the bijectivity
in~(5) follow from~\ref{DMtARAA}(ii'). Assume that~$M$  is free of
rank~$n$. Thanks to~\ref{reccllASdagger}, \ref{invariantsAA}
and~\ref{tAsdaggeretale} and \'etale descent we may pass to an
extension $\Sinfty\subset \Tinfty$ in~$\Rbar$ finite, \'etale and
Galois after inverting~$p$. By~[\AnBr, Thm.~4.40] there exists
such an extension $\Sinfty\subset \Tinfty$ so
that~$\cD^\dagger(M)\tensor_{\AA_S^\dagger}\AA_T^\dagger $ is a
free $\AA_T^\dagger$--module of rank~$n$. Fix a
basis~$\{e_1,\ldots,e_n\}$ and choose~$r\in\QQ_{>0}$ such that
these elements lie in~$M\tensor_{\ZZ_p}\AAtilde_\Rbar^{(0,r]}$.
By~\ref{DMtARAA}(iii') we have $M\tensor_{\ZZ_p}
\AAtilde^\dagger_\Rbar=\AAtilde^\dagger_\Rbar e_1 \dirsum \cdots
\AAtilde^\dagger_\Rbar e_n$.  For every~$s<\min\{r,r_T\}$,
see~\ref{recallAnBr}(vi)\&\ref{recallAnBr2},
let~$W_s:=M\tensor_{\ZZ_p} \AAtilde^{(0,s]}_\Rbar$.
Define~$\cDtilde^{(0,s]}(M):=W_s^{\cH_T}$
and~$\Dtilde^{(0,s]}(M):=W_s^{\H_T}$. Then,
$\cDtildedagger(M)=\cup_s\cDtilde^{(0,s]}(M)$ and $
\Dtildedagger(M)=\cup_s\Dtilde^{(0,s]}(M)$.

Note that~$\AAtilde^{(0,s]}_\Rbar$ satisfies~(TS1)--(TS4)
by~\ref{recallAnBr}(vi)\&\ref{recallAnBr2}. Hence, the
operators~$\tau^{(i)}_m$ (resp.~$t^{(i)}_m$)
preserve~$\cDtilde^{(0,s]}(M)$ (resp.~$\Dtilde^{(0,s]}(M)$) and we
further have decompositions
$$\cDtilde^{(0,s]}(M):=\cDtilde^{(0,s]}_m(M)\dirsum
\cDtilde^{(0,s],(0)}_m(M)\dirsum \cdots \dirsum
\cDtilde^{(0,s],(d)}_m(M)$$and
$$\Dtilde^{(0,s]}(M):=\Dtilde^{(0,s]}_m(M)\dirsum
\Dtilde^{(0,s],(1)}_m(M)\dirsum \cdots \dirsum
\Dtilde^{(0,s],(d)}_m(M)$$by~\ref{defDicDi}. This proves~(3) in
the overconvergent case. It follows from~\ref{decompletionofcoho}
that if~$\gamma_i^{p^n}\in \Gamma_S$ (resp.~$\tGamma_S$) then
$\gamma_i^{p^n}-1$ is bijective on $\cDtilde^{(0,s],(i)}_m(M)$
(resp.~$\Dtilde^{(0,s],(i)}_m(M)$ for~$m\gg 0$. We conclude that
the bijectivity in~(5) holds. Claim~(5) follows.

\

We deduce from~\ref{decompletionofcoho}, \ref{induced}
and~\ref{compareDicDi} the following theorem which summarizes the
results proven so far:

\label decompletecohomology. theorem\par\thm The natural maps
$$\H^n\bigl(\Gamma_S,\cD(M)\bigr)\llongrightarrow
\H^n\bigl(\Gamma_S,\cDtilde(M)\bigr) \llongrightarrow
\H^n\left(\cG_S,M\tensor_{\ZZ_p} \AAtilde_\Rbar\right),$$

$$\H^n\bigl(\tGamma_S,\D(M)\bigr)\llongrightarrow
\H^n\bigl(\tGamma_S,\Dtilde(M)\bigr) \llongrightarrow
\H^n\left(\G_S,M\tensor_{\ZZ_p} \AAtilde_\Rbar\right),$$

$$\H^n\bigl(\Gamma_S,\cDdagger(M)\bigr)\llongrightarrow
\H^n\bigl(\Gamma_S,\cDtildedagger(M)\bigr) \llongrightarrow
\H^n\left(\cG_S,M\tensor_{\ZZ_p}
\AAtilde_\Rbar^\dagger\right)$$and

$$\H^n\bigl(\Gamma_S,\D^\dagger(M)\bigr)\llongrightarrow
\H^n\bigl(\Gamma_S,\Dtilde^\dagger(M)\bigr) \llongrightarrow
\H^n\left(\G_S,M\tensor_{\ZZ_p} \AAtilde_\Rbar^\dagger\right)$$are
all isomorphisms.

\endthm

\endsection

\section Appendix II: Artin--Schreier theory\par The aim of this section is to prove the
following:

\label phiAA. proposition\par \prop The map $\varphi-1$
on~$\AAtilde_\Rbar$, $\AAtilde_\Rbar^\dagger$, $\AA_\Rbar$
and\/~$\AA_\Rbar^\dagger$ is surjective and its kernel is~$\ZZ_p$.
Furthermore, the exact sequence \labelf exactZZpAA\par$$
0\llongrightarrow \ZZ_p \llongrightarrow \AAtilde_\Rbar
\llongmaprighto{\phi-1} \AAtilde_\Rbar \llongrightarrow
0\eqno{{(\numfo)}}$$admits a continuous right splitting
$\sigma\colon \AAtilde_\Rbar \rightarrow\AAtilde_\Rbar$ (as sets)
so that $\sigma\bigl(\AA_\Rbar\bigr)\subset \AA_\Rbar$,
$\sigma\bigl(\AAtilde^\dagger_\Rbar\bigr)\subset
\AAtilde^\dagger_\Rbar$
and\/~$\sigma\bigl(\AA^\dagger_\Rbar\bigr)\subset
\AA^\dagger_\Rbar$.\endprop \Proof Note that by~[\AnBr, Prop.~4.3]
we have $\varphi\bigl(\AA_\Rbar^{(0,r]}\bigr) \subset
\AA_\Rbar^{(0,r/p]}$ and
$\varphi\bigl(\AAtilde_\Rbar^{(0,r]}\bigr) \subset
\AAtilde_\Rbar^{(0,r/p]}$ so
that~$(\phi-1)(\AA_\Rbar^\dagger)\subset \AA_\Rbar^\dagger$
and~$(\phi-1)(\AAtilde^\dagger_\Rbar)\subset
\AAtilde^\dagger_\Rbar$. We know from~\ref{AA} that~$p$ is a
regular element of~$\AA_\Rbar$ and~$\AAtilde_\Rbar$ and that
$\AA_\Rbar/p\AA_\Rbar=\EE_\Rbar$
and~$\AAtilde_\Rbar/p\AAtilde_\Rbar=\EEtilde_\Rbar$. In
particular, to prove that $\varphi-1$ on~$\AAtilde_\Rbar$
and\/~$\AA_\Rbar$ is surjective and its kernel is~$\ZZ_p$, it
suffices to prove that the kernel of $\varphi-1$
on~$\EEtilde_\Rbar$ is~$\FF_p$ and that $\varphi-1$ is surjective
on~$\EE_\Rbar$ and on~$\EEtilde_\Rbar$. Since $\EEtilde_\Rbar$ is
an integral domain by~\ref{EE}(5), the kernel of $\varphi-1$
is~$\FF_p$. The other claim follows from~\ref{phi1isonto}.

Since~$\AA_\Rbar^\dagger=\AA_\Rbar\cap \AAtilde_\Rbar^\dagger$, to
conclude that $\varphi-1$ is surjective on~$\AA_\Rbar^\dagger$ and
on~$\AAtilde_\Rbar^\dagger$, it suffices to prove that for
every~$x\in \AAtilde_\Rbar^\dagger$ the solutions
$y\in\AAtilde_\Rbar$ of~$(\varphi-1)(y)=x$ lie
in~$\AAtilde_\Rbar^\dagger$. Since any such solutions differ by an
element of~$\ZZ_p$ and the latter is contained
in~$\AAtilde_\Rbar^\dagger$, it suffices to show that $\varphi-1$
is surjective on~$\AAtilde_\Rbar^\dagger$.
Let~$z\in\AAtilde_\Rbar^\dagger$ and choose~$r\in\QQ_{>0}$ so
that~$z\in\AAtilde_\Rbar^{(0,r]}$. Write~$z =\sum_k [z_k] p^k$
with~$z_k\in\EEtilde_\Rbar$. Then, putting~$c=\min\{-1,w_r(z)\}$,
we have~$r \vv_\EE(z_k)+ k \geq c$ for every~$k\in\NN$ i.~e.,
$\vv_\EE(z_k)\geq {c-k\over r}$. By~\ref{phi1isonto}  there
exists~$y_k\in\EEtilde_\Rbar$ such that~$(\phi-1)(y_k)=z_k$ and
$\vv_\EE(y_k)=\vv_\EE(z_k)$ if $\vv_\EE(z_k)\geq 0$
or~$\vv_\EE(y_k)={\vv_\EE(z_k)\over p}$ if~$\vv_\EE(z_k)\leq 0$.
In any case, $\vv_\EE(y_k)\geq {c-k\over pr}$. Hence, $y:=\sum_k
p^k [y_k]$ lies in~$\AAtilde_\Rbar^{(0,pr]}$ and~$(\phi-1)(y)=z$.

\advance\ssnu by-1\slabel phi1isonto. lemma\par \slemma The map
$\varphi-1$ is surjective on~$\EE_\Rbar$, $\EE_\Rbar^+$,
$\EEtilde_\Rbar$ and\/~$\EEtilde_\Rbar^+$. Furthermore, given~$a$
and\/~$b\in\EEtilde_\Rbar$ such that~$a^p-a=b$ we have
$$\vv_\EE(a)=\cases{\vv_\EE(b) & if $\vv_\EE(b)\geq 0$; \cr {\vv_\EE(b)/ p} & if $\vv_\EE(b)\leq 0$. \cr}$$
\endslemma \Proof Recall that~$\EE_\Rbar:=\cup_\Sinfty \EE_S$
(resp.~$\EE_\Rbar^+:=\cup_\Sinfty \EE_S^+$) and the union is taken
over a maximal chain of finite normal extensions of~$\EE_R$
(resp.~$\EE^+_R$), \'etale after inverting~$\barpi$. Then,
$\H^1_\et\bigl(\EE_\Rbar,\ZZ/p\ZZ\bigr)=0$
and~$\H^1_\et\bigl(\EE^+_\Rbar,\ZZ/p\ZZ\bigr)=0$,
$\H^1_\et\bigl(\varphi^{-\infty}(\EE_\Rbar),\ZZ/p\ZZ\bigr)=0$
and~$\H^1_\et\bigl(\varphi^{-\infty}(\EE_\Rbar^+),\ZZ/p\ZZ\bigr)=0$.
By Artin--Schreier theory $\EE_\Rbar/(\varphi-1)\EE_\Rbar $
injects in $\H^1_\et\bigl(\EE_\Rbar,\ZZ/p\ZZ\bigr)$ and, hence, it
is zero. Analogously, $\EE^+_\Rbar/(\varphi-1)\EE^+_\Rbar=0 $.
This implies that
$\varphi^{-\infty}(\EE_\Rbar)/(\varphi-1)\bigl(\varphi^{-\infty}(\EE_\Rbar)\bigr)=0
$ and
$\varphi^{-\infty}(\EE^+_\Rbar)/(\varphi-1)\bigl(\varphi^{-\infty}(\EE^+_\Rbar)\bigr)=0
$. By~\ref{EE} the ring~$\EEtilde_\Rbar^+$ is the $\barpi$--adic
completion of~$\varphi^{-\infty}(\EE_\Rbar^+)$
and~$\EEtilde_\Rbar=\EEtilde_\Rbar^+\bigl[\barpi^{-1}\bigr]$. In
particular,
$\EEtilde_\Rbar=\varphi^{-\infty}(\EE_\Rbar)+\barpi\EEtilde_\Rbar^+
$ and we are left to prove that given a power
series~$b=\sum_{n=1}^{+\infty} b_n \barpi^n$ with~$\{b_n\}_n$ in~$
\varphi^{-\infty}(\EE_\Rbar^+)$, we can solve the equation
$(\varphi-1)(a)=b$. It suffices to find~$\{a_n\}_n$
in~$\varphi^{-\infty}(\EE^+_\Rbar)$ such that~$\barpi^{(p-1) n}
a_n^p-a_n=b_n$. Indeed, if we put~$a:=\sum_{n=1}^\infty
a_n\barpi^n$, then $(\varphi-1)(a)=b$. Given~$b_n\in
\varphi^{-\infty}(\EE^+_\Rbar)$ there exists~$\Sinfty$ and~$m$
such that~$b_n\in \varphi^{-m}\bigl(\EE_S^+\bigr)$. But~$\EE_S^+$
and~$\varphi^{-m}\bigl(\EE_S^+\bigr)$ are $\barpi$--adically
complete and separated, the equation~$\barpi^{(p-1) n} X^p-X=b_n$
in the variable~$X$ has~$1$ as derivative and admits~$b_n$ as
solution modulo~$\barpi$. By Hensel's lemma it admits a unique
solution~$a_n$ in~$\varphi^{-m}\bigl(\EE_S^+\bigr)$. The first
part of the lemma follows.

Assume that~$a^p-a=b$. Then, the properties of~$\vv_\EE^{\leq 1}$
recalled in~\ref{AA} imply that if $\vv_\EE(a)<0$ we have
$\vv_\EE(a^p)=p\vv_\EE(a)<\vv_\EE(a)$ and
$\vv_\EE(b)=\vv_\EE(a^p-a)=p\vv_\EE(a)$. On the other hand, if
$\vv_\EE(a)>0$ we have $\vv_\EE(a^p)=p\vv_\EE(a)>\vv_\EE(a)$ and
$\vv_\EE(b)=\vv_\EE(a^p-a)=\vv_\EE(a)$. The second claim follows.

\

\slabel lemmaspli1. lemma\par\slemma For every~$m$ and~$n\in\NN$
we have $(\phi-1)\bigl([\barpi]^n
\WW\bigl(\EEtilde^+_\Rbar\bigr)+p^m
\WW\bigl(\EEtilde^+_\Rbar\bigr)\bigr)=[\barpi]^n
\WW_m\bigl(\EEtilde^+_\Rbar\bigr)+p^m
\WW\bigl(\EEtilde^+_\Rbar\bigr)$ where~$[\barpi]$ is the
Teichm\"uller lift of~$\barpi$. In particular, the map
$\phi-1\colon \AAtilde_\Rbar\rightarrow \AAtilde_\Rbar$ is open
for the weak topology.\endslemma \Proof By construction
$\{[\barpi]^n \WW_m\bigl(\EEtilde^+_\Rbar\bigr)+p^m
\WW\bigl(\EEtilde^+_\Rbar\bigr)\}_{m,n}$ is a fundamental system
of neighborhoods for the weak topology on~$\AAtilde_\Rinfty$.
Since $\phi-1$ is linear, the first claim implies the second.
Since~$(\phi-1)(p^m a)=p^m(\phi-1)(a)$ for every~$a\in
\WW\bigl(\EEtilde_\Rbar\bigr)$,  since~$\phi-1$ is surjective
on~$\WW(\EEtilde_\Rbar^+)$ by~\ref{phi1isonto} and
since~$\WW_m\bigl(\EEtilde_\Rbar^+\bigr)=\WW\bigl(\EEtilde_\Rbar^+\bigr)/p^m
\WW\bigl(\EEtilde_\Rbar^+\bigr)$, it is enough to prove that for
every~$n$ we have $(\phi-1)\bigl([\barpi]^n
\WW_m\bigl(\EEtilde^+_\Rbar\bigr)\bigr)=[\barpi]^n
\WW_m\bigl(\EEtilde^+_\Rbar\bigr)$. Indeed,
$(\phi-1)\bigl([\barpi]^n \WW_m\bigl(\EEtilde^+_\Rbar\bigr)\bigr)
\subset \bigl([\barpi]^n \WW_m\bigl(\EEtilde^+_\Rbar\bigr)\bigr)$
remarking that $(\phi-1)([\barpi]^n a)=[\barpi]^{p n}
a^p-[\barpi]^n a=[\barpi]^n \bigl([\barpi]^{(p-1) n} a^p-a\bigr)$.
On the other hand, $[\barpi]^n \WW_m\bigl(\EEtilde^+_\Rbar\bigr)
\subset (\phi-1)\bigl([\barpi]^n
\WW_m\bigl(\EEtilde^+_\Rbar\bigr)\bigr)$ since for every~$b\in
\WW_m\bigl(\EEtilde^+_\Rbar\bigr)$ the equation~$[\barpi]^{(p-1)
n} X^p-X=b$ admits a solution modulo~$p$ (cf.~proof
of~\ref{phi1isonto}) and, hence,
in~$\WW_m\bigl(\EEtilde^+_\Rbar\bigr)$ by Hensel's lemma. The
lemma follows.

\slabel lemmaspli2. lemma\par\slemma There exists a  left
inverse~$\rho$ as $\ZZ_p$--modules of the inclusion~$\iota\colon
\ZZ_p \rightarrow \AAtilde_\Rbar$ of~(\ref{exactZZpAA}), which is
continuous  for the weak topology.\endslemma

\Proof Let~$R^*$ be the $p$--adic completion of the localization
of~$R$ at the generic point of~$R\tensor_V k$. We then have a
map~$\AAtilde_\Rbar\to \AAtilde_{\overline{R^\ast}}$, which is
continuous for the weak topology, so that it suffices to
construct~$\rho$ for~$R^\ast$. We may then assume that~$R=R^\ast$
is a complete discrete valuation ring with residue field~$L$. In
particular, $\EE_R$ is a discrete valuation field with valuation
ring~$\EE_R^+$ and~$\AA_R$ is a complete discrete valuation ring
with uniformizer~$p$ and residue field~$\EE_R$.

Recall that~$\EE_\Rbar$ is the union~$\cup_S \EE_S$ over all
finite normal extensions~$R \subset S\subset \Rbar$, \'etale after
inverting~$p$. Let~$R \subset S$ be any such. Since~$R$ is a
complete discrete valuation ring, also~$S$ is a complete discrete
valuation ring. Then, $\EE_S^+$ is a complete discrete valuation
ring. For~$S\subset T \subset \Rbar$ finite normal extensions,
\'etale after inverting~$p$ of degree~$n_{S,T}$, we get
that~$\EE_T^+$ is a finite and torsion free as $\EE_S^+$--module,
of rank~$n_{S,T}$; see~\ref{EE}. By loc.~cit.~the choice of a
$\EE_S^+$--basis of~$\EE_T^+$ determines a
$\phi^{-m}(\EE_S^+)$--basis of~$\phi^{-m}(\EE_T^+)$ for
every~$m\in\NN$ and maps
$$\barpi^{\ell_{S,S'}\over p^m}\EEtilde_\Tinfty^+ \to
\phi^{-m}(\EE_S^+)\tensor_{\phi^{-m}(\EE_S^+)} \EEtilde_\Sinfty^+
\cong \bigl(\EEtilde_{\Sinfty}^+\bigr)^{n_{S,T}} \to
\EEtilde_\Tinfty^+,$$where $\ell_{S,T}$ is a constant depending
on~$S\subset T$. We thus get an isomorphism
$\EEtilde_{\Sinfty}^{n_{S,T}} \to \EEtilde_\Tinfty$ as topological
groups (for the $\barpi$--adic topology). Since~$\EE_S^+$ is
integrally closed in~$\EE_T^+$, we may assume that the given
$\EE_S^+$--basis of~$\EE_T^+$ contains~$1$. Suppose furthermore
that~${\ell_{S,S'}\over p^m}< 1$. We then get a splitting of the
inclusion $\EEtilde_\Sinfty\subset \EEtilde_\Tinfty$ as
$\EEtilde_\Sinfty$--modules such that $\barpi\EEtilde_\Tinfty^+$
is mapped to~$\EEtilde_\Sinfty^+$. Consider the set~${\cal F}$ of
pairs $(A,t)$ where~$A$ is a normal sub--$
\EEtilde_\Rinfty$--algebra of~$\cup_S \EEtilde_\Sinfty$ and
$t\colon A\to \EEtilde_\Rinfty$ is a splitting of the inclusion
$\EEtilde_\Rinfty\subset A$ as $\EEtilde_\Rinfty$--modules such
that $ t\left(A\cap \bigl(\barpi\cdot \cup_S
\EEtilde_\Sinfty^+\bigr)\right)\subset \EEtilde_\Rinfty^+$. It is
an ordered set in which every chain has a maximal element.  Zorn's
lemma implies that~${\cal F}$ has a maximal element which, by the
discussion above, must coincide with~$\cup_S \EEtilde_\Sinfty$. We
conclude that there exists a left inverse $\xi$ as
$\EEtilde_\Rinfty$--modules of the inclusion~$\EEtilde_{\Rinfty}
\subset \cup_S \EEtilde_\Sinfty$ such that $\barpi
\EEtilde_\Sinfty^+$ is mapped to~$\EEtilde_\Rinfty^+$ for
every~$S$. Since~$\EEtilde_\Rbar$ (resp.~$\EEtilde_\Rbar^+$) is
the $\barpi$--adic completion of $\cup_S \EEtilde_\Sinfty$
(resp.~$\cup_S \EEtilde_\Sinfty^+$) and since~$\EEtilde_\Rinfty$
(resp.~$\EEtilde_\Rinfty^+$) is $\barpi$--adically complete and
separated, $\xi$ extends to a left inverse $\zeta$ as
$\EEtilde_\Rinfty$--modules of the inclusion~$\EEtilde_{\Rinfty}
\subset \EEtilde_\Rbar$ mapping $\barpi \EEtilde_\Rbar^+$
to~$\EEtilde_\Rinfty^+$. In particular, $\zeta$ is continuous for
the $\barpi$--adic topology.

On the other hand, recall from~\ref{EE} that $\EE_R^+=L\tensor_k
k_\infty[\![\pi_K]\!]$ and that $\EEtilde_\Rinfty^+$ is the
completion of $\cup_n \EE_R^+\left(\pi_K^{1\over p^n},x_1^{1\over
p^n},\ldots,x_d^{1\over p^n}\right)$ for the topology defined by
the fundamental system of neighborhoods $\left\{\pi_K^m
\left(\cup_n L\tensor_k \EE_V^+[\pi_K^{1\over
p^n}]\left(x_1^{1\over p^n},\ldots,x_d^{1\over
p^n}\right)\right)\right\}_m$. Define
$$\delta\colon \cup_n L\tensor_k k_\infty(\!(\pi_K)\!)\left(\pi_K^{1\over
p^n},x_1^{1\over p^n},\ldots,x_d^{1\over p^n}\right) \to
L\tensor_k k_\infty$$ as the $L$--linear map sending $\pi_K^{i_0}
x_1^{i_1} \cdots x_d^{i_d}$ to~$0$ for
every~$(i_0,i_1,\ldots,i_d)\in {\QQ}^{d+1}$ such
that~$(i_1,\ldots,i_d)$ is not equal to~$0$ in~$(\QQ/\ZZ)^d$. It
is well defined since~$\{\pi_K,x_1,\ldots,x_d \}$ is an absolute
$p$--basis of $\EE_R^+$. Furthermore, $\delta$ is continuous for
the $\pi_K$--topology and, hence, it extends to a continuous left
inverse~$\nu$ as $L$--modules of the inclusion~$L\tensor_k
k_\infty \subset \EEtilde_\Rinfty$ considering the $\barpi$--adic
topology on~$\EEtilde_\Rinfty$ and the discrete topology on~$L$.
Finally, choose a left splitting $\tau$ as $\FF_p$--vector spaces
of~$\FF_p\subset L\tensor_k k_\infty$.

Let~$\delta\colon  \AAtilde_\Rbar \to \ZZ_p$ be the map sending a
Witt vector $(a_0,\ldots,a_n,\ldots)$ of
$\AAtilde_\Rbar=\WW\bigl(\Etilde_\Rbar\bigr)$ to $\bigl(\tau\circ
\nu \circ \zeta(a_0),\ldots, \tau \circ \nu \circ
\zeta(a_n),\ldots\bigr) $. It is a left inverse of the inclusion
$\ZZ_p\subset \AAtilde_\Rbar$ and it is continuous for the weak
topology on~$\AAtilde_\Rbar$ and on~$\ZZ_p$. Note that the
topology induced on~$\ZZ_p$ from the weak topology
on~$\AAtilde_\Rbar$ is the $p$--adic topology. The lemma follows.

\

{\it End of the proof~\ref{phiAA}.} With the notations
of~\ref{lemmaspli2}, let $e:=\iota \circ \rho\colon \AAtilde_\Rbar
\rightarrow \AAtilde_\Rbar$. It is a continuous homomorphism of
$\ZZ_p$--modules and~$e^2=e$. Thus, if~$M:=\Ker(e)=\Im(e-1)$, we
have that~$M$ is closed in~$\AAtilde_\Rbar$
and~$\AAtilde_\Rbar=\ZZ_p\dirsum M$. Then, $(\phi-1)\vert_M\colon
M \rightarrow \AAtilde_\Rbar$ is bijective. It is open thanks
to~\ref{lemmaspli1}. Hence, its inverse is a continuous
homomorphism of $\ZZ_p$--modules. We let~$\sigma$ be the composite
of~$\bigl((\phi-1)\vert_M\bigr)^{-1}$ and the inclusion~$M\subset
\AAtilde_\Rbar$. It satisfies the requirements of~\ref{phiAA}.

\endsection

\

\

\centerline{{\sectionfont References.}}

\item{{[\AM]}} M.~F.~Atiyah, I.~G.~Macdonald: {\it Introduction to
commutative algebra.} Addison-Wesley, Reading,
Mass.~(1969).\spacing

\item{{[\An]}} F.~Andreatta: {\it Generalized ring of norms and
generalized $(\phi,\Gamma)$-modules}, to appear in Ann.~Sci.
E.N.S.\spacing

\item{{[\AnBr]}} F.~Andreatta,  O.~Brinon: {Surconvergence des
repr\'esentations $p$-adiques : le cas relatif}, preprint
(2005).\spacing

\item{{[\BL]}} L.~Berger: {\it Repr\'esentations $p$-adiques et
\'equations differentielles}, Invent.~Math.~{\bf 148}, 219-284,
(2002).\spacing

\item{{[\BG]}} S.~Bosch, U.~G\"ortz: {\it Coherent modules and
their descent on relative rigid spaces}, J.~Reine
Angew.~Math.~{\bf 595}, 119-134, (1998).

\spacing

\item{{[\BrinonSen]}} O.~Brinon: {Un g\'en\'eralisation de la
th\'eorie de Sen}. Math.~Ann.~{\bf 327}, 793-813, (223).\spacing

\item{{[\Ch]}} F.~Cherbonnier: {\it Repr\'esentations $p$-adiques
surconvergents}, thesis Orsay, 1996.\spacing

\item{{[\CCINV]}} F.~Cherbonnier, P.~Colmez: {\it
Repr\'esentations $p$-adiques surconvergents}, Invent.~Math. {\bf
133}, 581--611, (1998).\spacing

\item{{[\CCJAMS]}} F.~Cherbonnier, P.~Colmez: {\it Th\'eorie
d'Iwasawa des repr\'esentations $p$-adiques d'un corps local},
Journal of AMS {\bf 12}, 241-268, (1999).\spacing

\item{{[\ColmezCMP]}} P.~Colmez: {\it Les conjectures de
monodromie $p$-adique}, Sem. Bourbaki {\bf 897}, $54^{\hbox{\rm
\`eme}}$ ann\'ee (2001-2002).

\item{{[\Elkik]}} R.~Elkik: {\it Solutions d\'equations a
coefficients dans un anneau hens\'elien}, Ann.~Sci. E.N.S. {\bf
6}, 553--604, (1973).\spacing

\item{{[\FJAMS]}} G.~Faltings: {\it $p$-Adic Hodge Theory} Journal
AMS {\bf 1}, 255--299, (1988).\spacing

\item{{[\FALAN]}} G.~Faltings: {\it Crystalline cohomology and
$p$-adic Galois representations.} In ``Algebraic Analysis,
Geomtery and Number Theory" (J.I.Igusa ed.), John Hopkins
University Press, Baltimore, 25--80, (1998).

\spacing \item{[\FALAST]} G.~Faltings: {\it Almost \'etale
extensions.} In ``Cohomologies $p$-adiques et applications
arithm\'etiques,'' vol.~II. P.~Berthelot, J.-M.~Fontaine,
L.~Illusie, K.~Kato, M.~Rapoport eds. Ast\'erisque {\bf 279}
(2002), 185--270.

\spacing \item{[\Gabber]} O.~Gabber; {\it Affine analog of the
proper base change theorem}, Israel J.~of Math.~{\bf 87},
325--335, (1994).

\spacing\item{[\EGAIV]} A.~Grothendieck, J.~Dieudonn\'e: {\it
\'El\'ements de G\'eom\'etrie Alg\'ebrique.} Publ.~Math. IHES {\bf
24},  (1964--1965).

\spacing \item{{[\Fo]}} J.-M.~Fontaine: {\it Repr\'esentations
$p$-adiques des corps locaux I.} In ``The Grothendieck
Festschrift", vol 2, Progr.~Math.~{\bf 87}, Birkh\"auser, Boston,
249--309, (1990).\spacing

\item{{[\He]}} L.~Herr: {\it Sur la cohomologie galoisienne des
corps $p$-adiques}, Bull.~Soc.~Math.~France {\bf 126}, 563-600,
(1998).\spacing

\item{{[\IoSt]}} A.~Iovita, G.~Stevens: {\it A cohomological construction 
of $p$-adic families of modular forms}, in preparation.\spacing

\spacing

\item{{[\Jannsen]}} U.~Jannsen: {\it Continuous \'Etale
cohomology}, Math.~Ann.~{\bf 280},  (1988),  no. 2,
207--245.\spacing

\spacing

\item{{[\Sen]}} S.~Sen: {\it Continuous cohomology and $p$--adic
Galois representations.} Invent.~Math.~{\bf 62}, 89-116,
(1980).\spacing

\item{{[\Tate]}} J.~Tate: {\it $p$-divisible groups}. In ``
Proceedings of a conference on local fields", Springer Verlag,
158-183 (1967).




\bigskip
\bigskip

\indent Fabrizio Andreatta

\indent D{\eightpoint IPARTIMENTO} {\eightpoint DI} M{\eightpoint
ATEMATICA} P{\eightpoint URA} {\eightpoint ED} A{\eightpoint
PPLICATA}, U{\eightpoint NIVERSIT\'A} {\eightpoint DEGLI}\break
\indent S{\eightpoint TUDI} {\eightpoint DI} P{\eightpoint ADOVA},
{\eightpoint VIA} T{\eightpoint RIESTE} 63, P{\eightpoint ADOVA}
35121, I{\eightpoint TALIA}

\indent {\it E-mail address}: {fandreat@math.unipd.it}

\

\indent Adrian Iovita

\indent D{\eightpoint EPARTMENT} {\eightpoint OF} M{\eightpoint
ATHEMATICS} {\eightpoint AND} S{\eightpoint TATISTICS},
C{\eightpoint ONCORDIA} {\eightpoint UNIVERISTY}\break \indent
1455 {\eightpoint DE} M{\eightpoint AISONNENEUVE} B{\eightpoint
LV.} W{\eightpoint EST},

\indent H3G 1MB, M{\eightpoint ONTREAL}, Q{\eightpoint UEBEC},
C{\eightpoint ANADA}

\indent {\it E-mail address}: {iovita@mathstat.concordia.ca}

\bye